%% file: Superoptimal.tex
\documentclass[11pt, reqno]{amsart}

\makeatletter
\def\theequation{\@arabic\c@equation}
\usepackage{amsmath,amsthm,amscd,amssymb,latexsym,upref}
\usepackage{versions}
\usepackage{enumerate}
\usepackage{graphicx}
\usepackage[usenames,dvipsnames] {color}

\numberwithin{equation}{section}

\makeindex

\newcommand{\we}{\wedge}
\newcommand{\balpha}{\bar{\alpha}}
\newcommand{\bita}{\bar{\eta}}
\newcommand{\Dd}{\mathbb{D}}
\newcommand{\Tt}{\mathbb{T}}

\newcommand{\Cc}{\mathbb{C}}
\newcommand{\CCmn}{\mathbb{C}^{m \times n}}
\newcommand{\telwe}{\dot{\wedge}}

\newcommand{\eiu}{e^{i\theta}}
\newcommand{\hx}{\hat{x}}
\newcommand{\hy}{\hat{y}}

\newcommand\vmo{\mathrm{VMO}}
\DeclareMathOperator\ran{ran}

\DeclareMathOperator\p{P}
\DeclareMathOperator\clos{clos}

\DeclareMathOperator\Poc{POC}

\def\d{\mathbb{D}}
\def\c{\mathbb{C}}
\def\t{\mathbb{T}}

\def\cala{\mathcal{A}}
\def\calf{\mathcal{F}}

\def\be{\begin{equation}}
\def\ee{\end{equation}}

\def\s0{s_0}
\def\p0{p_0}

\newcommand\cmn{\mathbb{C}^{m\times n}}

 \newtheorem{theorem}{Theorem}[section]
 \newtheorem{corollary}[theorem]{Corollary}
 \newtheorem{lemma}[theorem]{Lemma}
 \newtheorem{proposition}[theorem]{Proposition}

\newtheorem{problem}[theorem]{Problem}

\newtheorem{definition}[theorem]{Definition}
\newtheorem{remark}[theorem]{Remark}

\newtheorem{fact*}{Fact}

\DeclareMathOperator\rank{rank}

\DeclareMathOperator\esssup{ess\;sup}
\newcommand\half{{\tfrac 12}}

\newcommand\id{{\mathrm{id}}}

\renewcommand\k{\mathcal{K}}

\newcommand\h{\mathcal{H}}

\newcommand{\T}{\mathbb{T}}

\newcommand{\D}{\mathbb{D}}
\newcommand{\C}{\mathbb{C}}

\newcommand{\inv}{^{-1}}

\newcommand{\ph}{\varphi}
\renewcommand\phi{\varphi}
\newcommand\al{\alpha}

\newcommand\la{\lambda}

\newcommand\beq{\begin{equation}}

\newcommand\eeq{\end{equation}}
\newcommand\df{\stackrel{\rm def}{=}}

\newcommand\nn{\nonumber}
\newcommand\bbm{\begin{bmatrix}}
\newcommand\ebm{\end{bmatrix}}
\newcommand\bpm{\begin{pmatrix}}
\newcommand\epm{\end{pmatrix}}
\numberwithin{equation}{section}

\newcommand\call{\mathcal{L}}

\theoremstyle{definition}

\begin{document}
\title[Exterior powers and superoptimal approximation] {Exterior products of operators and superoptimal analytic approximation}

\author{Dimitrios Chiotis}
\address{School of Mathematics,  Statistics and Physics, Newcastle University, Newcastle upon Tyne
 NE\textup{1} \textup{7}RU, U.K.}
\email{chiotisd@gmail.com}
\author{Zinaida A. Lykova}
\address{School of Mathematics,  Statistics and Physics, Newcastle University, Newcastle upon Tyne
 NE\textup{1} \textup{7}RU, U.K.}
\email{Zinaida.Lykova@ncl.ac.uk}

\author{N. J. Young}
\address{School of Mathematics, Statistics and Physics, Newcastle University, Newcastle upon Tyne NE1 7RU, U.K.
{\em and} School of Mathematics, Leeds University,  Leeds LS2 9JT, U.K.}
\email{Nicholas.Young@ncl.ac.uk}
\date{5th August 2021}



\thanks{Lykova and Young were supported by the Engineering and Physical Sciences Research Council grant EP/N03242X/1.   Chiotis was supported by a PhD studentship from the School of Mathematics, Statistics and Physics of  Newcastle University}

\begin{abstract} 
We give a new algorithm for the construction of the unique superoptimal analytic approximant of a given continuous
matrix-valued function on the unit circle, making use of exterior powers of operators in preference to spectral or {\em Wiener-Masani} factorizations.
\end{abstract} 

\maketitle
\sloppy
\fussy

\pagenumbering{arabic}
\setcounter{page}{1}
\maketitle
\input contents

\input intro

\input exterior_power

\input algorithm_state

\input orthonormal

\input algorithm_proofs

\input T_j-compact

\input application

\input references
\printindex
\end{document}

%% file: contents.tex
\section*{Contents} \label{contents}

\ref{intro}. Introduction \hfill  Page \pageref{intro}

\ref{history}. History and recent work \hfill   \pageref{history}

\ref{exterior}. Exterior powers of Hilbert spaces
\hfill \pageref{exterior}

\hspace*{1cm} \ref{ext_powers}. Exterior powers   \hfill \pageref{ext_powers}

\hspace*{1cm} \ref{point_w_p}. Pointwise wedge product and  pointwise creation operators
 \hfill \pageref{point_w_p}

\ref{statet_algorithm}.   Superoptimal analytic approximation \hfill \pageref{statet_algorithm}

\hspace*{1cm} \ref{known}. Known results \hfill \pageref{known}

\hspace*{1cm} \ref{Alg_statement}. {Algorithm for superoptimal analytic approximation \hfill \pageref{Alg_statement}

\ref{orthonomal}. Pointwise orthonormality of $\{\xi_i\}_{i=1}^j$ and $\{\bar{\eta}_i\}_{i=1}^j$ almost everywhere on $\Tt$ \hfill \pageref{orthonomal}

\ref{Xjsubset}. The closed subspace $X_{j+1}$ of $H^2(\Dd,\we^{j+2}\Cc^n)$
\hfill \pageref{Xjsubset}

\ref{Y_j-closed}. The closed subspace $Y_{j+1}$ of $H^2(\Dd,\we^{j+2}\Cc^m)^\perp$ \hfill \pageref{Y_j-closed}

\ref{Tj-well-def}. $T_j$ is a well-defined operator \hfill \pageref{Tj-well-def}

\ref{T_1-compact}. Compactness of the operators $T_1$ and $T_2$
\hfill \pageref{T_1-compact}

\ref{T_j-compact}. Compactness of the operator $T_j$
\hfill \pageref{T_j-compact}

\ref{applic}.  Application of the algorithm  \hfill \pageref{applic}

 References \hfill\pageref{bibliog}

%% file: intro.tex
\section{Introduction}\label{intro}

In this paper we put forward a new algorithm for the computation of the superoptimal analytic approximation of a continuous matrix-valued function on the circle, a notion that arises naturally in the context of the classical ``Nehari problem", and also in the ``robust stabilization problem" in control engineering.

To explain the term ``superoptimal", let us start from the elementary observation that a measure of the ``size" of a compact operator $T$ between Hilbert spaces is provided by the operator norm $\|T\|$ of $T$.  However, a single number can only ever provide a coarse measure of the size of a multi-dimensional object, and there is a well-developed classical theory \cite{GK} of {\em `s-numbers'} or {\em `singular values'} of an operator or matrix, which provides much more refined information about an operator than the operator norm.  Consider Hilbert spaces $\h,\k$ and an operator $T:\h\to\k$, and let $j\geq 0$. The quantity $s_j(T)$ is defined to be the distance, with respect to the operator norm, of $T$ from the set of operators of rank at most $j$: \index{$s_j(T)$}
\[
	s_j(T) \df \inf\{\|T-R\|: R\in\mathcal{L}(\h,\k), \rank R \leq j\}.
\]
Here, for  Hilbert spaces $H,K$, we define by $\mathcal{L}(H,K)$ the Banach space of bounded linear operators from $H$ to $K$ with the operator norm. We denote by $\mathcal{K}(H,K)$ the Banach space of compact linear operators from $H$ to $K$ with the operator norm.
In the setting of matrices $T$ (that is, in the case that $\h$ and $\k$ are finite-dimensional), $s_j(T)$ is often called the {\em $j$th singular value of $T$}. \index{singular values!of an operator} In this setting one can show that the singular values of $T$ are precisely the eigenvalues of $\sqrt{T^*T}$.
The largest singular value of $T$ is the spectral radius of $\sqrt{T^*T}$, that is, $\|T\|$, and so clearly the set of all singular values of $T$ contains much more information than the norm $\|T\|$ alone.  The use of $s$-numbers immediately gives rise to a measure of the error in an approximation of an operator- or matrix-valued function.  Consider, for example, an $m \times n$-matrix-valued function $G$ on the unit circle $\T$, and suppose we wish to approximate $G$ by a matrix-valued function $Q$ of a specified form (such as a rational function of a prescribed McMillan degree).   It is natural to regard the difference $G-Q$ as the ``error" in the approximation, and to regard the quantities
\[
     s_j^\infty (G-Q) \df  \esssup_{z\in\t} s_j(G(z)-Q(z))
\]
for $j\geq 0$ as measures of how good an approximation $Q$ is to $G$.  We set \index{ $s^\infty(F)$}
\[
s^\infty (G-Q) \df (s_0^\infty (G-Q), s_1^\infty(G-Q), \dots, s_j^\infty(G-Q), \dots),
\]
and say that $\cala G$ is a {\em superoptimal} approximation of $G$ in a given class $\calf$ of functions if $s^\infty(G-Q)$ attains its minimum with respect to  
 the lexicographic ordering of the set of sequences of non-negative real numbers  
 over  $\calf$  when $Q =\cala G$.

The notion of superoptimality pertains to matricial or operator-valued functions, and is therefore particularly relevant to control engineering and electrical networks more generally, since in these fields one must analyse  engineering constructs whose mathematical representations are typically matrix-valued functions on the circle or the real line.  In particular, a primary application is to
 the problem of designing  automatic controllers for  linear time-invariant plants with multiple inputs and outputs.
Such design problems are often formulated in the frequency domain, that is, in terms of the Laplace or $z-$transform of signals.  By this means
the problem becomes to construct an analytic matrix-valued function in a disc or half-plane, subject to various constraints.
An important requirement is usually to minimize, or at least to bound, some cost- or penalty-function.  In practical engineering
problems a wide variety of constraints and cost functions arise, and the engineer must take account of many complications,
such as the physical limitations of devices and the imprecision of models.  Engineers have developed numerous ways to cope with these
complications \cite{francis,DFT}.  One of them, developed in the 1980s, is $H^\infty$ control theory \cite{robust}.  It is a wide-ranging theory, that makes pleasing contact with some problems and results of classical analysis; a seminal role is played by Nehari's theorem on the best approximation of a bounded function on the circle by an analytic function in the disc.  Also important in the development of the theory was a series of deep papers by Adamyan, Arov and Krein \cite{aak1},\cite{aak2} which greatly extend Nehari's theorem and which apply to matrix-valued functions.

In this context the notion of a superoptimal analytic approximation arose very naturally.  Simple diagonal examples of a $2\times 2$-matrix-valued
function $G$ on $\mathbb T$ show that the set of best analytic approximants to $G$ in the $L^\infty$ norm typically comprises an entire infinite-dimensional ball of functions,
and so one is driven to ask for a stronger optimality criterion, and preferably one which will provide a unique optimum.   The very term ``superoptimal" was coined by engineers even before its existence had been proved in generality.  The paper \cite{superop} proved that the superoptimal approximant does indeed exist, and moreover is unique, as long as the approximand $G$ is the sum of a continuous function and an $H^\infty$ function on the circle.  In engineering examples $G$ is usually rational and so continuous on the circle.

Let us first provide some preliminary definitions and then formulate the problem. Throughout the paper, $\CCmn$ denotes the space of $m\times n$ complex matrices with the operator norm and $\Dd,\Tt$ denote the unit disc and the unit circle respectively.
 \index{$\Tt$}\index{$\Dd$} \index{$\CCmn$}
 
 \begin{definition}\label{1.5.2} Let $E$ be a Banach space.
 
 \noindent \index{$H^{\infty}(\mathbb{D},\Cc^{m\times n} )$} $H^{\infty}(\mathbb{D}, E)$ denotes the space of bounded analytic $E$-valued functions on the unit disk 
 with supremum norm: 
 
 $$\|Q\|_{H^{\infty}} \stackrel{\emph{def}}{=} \|Q\|_{\infty} \stackrel{\emph{def}}{=} \sup\limits_{z \in \mathbb{D}}\|Q(z)\|_{E}.$$
 
 \index{$L^{\infty}(\mathbb{T}, \Cc^{m\times n})$}\noindent$L^{\infty}(\mathbb{T},\Cc^{m\times n} )$ is the space of essentially bounded weakly measurable 
 $E$-valued functions on the unit circle with essential supremum norm
 $$\|f\|_{L^\infty}= \mathrm{ess} \sup\limits_{|z|=1}\|f(z)\|_{E},$$ 
and with functions equal almost everywhere identified.

\noindent Also, $C(\Tt, E )$ is the space of continuous $E$-valued functions from $\Tt$ to $E.$
 
 \end{definition}


Naturally engineers need to be able to {\em compute} the superoptimal approximant of $G$. 
\begin{problem}[The superoptimal analytic approximation problem]\label{mainproblem}
  Given a function $G\in L^{\infty}(\mathbb{T}, \Cc^{m\times n}),$ find a function $Q \in  H^{\infty}(\mathbb{D}, \mathbb{C}^{m\times n})$ such that the sequence  
 $s^{\infty}(G-Q)$ is minimized with respect to the lexicographic ordering.
 \end{problem}
In general, the superoptimal analytic approximant may not be unique. However, it has been proved that if the given function $G$ belongs to
 $H^\infty(\Dd,\Cc^{m\times n})+C(\Tt,\Cc^{m\times n}),$ then Problem \ref{mainproblem} has a unique solution. The following theorem, which was proved by V.V. Peller and N.J. Young in \cite{superop}, asserts what we have just stated.

\begin{theorem}[\cite{superop}, p. 303]\label{1.8}
 Let $G\in H^{\infty}(\Dd , \Cc^{m\times n})+C(\Tt, \Cc^{m\times n} ).$ Then the minimum with respect to the lexicographic ordering of 
 $s^\infty (G-Q)$ over all $Q\in H^\infty (\Dd, \CCmn ) $ is attained at a unique function $\cala G.$ Moreover, the singular values $s_j ( G(z) - \cala G(z))$ 
 are constant almost everywhere on $\Tt$ for $j\geq 0$ .   
\end{theorem}\vspace{3ex} 

The topic of this paper is not the existence and uniqueness of the function $\cala G$ described in Theorem \ref{1.8}, but rather the {\em construction} of $\cala G$.  In the proof of the validity of our construction, we have no compunction in making any use of results proved in \cite{superop}, such as the existence of some special matrix functions.  For example, to justify our algorithm we shall prove, using results of \cite{superop}, that certain operators that we introduce are unitarily equivalent to block Hankel operators, which fact enables us to make use of general properties of Schmidt vectors of Hankel operators, without the need to calculate the symbols of those Hankel operators.

The existence proof in \cite{superop} can in principle be
turned into an algorithm, but into a very computationally intensive one.  The construction is recursive, and at each step of the recursion 
one must augment a column-matrix function to a unitary matrix-valued function on the circle with some special properties.  Computationally this step requires a {\em spectral factorization} of a positive semi-definite matrix-valued function on the circle.  There are indeed algorithms
for this step, but they involve an iteration which may be slow to converge and badly conditioned, especially if some function values have eigenvalues
on or close to the unit circle.  

It is certainly desirable to avoid the matricial spectral factorization step if it is possible to do so.  Our aim in this project was to devise an algorithm
in which the iterative procedures are as few and as well-conditioned as possible.  Iteration cannot be completely avoided; even in the scalar case,
optimal error is the norm of a certain operator, and the best approximant is given by a simple formula involving the corresponding
Schmidt vectors.  Thus one has to perform a singular value decomposition.  In the case that the approximand $G$ is of type $m\times n$ one must expect to solve 
$\min(m,n)$ successive singular value problems.  However,  from the point of view of numerical linear algebra, singular value decomposition is regarded as a fast,
accurate and well-behaved operation.  In this paper we describe an algorithm that is, in a sense, parallel to the construction of \cite{Constr} and that in addition to the spectral factorisation of \emph{scalar} functions,
requires only rational arithmetic and singular-value decompositions.
Several engineers have developed alternative approaches \cite{statesapce},\cite{tsaigu} based on state-space methods.
These too are computationally intensive.

For practical purposes, before even looking for an algorithm for the construction of $\mathcal{A}G$, we need to know that the problem of superoptimal analytic approximation is well posed, in sense that arbitrarily small perturbations of $G$ do not result in large fluctuations in $\mathcal{A}G$. This issue arises even for scalar $G$, and in fact it is known \cite{pellkhr} that, for general continuous functions $G$, $\mathcal{A}G$ does {\em not} depend continuously on $G$. However, Peller and Khruschev have shown in \cite{pellkhr} that,  for $G$ in suitable subspaces $X$ of the continuous functions on $\T$, the best analytic 
approximation operator {\em is} continuous for $\|\cdot \|_X$, and so it makes sense to compute it. A similar assertion holds for matrix-valued functions $G$, as was shown by 
 Peller and Young in \cite{PY-cont}.

We believe that the present method, which makes use of exterior powers of Hilbert spaces and operators,
provides a conceptual approach to the construction of superoptimal approximants  which is a promising basis
for computation. The theoretical justification of the algorithm we present in this paper is lengthy and elaborate. However, the implementation of the algorithm should be straightforward.
It will be very interesting to see whether it leads to an efficient numerical method in the future.

For vector-valued $L^p$ spaces we use the terminology of \cite{NagyFoias}.
\begin{definition}\label{a.12}	
Let $E$ be a separable Hilbert space and let $1\leq p < \infty.$  Define
\begin{enumerate}
 \item[{\rm (i)}] $L^p (\Tt,E)$ to be the normed space of measurable 
(weakly or strongly, which amounts to the same thing, in view of the separability of $E$)
 $E$-valued maps  $f \colon \Tt \to E$ such that 
 $$\|f\|_p = \left(\displaystyle\frac{1}{2\pi}\int_{0}^{2\pi} \|f(e^{i\theta})\|_E^p d\theta\right)^{1/p} <\infty ;$$
      \item[{\rm (ii)}] 
     $H^p (\mathbb{D},E)$ to be the normed space of analytic $E$-valued maps $f\colon\mathbb{D} \to E $ such that
     $$ \|f\|_p = \sup\limits_{0<r<1} \left(\frac{1}{2\pi}\int_{0}^{2\pi} \|f(re^{i\theta})\|_E^p d\theta\right)^{1/p} < \infty,$$
the left hand side of this inequality defining a norm on $H^p(\D,E)$.
    \end{enumerate}
   \end{definition}
  
Our algorithm provides a solution $\mathcal{A}G$ to Problem \ref{mainproblem}. 
By computing the value of each $t_k$ at every step, we obtain each term $s_k^\infty(G-\mathcal{A}G)$ of the sequence $s^\infty(G-\mathcal{A}G).$   First we need the notion of a {\em Hankel operator} and the definitions of some long-established standard function spaces; for a more detailed account of these spaces see \cite[Chapter V]{NagyFoias}.

If $E$ is a separable Hilbert space, then every function $f\in H^2(\D,E)$ has a radial limit at almost every point of $\T$, by a theorem of Fatou \cite[Chapter V]{NagyFoias}, and the map that takes a function $f\in H^2(\D,E)$ to its radial limit function embeds $H^2(\D,E)$ isometrically in $L^2(\T,E)$.  In this paper we shall only envisage the case that $E$ is separable, and so we can always regard $H^2(\D,E)$ as a closed subspace of $L^2(\T,E)$.
The operators $P_+, P_-$ on $L^2(\T,E)$ are the operators of orthogonal projection onto the closed subspaces $H^2(\Dd, E)$ and 
${H^2(\Dd, E)}^\perp \df L^2(\T,E) \ominus H^2(\Dd,E)$.
\index{$P_{-}$}\index{$P_{+}$}\index{$H^{2}(\Dd, E)^{\perp}$}
 
\begin{definition}\label{defHankel} 
Let $E$ be a separable Hilbert space, and let $\ph$ be an essentially bounded measurable $\call(E)$-valued function on $\T$; then the {\em Hankel operator} $H_\ph$ is the operator from $H^2(\Dd, E)$ to 
${H^2(\Dd, E)}^\perp$ given by \index{Hankel operator}
\[
H_\ph x = P_-(\ph x) \quad \mbox{ for } x\in H^2(\Dd, E), \; \text{ where} \; (\ph x)(z)= \ph(z)x(z) \mbox{ for } \; z \in \T.
\]
\end{definition}

\begin{remark}\index{unitary operator}
	In this paper we call an operator $U\colon H \to K$ between Hilbert spaces $H,K$ a \emph{unitary operator} if $U$ is both isometric and surjective. Some authors restrict the name "unitary operator" to the case that $H=K.$ Such authors would use a terminology like "isometric isomorphism" for our "unitary operator" in the case that $H\neq K.$  
\end{remark}

For any vector $x$ in a Hilbert space $E$, we denote by $x^*$ the linear functional $ \langle \cdot, x \rangle_E $ on $E$. For an $E$-valued function $x$ on a set $S \subset \C$, we 
define $E^*$-valued function $x^*$ on $S$ by $x^*(z)= x(z)^*$  for all $z \in S$. We observe that if $ x \in L^p(\T,E)$, where $1 \le p \le\infty$, then $ x^* \in L^p(\T,E)$ and 
$\|x^*\|_p =\|x\|_p$.

If $ x, y\in E$, then $x y^*$ denotes the operator of rank one on $E$ defined by $x y^*(u)=
\langle u, y\rangle_E x $ for all $u \in E$. This operator is sometimes denoted by $x \otimes y$ (see, for example, \cite[equation (1.17)]{AMY}). If  $ x, y$ are $E$-valued functions on a set $S \subset \C$,
then $x y^*$ is the function from $S$ to ${\mathcal L}(E)$ given by $x y^*(z)= x(z) y(z)^*$  for all $z \in S$.

\begin{definition}[\cite{Young}, p. 206]
	 		 	Let $H,K$ be Hilbert spaces and let $T\colon H \to K$ be a compact operator. Suppose that $s$ is a singular value of $T.$ A \emph{Schmidt pair}
	 	for $T$ corresponding to $s$ is a pair $(x,y)$ of non-zero vectors, with $x\in H, \; y \in K,$ such that
	 	$$Tx=sy , \; T^*y=sx.$$\index{Schmidt pair}	 	
	 \end{definition}
The following lemma is elementary.
 \begin{lemma}\label{schmmax}
 	Let $T\in \mathcal{L}(H,K)$ be a compact operator and let $x\in H$, $y\in K$ be such that $(x,y)$ is a Schmidt pair for $T$ corresponding to $s=\|T\|.$ Then $x$ is a maximizing vector for $T,$ $y$ is a maximizing vector for $T^*,$ and $\|x\|_H = \|y\|_K.$ 
 \end{lemma}
 
 \begin{definition}[\cite{NagyFoias}, p. 190]
 {\rm (i)}. The matrix-valued bounded analytic function  $\Theta \in H^\infty(\Dd, \CCmn)$ is called  \emph{inner} if $\;\Theta(e^{it})$ is an isometry from $\Cc^n$ to $\Cc^m$ for almost every $e^{it}$ on $\Tt$. 
 	
 {\rm (ii)}.  An analytic $(m\times n)$-matrix-valued function $ \Phi$ on $\Dd$ is said to be \emph{outer} if  $$\Phi H^2(\Dd, \Cc^n)=\{ \Phi f : f \in H^2(\Dd,\Cc^n) \}$$ is a norm-dense subspace of $H^2(\Dd,\Cc^m)$, and 
 \emph{co-outer} if $$\Phi^TH^2(\Dd,\Cc^m)=\{\Phi^T g: g \in H^2(\Dd,\Cc^m)\}$$ is dense in $H^2(\Dd,\Cc^n)$.
 		\index{matrix-valued function!inner}\index{matrix-valued function! outer}\index{matrix-valued function!co-outer}\index{$\Theta^T$}
 	 \end{definition}  

The following is a brief summary of our algorithm. A full account of all the steps, with definitions and justifications will be given in Section \ref{statet_algorithm}.
Our method makes use of exterior powers $\we^pE$ of a finite-dimensional Hilbert space $E$ and of `pointwise wedge products' of $E$-valued functions $f,g$ on $\Dd$ or $\Tt$, defined by
\[
(f\telwe g)(z) = f(z) \we g(z) \quad \mbox{ for all } z\in\Dd \mbox{ or for all } z\in\Tt.
\]
These notions are explained more fully in Subsection \ref{point_w_p}.

 \index{Algorithm}
\begin{proof}[\emph{\textbf{Algorithm:}}]\let\qed\relax
		For a given $G \in H^\infty ( \Dd, \CCmn) + C(\Tt, \CCmn),$ the superoptimal analytic approximant $\mathcal{A}G \in H^\infty(\Dd,\CCmn)$ can be constructed as follows.   		

i) \textbf{ Step 0.} Let $T_0= H_G$ be the Hankel operator with symbol $G$. Let $ t_0 = \| H_G\| .$  If $t_0 = 0,$ then $H_G=0,$ which implies $G\in H^\infty ( \Dd, \CCmn).$ In this case, the algorithm terminates, we define $r$ to be zero and the superoptimal approximant $\mathcal{A}G$ is given by $\mathcal{A}G = G.$
			
			\noindent Let $t_0\neq 0.$ The Hankel operator $H_G$ is a compact operator and so there exists a Schmidt pair $(x_0 , y_0)$ corresponding to the singular value $t_0= \|H_G\|$ of $H_G.$ By the definition of a Schmidt pair $(x_0,y_0)$,			
			$$x_0 \in H^2(\Dd, \Cc^n ),\quad y_0 \in H^2 (\Dd, \Cc^m)^\perp $$
are non-zero vector-valued functions such that
			$$H_Gx_0 = t_0 y_0 , \quad H_G^* y_0 =t_0 x_0. $$
The functions $x_0 \in H^2(\Dd,\Cc^n)$ and $\bar{z}\bar{y}_0 \in H^2(\Dd, \Cc^m)$ admit the inner-outer factorizations 
			\begin{equation}\label{xi0eta01} x_0 = \xi_0 h_0 , \quad \bar{z}\bar{y}_0 = \eta_0 h_0 \end{equation} 
for some scalar outer factor $h_0 \in H^2(\Dd, \Cc)$ and column matrix inner functions $\xi_0\in H^\infty(\Dd, \Cc^n)$, $ \eta_0\in H^\infty(\Dd, \Cc^m). $ Then,
			\begin{equation}\label{equl1}  \|x_0 (z)\|_{\Cc^n} = |h_0(z)| = \|y_0 (z)\|_{\Cc^m} \;\; \text{almost everywhere on }\; \Tt.\end{equation}
\noindent We write equations (\ref{xi0eta01}) as 
			\begin{equation}\label{3112}
			\xi_0 = \frac{x_0}{h_0}, \quad \eta_0 = \frac{\bar{z}\bar{y}_0}{h_0}.
			\end{equation}
			\noindent Then \begin{equation}\label{xi0eta0=11}
			\| \xi_0 (z)  \|_{\Cc^n} =1= \| \eta_0(z) \|_{\Cc^m}\; \text{almost everywhere on}\; \Tt. \end{equation}
			
			There exists a function $Q_1 \in H^\infty(\Dd, \CCmn)$ which is at minimal distance from $G$; any such function satisfies  
			\begin{equation}\label{G-Q01} (G-Q_1)x_0 = t_0 y_0,\quad y_0^* (G-Q_1) = t_0 x_0^*.
			\end{equation}
Choose any function $Q_1 \in H^\infty(\Dd, \CCmn)$ which satisfies the equations \eqref{G-Q01}.

ii) \textbf{Step 1.} Let			
\begin{equation}\label{X_01} X_1 = \xi_0 \dot{\we} H^2(\Dd,\Cc^n ) \subset H^2(\Dd, \we^2 \Cc^n ), \vspace{2ex}\end{equation}
and let
			\begin{equation}\label{Y_01}   Y_1 = \bar{\eta}_0 \dot{\we} H^2(\Dd, \Cc^m)^\perp \subset H^2 (\Dd, \we^2 \Cc^m)^\perp.\end{equation}	
$X_1$ is a closed linear subspace of $H^2(\Dd,\we^{2}\Cc^n)$.
$Y_1$ is a closed linear subspace of $H^2 (\Dd, \we^2 \Cc^m)^\perp. $

			\noindent Define the operator $$T_1 : X_1 \to Y_1$$  by 
			\begin{equation}\label{T_01} T_1 ( \xi _0 \dot{\we} x ) = P_{Y_1} (\bar{\eta}_0\dot{\we} (G-Q_1)x) \; \text{for all} \; x\in H^2(\Dd, \Cc^n ),\end{equation} where $P_{Y_1}$ is the projection from $L^2(\Tt, \we^2 \Cc^m)$ on $Y_1.$ We show that $T_1$ is well-defined.

			If $T_1=0,$ then the algorithm terminates, we define $r$ to be $1$ and the superoptimal approximant $\mathcal{A}G$ is given by the formula
			$$G- \mathcal{A}G = \displaystyle\sum\limits_{i=0}^{r-1}\frac{t_i y_i x_i^*}{|h_i|^2}= \frac{t_0 y_0 x_0^*}{|h_0|^2}, $$
			and the solution is $$\mathcal{A}G =G - \frac{t_0 y_0 x_0^*}{|h_0|^2} .$$
			
			If $T_1\neq 0,$ let $t_{1} = \|T_1\| >0.$   $T_1$ is a compact operator and so
			there exist $v_1 \in H^2(\Dd, \Cc^n),\; w_1 \in H^2(\Dd, \Cc^m)^\perp$ such that
			$(\xi_0 \telwe v_1 , \bar{\eta}_0 \telwe w_1) $ is a Schmidt pair for $T_1$ corresponding to $t_1.$ Let $h_1$ be the scalar outer factor of
			$ \xi_0 \telwe v_1$ and let 
			\begin{equation}\label{x1y1eq1} x_{1} = ( I_{\Cc^n} - \xi_0 \xi_0^*)v_1, \;\; y_1= (I_{\Cc^m} - \bar{\eta}_0 \eta_0^T)w_1 ,\end{equation} where $\mathrm{I}_{\Cc^n}$ and $\mathrm{I}_{\Cc^m}$ are the identity operators in $\Cc^n$ and $\Cc^m$ respectively.\index{identity operator}
			Then
			\begin{equation}\label{x1y1h11} \|x_1 (z)\|_{\Cc^n} = |h_1(z)| = \|y_1 (z)\|_{\Cc^m}\;\text{almost everywhere on}\; \Tt.\end{equation} 
			There exists a function $Q_2 \in H^\infty(\Dd, \CCmn)$ such that both $s_0^\infty(G-Q_2)$ and $s_1^\infty(G-Q_2)$ are minimized and $$s_1^\infty(G-Q_2)=t_1.$$ Any such $Q_2$ satisfies 
			\begin{equation}\label{31111}\begin{aligned}(G-Q_2) x_0 = t_0 y_0 , \quad y_0^* (G-Q_2) = t_0 x_0^*  \\
			(G-Q_2)x_1 = t_1 y_1 , \quad y_1^* (G-Q_2)=t_1 x_1^*.\end{aligned}\end{equation} 
			Choose any function $Q_{2} \in H^\infty(\Dd, \CCmn)$ which satisfies the equations \eqref{31111}.			
			Define 
			\begin{equation}\label{311111}
			\xi_1 = \frac{x_1}{h_1}  , \quad \eta_1 = \frac{\bar{z}\bar{y}_1}{h_1}.
			\end{equation} 
                      Then $\| \xi_1(z)  \|_{\Cc^n} =1= \| \eta_1(z) \|_{\Cc^n}$ almost everywhere on $\Tt.$ 
\end{proof}
		\index{pointwise!orthonormal on $\Tt$} 
	\begin{definition}
		Let $E$ be a Hilbert space.  We say that a collection $\{\gamma_j\} $ of elements of $L^2(\Tt,E)$ is \emph{pointwise orthonormal on $\Tt$} if, for almost all $z \in \Tt$ with respect to Lebesgue measure, the collection of vectors $\{ \gamma_j(z)\}$ is orthonormal in $E.$
	\end{definition}		
iii) \textbf{Inductive step}. Suppose we have constructed 
			$$\begin{array}{clllllllll}
			&t_0 \geq t_1 \geq \cdots \geq t_j > 0\\
			&x_0, x_1, \cdots, x_j \in L^2 (\Tt, \Cc^n)\\
			&y_0 , y_1 , \cdots , y_j \in L^2(\Tt, \Cc^m) \\
			&h_0, h_1, \cdots, h_j \in H^2(\Dd,\Cc) \; \text{outer}\\
			& \xi_0,\xi_1, \cdots, \xi_j \in L^\infty(\Tt,\Cc^n)\; \text{pointwise orthonormal  on}\; \Tt \\
			& \eta_0, \eta_1, \cdots , \eta_j \in L^\infty(\Tt,\Cc^m) \;\text{pointwise orthonormal on}\; \Tt \\
			&X_0 = H^2(\Dd,\Cc^n),X_1, \cdots, X_j \\
			&Y_0 = H^2(\Dd,\Cc^m)^\perp, Y_1, \cdots, Y_j\\
			&T_0, T_1, \cdots, T_j \; \text{compact operators}.
			
			\end{array}$$
There exists a function $Q_{j+1} \in H^\infty(\Dd, \CCmn)$ such that
			$$\left(s_0^\infty(G-Q_{j+1}), s_1^\infty(G-Q_{j+1}), \cdots , s_j^\infty(G-Q_{j+1})\right)$$ is lexicographically minimized. Any such function $Q_{j+1}$ satisfies 
			
			\begin{equation}\label{g-qi1}(G-Q_{j+1})x_i = t_i y_i, \quad y_i^* (G-Q_{j+1}) = t_i x_i^*, \quad i=0, 1, \cdots, j. \end{equation}
Choose any function $Q_{j+1} \in H^\infty(\Dd, \CCmn)$ which satisfies the equations \eqref{g-qi1}.			
			
			Define \begin{equation}\label{X_j1} X_{j+1} = \xi_0 \dot{\we} \xi_1 \dot{\we} \cdots \dot{\we} \xi_j \dot{\we} H^2(\Dd,\Cc^n), \end{equation}
and let			
			\begin{equation}\label{Y_j1} Y_{j+1}= \bar{\eta}_0 \dot{\we} \bar{\eta}_1 \dot{\we} \cdots \dot{\we} \bar{\eta}_j \dot{\we} H^2 (\Dd, \Cc^m)^\perp .\end{equation}
$X_{j+1}$ is a closed subset of $H^2(\Dd,\we^{j+2}\Cc^n),$ and  $Y_{j+1}$ is a closed subspace of $H^2 (\Dd, \we^{j+2} \Cc^m)^\perp.$ 
			\noindent Consider the operator $$T_{j+1} : X_{j+1} \to Y_{j+1}$$ given, for all  $x \in H^2(\Dd,\Cc^n)$, by 
			\begin{equation}\label{T_j1}
			T_{j+1}(\xi_0 \dot{\we} \xi_1 \dot{\we} \cdots \dot{\we} \xi_j \dot{\we} x)= P_{Y_{j+1}} \left( \bar{\eta}_0 \dot{\we} \bar{\eta}_1 \dot{\we} \cdots \dot{\we} \bar{\eta}_j\dot{\we} (G-Q_{j+1})x \right).
\end{equation}
$T_{j+1}$ is well defined.
			
			If $T_{j+1}=0,$ then the algorithm terminates, we define $r$ to be $j+1,$ and the superoptimal approximant $\mathcal{A}G$ is given by the formula
			$$G- \mathcal{A}G = \sum\limits_{i=0}^{r-1} \frac{t_i y_i x_i^*}{|h_i|^2}= \sum\limits_{i=0}^{j} \frac{t_i y_i x_i^*}{|h_i|^2}.$$
			
			Otherwise, we define $t_{j+1} = \|T_{j+1}\| >0.$  Then $T_{j+1}$ is a compact operator and hence there exist $v_{j+1} \in H^2(\Dd,\Cc^n), \; w_{j+1} \in H^2(\Dd,\Cc^m)^\perp$ such that \begin{equation}\label{schmpairtj+11}(\xi_0 \dot{\we} \xi_1 \dot{\we} \cdots \dot{\we} \xi_j \dot{\we} v_{j+1},  \bar{\eta}_0 \dot{\we} \bar{\eta}_1 \dot{\we} \cdots \dot{\we} \bar{\eta}_j \dot{\we} w_{j+1})\end{equation}
			is a Schmidt pair for $T_{j+1}$ corresponding to the singular value $t_{j+1}.$ 
			
			\noindent Let $h_{j+1}$ be the scalar outer factor of $\xi_0 \dot{\we} \xi_1 \dot{\we} \cdots \dot{\we} \xi_j \dot{\we} v_{j+1},$ and let
			\begin{equation}\label{xj+1yj+11}x_{j+1} = (I_{\Cc^n} - \xi_0 \xi_0^* - \cdots - \xi_j \xi_j^*)v_{j+1}, \quad y_{j+1} = (I_{\Cc^m} - \bar{\eta}_0 \eta_0^T - \cdots- \bar{\eta}_j\eta_j^T) w_{j+1},\end{equation} and define 
			\begin{equation}\label{xij+1etaj+11}
			\xi_{j+1} = \frac{x_{j+1}}{h_{j+1}}, \quad \eta_{j+1}=\frac{\bar{z}\bar{y}_{j+1}}{h_{j+1}}.\end{equation}
			
			\noindent One can show that $\|\xi_{j+1}(z)\|_{\Cc^n}=1 $ and $\|\eta_{j+1}(z)\|_{\Cc^m}=1$ almost everywhere on $\Tt.$
			
			This completes the recursive step. The algorithm terminates after at most $\min(m,n)$ steps, so that, $r \leq \min(m,n)$ and the superoptimal approximant $\mathcal{A}G$ is given by the formula 
			
			$$G- \mathcal{A}G = \sum\limits_{i=0}^{r-1} \frac{ t_i y_i x_i^*}{|h_i|^2}.$$

				\begin{remark} {\em
				Observe that, in step $j$ of the algorithm, we define an operator $T_j$ in terms of any function $Q_j \in H^\infty(\Dd,\Cc^{m \times n})$ that satisfies the equations 
				\begin{equation}\label{g-q-j} (G-Q_{j})x_i = t_i y_i, \quad y_i^* (G-Q_{j}) = t_i x_i^*, \quad i=0, 1, \cdots, j-1. \end{equation}
				This constitutes a system of linear equations for $Q_j$ in terms of the computed quantities $x_i, t_i$ and $y_i$ for $i=0, \dots, j-1$, and we know, from Proposition \ref{g-qjj}, that the system has a solution for $Q_j$ in  $H^\infty(\Dd,\Cc^{m \times n})$.  By Proposition \ref{Twell} $T_j$ is independent of the choice of $Q_j$ that satisfies equations \eqref{g-q-j}.}
		\end{remark}
		\begin{remark}{\em 
				At each step we need to find $\|T_j\|$ and a Schmidt pair  
				\begin{equation}(\xi_0 \dot{\we} \xi_1 \dot{\we} \cdots \dot{\we} \xi_{j-1} \dot{\we} v_{j},  \bar{\eta}_0 \dot{\we} \bar{\eta}_1 \dot{\we} \cdots \dot{\we} \bar{\eta}_{j-1} \dot{\we} w_{j})
				\end{equation}
				for $T_{j}$ corresponding to the singular value $t_{j}.$ 
				Then we compute the scalar outer factor $h_{j}$ of $\xi_0 \dot{\we} \xi_1 \dot{\we} \cdots \dot{\we} \xi_{j-1} \dot{\we} v_{j}\in H^2(\Dd, \wedge^{j+1}\Cc^n)$.
				These are the only spectral factorisations needed in the algorithm.  Note that if $f\in H^2(\Dd, \Cc^n)$ has the inner-outer factorisation $f=hg$, with $h\in H^2(\Dd,\Cc)$ a scalar outer function and $g\in H^\infty(\Dd,\Cc^n)$ inner, then $(f^*f)(z)=|h(z)|^2 $ almost everywhere on $\Tt$, and so the calculation of $h$ requires us to find a spectral factorisation of the positive {\em scalar-valued} function $f^*f$ on the circle.			}
		\end{remark}
		\begin{remark}{\em
				In a numerical implementation of the algorithm one would need to find a way to compute the norms and Schmidt vectors of the compact operators $T_j$.
				For this purpose it would be natural to choose convenient orthonormal bases of the cokernel $X_j\ominus \ker T_j$ and the range $\ran T_j$.  It is safe to assume that in most applications $G$ will be a rational function, in which case the cokernel and range will be finite-dimensional.  At step 0, $T_0$ is a Hankel operator, and the calculation of the matrix of $T_0$ with respect to suitable orthonormal bases is a known task \cite{Young83}; we believe that similar methods will work for step $j$.
				
			}
		\end{remark}

	In Theorem \ref{mathcalAG} we arrive at the following conclusion about the superoptimal approximant $\mathcal{A}G.$
	 \begin{theorem} 
Let $G\in H^\infty(\Dd, \Cc^{m\times n})+C(\Tt, \Cc^{m\times n}).$ Let $T_i, x_i, y_i, h_i$, for $i\ge 0$, be defined by the algorithm above. Let $r$ be the least index $j \ge 0$ such that $T_j=0$. Then $r\leq \min(m,n)$ and the superoptimal approximant $\mathcal{A}G$ is given by the formula
	 	$$ G-\mathcal{A}G= \displaystyle \sum\limits_0^{r-1} \frac{t_i y_i x_i^*}{|h_i|^2} .$$
	 \end{theorem}

  Wedge products, and in particular pointwise wedge products, along with their properties are studied in detail in Section \ref{exterior}.

\section{History and recent work}\label{history}
The Nehari problem of approximating an essentially bounded Lebesgue measurable  function on the unit circle $\t$ by a bounded analytic function on the unit disk $\d$, has been attracting the interest of both pure mathematicians and engineers since the middle of the 20th century. The problem was first formulated and studied from the  viewpoint of scalar-valued functions, and, in the years that followed, from the operator-valued perspective also, which motivated research into the superoptimal approximation problem.

The Nehari problem in the scalar case first appeared in the paper of Nehari \cite{Neh1}. Given an essentially bounded complex valued function $g$ on $\t,$ one seeks its distance from $H^\infty$ with respect to the essential supremum norm, and wishes to determine for which elements of $H^\infty$ this distance is attained.  It is also of interest to know whether the distance is attained at a uniquely determined function. Such problems have been studied in detail by Nehari \cite{Neh1}, Sarason \cite{sarason} and Adamjan, Arov and Krein in \cite{aak1} and \cite{aak2}. These authors proved that the distance of $g$ from $H^\infty$ is equal to the norm of the Hankel operator  $H_g$ with symbol $g.$  Moreover, if $H_g$ has a maximizing vector in $H^2,$ then the bounded analytic complex-valued function $q$ that minimizes the essential supremum norm $\|g-q \|_{L^\infty}$ is uniquely determined and can be explicitly calculated (see, for example, \cite[p. 196]{Young}). Furthermore, if the essential norm $\|H_g\|_e$ is less than $\|H_g\|$, then $g$ has a unique best approximant.

 Pure mathematicians and engineers started seeking analogues of those results for matrix- and operator-valued functions. These generalizations are not only mathematically interesting, but are essential for applications in engineering, and especially in control theory.  There has accordingly been an explosion of research in this field since 1980, on the part of both pure mathematicians and engineers.

 Page \cite{page} and Treil \cite{treil1} gave various  extensions of the  results of Adamjan, Arov and Krein to operator-valued functions. Page proved that for operator-valued mappings $T\in L^\infty(\t, \mathcal{L}(E_1,E_2)),$  $\inf\{\|T-\Phi \|:\Phi \in H^\infty(\d,\mathcal{L}(E_1,E_2) \}$ is equal to $\|H_T\|$. Here $E_1,E_2$ are Hilbert spaces and $\mathcal{L}(E_1,E_2)$ denotes the Banach space of bounded linear operators from $E_1$ to $E_2.$ Treil extended the Adamjan, Arov and Krein theorem in \cite{aak2} to an operator-valued analogue. 

However, in the matrix-valued setting there are typically infinitely many functions that minimize the $L^\infty$ norm of the error function. This fact is simply illustrated by  the following example.
Let $G(z) = \mathrm{diag} \{ \bar{z} , 0 \}$, for $ z \in \t.$ 
The norm of $H_G$ in this case is easily seen to be $1,$ and hence \emph{all} matrix-valued functions $Q \in H^\infty(\d,\c^{2\times 2})$ of the form $Q(z) = \mathrm{diag}\{ 0 , q(z) \},$ where $q \in H^\infty$ and 
$\|q\|_{H^\infty} \leq 1,$  minimize the norm $\|G-Q\|_\infty$, yielding the error $1$. However, if one goes on to minimize in turn the essential suprema of both singular values of $G(z)-Q(z)$ over $Q \in H^\infty(\d,\c^{2\times 2})$, one finds that such a minimum occurs uniquely when $q(z)$ is equal to $0.$ This type of example suggests that the  enhanced approximation criterion based on successive singular values generates the ``very best" amongst the best approximants to $G$ by an element of $ H^\infty(\d,\c^{2\times 2})$.

Such reflections led to the formulation of a strengthened approximation problem, the superoptimal approximation problem \ref{mainproblem} as explained on above. In \cite{young2} N.J. Young introduced this strengthened notion of optimal analytic approximation, subsequently called {\em superoptimal} approximation. Given a $G$ as above, find a $Q \in H^\infty(\d,\cmn)$ such that the sequence $s^\infty(G-Q)$ is lexicographically minimized. This criterion obviously constitutes a considerable strengthening of the notion of optimality, as one needs to determine a $Q\in H^\infty(\d,\cmn)$ that not only minimizes $\|G-Q\|_{L^\infty},$ but minimizes the $L^\infty$ norm of all the subsequent singular values $s_j(G(z)-Q(z))$ for $j\geq 0$.  

A good starting point for the superoptimal approximation problem of matrix functions is \cite{superop}.  As we have said, the problem is to find, for a given $$G \in H^\infty(\d,\cmn)+C(\t,\cmn), $$ a function $Q \in H^\infty(\d,\cmn)$ such that the sequence $s^\infty(G-Q)$ is lexicographically minimized. Peller and Young proved some requisite preparatory results on ``thematic factorizations", on the analyticity of the minors of unitary completions of inner matrix columns and on the compactness of  some Hankel-type operators with matrix symbols. These results provided the foundation for their main theorem, namely that if $G$ belongs to $H^\infty(\d,\cmn)+C(\t,\cmn),$ then there exists a unique $Q \in H^\infty(\d,\cmn)$ such that the sequence $s^\infty(G-Q)$ is lexicographically minimized as $Q$ varies over $ H^\infty(\d,\cmn)$; moreover for this $Q$, the singular values $s_j(G(z)-Q(z)$ are constant almost everywhere for $z\in\t$, for $j=0,1,2,\dots\;.$ 

Later, in \cite{Constr} Peller and Young presented a conceptual algorithm for the computation of the superoptimal approximant. Their algorithm is based on the theory developed in \cite{superop}. Also in \cite{Constr}, the algorithm was applied to a concrete example of a rational $2\times 2$ matrix-valued function $G$ in $H^\infty(\d,\c^{2\times 2})+C(\t,\c^{2\times 2})$ and the superoptimal approximant $\mathcal{A}G$ was calculated by hand. 

Additionally, Peller and Young in \cite{PY} studied superoptimal approximation by {\em meromorphic} matrix-valued functions, that is, matrix-valued functions in $H^\infty$ that have at most $k$ poles for some prescribed integer $k$. They modified the  results of \cite{superop} and established a uniqueness criterion in the case that the given matrix-valued function $G$ is in $H^\infty + C$ and has at most $k$ poles. In addition, they provided an algorithm for the calculation of the superoptimal approximant. 

One can extend the above results to {\em operator}-valued functions on the circle; the operator-valued superoptimal approximation problem was studied by Peller in \cite{pellroper}. He generalized the notions of \cite{superop} and proved that there exists a unique superoptimal approximant in $H^\infty(\mathcal{B})$ for functions that belong to $H^\infty(\mathcal{B})+C(\mathcal{C}),$ where $\mathcal{B}$ denotes the space of bounded linear operators and $C(\mathcal{C})$ denotes the space of continuous functions on the circle taking values in the space of compact operators. 

Very badly approximable functions, that is, functions that have the zero function as a superoptimal approximant, were studied in the years that followed and a considerable amount of work was published. Peller and Young's paper \cite{superop} provided the motivation for the study of this problem, where they were able to algebraically characterise the very badly approximable matrix functions of class $H^\infty(\d,\cmn)+C(\t,\cmn). $ Their results were extended in \cite{fourblock} to the case of matrix functions $G$ for which $\|H_G\|_e$ is less than the smallest non-zero superoptimal singular value of $G.$ 
Very badly approximable matrix functions with entries in $H^\infty+C$ were completely characterised in \cite{pellerverybadly}.

Recent work in \cite{lpappr} by Baratchart, Nazarov and Peller explores the analytic approximation of matrix-valued functions in $L^p$ of the unit circle by matrix-valued functions from $H^p$ of the unit disk in the $L^p$ norm for $p \leq 2.$ They proved that if a given matrix-valued function $\Psi \in L^p(\t,\cmn)$ is a `respectable' matrix function, then its distance from $H^p(\d,\cmn)$ is equal to $\|H_\Psi\|,$ and they obtained a characterisation of that distance also in the case $\Psi$ is a `weird' matrix-valued function. Furthermore, they established the notion of $p$-superoptimal approximation and illustrated the fact that every $n\times n$ rational matrix function has a unique $p$-superoptimal approximant for $2\leq p <\infty.$ For the case $p=\infty$ they provided a counterexample. 

In a more recent paper of Condori \cite{condori}, the author considered the relation between the sum of the superoptimal singular values of admissible functions in $L^\infty(\t,\cmn)$ and the superoptimal analytic approximation problem in the space $L^\infty(\t,S_p^{m,n}),$ where $S_p^{m,n}$ denotes the space of $m\times n$ matrices endowed with the Schatten-von Neumann norm $\| \cdot \|_{S_p^{m,n}}.$ He illustrated the fact that if $\Phi \in L^\infty(\t, \c^{n\times n})$ is an admissible matrix function of order $k$, then $Q \in H^\infty(\d,\c^{n\times n})$ is a best approximant function under the $L^\infty(\t,S_1^{n,n})$-norm and the singular values $s_j((\phi-Q)(z))$ are constant almost everywhere on $\t$ for  $j=0,1,\dots, k-1$ if and only if $Q$ is a superoptimal approximant to $\Phi,$ 
$$\mathrm{ess}\sup_{z \in \t} s_j((\Phi-Q)(z))=0$$ for $j\geq k,$ and the sum of the superoptimal singular values of $\Phi$ is equal to 
$$\sup \left| \int_{\t} \mathrm{trace} (\Phi(\zeta)\Psi(\zeta))\;dm(\zeta) \right| ,$$ where $m,n>1,$ $1\leq k \leq \min(m,n)$ and the supremum is taken over all $\Psi \in H_0^1(\d,\c^{n\times m})$ for which $\|\Psi\|_{L^1(\t,\c^{n\times m})} \leq 1$ and $\mathrm{rank} \Psi(\zeta) \leq k$ almost everywhere on $\t.$

%% file: exterior_power.tex
\section{Exterior powers of Hilbert spaces}\label{exterior}

In this section  we recall the well-established notion of the wedge product of Hilbert spaces. One can find definitions and properties of wedge products in  \cite{depillis}, \cite{Greub}, \cite{pavan}, \cite{Wed} and \cite{Sim1,Sim2}.
Here we present a concise version of this theory which we need for the new 
superoptimal  algorithm.

\subsection{Exterior powers}\label{ext_powers}
In this subsection, we first present some results concerning the action of permutation operators on tensors, then we recall the definition of antisymmetric tensors and we define an inner product on the space of all antisymmetric tensors. In the following $E$ denotes a Hilbert space.  We shall assume known the notion of the algebraic tensor product of vector spaces, which is precisely explained in \cite{Tensor}.    For the Hilbert space tensor product of Hilbert spaces, see \cite{Dixmier}.
One can find proofs of many statements given below, in our paper \cite{YCL20-3}.
 
\begin{definition}
$\otimes^{p}E$ \index{$\otimes^{p}E$} is the $p$-fold algebraic
 tensor product 
of $E$ and is spanned by tensors of the form $x_1 \otimes x_2 \otimes \dots \otimes x_p,$ where $ x_j \in E$ for
$j=1, \dots,  p.$ \vspace{2ex}
\end{definition}

\begin{definition}\label{inner}
 An inner product on $\otimes^{p}E$ is defined on elementary tensors by
$$\langle x_1 \otimes x_2 \otimes \dots \otimes x_p, 
y_1 \otimes y_2 \otimes \dots \otimes y_p \rangle_{\otimes^{p}E}=p!\langle x_1 , y_1\rangle_E \cdots \langle x_p, y_p\rangle_E,$$ \vspace{3ex}
for any $x_1, \dots, x_p, y_1, \dots, y_p \in E$,
and is extended to 
$\otimes^p E$ by sesqui-linearity. 
\end{definition}

\begin{definition}
$\otimes_H^p E$ is the completion of $\otimes^p E$ with respect to the norm $\|u\| = \langle u, u\rangle^{1/2}_{\otimes^{p}E},$ for $u \in \otimes^p E.$

\end{definition}

\begin{definition}\label{sigma}
Let $\mathfrak{S}_p$ denote the \index{symmetric group} symmetric group on $\{1,\dots,p\},$ with the operation of composition.
 For $\sigma \in \mathfrak{S}_p$, we define 
\[
S_\sigma : \otimes^p E \to \otimes^p E
\]
 on elementary tensors by 
$$\displaystyle S_\sigma(x_1 \otimes x_2 \otimes \dots \otimes x_p)=x_{\sigma(1)} \otimes x_{\sigma(2)} \otimes \dots \otimes x_{\sigma(p)},$$
and we extend $S_\sigma$ to  $\otimes^p E$ by linearity, that is, for $u= \sum_{i=1}^n \la_i x_1^i \otimes \cdots \otimes x_p^i,$ we define 
$$S_\sigma (u) = \sum\limits_{i=1}^{n} \la_i S_\sigma (x_1^i \otimes \cdots \otimes x_p^i ).$$
for any $x_j^i \in E$ and $\la_i\in\C$.

\vspace{2ex}
\end{definition}

\begin{remark}
Clearly if $\sigma$ is a bijective self-map of $\{1, \dots, p\}$, then so is its inverse map $\sigma^{-1}$, and $(\mathfrak{S}_p, \circ)$ is a group under composition. Moreover,
 $$\sigma \circ \sigma^{-1} = \id = \sigma^{-1} \circ \sigma,$$ where $\id \in \mathfrak{S}_p$ is the identity map on $\{1,\dots,p\}.$ Then,  if $\epsilon_{\sigma}$ denotes the signature of the permutation $\sigma$,
$$\epsilon_{\sigma \circ \sigma^{-1}} = \epsilon_\sigma \epsilon_{\sigma^{-1}} =1, $$ hence $\epsilon_\sigma=  \epsilon_{\sigma^{-1}}.$
 
\end{remark}

\begin{proposition}\label{a.3}
Let $E$ be a Hilbert space, and let $p$ be any positive integer. Then, for any $\sigma\in \mathfrak{S}_p$,
 $S_\sigma$ is a linear operator on the normed space $(\otimes^p E , \| \cdot \|)$, which extends to an isometry $\mathbf{S}_\sigma$ on $(\otimes^p_H E,\| \cdot \|) $. Furthermore, $\mathbf{S}_\sigma$ is a unitary operator on $\otimes_H^p E $,  $\mathbf{S}_\sigma^*= \mathbf{S}_{\sigma\inv}$, and therefore
\[
\mathbf{S}_\sigma^*\mathbf{S}_\sigma = \mathbf{S}_{\sigma\inv} \mathbf{S}_\sigma = I
\]
is the identity operator on $\otimes^p_H E$. 
\end{proposition}

Henceforth we shall denote the extended operator $\mathbf{S}_\sigma$ by $S_\sigma$.

\begin{definition}\label{a.4}
A tensor $u \in \otimes_H^{p}E$ is said to be \emph{symmetric} if $S_\sigma(u)=u$ for all $\sigma \in \mathfrak{S}_p.$ 
\index{symmetric tensor}
A tensor $u \in \otimes_H^{p}E$ is said to be \emph{antisymmetric} if $u=\epsilon_{\sigma}S_\sigma u$ for all $\sigma \in \mathfrak{S}_p$
\index{antisymmetric tensor}
where $\epsilon_{\sigma}$ is the signature of $\sigma.$
\end{definition}

\begin{definition}
The space of all antisymmetric tensors in $\otimes_H^{p}E$ will be denoted by $\wedge^p E$. \index{$\wedge^p E$}
\end{definition}

\begin{theorem}\label{a.10}
Let $E$ be a Hilbert space. Then $\bigwedge^{p}E$ is a closed linear subspace of the Hilbert space $\otimes_H^{p}E$ for any $p \geq 2.$
\end{theorem}\vspace{2ex}

\begin{proof}
For $\sigma \in \mathfrak{S}_p$ define the operator 
 $$
f_\sigma \stackrel{\emph{def}}{=} S_\sigma - \epsilon_\sigma  I\;\colon \otimes_H^p E \to \otimes_H^p E,
$$
where $I$ denotes the identity operator on $ \otimes_H^p E $.  Since $S_\sigma$ is a continuous linear operator on $\otimes_H^pE$, $f_\sigma$ is a continuous linear operator. 
 The kernel of the operator $f_\sigma$ is 
 $$
 \begin{array}{cllllllllll}
 \ker f_\sigma &=\{ u \in \otimes_H^p E \colon (S_\sigma - \epsilon_\sigma I)(u)= 0\} \\
 &= \{ u \in \otimes_H^pE \colon \epsilon_\sigma S_\sigma(u)= u\}.
  \end{array}
$$
 
\noindent Since $f_\sigma$ is a continuous linear operator on $\otimes_H^p E,$ $\ker f_\sigma$ is a closed linear subspace of $\otimes_H^p E.$ Thus $\we^p E$ is a closed linear subspace of $\otimes_H^p E,$
 since 
$$
\we^p E = \{u \in \otimes_H^p E \; \colon\epsilon_\sigma S_\sigma(u) = u \mbox{ for all } \sigma \in \mathfrak{S}_p\} = \bigcap\limits_{\sigma \in \mathfrak{S}_p} \ker f_\sigma. 
$$
\end{proof}
Theorem \ref{a.10} implies that the orthogonal projection from $\otimes^p_H E$ onto $\we^p E$
is well defined on $\otimes^p_H E$.

 \begin{definition}\label{a.5}
Let $E$ be a Hilbert space.  For $x_1, \dots, x_p \in E,$ define $x_1 \wedge x_2 \wedge \dots \wedge x_p$ to be the orthogonal
  projection  of the \index{elementary tensor}elementary tensor $x_1 \otimes x_2 \otimes \dots \otimes x_p$ onto $\wedge^p E$, that is
  $$x_1 \wedge x_2 \wedge \dots \wedge x_p= P_{\we^p E}( x_1 \otimes \cdots \otimes x_p).$$
   \end{definition}

One can find a proof of the following  statement in \cite{YCL20-3}.

 \begin{theorem}\label{a.6} Let $E$ be a Hilbert space.  
  For all $u \in \otimes_H^p E,$  
  $$P_{\wedge^p E}(u) = \displaystyle\frac{1}{p!} \sum\limits_{\sigma \in \mathfrak{S}_p}\epsilon_{\sigma}S_\sigma (u).$$
 \end{theorem}

\begin{proposition}\label{we}
Let $E$ be a Hilbert space.
 The inner product in $\wedge^p_H E$ is given by :
 
 $$\langle x_1 \wedge \dots \wedge x_p, y_1 \wedge \dots \wedge y_p \rangle_{\wedge^p_H E} = \det \begin{pmatrix}
                                                                                   
                                                                                   \langle x_1, y_1 \rangle_E &\dots & \langle x_1, y_p \rangle_E\\
                                                                                   					\vdots &\ddots & \vdots \\
                                                                                   \langle x_p, y_1  \rangle_E & \dots &\langle x_p, y_p \rangle_E

                                                                                  \end{pmatrix},$$
for all $x_1, \dots, x_p, y_1, \dots, y_p \in E$.
\end{proposition}

\begin{proof}
 By Theorem \ref{a.6}, we have
 $$
 \begin{array}{cllllllllllllllllll}
 &\langle x_1 \wedge \dots \wedge x_p, y_1 \wedge \dots \wedge y_p \rangle_{\wedge^p E} =\\ [5ex]
 &=\left\langle \displaystyle\frac{1}{p!}\sum\limits_{\sigma \in \mathfrak{S}_p}\epsilon_{\sigma}S_\sigma (x_1 \otimes x_2 \otimes \dots \otimes x_p) , \displaystyle\frac{1}{p!}\sum\limits_{\tau \in \mathfrak{S}_p}\epsilon_{\tau}S_\tau 
 (y_1 \otimes y_2 \otimes \dots \otimes y_p) \right\rangle_{\otimes_H^p E} \vspace{2ex}\\ 
 &=\;\displaystyle \frac{1}{p!} \sum\limits_{\sigma' \in \mathfrak{S}_p}\epsilon_{\sigma'}  \langle  x_1 \otimes x_2 \otimes \dots \otimes x_p,S_{\sigma'} (y_1 \otimes y_2 \otimes \dots \otimes y_p )  \rangle_ {\otimes_H^p E}\vspace{2ex} \\ 
  &= \;\displaystyle  \sum\limits_{\sigma'\in \mathfrak{S}_p}\epsilon_{\sigma'} \prod\limits_{i=1}^p \langle  x_{i}, y_{\sigma'(i)}\rangle_E \vspace{2ex} \\
  
  &= \det \begin{pmatrix}
           \langle x_1 , y_1 \rangle_E & \cdots & \langle x_1 , y_p\rangle_E\\
           \vdots &\ddots & \vdots\\
           \langle x_p , y_1 \rangle_E & \cdots & \langle x_p , y_p \rangle_E
           
          \end{pmatrix}, \;\; \text{by Leibniz' formula}.

  \end{array}$$
\end{proof}
See \cite[Proposition 2.14]{YCL20-3} for slightly more detail.

\begin{corollary}\label{lin-depend}
	Let $E$ be a Hilbert space and let $x_1,\dots,x_p\in E$. Then $x_1\we \cdots \we x_p =0$ if and only if $x_1,\dots,x_p$ are linearly dependent. 
	\end{corollary}

\begin{proof} Note that 
	$x_1\we \cdots \we x_p =0$ if and only if $\| x_1\we \cdots \we x_p\|_{\we^{p}E}^2=0, $ which, by Proposition \ref{we}, holds if and only if
	$$\det [\langle x_i, x_j \rangle]_{i,j=1}^p=0 .$$ 

Thus $x_1\we \cdots \we x_p =0$ if and only if there exist complex numbers $\lambda_1,\dots,\lambda_p,$ which are not all zero, such that
	$$\begin{pmatrix}
	
	\langle x_1, y_1 \rangle_E &\dots & \langle x_1, y_p \rangle_E\\
	\vdots &\ddots & \vdots \\
	\langle x_p, y_1  \rangle_E & \dots &\langle x_p, y_p \rangle_E
		
\end{pmatrix} \begin{pmatrix}
\bar{\lambda}_1 \\ \vdots \\ \bar{\lambda}_p
\end{pmatrix}=0 .$$
 This holds if and only if there exist complex numbers $\lambda_1,\dots,\lambda_p,$ which are not all zero, such that
$$ \langle x_i, \sum_{j=1}^p \lambda_j x_j \rangle_E=0 \quad \text{for}\; i=1,\dots,p .$$
The latter statement is equivalent to the assertion that there exist complex numbers $\lambda_1,\dots,\lambda_p,$ which are not all zero, such that
$$ \langle \sum_{i=1}^p \lambda_i x_i, \sum_{j=1}^p \lambda_j x_j \rangle_E=0 ,$$ which in turn is equivalent to the condition that there exist complex numbers $\lambda_1,\dots,\lambda_p,$ not all zero, such that
$$ \sum_{j=1}^p \lambda_j x_j =0.$$ The latter statement is equivalent to the linear dependence of $x_1,\dots,x_p$ as required.
\end{proof}

\begin{lemma}\label{weon}
	Suppose $\{u_1, \cdots, u_n\}$ is an orthonormal set in $E.$ Then, for  $j=1, \dots,n-1$ and for every $x \in E,$
	$$\| u_1  \we \cdots \we u_j \we x\|_{\we^{j+1}E} =\| x - \displaystyle\sum\limits_{i=1}^j \langle x, u_i  \rangle u_i\|_{E}.   $$
\end{lemma}
See \cite[Lemma 2.15]{YCL20-3}.

\begin{definition}
	Let $(E, \| \cdot\|_E) $ be a Hilbert space. The \emph{$p$-fold Cartesian product of $E$} is defined to be the set 
	$$ \underbrace{E\times \dots \times E}_{p-times} =\{ (x_1,\dots, x_p): x_i \in E \}.$$ Moreover, we define a norm on  $\underbrace{E\times \dots \times E}_{p-times}$ by 
	$$\|(x_1,\dots, x_p) \|=\{\sum_{i=1}^p \|x_i\|_E^2\}^\half. $$
\end{definition}

\begin{definition}
Let $E$ be a Hilbert space. We define the multilinear operator 
$$\Lambda \colon \underbrace{E\times \dots \times E}_{p-times} \to \we^p E$$ 
by 
$$\Lambda(x_1,\dots,x_p)= x_1 \we\dots \we x_p \quad \text{for all}\quad x_1,\dots,x_p \in E. $$

\end{definition}

\begin{proposition}\label{Hadam}{\rm [Hadamard's inequality, \cite{Sing}, p. 477]}
	For any matrix $$A=(a_{ij})\in \Cc^{n\times n},$$
	$$|\det(A)| \leq \prod\limits_{j=1}^{n}\left( \sum\limits_{i=1}^n |a_{ij} |^2  \right)^{1/2} \quad \text{and} \quad |\det(A)| \leq \prod\limits_{i=1}^{n}\left( \sum\limits_{j=1}^n |a_{ij} |^2  \right)^{1/2}. $$
\end{proposition}

\begin{proposition}\label{weopiscontinuous1}
Let $E$ be a Hilbert space. Then the multilinear mapping $$\Lambda \colon \underbrace{E\times \dots \times E}_{p-times} \to \we^p E$$ is bounded and
\begin{equation}\label{Had-C-S-inq}
\| \Lambda(x_1,\dots,x_p )\|_{\we^p E}^2 \leq  \prod\limits_{j=1}^p \|x_j\|_E \left( \sum\limits_{i=1}^p \|x_i\|_E^2 \right)^{1/2}.
\end{equation}
\end{proposition}
See \cite[Proposition 2.19]{YCL20-3} for more detail.

\subsection{Pointwise wedge products and pointwise creation operators}\label{point_w_p}

For the purposes of this paper we need to consider the wedge product of mappings defined on the unit circle or in the unit disk that take values in Hilbert spaces. To this end we introduce a notion of pointwise wedge product and we study its properties.

\begin{definition}\label{a.17}
\index{ pointwise wedge product on $\Tt$} 
Let $E$ be a Hilbert space and let $f,g\colon \Dd \to E$  $\mathrm{(} f,g\colon \Tt \to E \mathrm{)}$ be $E$-valued maps. We define the \emph{pointwise wedge product of $f$ and $g,$}
$$f\telwe g \colon \Dd \to \we^2E \quad \mathrm{(} f\telwe g \colon \Tt \to \we^2E\mathrm{)}$$ 
by $$(f\telwe g) (z) = f(z) \we g(z) \quad \text{for all}\; z \in \Dd \quad \mathrm{(} \text{for almost all}\;  z \in \Tt\mathrm{)}.$$

\end{definition}\vspace{3ex}

 \begin{definition}\label{pointwiseld}
   Let $E$ be a Hilbert space and let $\chi_1,\dots,\chi_n \colon \mathbb{D} \to E$ 
\newline
$\mathrm{(}\chi_1,\dots,\chi_n \colon \mathbb{T} \to E \mathrm{)}$ be $E$-valued maps. We call $\chi_1,\dots \chi_n$ \emph{pointwise linearly dependent} on $\Dd$ \emph{(}or  on $\Tt$\emph{)} if  for all $z\in\Dd$ 
\rm{(}for almost all $z\in\Tt$ respectively\rm{)} the vectors $\chi_1(z), \dots,\chi_n(z)$ are linearly dependent in $E$.
\index{pointwise linearly dependent on $\Tt$}
   \end{definition}

\begin{remark}
 If $x_1, \dots, x_n$ are pointwise linearly dependent on $\Tt$, then 
 $$(x_1\telwe \dots \telwe x_n)(z)=0$$ for almost all $z \in \Tt.$
\end{remark}

\begin{proposition}\label{wejanalytic}
Let $E$ be a Hilbert space and $x_1, x_2, \dots, x_n \colon \Dd \to E$ be  analytic $E$-valued maps on $\Dd.$
Then, $$x_1 \telwe x_2\telwe \dots \telwe x_n\colon \Dd \to \we^n E $$ is also analytic on $\Dd$ and 
$$(x_1 \telwe x_2 \telwe \dots \telwe x_n)'(z) = x_1'(z)\we x_2(z) \we \dots \we x_n(z)  + \dots + x_1(z) \we x_2(z)\we \dots \we x_n'(z)$$ for all $z\in \Dd.$
\end{proposition}
 The proof is straightforward. It follows from Proposition \ref{we}, continuity of $\Lambda$ (Proposition \ref{weopiscontinuous1}) and Hadamard's inequalities \ref{Had-C-S-inq}.

\begin{definition}
Let $E$ be a separable Hilbert space.
 If $x\in L^2(\Tt,E)$ and $y \in L^{\infty}(\Tt, E)$, then $y^*x \in L^2(\Tt ,E)$ is given by $(y^*x)(z)=\langle x(z),y(z) \rangle_E$ almost everywhere on $\mathbb{T}$.
\end{definition}

 \begin{proposition}\label{xweyh2}
 Let $E$ be a separable Hilbert space, let $x\in H^2(\Dd,E)$ and let $y \in H^{\infty}(\Dd, E).$ Then
  $$x \dot{\wedge} y \in H^2 (\Dd, \wedge^2 E).$$
 \end{proposition}
 See
\cite[Proposition 3.8]{YCL20-3}. 

\begin{definition} Let $E$ be a Hilbert space.
We say that a family  $\{ f_\lambda \}_{ \lambda \in \Lambda}$
of maps from $\Tt$ to $E$ is {\em pointwise orthonormal} on $\Tt$, if for all $z$ in a set of full measure in $\Tt$, 
the family of vectors $\{ f_\lambda(z) \}_{ \lambda \in \Lambda}$ is
orthonormal in $E$.
\end{definition}

 \begin{proposition}\label{wel2conv}  Let $E$ be a separable Hilbert space, and let 
$\; \xi_0,\xi_1, \cdots, \xi_j \in L^\infty(\Tt,E) $ be a pointwise orthonormal set on $\Tt$, and let $x\in L^2(\Tt,E)$.  Then
  $$ \xi_0 \dot{\we} \xi_1 \dot{\we} \cdots \dot{\we} \xi_j \dot{\we} x \in L^2 (\Tt, \wedge^{j+2} E),$$
and 
\[
\| \xi_0 \dot{\we} \xi_1 \dot{\we} \cdots \dot{\we} \xi_j \dot{\we} x \|_{L^2 (\Tt, \wedge^{j+2} E)} \le
\| x \|_{L^2(\Tt,E)}.
\]
 \end{proposition}
It follows from Lemma \ref{weon}. See \cite[Proposition 3.11]{YCL20-3} for a proof of this proposition.

\begin{definition}\label{pwcre} Let $E$ be a separable Hilbert space. 
	Let $\xi \in H^\infty(\Dd, E)$. We define the \emph{pointwise creation operator} $$C_\xi :H^2(\Dd,E) \to H^2(\Dd, \we^2E)$$ by 
	$$ C_\xi f = \xi \telwe f\; \text{for} \; f \in H^2(\Dd,E).$$  
\end{definition}\index{pointwise creation operator} \index{$C_\xi$}


\begin{remark}\label{genfatouwe}  Let $E$ be a separable Hilbert space. Let $\xi \in H^\infty(\Dd, E)$ and let $f \in H^2(\Dd, E)$.
 By the generalized Fatou's Theorem \cite[Chapter V]{NagyFoias}, the radial limits 
	$$\lim_{r\to 1}\xi(r\eiu)\underset{\|\cdot\|_E}{=} \tilde{\xi}(\eiu), \quad \lim_{r\to 1} f(r\eiu)\underset{\|\cdot\|_E}{=} \tilde{f}(\eiu) \;\;
(0<r<1)
$$ exist almost everywhere on $\Tt$ and define functions $\tilde{\xi} \in L^\infty(\Tt,E)$ and $\tilde{f}\in L^2(\Tt,E)$ respectively, which satisfy the relations 
	$$ \lim_{r \to 1} \| \xi(r\eiu) - \tilde{\xi}(\eiu) \|_E=0,\quad \lim_{r \to 1} \| f(r\eiu) - \tilde{f}(\eiu) \|_{E}=0 \;\;
(0<r<1) $$ 
for almost all $ \eiu  \in \Tt$.
\end{remark}\index{$\tilde{f}$}

\begin{lemma}\label{xitelweh21} Let $E$ be a separable Hilbert space.
 Let $\xi\in H^\infty(\Dd, E)$ and let $f \in H^2(\Dd, E).$ Then the radial limits $\lim_{r \to 1} (\xi(r\eiu) \we f(r\eiu))$ exist for almost all $\eiu \in \Tt$ and define functions in $L^2(\Tt,\we^2 E).$
\end{lemma}

One can find a proof of this statement in \cite[Lemma 4.3]{YCL20-3}.

\begin{remark}\label{H2subsetL2}  Let $E$ be a separable Hilbert space.
By \cite[Chapter 5, Section 1]{NagyFoias}, for any separable Hilbert space $E$, the map 
$f \mapsto \tilde{f}$ is an isometric embedding of $H^2(\Dd,E)$ in $L^2(\Tt,E)$, where 
$ \tilde{f}(\eiu)= \lim_{r \to 1} f(r\eiu).$ Since $H^2(\Dd,E)$ is complete and the embedding is isometric, the image of the embedding is complete, and therefore is closed in $L^2(\Tt, E).$ Therefore, the space $H^2(\Dd,E)$ is identified isometrically with a closed linear subspace of $L^2(\Tt, E).$ 
In future we shall use the same notation for $f$ and $\tilde{f}.$
\end{remark}

\begin{definition}\label{pls} Let $E$ be a separable Hilbert space.
Let $F$ be a subspace of $L^2(\Tt, E)$ and let $X$ be a subset of $L^2(\Tt,E).$ 

We define the \emph{pointwise orthogonal complement} of $X$ in $F$ to be the set
$$\Poc(X,F) = \{ f \in F: f(z)\perp \{x(z):x \in X\}\; \text{for almost all}\;z \in \Tt\}.$$

\index{pointwise orthogonal complement} \index{$\Poc(X,F)$}

\end{definition}

\begin{proposition}\label{vclosed} Let $E$ be a separable Hilbert space.
Let $\eta \in L^2(\Dd,E)$. Then
\begin{enumerate}
	\item[(i)] The space $V= \{ f \in H^2 (\Dd,E):  \langle f(z) , \eta(z) \rangle_{E} =0 \; \text{for almost all}\; z \in \Tt \}$ is a closed subspace of $H^2(\Dd, E).$ 
\item[(ii)] The space  $V=\{ f \in L^2(\Tt,E): \langle f(z), \eta(z) \rangle_{E} =0 \;\text{for almost all}\; z\in \Tt  \} $ is a closed subspace of $L^2(\Tt,E).$ 
\end{enumerate}
\end{proposition}
One can find a proof of this statement in \cite[Proposition 4.6]{YCL20-3}.

%% file: algorithm_state.tex
\section{Superoptimal analytic approximation}\label{statet_algorithm}

	In this section we present our main result, which is an algorithm for the superoptimal analytic approximation of a matrix-valued function on the circle. In Subsection \ref{known} we recall certain known results and Peller and Young's algorithm (Theorem \ref{superopconstruct}). In Subsection \ref{Alg_statement}
we present an alternative algorithm for the superoptimal approximant, based on exterior powers of Hilbert spaces. The proof of the validity of the new algorithm relies on the cited work given in Subsection \ref{known}.

\subsection{Known results}\label{known}

\begin{theorem}[Hartman's Theorem, \cite{Peller}, p. 74]\label{2.04}
Let $E,F$ be separable Hilbert spaces and let $\Phi \in L^\infty(\Tt,  \mathcal{L}(E,F)).$ The following statements are equivalent: 
\begin{enumerate}
	\item[i)] The Hankel operator $H_\Phi$ is compact on $H^2(\Dd,E)$;
	\item[ii)] $\Phi \in H^\infty(\Dd, \mathcal{L}(E,F))+ C (\Tt,\mathcal{K}(E,F))$;
	\item[iii)] there exists a function $\Psi \in C (\Tt,\mathcal{K}(E,F))$ such that $\hat{\Phi}(n) = \hat{\Psi}(n)$ for $n<0.$
 \end{enumerate}
\end{theorem}
\index{Hartman's Theorem} \index{compact Hankel operator}

\begin{theorem}[\cite{page}]\label{nehtmatr}
	For any matrix-valued $\phi \in L^{\infty}(\mathbb{T}, \Cc^{m\times n} ),$	
	$$\inf\limits_{Q\in H^{\infty}(\mathbb{D}, \Cc^{m\times n})} \|\phi-Q\|_\infty=\|H_\phi \|$$and the infimum is attained.	
\end{theorem}\vspace{4ex}

\begin{definition}[\cite{superop}, p. 306]

The class of \emph{quasi-continuous} functions is defined by
$$QC=(H^\infty(\Dd, \CCmn) + C(\Tt, \CCmn)) \cap \overline{(H^\infty(\Dd, \CCmn) + C(\Tt, \CCmn))}.$$ In other words this class consists of functions on the circle which belong to $H^\infty + C$ and have the property that their complex conjugates belong to $H^\infty + C$ as well. 
 
 \index{quasi-continuous function}
\end{definition}
We shall also need the class of functions of vanishing mean oscillation, as described, for example, in \cite[Appendix 2, Section 5]{Peller}.
\begin{definition}
  For any function $f\in L^1(\T)$ and any arc $I$ in $\T$ let
\[
f_I \df \frac{1}{m(I)} \int_I f dm
\]
where $m$ is Lebesgue measure on $\T$.  Thus, $f_I$ is the mean of $f$ over $I$.  The function $f$ is said to have {\em vanishing mean oscillation} if
\[
\lim_{m(I) \to 0} \frac{1}{m(I)} \int_I |f-f_I| \, dm =0.
\]
The space of functions of vanishing mean oscillation on $\T$ is denoted by $\vmo$.
\end{definition}
\index{$\vmo$}
$\vmo$ is also related to the compactness of Hankel operators.  The following is \cite[Theorem 5.8]{Peller}.
\begin{theorem}
Let $\ph\in L^2$.  Then $H_\ph$ is compact if and only if $P_-\ph\in \vmo$.
\end{theorem}
It is therefore not surprising that the spaces $QC$ and $\vmo$ are closely related.  In fact
\begin{theorem}\label{vmoqc}\cite[Page 729]{Peller}
\[
QC=\vmo \cap L^\infty.
\]
\end{theorem}
\noindent It follows from another characterization of $\vmo$, to wit
\[
\vmo \ =\ \{f+\tilde g: f,g \in C(\T)\},
\]
where $\tilde g$ denotes the harmonic conjugate of $g$ \cite[Theorem A2.8]{Peller}.

\begin{remark}
For $G \in H^\infty(\Dd, \CCmn) + C(\Tt, \CCmn)$ we will say that a function $Q \in H^\infty(\Dd, \CCmn)$ which minimizes the norm 
 $\|G-Q\|_{L^\infty}$ is a function {\em at minimal distance from $G.$} By Nehari's Theorem, all such functions $Q$ satisfy  $ \|G-Q\|_{L^\infty}= \|H_G\|.$
 \index{function at minimal distance from $G$}
\end{remark}

\begin{definition}[\cite{NagyFoias}, p. 190]
	For a separable Hilbert space $E,$ a function $\xi \in H^\infty(\Dd, E)$ will be called \emph{inner} if for almost every $z \in \Tt,$
	$$ \|\xi(z)\|_E=1.$$
\end{definition}

\begin{theorem}[\cite{superop}, Theorem 1.1]\label{coouter}
	Let $\varphi$ be an $n\times 1 $ inner matrix function. There exists a co-outer function $\varphi_c \in H^\infty(\Dd,\Cc^{n\times (n-1)})$ such that 
	$$ \Phi = \begin{pmatrix} \varphi & \bar{\varphi}_c\end{pmatrix}$$
is unitary-valued on $\Tt$ and all minors of $\Phi$ on the first column are in $H^\infty.$

\end{theorem}

For a function $G: \T \to \Cc^{m\times n}$ and a space $X$ of scalar functions on $\T$, we write $G \in X$ to mean that each entry of $G$ belongs to $X$.

Next we describe some properties that a space $X$ of equivalence classes of scalar functions on the circle may possess \cite[Page 330]{superop}.
Define the non-linear operator $\cala = \cala^{(m,n)}$ on the space of $m\times n$ functions 
$G \in H^\infty(\Dd, \CCmn) + C(\Tt, \CCmn)$ by saying that $\cala^{(m,n)} G$ is the unique superoptimal approximation in $H^\infty(\Dd, \CCmn)$ to $G$.

We say that $X$ is {\em hereditary for} $\cala$ if, for every scalar function $g\in X$, the best analytic approximation $\cala g$ of $g$ belongs to $X$.

Consider the following conditions on $X$ from 
\cite{pellkhr}:
\begin{enumerate}
\item[($\al$1)] $X$ contains all polynomial functions and $X\subset \vmo$;
\item[($\al$2)]  $X$ is hereditary for $\cala$;
\item[($\al$3)] if $f\in X$ then $\bar z \bar f \in X$ and $P_+f \in X$;
\item[($\al$4)]  if $f,g \in X \cap L^\infty$ then $fg \in X\cap L^\infty$;
\item[($\al$5)]  if $f\in X \cap H^2$ and $h\in H^\infty$ then $T_{\bar h}f\in X \cap H^2$.
\end{enumerate}
The relevance of these properties is contained in the following statement, which is \cite[Lemma 5.3]{superop}.
Recall that a function $f\in L^\infty$ is said to be {\em badly approximable} if the best analytic approximant to $f$ is the zero function.  In view of Nehari's Theorem,
$f$ is badly approximable if and only if $\|f\|_\infty = \|H_f\|$. \index{badly approximable}
\begin{lemma}\label{lem5.3}
Let $X$ satisfy $(\al 1)$ to $(\al 5)$ and let $\ph\in X$ be an $n\times 1$ inner function.  Let $\ph_c$ be an $n \times (n-1)$ function in $H^\infty$ such that 
$[\ph \ \bar\ph_c]$ is unitary-valued a.e. on $\T$ and has all its minors on the first column in $H^\infty$.  Then each entry of $\ph_c$ belongs to $X$.
\end{lemma}

Below we shall use a modified version of \cite[Theorem 0.2]{superop}.


\begin{theorem}\label{1.7}\emph{\cite[Theorem 0.2]{superop}}
	Let $\phi \in  L^{\infty}(\Tt, \Cc^{m\times n})$ be such that $H_\phi$ has a Schmidt pair $(v,w)$ corresponding to the singular value $t=\|H_\phi\|$.
	Let $Q$ be a function in $H^{\infty}(\Dd, \Cc^{m \times n})$ at minimal $L^\infty$-distance from $\phi$. Then
$$(\phi - Q) v = tw$$  and 
\[
(\phi-Q)^*w= tv.
\]
Moreover
\begin{equation}\label{normvz}
	\|w(z)\|_{\mathbb{Cc}^m} =  \|v(z)\|_{\mathbb{Cc}^n} ~ \text{almost everywhere on}~ \mathbb{T}
\end{equation}	
and $$\| \phi (z) - Q(z)\| = t ~\text{almost everywhere on}~ \mathbb{T}.$$
\end{theorem}
\begin{proof}
	By Nehari's Theorem, $\|\phi-Q\|_{L^\infty}= t$ and, by hypothesis,
\[
H_\phi v=tw, \qquad H_\phi^* w=tv.
\]
If $t=0$ then $\phi\in  H^\infty(\Dd,\Cc^{m\times n})$, so that $\phi=Q$ and the statement of the theorem is trivially true.  We may therefore assume $t>0$.
Thus $H_\phi^*H_\phi v = t^2 v$, and so $v$ is a maximising vector for $H_\phi$.  We can assume that $v$ is a unit vector in $H^2(\Dd, \Cc^n)$, and then $w$ is a unit vector in  $H^2(\Dd, \Cc^m)^\perp$ and is a maximising vector for $H_\phi^*$.   We have
\[
t = \|H_\phi v \| = \|H_{\phi - Q} v \|= \| P_{-} (\phi - Q) v\| \leq \| (\phi - Q) v\| \leq \|\phi - Q\|_{L^\infty} = t.
\]
The inequalities must hold with equality throughout, and therefore	
	
	$\| P_{-} (\phi - Q) v\| = \| (\phi - Q) v\|$, which implies that  $(\phi - Q) v \perp H^2$ and so
	
	$$H_\phi v= P_{-} (\phi -Q)v = (\phi -Q)v.$$
	
	Furthermore $\| (\phi - Q) v\|$ = $\| (\phi - Q) \|_{L^\infty} \|v\|$ and since $v(z)$ is therefore a maximizing vector for $\phi(z)-Q(z)$ for almost all $z$, 
	we have $\| \phi(z)-Q(z) \| = \|H_\phi \|.$

Likewise,
\begin{align*}
 t=\|H_\phi^*\|= \|H_{\phi-Q}^*\|&= \|H_{\phi-Q}^* w\|= \|P_+(\phi-Q)^*w\|_{L^2} \leq \|(\phi-Q)^*w\|_{L^2}\\ 
& \leq \|(\phi-Q)^*\|_{L^\infty}\|w\|_{L^2}
 = \|(\phi-Q)^*\|_{L^\infty} = t.
\end{align*}
Again, the inequalities hold with equality throughout, and in particular
\[
\|P_+(\phi-Q)^*w\|_{L^2} = \|(\phi-Q)^*w\|_{L^2},
\]
so that $(\phi-Q)^*w \in H^2$ and
\[
(\phi-Q)^*w= H_\phi^*w = tv.
\]
\end{proof}

\begin{lemma}\cite[p. 315-316]{superop}  \label{2.2}
Let $m,n>1,$ let $G \in H^\infty ( \Dd, \CCmn) +C (\Tt, \CCmn)$ and 
\noindent $t_0 = \| H_G\| \neq 0.$ Suppose that $v$ is a maximizing vector of $H_G$ and let  
\begin{equation}\label{eqhgv} H_{G} v = t_{0} w. 
\end{equation}

\noindent Then $v, \bar{z}\bar{w} \in H^2(\Dd, \Cc^n)$ have the factorizations 
\begin{equation} \label{eq212} v= v_0 h , \; \; \bar{z}\bar{w} = \phi w_0 h \end{equation}

\noindent for some scalar outer function $h,$ some scalar inner $\phi,$ and column-matrix inner functions $v_0 , w_0.$ Moreover there exist unitary-valued functions $V, W$ of types $n \times n , \; m \times m$ respectively, of the form 
\begin{equation}\label{eq222} V= \begin{pmatrix}
v_0 & \bar{\alpha}
\end{pmatrix}, \; W^T = \begin{pmatrix}
w_0 & \bar{\beta}
\end{pmatrix},\end{equation}
              
\noindent where $\alpha, \beta $ are inner, co-outer functions, quasi-continuous functions of types $n\times (n-1), \; m\times (m-1)$ respectively, and all minors on the first columns of $V, W^T$ are in $H^\infty.$  

\noindent Furthermore every $Q \in H^\infty (\Dd, \CCmn)$ which is at minimal distance from $G$ satisfies 
\begin{equation}\label{eq2223}W(G-Q)V = \begin{pmatrix}
t_0 u_0 &~& 0\\
0 &~& F
\end{pmatrix}\end{equation} for some $F \in H^\infty (\Dd, \Cc^{(m-1)\times(n-1)}) + C(\Tt, \Cc^{(m-1)\times(n-1)})$ and some quasi-continuous function $u_0$ given by
\begin{equation}\label{eq223}u_0 = \frac{\bar{z}\bar{\phi}\bar{h}}{h}\end{equation} with $|u_0 (z)|=1$ almost everywhere on $\Tt.$

\end{lemma}
In the statement of the lemma, in saying that an $m \times n$ matrix-valued function $\al$ is {\em co-outer}
\index{co-outer}
we mean that each column of $\al$ is in $H^{\infty}_m$ and $\al^TH^2_m$ is dense in $H^2_n$. (In \cite[Page 190]{NagyFoias}, such a function $\al$ is said to be {\em *-outer}).

\begin{proof}
First we construct $V$ and $W$ with the properties \eqref{eqhgv} to \eqref{eq2223}.
 By equation \eqref{normvz}, $\|v(z) \| = \|w(z)\|$ almost everywhere, and so the column-vector functions $v, \bar{z}\bar{w}$ in $H^2$ have the same (scalar) outer factor $h.$
 This property yields the inner-outer factorizations (\ref{eq212}) for some column inner functions $v_0, w_0$.  By Theorem \ref{coouter}, there exists an inner co-outer function $\al$ of type $n\times (n-1)$
such that $V\df [v_0 \  \bar\al]$ is unitary-valued almost everywhere on $\T$ and all minors on the first column of $V$ are in $H^\infty$.  Likewise there exists an inner co-outer function $\beta$ of type $m \times (m-1)$ such that $W\df [w_0 \ \bar\beta]^T$ is unitary-valued almost everywhere on $\T$ and all minors on the first column of $W^T$ are in $H^\infty$.
 
Next we show that $u_0$ given by equation \eqref{eq223} is quasicontinuous.
Let $Q \in H^\infty(\Dd, \CCmn)$ be at minimal distance from $G.$ Then 
 $$\|G-Q\|_\infty = \|H_G\|= t_0.$$
 
\noindent By Theorem \ref{1.7}, 
 $$(G-Q)v = t_0 w$$ and by the factorizations (\ref{eq212}) we have
 $$(G-Q) v_0 h = t_0 \bar{z} \bar{\phi} \bar{h} \bar{w}_0$$ and by equations (\ref{eq222}) and (\ref{eq223})
 
 $$(G-Q) V \begin{pmatrix} 1 & 0 & \cdots & 0 \end{pmatrix}^T = W^* \begin{pmatrix} t_0 u_0 & 0 & \cdots & 0 \end{pmatrix}^T.$$
 
\noindent Thus $$W(G-Q)V = \begin{pmatrix} t_0 u_0 & f\\ 0 & F \end{pmatrix}$$ for some $f \in L^\infty(\Tt, \Cc^{1\times (n-1)}), \; F \in L^\infty(\Tt, \Cc^{(m-1)\times (n-1)}).$
 
\noindent Because $t_0 = \|H_G\|,$ it follows that $|u_0|=1$ almost everywhere, and from Nehari's Theorem 
 $$\|W(G-Q)V\|_\infty = \|G-Q\| = \|H_G\| = t_0,$$
and we have $f=0.$ So, $W(G-Q)V$ has the form (\ref{eq2223}).
Now, $\|H_{u_0}\| \leq \|u_0\|_\infty =1$ and $\|H_{u_0}h\|=\|\bar{z}\bar{\phi} \bar{h}\| = \|h\|.$ Hence
$$\|H_{u_0}\|=1 = \|u_0\|_\infty,$$
which implies that $u_0$ is badly approximable.
The $(1,1)$ entries of equation (\ref{eq2223}) are 
$$w_0^T (G-Q)v_0 = t_0 u_0.$$
Since $v_0 \in H^\infty(\Dd, \Cc^n ) , \; w_0 \in H^\infty(\Dd, \Cc^m)$  and $H^\infty(\Dd,\C) + C(\Tt,\C)$ is an algebra, \\ 
$u_0 \in H^\infty + C.$ 
By a result in \cite[Section 3.1]{pellkhr}, if $u_0 \in H^\infty + C$ and $u_0$ is badly approximable then $\bar{u}_0 \in H^\infty +C.$ Thus $u_0$ is quasi-continuous. 

Now we show that $v_0, w_0 \in QC$.  It follows from Nehari's Theorem that
\[
(G-Q)^* w = t_0 v, 
\]
much as in the proof of \cite[Theorem 0.2]{superop}.  Indeed, since $H_G^* w = t_0 v$ and $H_G^* = P_+M_{(G-Q)^*} | {H^2}^\perp$,
we have (assuming, as we may, that $v$ and $w$ are unit vectors),
\begin{align*}
 t_0 = \|H_G^*w\| &= \|P_+(G-Q)^*w\| \\
	&\leq \|(G-Q)^* w\| \leq \|G-Q\|_\infty \|w\| =t_0.
\end{align*}
It follows that the inequalities hold with equality, and so
\[
\|P_+(G-Q)^*w\| =  \|(G-Q)^* w\|,
\]
whence
\[
P_+(G-Q)^*w =  (G-Q)^* w,
\]
and so 
\be\label{usethis}
 (G-Q)^* w = H_G^*w = t_0 v,
\ee
as claimed.

Taking complex conjugates in the last equation we have
\[
(G-Q)^T \bar w = t_0 \bar v.
\]
Thus, by equation \eqref{eq212},
\[
(G-Q)^Tz\ph w_0 h = t_0 \bar h \bar{v}_0
\]
for some outer function $h$ and scalar inner $\ph$.  Therefore
\[
\bar{v}_0 = \frac{(G-Q)^Tz\ph w_0 h}{t_0 \bar h}.
\]
Recall that $u_0 = \bar z \bar \ph \bar h/h$, and so
\[
\bar{v}_0 = \frac{1}{t_0}(G-Q)^T\bar{u}_0w_0.
\]
Since $u_0\in QC, \ G-Q\in H^\infty+C$ and $w_0 \in H^\infty$, it follows that $\bar{v}_0 \in H^\infty +C$.  Since also $v_0\in H^\infty$, we have $v_0\in QC$.

To complete the proof of Lemma \ref{2.2}, all that remains is to show that $\al,\beta $ are quasicontinuous and $F\in H^\infty + C$.  This will follow from Lemma \ref{lem5.3} above.

The space $\vmo$ satisfies conditions $(\al1)$ to $(\al5)$,  as stated on \cite[Page 335]{superop}, and we have $v_0 \in QC \subset \vmo$.  Hence we may apply Lemma \ref{lem5.3} with $\ph=v_0$ to deduce that $\al \in \vmo$.  Since also $\al\in L^\infty$, it follows from Theorem \ref{vmoqc} that $\al\in QC$.  Likewise, $\beta\in QC$.

 To show that $F\in H^\infty + C$, for $1<i \leq m, \; 1<j\leq n$ consider the $2\times 2$ minor of equation (\ref{eq2223})  with indices $1i , 1j :$ 
\begin{equation}\label{eqsum}
\sum\limits_{r<s,\; k<l} W_{1i , rs} (G-Q)_{rs, kl} V_{kl , 1j} = t_0 u_0 F_{i-1 , j-1}.
\end{equation}

\noindent By the analytic minors property of $W,V,$
$$V_{kl,1j}, W_{1i,rs} \in H^\infty.$$

Since $G-Q \in H^\infty(\Dd, \CCmn) + C(\Tt, \CCmn) ,$ the left hand side of equation (\ref{eqsum}) is in $H^\infty (\Dd, \Cc^{2\times 2}) + C(\Tt,\Cc^{2\times 2})$ and hence $u_0F \in H^\infty(\Dd, \C^{(m-1)\times(n-1)} + C(\Tt, \C^{(m-1)\times(n-1)}).$
Thus 
$$
F = \bar{u}_0 (u_0 F) \in H^\infty(\Dd, \C^{(m-1)\times(n-1)}) + C(\Tt, \C^{(m-1)\times(n-1)}).
$$

\end{proof}


\begin{definition}\label{thematicdef}
	We say that a unitary-matrix-valued function $V$ is a  \emph{thematic completion} of a column-matrix inner function $v_0\in H^\infty(\Dd,\Cc^n)$  if $V = \begin{pmatrix}
	v_0 & \bar{\alpha}
	\end{pmatrix}$, for some co-outer function $\alpha \in H^\infty(\Dd,\Cc^{n\times(n-1)})$ such that $V(z)$ is a unitary matrix for almost all $z \in \Tt$ and such that all minors on the first column of $V$ are analytic. 
\end{definition}\index{thematic completion}
\begin{remark} \rm 
	By Theorem $1.1$ of \cite{superop}, every column-matrix inner function has a thematic completion. Thematic completions are not unique, for if $V = \begin{pmatrix}
	v_0 & \bar{\alpha}
	\end{pmatrix}$ is a thematic completion of $v_0$, then so is $\begin{pmatrix}
		v_0 & \bar{\alpha}U\end{pmatrix}$ for any constant $(n-1)$-square unitary matrix $U.$ However, by Corollary $1.6$ of \cite{superop}, the thematic completion of $v_0$ \emph{is} unique up to multiplication on the right by a constant unitary matrix of the form $\mathrm{diag}\{1,U\}$ for some constant $(n-1)$- square matrix $U,$ and so it is permissible to speak of  ``{\em the} thematic completion of $v_0$".

Furthermore, by Theorem $1.2$ of \cite{superop}, thematic completions have constant determinants almost everywhere on $\Tt,$ and hence $\alpha, \beta$ are inner matrix functions. Observe that, as we showed above, if the column $v_0$ belongs to $\mathrm{VMO},$ then the thematic completion of $v_0$ is quasi-continuous. Similarly, if the column $w_0$ belongs to $\mathrm{VMO}$, then the thematic completion of $w_0$ is quasi-continuous. Thus $\alpha,\beta$ are inner, co-outer, quasi-continuous functions of types $n\times(n-1)$ and $m\times (m-1)$ respectively.  	
\end{remark}

\begin{lemma}[\cite{superop}, p. 316]\label{f+hinfty}
 Let $m,n>1,$ let $G \in H^\infty(\Dd, \CCmn)+C(\Tt, \CCmn),$ let $\|H_G\|=t_0$ and let $Q_1 \in H^\infty(\Dd, \CCmn)$ be at minimal distance from $G,$ so that in the notation of Lemma \ref{2.2}, \begin{equation}\label{g-q1vf}
 W(G-Q_1)V= \begin{pmatrix}
 t_0 u_0 & 0 \\ 0 & F
 \end{pmatrix}\end{equation} for some $F \in H^\infty(\Dd, \Cc^{m-1\times n-1})+ C(\Tt, \Cc^{m-1\times n-1}).$
 Let $$\mathcal{E} = \{ G-Q : Q \in H^\infty(\Dd, \CCmn), \|G-Q\|_\infty=t_0\}.$$ Then 
 \begin{equation}\label{wev}W \mathcal{E} V = \begin{pmatrix}
 t_0 u_0 & 0 \\ 0 & F +H^\infty(\Dd, \Cc^{m-1\times n-1})
 \end{pmatrix}\cap B(t_0),\end{equation} where $B(t_0)$ is the closed  ball of radius $t_0$ in $L^\infty(\Tt, \CCmn).$
\end{lemma}
\begin{proof}
Let $E_1 = G-Q_1 \in \mathcal{E}$ and for $Q \in H^\infty(\Dd,\CCmn)$ consider 
$$E = E_1 -Q.$$
By Lemma \ref{2.2}, there exists a function $g \in L^\infty(\Dd,\Cc^{m-1\times n-1})$ such that 
$$ WEV = \begin{pmatrix}
t_0 u_0 & 0 \\ 0 & g 
\end{pmatrix}.	$$ The latter equation combined with equation (\ref{g-q1vf}), yields 
$$\begin{array}{clllll} WQV&= W(G-Q_1)V-WEV\vspace{2ex}\\&= \begin{pmatrix}
0 & 0 \\ 0 & F-g
\end{pmatrix}\in \begin{pmatrix}
0 & 0 \\ 
0 & L^\infty(\Tt, \Cc^{m-1\times n-1})
\end{pmatrix}\cap WH^\infty(\Dd,\CCmn)V.\end{array}$$
By (\cite{superop}, Lemma 1.5),  $WQV \in H^\infty(\Dd,\CCmn),$ say $F-g=q \in H^\infty(\Dd,\CCmn).$ Then 
$$WEV= \begin{pmatrix}
t_0 u_0 & 0 \\ 0 & g
\end{pmatrix}= \begin{pmatrix}
t_0 u_0 & 0 \\ 0 & F-q
\end{pmatrix},  $$which proves the inclusion ($\subseteq$) in equation \eqref{wev}. 

\noindent Conversely, suppose $q \in H^\infty(\Dd,\Cc^{m-1\times n-1} )$ and $$\|F-q\|_\infty \leq t_0.$$
By (\cite{superop}, Lemma 1.5), there exists a function $Q\in H^\infty(\Dd,\CCmn)$ such that 
$$ WQV = \begin{pmatrix}
0 & 0 \\ 0 & q
\end{pmatrix}$$ and thus 
$$ W(E_1 - Q)V = \begin{pmatrix}
t_0 u_0 & 0 \\ 0 & F-q
\end{pmatrix}. $$Then 
$$ E_1 - Q = G-(Q_1-Q) \in \mathcal{E},$$ and so 
$$\begin{pmatrix}
t_0 u_0 & 0 \\ 0 & F-q
\end{pmatrix}\in W\mathcal{E}V. $$ Hence equality holds in equation (\ref{wev}).
\end{proof}

\begin{lemma}[\cite{Constr}, p. 16]\label{3.1constr}
	Let $G\in H^\infty(\Dd,\CCmn)+C(\Tt,\CCmn)$ and let $(x_0, y_0)$ be a Schmidt pair for the Hankel operator $H_G$ corresponding to the singular value $t_0 = \|H_G\|.$	Let $x_0 = \xi_0 h_0$ be the inner-outer factorization of $x_0,$ where $\xi_0\in H^\infty(\Dd,\Cc^n)$ is the inner and $h_0\in H^2(\Dd,\Cc)$ is the scalar outer factor of $x_0 \in H^2(\Dd,\Cc^n),$ and let 
	$$ V_0 = \begin{pmatrix} \xi_0 & \overline{\alpha_0} \end{pmatrix}$$ be a unitary-valued function on $\Tt$, where $\alpha_0 \in H^\infty(\Dd,\Cc^{n \times (n-1)})$ is co-outer. Then $$V_0 \begin{pmatrix}
	0 & H^2(\Dd,\Cc^{n-1})
	\end{pmatrix}^T$$ is the orthogonal projection of $H^2(\Dd,\Cc^n)$ onto the pointwise orthogonal complement of $x_0$ in $L^2(\Tt,\Cc^n).$
	Similarly $$V_0^* \begin{pmatrix}
	0 & H^2(\Dd,\Cc^{n-1})^\perp 
	\end{pmatrix}$$ is the orthogonal projection of $H^2(\Dd,\Cc^n)^\perp$ onto the pointwise orthogonal complement of $x_0$ in $L^2(\Tt,\Cc^n).$

\end{lemma}

\begin{lemma}[\cite{Constr}, p. 16]\label{3.2constr}
	Let $G,x_0, y_0$ be defined as in Lemma \ref{3.1constr} and let $\mathcal{K}, \mathcal{L}$ be the projections of $H^2(\Dd,\Cc^n),H^2(\Dd,\Cc^m)^\perp$ onto the pointwise orthogonal complements of $x_0, y_0$ in $L^2(\Tt,\Cc^n),L^2(\Tt,\Cc^m)$ respectively. Let $Q_0 \in H^\infty(\Dd,\CCmn)$ be at minimal distance from $G,$ let $F$ be the $(2,2)$ block of $W_0(G-Q_0)V_0,$ as in Lemma \ref{2.2}, that is,   \begin{equation}\label{V0W0}V_0=\begin{pmatrix} \xi_0 & \overline{\alpha_0} \end{pmatrix},\quad W_0 = \begin{pmatrix}
	\eta_0 & \overline{\beta_0} 
	\end{pmatrix}^T\end{equation} are unitary-valued functions on $\Tt,$ $\alpha_0, \beta_0$ are co-outer functions of size $n\times (n-1) , m\times (m-1)$ respectively and all minors on the first columns of $V_0, W_0^T$ are in $H^\infty.$ Let $Q \in H^\infty(\Dd,\CCmn)$ satisfy 
	\begin{equation}\label{g-q3.1}(G-Q)x_0 = \|H_G\| y_0 , \quad y_0^* (G-Q) = \|H_G\| x_0^*.$$ Then $H_F$ is a unitary multiple of the operator
	$$ \Gamma := P_\mathcal{L} M_{G-Q}|\mathcal{K},\end{equation} where 
	$M_{G-Q}:L^2(\Tt,\Cc^n) \to L^2(\Tt,\Cc^m)$ is the operator of multiplication by $G-Q.$ More explicitly, if $U_1 : H^2(\Dd,\Cc^{n-1}) \to \mathcal{K},$ $U_2 : H^2(\Dd,\Cc^{m-1})^\perp \to \mathcal{L}$ are defined by 
	\begin{equation}
	\label{u1u2}
	U_1 \chi = V_0 \begin{pmatrix}
	0 \\ \chi 
	\end{pmatrix}, \quad U_2 \psi = W_0^* \begin{pmatrix}
	0 \\ \psi
	\end{pmatrix}\quad \text{for all}\quad \chi \in H^2(\Dd,\Cc^{n-1}), \psi \in H^2(\Dd,\Cc^{m-1}),\end{equation} then $U_1 , U_2$ are unitaries and 
	$$H_F = U_2^* \Gamma U_1.$$

\end{lemma}

\begin{lemma}[\cite{superop}, p. 337]\label{L6.2}
	Let $\alpha\in \mathrm{QC}$ of type $m\times n,$ where $m\geq n,$ be inner and co-outer. There exists $A\in H^\infty(\Dd,\Cc^{n\times m})$ such that $A\alpha=I_n.$ Here $I_n$ denotes the $n\times n$ identity matrix. 
	
\end{lemma}

Theorem \ref{superopconstruct} gives the algorithm for the superoptimal analytic approximant constructed in \cite{Constr}.

\begin{theorem}[\cite{Constr}, p. 17]\label{superopconstruct}
	Let $G \in H^\infty(\Dd,\CCmn)+ C(\Tt,\CCmn).$ The superoptimal approximant $\mathcal{A}G$ to $G$ is given by the following formula. 
	
	\noindent If $H_G=0,$ then $\mathcal{A}G=G.$ Otherwise define spaces $K_j\subset L^2(\Tt,\Cc^n), N_j\subset L^2(\Tt,\Cc^m),$ vectors $\chi_j \in K_j,$ $\psi_j \in N_j,$ $H^\infty$ functions $Q_j,$ operators $\Gamma_j$ and positive $\lambda_j$ as follows.
	
	\noindent Let $$K_0  = H^2(\Dd,\Cc^n),\quad N_0= H^2(\Dd,\Cc^m)^\perp,\quad Q_0=0.$$ Let    $$ \Gamma_j = P_{N_j} M_{G-Q_j} |K_j : K_j \to N_j,\quad \lambda_j=\|\Gamma_j\|,$$where $P_{N_j}$ is the orthogonal projection onto $N_j$. If $\lambda_j=0$ set $r=j$ and terminate the construction. Otherwise let $\chi_j, \psi_j$ be a Schmidt pair for $\Gamma_j$ corresponding to the singular value $\lambda_j.$ Let $K_{j+1}$ be the range of the orthogonal projection of $K_j$ onto the pointwise orthogonal complement of $\chi_0, \cdots, \chi_j$ in $L^2(\Tt,\Cc^n).$ Let $N_{j+1}$ be the projection of $N_j$ onto the pointwise orthogonal complement of $\psi_0, \cdots, \psi_j$ in $L^2(\Tt,\Cc^m).$ Let $Q_{j+1} \in H^\infty(\Dd,\CCmn)$ be chosen to satisfy, for $0\leq k \leq j,$
	
	\begin{equation}\label{32constr} Q_{j+1} \chi_k = G\chi_k - t_k\psi_k, \quad \psi_k^* Q_{j+1}=\psi_k^* G - t_k \chi_k^*  .\end{equation}
	
	\noindent Then each $\Gamma_j$ is a compact operator, $Q_j$ with the above properties does exist, the construction terminates with $r\leq \min(m,n)$ and 
	
	\begin{equation}\label{g-ag=sum} G-\mathcal{A}G= \displaystyle \sum\limits_{j=0}^{r-1} \frac{\lambda_j \psi_j \chi_j^*}{|h_j|^2}.\end{equation}
	
\end{theorem}

\noindent We shall derive a similar formula for the superoptimal analytic approximant $\mathcal{A}G,$ by making use of exterior products of Hilbert spaces.

\subsection{Algorithm for superoptimal analytic approximation} \label{Alg_statement}

In this section we consider the superoptimal analytic approximation problem  for a function $G  \in H^\infty ( \Dd, \CCmn) + C(\Tt, \CCmn).$ 
We first state the algorithm for the solution of Problem \ref{mainproblem}; later we shall prove the claims that are made in this description of the algorithm. We will assume here the result of Peller and Young \cite{superop} that   Problem \ref{mainproblem} has a unique solution (see Theorem \ref{1.8}).  For convenience, we give citations of the steps in this paper where the corresponding claims are proved.

\begin{proof} [\emph{\textbf{Algorithm:}}]
Let $G\in H^\infty(\Dd,\CCmn)+C(\Tt,\CCmn).$ In this subsection we shall give a fuller and more precise statement of the algorithm for $\mathcal{A}G$ outlined in the Introduction, Section \ref{intro}, in preparation for a subsequent formal proof of Theorem \ref{mathcalAG}, which asserts that if entities $r,t_i,x_i,y_i,h_i$ for $i=0,\dots,r-1,$ are generated by the algorithm, then the superoptimal approximant is given by equation
 $$\mathcal{A}G = \displaystyle G- \sum\limits_{i=0}^{r-1} \frac{t_i y_i x_i^*}{|h_i|^2}.  $$
 
 The proof will be by induction on $r,$ which is the least index $j\geq 0$ such that $T_j =0,$ where $T_0 = H_G, T_1,T_2,\dots$ is a sequence of operators recursively generated by the algorithm.    \\  

\textbf{Step 0.} 
\noindent Let $ t_0 = \| H_G\| .$ 
\noindent If $t_0 = 0,$ then $H_G=0,$ which implies $G\in H^\infty ( \Dd, \CCmn).$ In this case, the algorithm terminates, we define $r$ to be zero and the superoptimal approximant $\mathcal{A}G$ is given by $\mathcal{A}G = G$, in agreement with the formula
\begin{equation}\label{formAG}
 G-\mathcal{A}G = \sum\limits_{j=0}^{r-1} \frac{t_j y_j x_j^*}{x_j^*x_j}  
\end{equation}
 contained in the statement of Theorem \ref{mathcalAG} (since the sum on the right hand side of equation \eqref{formAG} is empty, and therefore by convention is interpreted as being zero).

Otherwise, $t_0 > 0$.   By Theorem \ref{2.04} and Lemma \ref{2.2}, $H_G$ is a compact operator and so there exists a Schmidt pair $(x_0 , y_0)$ corresponding to the singular value $t_0$ of $H_G.$ By the definition of the Schmidt pair $(x_0, y_0)$ corresponding to $t_0$ for the Hankel operator
$H_G : H^2 (\Dd, \Cc^n ) \to H^2 (\Dd, \Cc^m)^\perp,$
$$x_0 \in H^2(\Dd, \Cc^n ),\quad y_0 \in H^2 (\Dd, \Cc^m)^\perp $$
are non-zero vector-valued functions such that
$$H_Gx_0 = t_0 y_0\quad \mbox{ and } \quad H_G^* y_0 =t_0 x_0. $$

By Lemma \ref{2.2}, $x_0 \in H^2(\Dd,\Cc^n)$ and $\bar{z}\bar{y}_0 \in H^2(\Dd, \Cc^m)$ admit inner-outer factorizations 
\begin{equation}\label{xi0eta0} x_0 = \xi_0 h_0 , \quad \bar{z}\bar{y}_0 = \eta_0 h_0 \end{equation} for some scalar outer factor $h_0 \in H^2(\Dd, \Cc)$ and column matrix inner functions $\xi_0\in H^\infty(\Dd, \Cc^n)$, $ \eta_0\in H^\infty(\Dd, \Cc^m). $ Then
\begin{equation}\label{equl}  \|x_0 (z)\|_{\Cc^n} = |h_0(z)| = \|y_0 (z)\|_{\Cc^m} \;\; \text{almost everywhere on }\; \Tt.\end{equation}
                                                                                                              
\noindent We write equations (\ref{xi0eta0}) as 
\begin{equation}\label{31}
\xi_0 = \frac{x_0}{h_0} \; , \; \eta_0 = \frac{\bar{z}\bar{y}_0}{h_0}.
\end{equation}
By equations (\ref{equl}) and (\ref{31}), \begin{equation}\label{xi0eta0=1}
                                           \| \xi_0 (z)  \|_{\Cc^n} =1= \| \eta_0(z) \|_{\Cc^m}\; \text{almost everywhere on}\; \Tt. \end{equation}

\noindent Let $t_0\neq0.$ By Lemma \ref{2.2}, every function $Q_1 \in H^\infty(\Dd, \CCmn)$ which is at minimal distance from $G$ satisfies  
\begin{equation}\label{G-Q0}(G-Q_1)x_0 = t_0 y_0,\quad y_0^* (G-Q_1) = t_0 x_0^*.
\end{equation}

\textbf{Step 1.} Let
\begin{equation}\label{X_0} X_1 \stackrel{\text{def}}{=} \xi_0 \dot{\we} H^2(\Dd,\Cc^n ) \vspace{2ex}.\end{equation}
By Proposition \ref{xjsubseth2}, 
$ X_1$ is a closed subspace of  $H^2(\Dd, \we^2 \Cc^n )$.
\noindent Similarly, $$\eta_0 \telwe zH^2(\Dd,\Cc^m) \subset zH^2(\Dd,\we^2\Cc^m) $$ and therefore
$$\bar{\eta}_0 \telwe \overline{z H^2(\Dd,\Cc^m)} \subset \bar{z} \overline{H^2(\Dd,\we^2\Cc^m)}, $$ that is, if 
\begin{equation}\label{Y_0}   Y_1 \stackrel{\text{def}}{=} \bar{\eta}_0 \dot{\we} H^2(\Dd, \Cc^m)^\perp,\end{equation}
then
$$  Y_1  \subset H^2 (\Dd, \we^2 \Cc^m)^\perp.$$
By Proposition \ref{clwe},
$ Y_1$ is a closed subspace of  $H^2 (\Dd, \we^2 \Cc^m)^\perp .$

\noindent	Choose any function $Q_{1} \in H^\infty(\Dd, \CCmn)$ which satisfies the equations \eqref{G-Q0}.	
 Consider the operator $T_1 : X_1 \to Y_1$ defined by 
\begin{equation}\label{T_0} T_1 ( \xi _0 \dot{\we} x ) = P_{Y_1} (\bar{\eta}_0\dot{\we} (G-Q_1)x) \; \text{for all} \; x\in H^2(\Dd, \Cc^n ),\end{equation}where $P_{Y_1}$ is the projection from $L^2(\Tt, \we^2 \Cc^m)$ on $Y_1.$ 

By Corollary \ref{projwell} and Proposition \ref{Twell}, $T_1$ is well-defined.

  If $T_1=0,$ then the algorithm terminates, we define $r$ to be $1$ and, in agreement with Theorem \ref{mathcalAG}, the superoptimal approximant $\mathcal{A}G$ is given by the formula
$$G- \mathcal{A}G = \displaystyle\sum\limits_{i=0}^{r-1}\frac{t_i y_i x_i^*}{|h_i|^2}= \frac{t_0 y_0 x_0^*}{|h_0|^2}, $$
and the solution is $$\mathcal{A}G =G - \frac{t_0 y_0 x_0^*}{|h_0|^2} .$$

If $T_1\neq 0,$ let $t_{1} = \|T_1\| >0.$ By Theorem \ref{T0compact}, $T_1$ is a compact operator and so
 there exist $v_1 \in H^2(\Dd, \Cc^n),\; w_1 \in H^2(\Dd, \Cc^m)^\perp$ such that
 $(\xi_0 \telwe v_1 , \bar{\eta}_0 \telwe w_1) $ is a Schmidt pair for $T_1$ corresponding to $t_1.$ Let $h_1$ be the scalar outer factor of
 $ \xi_0 \telwe v_1$ and let 
 \begin{equation}
x_{1} = ( I_{\Cc^m} - \xi_0 \xi_0^*)v_1, \;\; y_1= (I_{\Cc^m} - \bar{\eta}_0 \eta_0^T)w_1,
\end{equation} 
where $\mathrm{I}_{\Cc^n}$ and $\mathrm{I}_{\Cc^m}$ are the identity operators in $\Cc^n$ and $\Cc^m$ respectively.\index{identity operator}
Then, by Proposition \ref{onxi},
  \begin{equation}\label{x1y1h1} \|x_1 (z)\|_{\Cc^n} = |h_1(z)| = \|y_1 (z)\|_{\Cc^m}\;\text{almost everywhere on}\; \Tt.\end{equation} 
  
 By Theorem \ref{1.8}, there exists a function $Q_2 \in H^\infty(\Dd, \CCmn)$ such that both\newline $s_0^\infty(G-Q_2)$ and $s_1^\infty(G-Q_2)$ are minimized and 
$$
s_1^\infty(G-Q_2)=t_1.
$$ 
By Proposition \ref{g-q1y1t1x1}, any   such $Q_2$ satisfies 
 \begin{equation}\label{3111}
\begin{aligned}(G-Q_2) x_0 = t_0 y_0 , \quad y_0^* (G-Q_2) = t_0 x_0^*  \\
  (G-Q_2)x_1 = t_1 y_1 , \quad y_1^* (G-Q_2)=t_1 x_1^*.\end{aligned}
\end{equation} 
Choose any function $Q_{2} \in H^\infty(\Dd, \CCmn)$ which satisfies the equations \eqref{3111}.	

Define 
\begin{equation}\label{311}
\xi_1 = \frac{x_1}{h_1}  , \quad \eta_1 = \frac{\bar{z}\bar{y}_1}{h_1}.
\end{equation} By equations (\ref{x1y1h1}) and (\ref{311}),  $\| \xi_1(z)  \|_{\Cc^n} =1= \| \eta_1(z) \|_{\Cc^n}$ almost everywhere on $\Tt.$\\

\textbf{Step 2.} Define
                                                                                          
$$\begin{array}{cllll}
   X_2 &\stackrel{\text{def}}{=} \xi_0 \dot{\we} \xi_1 \dot{\we} H^2(\Dd, \Cc^n)\vspace{2ex}\\
   Y_2 &\stackrel{\text{def}}{=} \bar{\eta}_0\dot{\we} \bar{\eta}_1 \dot{\we} H^2(\Dd, \Cc^m)^\perp .
  \end{array}$$

 Note that, by Proposition \ref{xjclosed}, $X_2$ is a closed linear subspace of $H^2(\Dd, \we^3 \Cc^n),$ and, by Proposition \ref{clwegen}, $Y_2$ is a closed linear subspace of $H^2 (\Dd, \we^3 \Cc^m)^\perp.$

\noindent Now consider the operator $T_2 : X_2 \to Y_2$ given by 
\begin{equation}\label{T1} T_2 ( \xi _0 \dot{\we} \xi_1\dot{\we} x ) = P_{Y_2} (\bar{\eta}_0\dot{\we} \bar{\eta}_1\dot{\we} (G-Q_2)x),\end{equation} where $P_{Y_2}$ is the projection from $L^2(\Tt,\Cc^m)$ on $Y_2.$ 
  
\noindent By Corollary \ref{projwellgen} and Proposition \ref{Twell}, $T_2$ is well defined, that is, it does not depend on the choice of $Q_2 \in H^\infty(\Dd,\CCmn)$ satisfying equations \eqref{3111}.
If $T_2 =0,$ then the algorithm terminates, we define $r$ to be $2$ and, according to Theorem \ref{mathcalAG}, the superoptimal approximant $\mathcal{A}G$ is given by the formula 
$$G- \mathcal{A}G = \sum\limits_{i=0}^{r-1} \frac{t_i y_i x_i^*}{|h_i|^2}= 
\frac{t_0 y_0 x_0^*}{|h_0|^2}+ \frac{t_1 y_1 x_1^*}{|h_1|^2}.$$ 
If $T_2 \neq 0$, then let $t_2= \|T_2\|.$ By Theorem \ref{T2compactt}, $T_2$ is a compact operator and hence there exist $v_2 \in H^2(\Dd,\Cc^n), \; w_2 \in H^2(\Dd,\Cc^m)^\perp$ such that 
$$( \xi _0 \dot{\we} \xi_1\dot{\we} v_2 ,\; \bar{\eta}_0\dot{\we} \bar{\eta}_1\dot{\we} w_2)$$is a Schmidt pair for $T_2$ corresponding to $\|T_2\|=t_2.$ \\ 

Let $h_2$ be the scalar outer factor of
 $ \xi _0 \dot{\we} \xi_1\dot{\we} v_2$.
Note that, by Proposition \ref{xjclosed}, $ \xi _0 \dot{\we} \xi_1\dot{\we} v_2 \in H^2(\Dd, \we^3 \Cc^n).$
Let 
 \begin{equation}
 x_{2} = ( I_{\Cc^m} - \xi_0 \xi_0^*-\xi_1 \xi_1^* )v_2, \;\; y_2= (I_{\Cc^m} - \bar{\eta}_0 \eta_0^T - \bar{\eta}_1 \eta_1^T)w_2.
\end{equation} 
 Then, by Proposition \ref{onxi},
  \begin{equation}\label{x2y2h2} 
\|x_2 (z)\|_{\Cc^n} = |h_2(z)| = \|y_2 (z)\|_{\Cc^m}\;\text{almost everywhere on}\; \Tt.\end{equation} 
Define \begin{equation}\label{xisetais}\xi_2 = \frac{x_2}{h_2} , \; \eta_2 = \frac{\bar{z} \bar{y}_2}{h_2}.  
\end{equation}
\noindent Clearly $\|\xi_2(z)\|_{\Cc^n}=1 $ and $\|\eta_2(z)\|_{\Cc^m}=1$ almost everywhere on $\Tt.$\\

\textbf{Recursive step.} Suppose that, for $j \le \min(m,n) -2$, we have constructed 
\begin{equation}\label{rec_step}\begin{aligned}	&t_0 \geq t_1 \geq \cdots \geq t_j > 0\\
&x_0, x_1, \cdots, x_j \in L^2 (\Tt, \Cc^n)\\
&y_0 , y_1 , \cdots , y_j \in L^2(\Tt, \Cc^m) \\
&h_0, h_1, \cdots, h_j \in H^2(\Dd,\Cc) \; \text{outer}\\
& \xi_0,\xi_1, \cdots, \xi_j \in L^\infty(\Tt,\Cc^n)\; \text{pointwise orthonormal on}\;\Tt\\
& \eta_0, \eta_1, \cdots , \eta_j \in L^\infty(\Tt,\Cc^m)\; \text{pointwise orthonormal on}\;\Tt\\
&X_0 = H^2(\Dd,\Cc^n),X_1, \cdots, X_j \\
&Y_0 = H^2(\Dd,\Cc^m)^\perp, Y_1, \cdots, Y_j\\
&T_0, T_1, \cdots, T_j \; \text{compact operators}.\end{aligned}\end{equation}

 By Theorem \ref{1.8}, there exists a function $Q_{j+1} \in H^\infty(\Dd, \CCmn)$ such that
 $$\left(s_0^\infty(G-Q_{j+1}), s_1^\infty(G-Q_{j+1}), \cdots , s_{j+1}^\infty(G-Q_{j+1})\right)$$ is lexicographically minimized. By Proposition \ref{g-qjj}, any such function $Q_{j+1}$ satisfies 
\begin{equation}\label{g-qi}(G-Q_{j+1})x_i = t_i y_i, \quad y_i^* (G-Q_{j+1}) = t_i x_i^*, \quad i=0, 1, \cdots, j. \end{equation}
\index{$X_{j}$} 
Choose any function $Q_{j+1} \in H^\infty(\Dd, \CCmn)$ which satisfies the equations \eqref{g-qi}.	
Define \begin{equation}\label{X_j} X_{j+1} \stackrel{\text{def}}{=} \xi_0 \dot{\we} \xi_1 \dot{\we} \cdots \dot{\we} \xi_j \dot{\we} H^2(\Dd,\Cc^n) \end{equation}
\index{$Y_{j}$} 
\begin{equation}\label{Y_j} Y_{j+1} \stackrel{\text{def}}{=} \bar{\eta}_0 \dot{\we} \bar{\eta}_1 \dot{\we} \cdots \dot{\we} \bar{\eta}_j \dot{\we} H^2 (\Dd, \Cc^m)^\perp .\end{equation}
Note that, by Proposition \ref{xjsubseth2}, $X_{j+1}$ is a subset of $H^2(\Dd,\we^{j+2}\Cc^n),$ and, by Proposition \ref{etatelwelj},  $Y_{j+1}$ is a closed subspace of $H^2 (\Dd, \we^{j+2} \Cc^m)^\perp.$ 
\index{$T_{j}$}
\noindent Consider the operator $$T_{j+1} : X_{j+1} \to Y_{j+1}$$ given by 
\begin{equation}\label{T_j}
T_{j+1}(\xi_0 \dot{\we} \xi_1 \dot{\we} \cdots \dot{\we} \xi_j \dot{\we} x)= P_{Y_{j+1}} \left( \bar{\eta}_0 \dot{\we} \bar{\eta}_1 \dot{\we} \cdots \dot{\we} \bar{\eta}_j\dot{\we} (G-Q_{j+1})x \right).\end{equation}
By Corollary \ref{projwellgen} and Proposition \ref{Twell}, $T_{j+1}$ is well-defined and does not depend on the choice of $Q_{j+1}$ subject to equations \eqref{g-qi}.

If $T_{j+1}=0,$ then the algorithm terminates, we define $r$ to be $j+1,$ and, according to Theorem \ref{mathcalAG}, the superoptimal approximant $\mathcal{A}G$ is given by the formula
$$G- \mathcal{A}G = \sum\limits_{i=0}^{r-1} \frac{t_i y_i x_i^*}{|h_i|^2}.$$

Otherwise, we define $t_{j+1} = \|T_{j+1}\| >0.$ By Theorem \ref{T0compact}, $T_{j+1}$ is a compact operator and hence there exist $v_{j+1} \in H^2(\Dd,\Cc^n), \; w_{j+1} \in H^2(\Dd,\Cc^m)^\perp$ such that \begin{equation}\label{schmpairtj+1}(\xi_0 \dot{\we} \xi_1 \dot{\we} \cdots \dot{\we} \xi_j \dot{\we} v_{j+1},  \bar{\eta}_0 \dot{\we} \bar{\eta}_1 \dot{\we} \cdots \dot{\we} \bar{\eta}_j \dot{\we} w_{j+1})\end{equation}
is a Schmidt pair for $T_{j+1}$ corresponding to the singular value $t_{j+1}.$ 

\noindent Let $h_{j+1}$ be the scalar outer factor of $\xi_0 \dot{\we} \xi_1 \dot{\we} \cdots \dot{\we} \xi_j \dot{\we} v_{j+1},$ and let
\begin{equation}\label{xj+1yj+1}x_{j+1} = (I_{\Cc^n} - \xi_0 \xi_0^* - \cdots - \xi_j \xi_j^*)v_{j+1}, \quad y_{j+1} = (I_{\Cc^m} - \bar{\eta}_0 \eta_0^T - \cdots- \bar{\eta}_j\eta_j^T) w_{j+1}.\end{equation} 
 Then, by Proposition \ref{onxi},
  \begin{equation}\label{xjyjhj} \index{$x_{j}$} \index{$y_{j}$}
\|x_{j+1} (z)\|_{\Cc^n} = |h_{j+1}(z)| = \|y_{j+1} (z)\|_{\Cc^m}\;\text{almost everywhere on}\; \Tt.\end{equation} 
We define 
\begin{equation}\label{xij+1etaj+1}
\xi_{j+1} = \frac{x_{j+1}}{h_{j+1}}, \quad \eta_{j+1}=\frac{\bar{z}\overline{y_{j+1}}}{h_{j+1}}.\end{equation}

\noindent Clearly $\|\xi_{j+1}(z)\|_{\Cc^n}=1 $ and $\|\eta_{j+1}(z)\|_{\Cc^m}=1$ almost everywhere on $\Tt.$\\

This completes the recursive step. The algorithm terminates after at most $\min(m,n)$ steps, so that $r \leq \min(m,n)$ and, in accordance with Theorem \ref{mathcalAG} the superoptimal approximant $\mathcal{A}G$ is given by the formula 

$$G - \mathcal{A}G = \sum\limits_{i=0}^{r-1} \frac{ t_i y_i x_i^*}{|h_i|^2}.$$
 \end{proof}

%% file: orthonormal.tex
 \section{Pointwise orthonormality of $\{\xi_i\}_{i=1}^j$ and $\{\bar{\eta}_i\}_{i=1}^j$ almost everywhere on $\Tt$}\label{orthonomal}
These orthonormality properties will be needed for the justification of the main algorithm.

\begin{proposition}\label{onxi}
 Let $m,n$ be positive integers with $\min(m,n) \geq 2,$ let \newline $G\in H^\infty(\Dd,\CCmn)+C(\Tt,\CCmn)$ and let $0\leq j\leq \min(m,n)-2.$ Suppose we have applied steps $0,\dots,j$ of the superoptimal analytic approximation algorithm from Subsection \ref{Alg_statement} to $G$ and we have obtained $x_i,y_i$ as in equations \eqref{xj+1yj+1}, and $\xi_i, \eta_i$ as in equations \eqref{xij+1etaj+1} for $i=0,\cdots, j.$
Then
\begin{enumerate}
	\item[\emph{(i)}] $\xi_0 \telwe v_1 = \xi_0 \telwe x_1, \quad \xi_0\telwe \cdots \telwe \xi_{j-1} \telwe v_j =\xi_0\telwe \cdots \telwe \xi_{j-1} \telwe x_j,  $  
	$\bar{\eta}_0 \telwe w_1 = \bar{\eta}_0 \telwe y_1,$ \newline and $\bar{\eta}_0 \telwe\cdots \telwe \bar{\eta}_{j-1}\telwe w_j =  \bar{\eta}_0 \telwe\cdots \telwe \bar{\eta}_{j-1}\telwe y_j ;$
	\item[\emph{(ii)}] $\|x_j(z)\|_{\Cc^n} = \|y_j(z)\|_{\Cc^m}=|h_j(z)| $ almost everywhere on $\Tt;$
\item[\emph{(iii)}] The sets $\{\xi_i(z)\}_{i=0}^{j}$ and $\{\bar{\eta}_i(z) \}_{i=0}^j $   are orthonormal in $\Cc^n$ and $\Cc^m$ respectively for almost every $z \in \Tt.$ 
 
\end{enumerate}
 
\end{proposition}

\begin{proof}
We will prove statement (ii) in Propositions \ref{x0wev1eta1wew1} and \ref{x0wev2eta2wew2}. Statement (i) is proven below in equations \eqref{xi0telx1}, \eqref{xi0xijisxi0xj}, \eqref{eta-wi-yi}. 
Let us prove assertion (iii).

 Since the function $G$ belongs to $ H^\infty(\Dd,\CCmn)+C(\Tt,\CCmn),$ by Hartman's theorem, the Hankel operator with symbol $G$, denoted by $H_G,$ is a compact operator, and so there exist functions $$x_0 \in H^2(\Dd,\Cc^n),\quad y_0 \in H^2(\Dd,\Cc^m)^\perp$$ such that $(x_0,y_0)$ is a Schmidt pair corresponding to the singular value $t_0= \|H_G\|\neq 0.$ By Lemma  \ref{2.2}, $x_0 , \bar{z}\bar{y}_0 $ admit the inner-outer factorizations 
 $$ x_0 = \xi_0 h_0, \quad \bar{z}\bar{y}_0 = \eta_0 h_0$$ for column matrix inner functions $\xi_0\in H^\infty(\Dd, \Cc^n)$, $ \eta_0\in H^\infty(\Dd, \Cc^m) $ and some scalar outer factor $h_0 \in H^2(\Dd, \Cc).$ By Theorem \ref{1.7},
 \begin{equation}\label{normxoyo} \|x_0 (z)\|_{\Cc^n} = |h_0(z)| = \|y_0 (z)\|_{\Cc^m} \;\; \text{almost everywhere on }\; \Tt.\end{equation}
Thus \begin{equation}
\label{xi01}\|\xi_{0}(z)\|_{\Cc^n} =1 \;\text{almost everywhere on}\; \Tt.\end{equation} Hence (iii) of Proposition \ref{onxi} holds for $\{\xi_i(z)\}_{i=0}^j$ in the case that $j=0.$ 
 
 \noindent Let $T_1$ be given by equation (\ref{T_0}). By the hypothesis \eqref{rec_step}, $T_1$ is a compact operator, and if $T_1 \neq 0,$ then there exist $v_1\in H^2(\Dd,\Cc^n)$ and $w_1\in H^2(\Dd,\Cc^m)^\perp$ such that $(\xi_0 \telwe v_1, \bar{\eta}_0\telwe w_1)$ is a Schmidt pair corresponding to $\|T_1\|=t_1.$ By Proposition \ref{xweyh2}, $\xi_0 \telwe v_1 \in H^2(\Dd,\we^2\Cc^n).$ Let $h_1$ be the scalar outer factor of $\xi_0\telwe v_1.$ We define
 \begin{equation}\label{x1} x_1= (I_{n} - \xi_0 \xi_0^*) v_1 \end{equation} and  \begin{equation}\label{xi11} \xi_1= \frac{x_1}{h_1}. \end{equation}
Then, for $z\in \Dd,$
 $$ \xi_1 (z) =\displaystyle\frac{1}{h_1(z)} v_1(z)-\frac{1}{h_1(z)} \xi_0 (z) \xi_0 (z)^* v_1(z). $$
Note that by equation (\ref{xi01}),
$$ \xi_0^*(z)\xi_0(z) = \langle \xi_0 (z), \xi_0 (z)\rangle_{\Cc^n} =1\quad\text{almost everywhere on}\; \Tt,$$hence$$ \langle  \xi_1(z), \xi_0(z)\rangle_{\Cc^n} = \xi_0^*(z)  \xi_1(z)=\displaystyle\frac{1}{h_1(z)}\xi_0(z)^* v_1(z) - \displaystyle\frac{1}{h_1(z)}\xi_0 (z)^* \xi_0 (z) \xi_0 (z)^* v_1(z) =0$$ almost everywhere on $\Tt.$ Note that, by equation \eqref{x1}, for almost every $z \in \Tt,$
\begin{align}  \xi_0 (z) \we v_1(z) &=\xi_0 (z) \we (x_1(z) + \xi_{0}(z) \xi_0(z)^* v_1(z) )\nonumber\vspace{2ex} \\&= \xi_0 (z) \we x_1(z) + \xi_0 (z) \we \xi_{0}(z) \xi_0(z)^* v_1(z) )\nonumber\vspace{2ex} \\ &= \xi_0(z) \we x_1(z)\label{xi0telx1},\end{align} the last equality following from the pointwise linear dependence of the vectors $\xi_0$ and 
$z\mapsto \xi_0(z) \langle v_1(z), \xi_0(z)\rangle_{\Cc^n}$ almost everywhere on $\Tt.$ 

 Moreover,
 since $h_1$ is the scalar outer factor of $\xi_0 \telwe v_1,$  for almost every $z \in \Tt,$ we have $$\begin{array}{clll}|h_1(z)|&= \| \xi_0 (z) \we v_1(z)\|_{\we^2{\Cc^n}} = \| \xi_0(z) \we x_1(z)\|_{\we^2\Cc^n},\end{array}$$ By Lemma \ref{weon}, 
 
 $$  \| \xi_0(z) \we x_1(z)\|_{\we^2\Cc^n} = \| x_1(z)- \langle x_1(z) , \xi_0(z) \rangle_{\Cc^n} \xi_0(z)\|_{\Cc^n} = \| x_1(z) \|_{\Cc^n}  $$ almost everywhere on $\Tt.$    Hence, for almost every $z\in \Tt,$ \begin{equation}\label{h1=xi} |h_1(z)| = \|x_1(z)\|_{\Cc^n} \end{equation} and thus
 $$\| \xi_1(z) \|_{\Cc^n} = \displaystyle\frac{\|x_1(z)\|_{\Cc^n}}{|h_1(z)|} =1  \quad \text{almost everywhere on} \; \Tt.$$ 

\noindent Consequently, $\{ \xi_0(z) , \xi_1(z)\}$ is an orthonormal set in $\Cc^n$ for almost every $z \in \Tt.$ Hence (iii) of Proposition \ref{onxi} holds for $\{\xi_i(z)\}_{i=0}^j$ in the case that $j=1.$ 

\textbf{ Recursive step:} Suppose the entities in equations \eqref{rec_step} have been constructed and have the stated properties. Since by the inductive hypothesis $T_j$ is a compact operator,
there exist $$v_j \in H^2(\Dd,\Cc^n),\quad w_j \in H^2(\Dd,\Cc^m)^\perp$$ such that 
$$ (\xi_0 \telwe \xi_1 \telwe \cdots \telwe \xi_{j-1} \telwe v_j, \overline{ \eta}_0\telwe \overline{ \eta}_1\telwe \cdots \telwe \overline{ \eta}_{j-1}\telwe w_j) $$is a Schmidt pair for $T_j$ corresponding to $\|T_j\|=t_j.$  Let us first prove that $\xi_0\telwe \xi_1 \telwe  \cdots \telwe \xi_{j-1}\telwe v_j$ is an element of $ H^2(\Dd,\we^{j+1}\Cc^n).$ By hypothesis,
\[x_i = (I_{n}-\xi_0 \xi_0^* -\dots - \xi_{i-1}\xi_{i-1}^*)v_i\quad \text{and}\quad \xi_i=\frac{x_i}{h_i}\]for $i=0,\dots,j-1.$ Then, for all $z\in \Dd,$
\begin{align*}
\left( \xi_0\telwe\xi_1 \telwe  \cdots \telwe \xi_{j-1} \telwe v_j\right) (z) &= \left(\xi_0\telwe\frac{x_1}{h_1}\telwe  \cdots \telwe \frac{x_{j-1}}{h_{j-1}}\telwe v_j\right)(z)\\
&= \left(\frac{1}{h_1}\cdots \frac{1}{h_{j-1}} \xi_0 \telwe x_1 \telwe \cdots \telwe x_{j-1}\telwe v_j\right)(z).
\end{align*} We obtain
\[ (\xi_0\telwe\xi_1 \telwe  \cdots \telwe \xi_{j-1} \telwe v_j)(z) = \left(\frac{1}{h_1}\cdots \frac{1}{h_{j-1}} \xi_0 \telwe v_1 \telwe \cdots \telwe v_{j-1}\telwe v_j\right)(z),\; \text{for all}\;z\in \Dd, \]due to pointwise linear dependence of $\xi_k$ and 
$z\mapsto \xi_k(z) \langle v_{j}(z), \xi_k(z)\rangle_{\Cc^n}$ 
 on $\Dd,$ for all $k=0,\dots,j-1.$ By Proposition \ref{wejanalytic}, 
\[\frac{1}{h_1}\cdots \frac{1}{h_{j-1}} \xi_0 \telwe v_1 \telwe \cdots \telwe v_{j-1}\telwe v_j\] is analytic on $\Dd.$ Moreover, by Proposition \ref{wel2conv}, since $\xi_0,\xi_1,\dots,\xi_{j-1}$ are pointwise orthogonal on $\Tt,$ \[ \|\xi_0\telwe \xi_1 \telwe \cdots\telwe \xi_{j-1}\telwe v_j\|_{L^2(\Tt,\we^{j+1}\Cc^n)} < \infty. \]Therefore \[ \xi_0\telwe \xi_1 \telwe \cdots\telwe \xi_{j-1}\telwe v_j \in H^2(\Dd,\we^{j+1}\Cc^n). \]

 Let $h_{j}$ be the scalar outer factor of $\xi_1 \telwe \xi_2 \telwe \cdots \telwe \xi_{j-1} \telwe v_j.$ We define
\begin{equation}\label{xj} x_{j} = (I_{n} - \xi_0 \xi_0^* - \cdots - \xi_{j-1} \xi_{j-1}^*)v_j\end{equation} and 
\begin{equation}\label{xijj} \xi_{j}=\frac{x_{j}}{ h_{j}}.\end{equation}

Let us show that $\{  \xi_0(z), \cdots, \xi_{j-1}(z), \xi_{j}(z)\}$ is an orthonormal set in $\Cc^n$ almost everywhere on $\Tt.$ 
We have
$$ \xi_{j} = \displaystyle\frac{1}{h_{j}} v_j - \displaystyle\frac{1}{h_{j} } \xi_0  \xi_0^*  v_j - \cdots - \displaystyle\frac{1}{h_{j} }\xi_{j-1}  \xi_{j-1}^*  v_j ,$$ and so, for $i=0,\dots, j-1,$
$$
\begin{array}{cll}

&\langle \xi_{j}(z), \xi_{i}(z) \rangle_{\Cc^n}= \displaystyle\frac{1}{h_j(z)}\xi_{i}^*(z) v_j(z) - \displaystyle\frac{1}{h_j(z)}\xi_{i}^*(z) \xi_0(z) \xi_0^*(z) v_1(z)- \cdots \vspace{2ex}\\
&\hspace{37ex}-\displaystyle\frac{1}{h_j(z)}\xi_{i}^*(z) \xi_{j-1}(z) \xi_{j-1}^*(z) v_1(z) \end{array}$$ \text{almost everywhere on} $\Tt.$
Note that by the inductive hypothesis, for $i,k=0,1,\cdots, j-1$ and for almost all $z \in \Tt$,  
$$  \xi_{i}^*(z) \xi_k (z)= \left\{ \begin{array}{ll} 0,\quad \text{for} \; i\neq k\\ 1,\quad \text{for} \;i=k\end{array}.\right.$$
Thus, for $i=0,\dots,j-1,$ 
$$\langle \xi_{j}(z), \xi_{i}(z) \rangle_{\Cc^n}= \displaystyle\frac{1}{h_{j}(z)}\xi_{i}^*(z) v_j(z) - \displaystyle\frac{1}{h_{j}(z)}\xi_{i}^*(z) \xi_i(z) \xi_i^* v_j(z)=0 $$almost everywhere on $\Tt,$ and hence, by induction on $j$ and for all integers $j=0,\dots, r-1,$
 $\{  \xi_0(z), \cdots, \xi_{j-1}(z), \xi_{j}(z)\}$ is an orthogonal set in $\Cc^n$ for almost all $z \in \Tt.$ 

Let us show that $$\xi_0(z) \we \cdots \we \xi_{j-1}(z)\we v_j(z) =\xi_0(z) \we \cdots \we \xi_{j-1}(z)\we  x_{j}(z) $$almost everywhere on $\Tt.$ Equation \eqref{xj} yields $$
\begin{array}{cllll}

&\xi_0(z) \we \cdots \we \xi_{j-1}(z)\we v_j(z) \vspace{2ex}\\&= \xi_0(z) \we \cdots \we \xi_{j-1}(z)\we (x_{j}(z)+ \xi_0 (z) \xi_0^* (z) v_j(z) \\
&\hspace{32ex}+ \cdots + \xi_{j-1}(z)\xi_{j-1}^*(z)v_j(z))  \\  
&=  \xi_0(z) \we \cdots \we \xi_{j-1}(z)\we  (x_{j}(z) + \xi_0 (z) \langle v_j(z), \xi_0(z)\rangle_{\Cc^n} + \cdots  \\
&\hspace{32ex}+ \cdots + \xi_{j-1}(z) \langle v_j(z) ,\xi_{j-1}(z)\rangle_{\Cc^n} ) \end{array}$$almost everywhere on $\Tt.$ 
\noindent Notice that, for  $i=0,\cdots, j-1,$ 
the vectors $\xi_i$ and $z \mapsto \xi_i(z)\langle v_j(z), \xi_i(z)\rangle_{\Cc^n}$ are pointwise linearly dependent almost everywhere on $\Tt.$ Thus for all $i=0,\cdots,j-1,$
$$ \xi_0 (z) \we \cdots \we \xi_{j-1}(z) \we \xi_i (z) \langle v_{i+1}(z), \xi_i(z) \rangle_{\Cc^n} =0   $$ almost everywhere on $\Tt.$  

\noindent Hence 
\begin{equation}\label{xi0xijisxi0xj}\xi_0(z) \we \cdots \we \xi_{j-1}(z)\we  v_{j}(z) =\xi_0(z) \we \cdots \we \xi_{j-1}(z)\we  x_{j}(z) \quad \text{almost everywhere on} \; \Tt. \end{equation}

\noindent Next, we shall show that $\|\xi_j(z)\|_{\Cc^n}=1$ for almost all $z \in \Tt.$ Recall that $h_j$ is the scalar outer factor of $\xi_1 \telwe \xi_2 \telwe \cdots \telwe \xi_{j-1} \telwe v_j,$ and therefore 

$$
\begin{array}{cllll}

|h_{j}(z)|&=\|\xi_0(z) \we \cdots \we \xi_{j-1}(z)\we v_j(z)\|_{\we^{j+1}\Cc^n} =  \| \xi_0(z) \we \cdots \we \xi_{j-1}(z)\we  x_{j}(z)\|_{\we^{j+1}\Cc^n}\\  
 \end{array}$$almost everywhere on $\Tt.$
 
\noindent By the inductive hypothesis, $\{\xi_0(z) , \cdots, \xi_{j-1}(z)\}$ is an orthonormal set in $\Cc^n$ for almost all $z \in \Tt,$  hence, by Lemma \ref{weon},  
\begin{align}
|h_j(z)|&= \| \xi_0(z) \we \cdots \we \xi_{j-1}(z)\we  x_{j}(z)\|_{\we^{j+1}\Cc^n}\nonumber \vspace{2ex}\\
&= \| x_j(z) - \sum\limits_{i=0}^{j-1} \langle x_{j}(z), \xi_i(z)\rangle \xi_i(z) \|_{\Cc^n}\nonumber \vspace{2ex} \\
&= \|x_{j}(z)\|_{\Cc^n} \;\text{almost everywhere on}\;\Tt.\label{hj=xj} \end{align}
Thus $$\| \xi_{j}(z)\|_{\Cc^n} =\frac{\| x_j(z)\|_{\Cc^n}}{|h_j(z)|}  =1$$ almost everywhere on $\Tt,$ and hence, by induction on $j,$
$\{  \xi_0(z), \cdots, \xi_{j-1}(z), \xi_{j}(z)\}$ is an orthonormal set in $\Cc^n$ for almost all $z \in \Tt,$ and for all integers $j=0,\dots,r-1.$

\noindent Next, we will prove inductively that the set $\{\bar{\eta}_i(z) \}_{i=0}^j, $ defined in equations (\ref{xij+1etaj+1}), is orthonormal. For $i=0,$ by equation (\ref{normxoyo}), we have
 \begin{equation}
	\label{eta01}\|\bar{\eta}_{0}(z)\|_{\Cc^m} =1 \;\text{almost everywhere on}\; \Tt.\end{equation} 
\noindent 	Let $T_1$ be given by equation (\ref{T_0}). $T_1$ is assumed to be a compact operator, and if $T_1 \neq 0,$ there exist $v_1\in H^2(\Dd,\Cc^n)$ and $w_1\in H^2(\Dd,\Cc^m)^\perp$ such that $(\xi_0 \telwe v_1, \bar{ \eta}_0\telwe w_1)$ is a Schmidt pair corresponding to $\|T_1\|=t_1.$ Suppose $h_1$ is the scalar outer factor of $\xi_0\telwe v_1.$  Let
	\begin{equation}\label{y1} y_1= (I_{m} - \bar{ \eta}_0 \eta_0^T) w_1  = w_1 -\bar{ \eta}_0  \eta_0^T   w_1 \end{equation} and let $$\eta_1(z) = \frac{\bar{z} \bar{y}_1(z)}{h_1(z)} \quad \text{almost everywhere on} \; \Tt. $$ Then,
$$ \bar{\eta}_1(z)= \frac{zy_0(z)}{\bar{h}_1(z)} = \frac{z w_1(z) }{\bar{h}_1(z)} - \frac{z \bar{\eta}_0(z) \eta_0^T(z) w_1(z)}{\bar{h}_1(z)} \quad \text{almost everywhere on} \; \Tt.$$
By equation (\ref{eta01}), $\left\|\bar{ \eta}_0(z)\right\|_{\Cc^m}=1$ almost everywhere on $\Tt.$
Hence
$$
\begin{array}{cllllll}

\langle \bar{\eta}_1(z) , \bar{\eta}_0(z) \rangle_{\Cc^m} &= \eta_0^T (z) \bar{\eta}_1(z)\vspace{2ex}\\
&=\displaystyle \frac{z \eta_0^T (z) w_1(z) }{\bar{h}_1(z)} - \frac{z \eta_0^T (z)\bar{\eta}_0(z) \eta_0^T(z) w_1(z)}{\bar{h}_1(z)}\vspace{2ex}\\
&= \displaystyle\frac{z \eta_0^T (z) w_1(z) }{\bar{h}_1(z)} - \frac{z \langle \bar{\eta}_0(z), \bar{\eta}_0(z)\rangle_{\Cc^m}\eta_0^T (z)w_1(z)}{\bar{h}_1(z)}\vspace{2ex}\\
&= \displaystyle\frac{z \eta_0^T (z) w_1(z) }{\bar{h}_1(z)} - \frac{z \eta_0^T (z) w_1(z) }{\bar{h}_1(z)}\vspace{2ex} \\ &= 0\quad \text{almost everywhere on} \; \Tt.

\end{array}$$
Recall that $h_1$ is the scalar outer factor of $\xi_0 \telwe v_1.$ 
By equation \eqref{h1=xi} and Proposition \ref{x0wev1eta1wew1}, $$\|x_1(z)\|_{\Cc^n} = \|y_1(z)\|_{\Cc^m} = |h_1(z)| $$almost everywhere on $\Tt,$ thus $$\|\bar{ \eta}_1(z) \|_{\Cc^m} =\frac{\|z y_1(z) \|_{\Cc^m}}{|\bar{h}_1(z)|}=1 \quad \text{almost everywhere on} \; \Tt. $$ Consequently, $\{ \bar{ \eta}_0(z) , \bar{ \eta}_1(z)\}$ is an orthonormal set in $\Cc^m$ for almost every $z \in \Tt.$ Hence (iii) of Proposition \ref{onxi} holds for $\{\bar{\eta}_i\}_{i=0}^j$ in the case $j=1.$ 

\textbf{Recursive step:} Suppose the entities in equations \eqref{rec_step} have been constructed and have the stated properties. Since by the inductive hypothesis $T_j$ is a compact operator,
there exist $$v_j \in H^2(\Dd,\Cc^n), \quad w_j \in H^2(\Dd,\Cc^m)^\perp$$ such that 
$$ (\xi_1 \telwe \xi_2 \telwe \cdots \telwe \xi_{j-1} \telwe v_j, \bar{ \eta}_0\telwe \bar{ \eta}_1\telwe \cdots \telwe \bar{ \eta}_{j-1}\telwe w_j) $$is a Schmidt pair for $T_j$ corresponding to $\|T_j\|=t_j.$ We have proved above that $$\xi_0 \telwe \cdots\telwe \xi_{j-1} \telwe v_j \in H^2(\Dd,\we^{j+1}\Cc^n).$$
Let $h_{j}$ be the scalar outer factor of $\xi_0 \telwe \xi_1 \telwe \cdots \telwe \xi_{j-1} \telwe v_j.$ We define
$$y_{j} =( I_{m} - \bar{\eta}_0 \eta_0^T - \cdots - \bar{\eta}_{j-1} \eta_{j-1}^T)w_j$$ and 
\begin{equation}\label{etajj} \bar{\eta}_{j}=\frac{zy_{j}}{ \bar{h}_{j}}\end{equation} 
Let us show that $\{ \bar{ \eta}_0 (z), \dots, \bar{ \eta}_{j}(z) \}$ is an orthonormal set in $\Cc^m$ almost everywhere on $\Tt.$
 We have 
$$ \bar{\eta}_{j} =  \frac{z w_j  }{\bar{h}_{j} } -\cdots -\frac{{z} \bar{\eta}_{j-1} \eta_{j-1}^T w_j }{\bar{h}_{j} } $$ 
and so, for  $i=0,\dots, j-1,$

$$
\begin{array}{cllllll}
\langle \bar{\eta}_{j}(z), \bar{\eta}_{i}(z) \rangle_{\Cc^m}&= \eta_i^T(z) \bar{\eta}_{j}(z)\vspace{2ex} \\
&= \displaystyle \frac{z \eta_i^T (z) w_j(z) }{\bar{h}_{j}(z)} -\cdots -\frac{z\eta_i^T(z) \bar{\eta}_{j}(z)\eta_j^T(z)w_j(z)}{\bar{h}_{j}(z)} \end{array}$$ 
almost everywhere on $\Tt.$

\noindent Notice that, by the inductive hypothesis, for $i,k= 0 ,\dots, j-1$ and for almost all $z \in \Tt,$ 

$$  \eta_{i}^T(z) \bar{\eta}_k (z)= \left\{ \begin{array}{ll} 0,\quad \text{for} \; i\neq k\\ 1,\quad \text{for} \;i=k\end{array}.\right.$$ Hence, for  $i=0,\dots,j-1,$

$$\langle \bar{\eta}_{j}(z), \bar{\eta}_{i}(z) \rangle_{\Cc^m}= \displaystyle \frac{z \eta_i^T(z)  w_j(z) }{\bar{h}_{j}(z)} - \displaystyle \frac{z \eta_i^T (z) w_j(z) }{\bar{h}_{j}(z)}=0$$ almost everywhere on $\Tt.$ Thus by induction on $j$, for all integers $j=0,\dots,r-1,$ $\{ \bar{ \eta}_0 (z), \dots, \bar{ \eta}_{j}(z) \}$ is an orthogonal set in $\Cc^m$ almost everywhere on $\Tt.$

\noindent To complete the proof, we have to prove that 
$\left\|\bar{\eta}_{j}(z)\right\|_{\Cc^m}=1$
for almost all $z \in \Tt.$ Recall that $h_j$ is the scalar outer factor of  $\xi_0 \telwe \xi_1 \telwe \cdots \telwe \xi_{j-1} \telwe v_j.$ By Proposition \ref{xjwevjetajwewj}, $$|h_j(z)|= \| x_j(z)\|_{\Cc^n} = \|y_j(z)\|_{\Cc^m} $$ almost everywhere on $\Tt,$ thus 
$$
\left\|\bar{\eta}_{j}(z)\right\|_{\Cc^m}=\displaystyle \|\frac{zy_j(z)}{\bar{h}_j(z)}\|_{\Cc^m}=1
$$ 
almost everywhere on $\Tt,$ and hence, $\{ \bar{ \eta}_0 (z), \dots, \bar{ \eta}_{j}(z) \}$ is an orthonormal set in $\Cc^m$ almost everywhere on $\Tt.$

Note that, for $j= 1, \dots, r-1 $,
\begin{align}\label{eta-wi-yi} \bar{\eta}_0 \telwe \cdots \telwe \bita_{j-1} \telwe y_{j} &=\bar{\eta}_0 \telwe \cdots \telwe \bita_{j-1}\telwe  (I_{m}-\bar{\eta}_0\eta_0^T-\dots- \bita_{j-1}\eta_{j-1}^T)w_{j}\nonumber\vspace{2ex}\\ &=\bar{\eta}_0 \telwe \cdots \telwe \bita_{j-1} \telwe w_{j} - \sum\limits_{k=0}^{j-1} \bar{\eta}_0 \telwe \cdots \telwe \bita_{j-1}\telwe  \bita_k \eta_k^Tw_{j}\nonumber\vspace{2ex}\\&= \bar{\eta}_0 \telwe \cdots \telwe \bita_{j-1} \telwe w_{j}  \end{align}
on account of the pointwise linear dependence of $\bita_{k}$ and $z\mapsto \bita_k(z) \langle w_{j}(z), \bita_k(z)\rangle_{\Cc^m}$ almost everywhere on $\Tt$. 
\end{proof}

\section{The closed subspace $X_{j+1}$ of $H^2(\Dd,\we^{j+2}\Cc^n)$}
\label{Xjsubset}

\noindent Notice that, although $x_0 \in H^2(\Dd, \Cc^n)$ and $\xi_0$ is inner,  $x_i$ and $\xi_i$ might not be in $H^2(\Dd, \Cc^n)$ in general for $i=1,\cdots, \min(m,n)-1.$ However, for every $x \in H^2(\Dd,\Cc^n),$ the pointwise wedge product $$\xi_0\telwe\cdots\telwe\xi_j\telwe x $$ \emph{is} an element of $H^2(\Dd,\we^{j+2}\Cc^n)$ as the following proposition asserts.

\begin{proposition}\label{xjsubseth2}
	Let $G \in H^\infty(\Dd, \CCmn)+C(\Tt,\CCmn),$ and let $j\leq n-1. $ Let the vector-valued functions $\xi_0, \xi_1, \cdots, \xi_j$ be constructed after applying steps $0,\dots,j$ of the algorithm above and be given by equations \eqref{xij+1etaj+1}. Then 
	$$ \xi_0 \telwe \dots \telwe \xi_j \telwe H^2(\Dd,\Cc^n)$$
	is a subset of $H^2(\Dd,\we^{j+2}\Cc^n).$
\end{proposition}
\begin{proof} For $j=0,$ since $G\in H^\infty(\Dd,\CCmn)+C(\Tt,\CCmn),$ the Hankel operator $H_G$ is compact. There exist $x_0 \in H^2(\Dd,\Cc^n), y_0 \in H^2(\Dd,\Cc^m)^\perp$ such that $(x_0,y_0)$ is a Schmidt pair for the Hankel operator $H_G$ corresponding to the singular value $\|H_G\|.$ By Lemma $\ref{2.2},$ $x_0, y_0$ admit the inner-outer factorizations $$x_0 = \xi_0 h_0 ,\quad  \bar{z} \bar{y}_0 = \eta_0 h_0$$ for some inner $\xi_0 \in H^\infty(\Dd,\Cc^n), \eta_0 \in H^\infty(\Dd,\Cc^m)$ and some scalar outer $h_0\in H^2(\Dd,\Cc).$ 
	
	\noindent Then, by Proposition \ref{xweyh2}, $\xi_0 \telwe H^2(\Dd,\Cc^n) \subset H^2(\Dd,\we^2\Cc^n).$ 
	
	Let us now consider the case where $j=1.$ By definition, $$X_1 = \xi_0 \telwe H^2(\Dd,\Cc^n), \quad Y_1 = \bar{\eta}_0 \telwe H^2(\Dd,\Cc^m)^\perp$$ and, by the inductive hypothesis, $T_1 \colon X_1 \to Y_1$ given by equation (\ref{T_0}) is a compact operator. Suppose $\|T_1\|\neq 0$ and let $(\xi_0 \telwe v_1 , \bar{\eta}_0 \telwe w_1)$ be a Schmidt pair corresponding to $\|T_1\|,$ where $v_1 \in H^2(\Dd,\Cc^n)$ and $w_1 \in H^2(\Dd,\Cc^m)^\perp.$ We define 
	$$x_1 = (I_{n}-\xi_0\xi_0^*)v_1.$$ Note that, by Proposition \ref{xweyh2}, $\xi_0 \telwe v_1 \in H^2(\Dd,\we^2\Cc^n).$ Let $h_1 \in H^2(\Dd,\Cc)$ be the scalar outer factor of $\xi_0 \telwe v_1 \in H^2(\Dd,\we^2\Cc^n).$ 
	Then we define $$\xi_1 = \frac{x_1}{h_1}.$$ Note that $\xi_0$ and 
	$z \mapsto \xi_0(z)\langle v_1(z), \xi_0(z)\rangle_{\Cc^n}$
	 are pointwise linearly dependent on $\Dd,$ since $\xi_0^*v_1$ is a mapping from $\Dd$ to $\Cc.$ Thus, for all $x \in H^2(\Dd,\Cc^n)$ and $z \in \Dd,$ we have 
	
	$$	(\xi_0 \telwe \xi_1 \telwe x)(z)
	= \displaystyle \xi_0(z) \we \xi_1(z)\we x(z)\vspace{2ex}\\
	= \displaystyle \xi_0(z) \we \frac{x_1(z)}{h_1(z)}\we x(z),$$
	and by substituting the value of $x_1$, we find

	$$\begin{array}{cllll}
	&\displaystyle \xi_0(z) \we \frac{x_1(z)}{h_1(z)}\we x(z)\vspace{2ex}\\&=\displaystyle \frac{1}{h_1(z)} \xi_0(z) \we (v_1(z)- \xi_0(z) \xi_0(z)^* v_1(z))\we x(z)\vspace{2ex}\\
	&= \displaystyle\frac{1}{h_1(z)} \xi_0(z) \we v_1(z)\we x(z) - \displaystyle\frac{1}{h_1(z)} \xi_0(z) \we \xi_0(z) \xi_0(z)^* v_1(z) \we x(z) \vspace{2ex} \\
	&=  \left(\displaystyle\frac{1}{h_1} \xi_0 \telwe v_1 \telwe x\right)(z). \end{array}$$ Note that $v_1\in H^2(\Dd,\Cc^n),$ $\xi_0\in H^\infty(\Dd,\Cc^n)$ and $h_1 \in H^2(\Dd,\Cc)$ is the scalar outer factor of $\xi_0 \telwe v_1$. By Proposition \ref{wejanalytic}, for every $x \in H^2(\Dd,\Cc^n),$ 
	$$\frac{1}{h_1} \xi_0 \telwe v_1 \telwe x $$ 
is analytic on $\Dd.$ By Proposition \ref{wel2conv}, since $\xi_0$ and $\xi_1$ are pointwise orthogonal almost everywhere on $\Tt,$
	$$\| \xi_0 \telwe \xi_1 \telwe x\|_{L^2(\Tt,\we^3\Cc^n)} < \infty. $$ Hence, 
	$$\xi_0 \telwe \xi_1 \telwe x \in \frac{1}{h_1} \xi_0 \telwe v_1 \telwe H^2(\Dd,\Cc^n) \subset H^2(\Dd,\we^3\Cc^n) .$$
	
\textbf{Recursive step:} suppose we have constructed vector-valued functions $\xi_0, \dots , \xi_{j-1},$ $\eta_0 , \dots , \eta_{j-1},$ spaces $X_j,Y_j$ and a compact operator $T_j \colon X_j \to Y_j$ after applying steps $0,\dots,j$ of the algorithm from Subsection \ref{Alg_statement} satisfying 
\begin{equation}\label{3.2.3sub}
\xi_0 \telwe \cdots \telwe \xi_{j-1} \telwe H^2(\Dd,\Cc^n) \subset H^2(\Dd,\we^{j+1}\Cc^n).\end{equation} 
Since $T_j$ is a compact operator, there exist vector-valued functions 
$v_{j}\in H^2(\Dd,\Cc^n), w_{j} \in H^2(\Dd,\Cc^m)^\perp$ such that 
$$(\xi_0 \telwe \cdots \telwe \xi_{j-1} \telwe v_{j}, \overline{\eta}_0 \telwe \dots \telwe \overline{\eta}_{j-1}\telwe w_{j}) $$ is a Schmidt pair for $T_j$ corresponding to $\|T_{j}\|.$
	Define \begin{equation}
	\label{xj321}
	 x_{j} = (I_{n} - \xi_0 \xi_0^* - \dots - \xi_{j-1}\xi_{j-1}^*)v_{j}.
\end{equation} 
By assumption, $\xi_0 \telwe \cdots \telwe \xi_{j-1} \telwe v_{j} $ lies in $H^2(\Dd,\we^{j+1} \Cc^n).$ Let $h_j \in H^2(\Dd,\Cc)$ be the scalar outer factor of  $\xi_0 \telwe \cdots \telwe \xi_{j-1} \telwe v_{j}.$ Define $\xi_{j} = \frac{x_j}{h_j}.$ Note that 
$\xi_i(z)$ and  $z \mapsto \xi_i(z)\langle v_j(z), \xi_i(z)\rangle_{\Cc^n}$
 are pointwise linearly dependent on $\Dd$ for  $i=0,\dots,j-1.$ Thus, for all $x\in H^2(\Dd,\Cc^n)$ and all $z \in \Dd,$ 
	
	\begin{align}
	(\xi_0 \telwe \cdots \telwe \xi_{j-1} \telwe \xi_j \telwe x)(z) &=  (\xi_0 \telwe \cdots \telwe \xi_{j-1} \telwe \frac{x_j}{h_j}\telwe x)(z)\nonumber\vspace{2ex}\\
	&=  \xi_0(z) \we \cdots \we \xi_{j-1}(z) \we \frac{1}{h_j(z)} \bigg( v_{j}(z)- \xi_0(z)\xi_0^*(z)v_j(z) -\cdots \bigg.\nonumber\\ &\hspace{30ex}- \bigg. \xi_{j-1}(z)\xi_{j-1}^*(z) v_{j}(z) \bigg) \we x(z)\nonumber\vspace{2ex}\\
	&=  \xi_0(z) \we \cdots \we \xi_{j-1}(z) \we \frac{1}{h_j(z)} v_{j}(z) \we x(z) \nonumber\\ &- \sum\limits_{i=0}^{j-1}\xi_0(z)\we \dots\we \xi_{j-1}(z) \we \xi_i(z)\xi_i^*(z)v_j(z)\we x(z)\nonumber \vspace{2ex}\\
	&= \left(\frac{1}{h_j}\xi_0 \telwe \cdots \telwe \xi_{j-1} \telwe v_{j} \telwe x\right)(z)\label{1hjvj}.\end{align} Recall that, for $i=0,\dots,j-1,$ by the algorithm from Subsection $3.2.1,$ $$x_i = (I_{n} -\xi_0\xi_0^* - \dots-\xi_{i-1}\xi_{i-1}^*)v_i  $$ and 
	$$ \xi_i = \frac{x_i}{h_i}.$$ By equation \eqref{1hjvj}, for all $z \in \Dd,$ 
	
$$(\xi_0 \telwe \cdots \telwe \xi_{j-1} \telwe \xi_j \telwe x)(z) =\left(\frac{1}{h_j}\xi_0 \telwe \cdots \telwe \xi_{j-1} \telwe v_{j} \telwe x\right)(z) . $$ Substituting $\frac{x_i}{h_i}$ for $\xi_i$ in the latter equation, where $x_i$ are given by equation \eqref{xj321} for  $i=1,\dots,j-1,$  we obtain 	
	
$$(\xi_0 \telwe \cdots \telwe \xi_{j-1} \telwe \xi_j \telwe x)(z)=\left(\frac{1}{h_1} \frac{1}{h_2} \cdots \frac{1}{h_j}\xi_0 \telwe v_1 \telwe \cdots \telwe v_{j-1} \telwe v_{j} \telwe x\right)(z),\; z\in\Dd	$$
on account of the pointwise linear dependence of $\xi_k$ and
$z \mapsto \langle v_k(z),\xi_i(z) \rangle_{\Cc^n} \xi_k(z) $ 
 on $\Dd,$ for $k=0,\dots,j.$ By Proposition \ref{wejanalytic}, for every $x \in H^2(\Dd,\Cc^n),$ $$\frac{1}{h_1} \frac{1}{h_2} \cdots \frac{1}{h_j} \xi_0 \telwe v_1 \telwe\cdots \telwe v_j \telwe x$$ is analytic on $\Dd.$ By Proposition \ref{wel2conv}, since $\xi_0,\xi_1,\dots,\xi_j$ are pointwise orthogonal on $\Tt,$
	$$\| \xi_0 \telwe \xi_1 \telwe \cdots \xi_j \telwe x\|_{L^2(\Tt,\we^{j+2}\Cc^n)} < \infty. $$ Thus, for every $x \in H^2(\Dd,\Cc^n),$
	$$\xi_0 \telwe \xi_1 \telwe\cdots \telwe \xi_j \telwe x \in H^2(\Dd,\we^{j+2}\Cc^n) $$ and the claim has been proved.

\end{proof}
	
	\begin{proposition}\label{xjclosed}
		In the notation of Proposition \ref{xjsubseth2}, $$ \xi_0 \telwe \dots \telwe \xi_j \telwe H^2(\Dd,\Cc^n)$$
		is a closed subspace of $H^2(\Dd,\we^{j+2}\Cc^n).$
	\end{proposition}
\begin{proof}
 Let us first show that $\xi_0 \telwe H^2(\Dd,\Cc^n)$ is a closed subspace of $H^2(\Dd,\we^2\Cc^n).$ Observe that, by Proposition \ref{xweyh2}, $\xi_0 \telwe H^2(\Dd,\Cc^n)\subset H^2(\Dd,\we^2\Cc^n).$ Let $$\Xi_0 =\{ f \in H^2(\Dd,\Cc^n): \langle f(z), \xi_0 (z) \rangle_{\Cc^n}=0\quad \text{almost everywhere on} \; \Tt \}.$$ Consider a vector-valued function $w\in H^2(\Dd,\Cc^n).$ For all $z\in \Dd,$ we may write $w$ as 
	$$w(z)= w(z) -\langle w(z),\xi_0(z)\rangle_{\Cc^n}\xi_0(z)+\langle w(z),\xi_0(z)\rangle_{\Cc^n}\xi_0(z).$$ Then, for all $w\in H^2(\Dd,\Cc^n)$ and for all $z\in \Dd,$
	
	$$\begin{array}{clllll}(\xi_0 \telwe w)(z)& = \xi_0(z)\we \big(w(z) -\langle w(z),\xi_0(z)\rangle_{\Cc^n}\xi_0(z)+\langle w(z),\xi_0(z)\rangle_{\Cc^n}\xi_0(z)\big)\vspace{2ex}\\&=\xi_0(z)\we \big(w(z) -\langle w(z),\xi_0(z)\rangle_{\Cc^n}\xi_0(z) \big)  \end{array}$$
on account of the  the pointwise linear dependence of $\xi_0$ and 
	$z \mapsto \langle w(z),\xi_0(z) \rangle_{\Cc^n} \xi_0(z) $  on $\Dd.$ 

Note that $$w(z)- \langle w(z),\xi_0(z)\rangle_{\Cc^n}\xi_0(z)\in \Xi_0,$$ thus 
	$$\xi_0 \telwe H^2(\Dd,\Cc^n) \subset \xi_0 \telwe \Xi_0. $$ 
By Corollary \ref{vclosed}, $\Xi_0$ is a closed subspace of $H^2(\Dd,\Cc^n),$ hence
	$$ \xi_0 \telwe H^2(\Dd,\Cc^n) \supset \xi_0 \telwe \Xi_0,$$and so,
	$$\xi_0 \telwe H^2(\Dd,\Cc^n) = \xi_0 \telwe \Xi_0 .$$
	Consider the mapping 
	$$C_{\xi_0}\colon \Xi_0 \to \xi_0 \telwe \Xi_0 $$
given by
	$$C_{\xi_0} w = \xi_0 \telwe w $$for all $w\in \Xi_0.$  Notice that, by Proposition \ref{onxi}, $\|\xi_0(e^{i\theta})\|_{\Cc^n}^2=1$ for almost every $\eiu \in  \Tt.$ Therefore, for any $w \in \Xi_0,$ we have 
	$$\begin{array}{cllllll}
	\|\xi_0 \telwe w\|_{L^2(\Tt,\we^{2}\Cc^n)}^2 &=\displaystyle \frac{1}{2\pi} \int\limits_0^{2\pi} \langle \xi_0 \telwe w, \xi_0 \telwe w \rangle (e^{i\theta})d\theta \vspace{3ex} \\
	&=  \displaystyle \frac{1}{2\pi} \int\limits_0^{2\pi}\left( \|\xi_0(e^{i\theta})\|_{\Cc^n}^2 \|w(e^{i\theta})\|_{\Cc^n}^2 - |\langle w(e^{i\theta}), \xi_0(e^{i\theta})\rangle|^2\right)\; d\theta \vspace{3ex}\\
	&= \|w\|_{L^2(\Tt, \Cc^n)}^2, 
	\end{array}$$
since	$w$ is pointwise orthogonal to $\xi_0$ almost everywhere on $\Tt.$  Thus the mapping $$C_{\xi_0}\colon  \Xi_0 \to \xi_0 \telwe \Xi_0$$ is an isometry. Furthermore, $C_{\xi_0}\colon  \Xi_0 \to \xi_0 \telwe \Xi_0$ is a surjective mapping, 
	thus $\Xi_0$ and $\xi_0 \telwe \Xi_0$ are isometrically isomorphic. Therefore, since $\Xi_0$ is a closed subspace of $H^2 (\Dd, \Cc^n),$ the space $\xi_0 \telwe \Xi_0$ is complete, hence is a closed subspace of $H^2(\Dd,\we^2\Cc^n)$. Hence $\xi_0 \telwe H^2(\Dd,\Cc^n)$ is a closed subspace of $H^2(\Dd,\we^2\Cc^n).$
	
	To prove that 
	$\xi_0 \telwe \dots \telwe \xi_j \telwe H^2(\Dd,\Cc^n) $ is a closed subspace of $H^2(\Dd,\we^{j+2}\Cc^n)$, let us consider
	$$\Xi_j =\{f \in H^2(\Dd,\Cc^n):\langle f(z), \xi_i (z)\rangle_{\Cc^n} =0, \; \text{for}\; i=0,\cdots,j \} $$
which is the pointwise orthogonal complement of $\xi_0, \dots,\xi_j$ in $H^2(\Dd,\Cc^n).$ Let $\psi \in H^2(\Dd,\Cc^n).$ We may write $\psi$ as 
	$$\psi(z) = \psi(z) - \sum\limits_{i=0}^j \langle \psi(z), \xi_i(z)\rangle_{\Cc^n}\xi_i(z) +\sum\limits_{i=0}^j \langle \psi(z), \xi_i(z)\rangle_{\Cc^n}\xi_i(z).  $$
	Then, for all $\psi \in H^2(\Dd,\Cc^n)$ and for almost all $z \in \Tt,$ 
	
	$$(\xi_0 \telwe \cdots \telwe \xi_j \telwe \psi)(z) = \xi_0(z)\we\cdots\we \left( \psi(z)- \sum\limits_{i=0}^j \langle \psi(z), \xi_i(z)\rangle_{\Cc^n}\xi_i(z)  \right)$$
due to the pointwise linear dependence of $\xi_k$ and 
	$z \mapsto \xi_k(z) \langle \psi(z), \xi_k(z)\rangle_{\Cc^n} $ almost everywhere 
on $\Tt.$
	
	\noindent Notice that  $\left( \psi(z)- \sum\limits_{i=0}^j \langle \psi(z), \xi_i(z)\rangle_{\Cc^n}\xi_i(z)  \right)$ is in $\Xi_j,$ thus
	$$\xi_0 \telwe \cdots \telwe \xi_j \telwe H^2(\Dd,\Cc^n) \subset \xi_0 \telwe \cdots \telwe \xi_j \telwe \Xi_j. $$The reverse inclusion holds by the definition of $\Xi_j,$ hence 
	
	$$\xi_0 \telwe \cdots \telwe \xi_j \telwe H^2(\Dd,\Cc^n) = \xi_0 \telwe \cdots \telwe \xi_j \telwe \Xi_j.$$
	
	Consequently, in order to prove the proposition it suffices to show that $\xi_0 \telwe \cdots \telwe \xi_j \telwe \Xi_j $ is a closed subspace of $H^2(\Dd, \we^{j+2} \Cc^n).$
	By Corollary \ref{vclosed}, $\Xi_j$ is a closed subspace of $H^2(\Dd, \Cc^n),$ being an  intersection of closed subspaces. For any $f \in \Xi_j,$   	$$ \begin{array}{clllll}	
	&\|\xi_0 \telwe \xi_1 \telwe\cdots\telwe \xi_j \telwe f\|_{L^2 ( \Tt, \we^{j+2} \Cc^n)}^2 \vspace{2ex} \\
	&=\displaystyle\frac{1}{2\pi} 
	\int_0^{2\pi} \det \begin{pmatrix}
\|\xi_0 (e^{i\theta})\|_{\Cc^n}^2 & \cdots  & \cdots & \langle \xi_0 (e^{i\theta}) , f(e^{i\theta})\rangle_{\Cc^n}\\
\langle \xi_1 (e^{i\theta}) , \xi_0 (e^{i\theta}) \rangle_{\Cc^n} &  \|\xi_1 (e^{i\theta})\|_{\Cc^n}^2 & \cdots & \langle \xi_1 (e^{i\theta}) , f(e^{i\theta})\rangle_{\Cc^n}\\
	\vdots  & \vdots   &      \ddots & \vdots\\
	\langle f(e^{i\theta}) , \xi_0 (e^{i\theta})\rangle_{\Cc^n} & \cdots & \cdots &\|f(e^{i\theta})\|_{\Cc^n}^2   \end{pmatrix}
	d\theta.\end{array}$$ 
Note that $f$ and $\xi_i$ are pointwise orthogonal almost everywhere on $\Tt,$ and, by Proposition \ref{onxi}, $\{\xi_0(z), \dots, \xi_j(z)\}$ is an orthonormal set for almost every $z \in \Tt.$
	Hence   
	$$\begin{array}{clll}
	&\|\xi_0 \telwe \xi_1 \telwe\cdots\telwe \xi_j \telwe f\|_{L^2 ( \Tt, \we^{j+2} \Cc^n)}^2  \vspace{2ex}\\
	&=  \displaystyle\frac{1}{2\pi} 
	\int_0^{2\pi} \det \begin{pmatrix}
	1&    0&   \cdots    &       0\\
	0&    1&    \cdots   &       0\\
  \vdots &  \vdots&   \ddots   &  \vdots \\
	0 &   0   & \cdots  &  \|f(e^{i\theta})\|_{\Cc^n}^2 \end{pmatrix}d\theta \vspace{2ex}\\
	&= \|f\|_{L^2 ( \Tt, \Cc^n)}^2.\end{array}$$
	Thus $$\xi_0 \telwe \xi_1 \telwe\cdots\telwe \xi_j \telwe \cdot \colon \Xi_j \to  \xi_0 \telwe \xi_1 \telwe\cdots\telwe \xi_j \telwe \Xi_j $$ 
	is an isometry. 
	Furthermore  $$(\xi_0 \telwe \xi_1 \telwe\cdots\telwe \xi_j \telwe \cdot)\colon \Xi_j \to \xi_0 \telwe \xi_1 \telwe\cdots\telwe \xi_j \telwe \Xi_j $$ is a surjective mapping, thus $\Xi_j$ and $\xi_0 \telwe \cdots \telwe \xi_j \telwe \Xi_j$ are isometrically isomorphic. Therefore, since $\Xi_j$ is a closed subspace of $H^2(\Dd,\Cc^n),$ the space $\xi_0 \telwe \cdots \telwe \xi_j \telwe \Xi_j $ is a closed subspace of $H^2(\Dd,\we^{j+2}\Cc^n).$ Hence 
	$$\xi_0 \telwe \cdots \telwe \xi_j \telwe H^2(\Dd,\Cc^n) $$ is a closed subspace of $H^2(\Dd,\we^{j+2}\Cc^n).$
\end{proof}

\section{The closed subspace $Y_{j+1}$ of $H^2(\Dd,\we^{j+2}\Cc^m)^\perp$}
\label{Y_j-closed}

\begin{proposition}\label{clwe}
Given $\bar{\eta}_0 = \frac{zy_0}{\overline{h_0}}$ as constructed in the algorithm in Subsection \ref{Alg_statement} the space
$\bar{\eta}_0 \telwe H^2(\Dd, \Cc^m)^\perp$ is a closed subspace of $H^2(\Dd, \we^2 \Cc^m)^\perp.$

\end{proposition}

\begin{proof}
As in Proposition \ref{xjsubseth2}, one can show that $$\eta_0\telwe zH^2(\Dd,\Cc^m) \subset zH^2(\Dd,\we^2\Cc^m) $$and therefore 

$$\bar{ \eta}_0 \telwe \bar{z} \overline{ H^2(\Dd,\Cc^m)} \subset \bar{z} \overline{ H^2(\Dd,\we^2\Cc^m)}. $$Hence 
$$\bar{ \eta}_0 \telwe H^2(\Dd,\Cc^m)^\perp \subset H^2(\Dd,\we^2\Cc^m)^\perp. $$
By virtue of the fact that complex conjugation is a unitary operator on $L^2(\Tt,\Cc^m),$ an equivalent statement to Proposition \ref{clwe} is that $\eta_0 \telwe zH^2(\Dd,\Cc^m)$ is a closed subspace of $zH^2(\Dd,\we^2 \Cc^m).$  Let 
$$
V= \{ f \in zH^2 (\Dd,\Cc^m) \; : \; \langle f(z) , \eta_0 (z) \rangle_{\Cc^m} =0\quad \text{for almost all}\; z \in \Tt \}
$$ 
be the pointwise orthogonal complement of $\eta_0$ in $zH^2(\Dd,\Cc^m).$

\noindent Consider $g \in zH^2 (\Dd,\Cc^m).$ We may write $g$ as 
$$g(z) = g(z) - \langle g(z) , \eta_0 (z) \rangle_{\Cc^m}\cdot \eta_0 (z) +  \langle g(z) , \eta_0 (z) \rangle_{\Cc^m} \cdot\eta_0 (z)$$ for every $z \in \Dd.$ Then, for all $g \in zH^2(\Dd, \Cc^m)$ and for all $z \in \Dd,$

$$(\eta_0  \telwe g) (z) = \eta_0 (z) \we [g(z) - \langle g(z) , \eta_0 (z) \rangle_{\Cc^m}  \eta_0(z)]$$
on account of the pointwise linear dependence of $\eta_0$ and $ z \mapsto  \langle g(z) , \eta_0(z) \rangle_{\Cc^m}  \eta_0(z)$ on $\Dd.$

\noindent Note that $g(z)-\langle g(z) , \eta_0 (z) \rangle_{\Cc^m}\eta_0 (z) \in V,$
thus $$\eta_0 \telwe zH^2(\Dd, \Cc^m) \subset \eta_0 \telwe V.$$
The reverse inclusion is obvious, hence $$\eta_0 \telwe zH^2(\Dd, \Cc^m) = \eta_0 \telwe V.$$ 

\noindent To prove the proposition, it suffices to show that  $\eta_0 \telwe V$ is a closed subspace of $zH^2 ( \Dd, \we^2 \Cc^m).$

\noindent Consider the mapping $$C_{\eta_0}\colon V \to \eta_0 \telwe V$$ defined by 
$$ C_{\eta_0} \nu = \eta_0 \telwe \nu$$ for all $\nu \in V.$ Notice that, by Proposition \ref{onxi}, $\|\eta_0(e^{i\theta})\|_{\Cc^m}^2=1$ for almost every $\eiu \in  \Tt.$ Then, for any $\upsilon \in V,$ we have 
$$\begin{array}{cllllll}
\|\eta_0 \telwe \upsilon\|_{L^2(\Tt,\we^{2}\Cc^m)}^2 &=\displaystyle \frac{1}{2\pi} \int\limits_0^{2\pi} \langle \eta_0 \telwe \upsilon, \eta_0 \telwe \upsilon \rangle (e^{i\theta})d\theta \vspace{3ex} \\
&=  \displaystyle \frac{1}{2\pi} \int\limits_0^{2\pi}\left( \|\eta_0(e^{i\theta})\|_{\Cc^m}^2 \|\upsilon(e^{i\theta})\|_{\Cc^m}^2 - |\langle \upsilon(e^{i\theta}), \eta_0(e^{i\theta})\rangle|^2\right)\; d\theta \vspace{3ex}\\
&= \|\upsilon\|_{L^2(\Tt, \Cc^m)}^2, 
\end{array}$$ since
$\upsilon$ is pointwise orthogonal to $\eta_0$ almost everywhere on $\Tt.$ Thus the mapping $C_{\eta_0}\colon  V \to \eta_0 \telwe V$ is an isometry.

\noindent  Note that by Corollary \ref{vclosed}, $V$ is a closed subspace of $zH^2 (\Dd, \Cc^m).$ Furthermore, $$C_{\eta_0}\colon  V \to \eta_0 \telwe V$$ is a surjective mapping, 
thus $V$ and $\eta_0 \telwe V$ are isometrically isomorphic. Therefore, since $V$ is a closed subspace of $zH^2(\Dd,\Cc^m),$ the space $\eta_0 \telwe V$ is complete and therefore a closed subspace of $zH^2(\Dd,\we^2\Cc^m)$. Hence $\bar{\eta}_0 \telwe H^2(\Dd,\Cc^m)^\perp$ is complete and therefore a closed subspace of $H^2(\Dd,\we^2\Cc^m)^\perp.$

\end{proof}

\begin{corollary}\label{projwell}
The orthogonal projection $P_{Y_1}$ from $L^2(\Tt, \we^2\Cc^m)$ onto $\bar{\eta}_0 \telwe H^2(\Dd, \Cc^m)^\perp $ is well defined.
\end{corollary}

\begin{proof}
By Proposition \ref{H2subsetL2}, $H^2(\Dd,\we^2\Cc^m)$ can be identified with a closed subspace of $L^2(\Tt,\we^2\Cc^m),$ thus we have $$ H^2(\Dd,\we^2\Cc^m)^\perp= L^2(\Tt,\we^2\Cc^m)\ominus H^2(\Dd,\we^2\Cc^m).$$ Now the assertion follows immediately from Proposition \ref{clwe}. 
\end{proof}

\begin{proposition}\label{clwegen}
Let $0\leq j \leq m-2.$ Let the functions $\bar{\eta}_i$ be given by equations \eqref{xij+1etaj+1} in the algorithm from Subsection \ref{Alg_statement}, that is, $\bar{\eta}_i= \displaystyle \frac{zy_i}{\overline{h}_i}$ for  $i =0, \cdots, j.$ Then, the space
$$\bar{\eta}_0 \telwe \bar{\eta}_1\telwe\cdots\telwe  \bar{\eta}_j \telwe  H^2(\Dd, \Cc^m)^\perp$$ is a closed linear subspace of $H^2(\Dd,\we^{j+2}\Cc^m)^\perp.$

\end{proposition}

\begin{proof}
First let us show that, for every $x\in H^2(\Dd,\Cc^m),$ 
$$\eta_0 \telwe \eta_1 \telwe \cdots \telwe \eta_j \telwe z x \in zH^2(\Dd,\we^{j+2}\Cc^m). $$
Recall that
$$ y_{j} = (I_{m} -\bar{\eta}_0\eta_0^T-\dots - \bar{\eta}_{j-1}\eta_{j-1}^T )w_{j}$$and 
\begin{equation}\label{etajwj1} \eta_0 \telwe \cdots \telwe \eta_{j-1} \telwe \bar{z}\bar{y}_{j}  =\eta_0 \telwe \cdots \telwe \eta_{j-1} \telwe (\bar{z}\bar{w}_{j} - \sum\limits_{i=0}^{j-1}{\eta}_i\eta_i^* \bar{z} \bar{w}_{j} )=  \eta_0 \telwe \cdots \telwe \eta_{j-1} \telwe \bar{z} \bar{w}_{j} \end{equation} because of the pointwise linear dependence of $\eta_i$ and 
$ z \mapsto  \langle \bar{z} \bar{w}_{j+1}(z), \eta_i(z) \rangle_{\Cc^m}  \eta_i(z)$ on $\Dd.$

\noindent By Proposition \ref{xjwevjetajwewj},
$$|h_i(z)| = \|y_i(z)\|_{\Cc^m} $$almost everywhere on $\Tt.$ 

\noindent Substituting $\eta_i = \frac{\bar{z}\bar{y}_i}{h_i} $ for all $i=0,\dots,j-1$ in equation \eqref{etajwj1}, we obtain  $$\eta_0 \telwe \cdots \telwe \eta_{j-1} \telwe \bar{z}\bar{y}_{j}= \frac{1}{h_0} \frac{1}{h_1} \cdots \frac{1}{h_j} \eta_0 \telwe \bar{z}\bar{w_1} \telwe \cdots \telwe \bar{z}\bar{w_j}.$$ Observe that, by Proposition \ref{wejanalytic}, for every $x \in H^2(\Dd,\Cc^m),$
$$ \frac{1}{h_0} \frac{1}{h_1} \cdots \frac{1}{h_j} \eta_0 \telwe \bar{z}\bar{w_1} \telwe \cdots \telwe \bar{z}\bar{w_j}\telwe zx $$is analytic on $\Dd.$ By Proposition \ref{wel2conv}, for all $x \in H^2(\Dd,\Cc^m),$ since $\eta_0, \cdots,\eta_j$ are pointwise orthogonal on $\Tt,$ $$\|\eta_0 \telwe \eta_1 \telwe \cdots \telwe \eta_j \telwe zx  \|_{L^2(\Tt,\we^{j+2}\Cc^m)} < \infty. $$ Hence, for every $x\in H^2(\Dd,\Cc^m),$ $$\eta_0 \telwe \eta_1 \telwe \cdots \telwe \eta_j \telwe zx = z\frac{1}{h_0} \frac{1}{h_1} \cdots \frac{1}{h_j} \eta_0 \telwe \bar{z}\bar{w}_1 \telwe \cdots \telwe \bar{z}\bar{w}_j \telwe x$$ is in $ zH^2(\Dd,\we^{j+2}\Cc^m) .$

\noindent Taking complex conjugates, we infer that
$$Y_{j+1}\stackrel{\text{def}}{=} \bar{\eta}_0 \dot{\we}  \cdots \telwe \bar{\eta}_{j-1} \dot{\we} \bar{\eta}_j\dot{\we} H^2 (\Dd, \Cc^m)^\perp \subset  H^2(\Dd,\we^{j+2}\Cc^m)^\perp.$$

Let us prove that $Y_{j+1}$ is a closed linear subspace of $H^2(\Dd,\we^{j+2}\Cc^m)^\perp.$ Since complex conjugation is a unitary operator on $L^2(\Tt,\Cc^m),$ an equivalent statement to the above is that  
 $$\eta_0 \telwe \eta_1 \telwe\cdots\telwe \eta_j \telwe  zH^2(\Dd, \Cc^m) $$ is a closed linear subspace of $zH^2(\Dd, \we^{j+2} \Cc^m).$
 
\noindent Let $$V_j =\{ \varphi \in  zH^2(\Dd, \Cc^m) \; : \; \langle \varphi(z) , \eta_i (z) \rangle_{\Cc^m} = 0 , \;\text{for}\; i=0, \cdots, j\}$$ be the pointwise orthogonal complement of
 $\eta_0 , \cdots , \eta_j$ in $zH^2(\Dd, \Cc^m).$ Consider $f \in zH^2(\Dd, \Cc^m).$ We may write $f$ as 
 $$ f(z) = f(z) - \sum\limits_{i=0}^j \langle f(z) , \eta_i (z) \rangle\eta_i(z) + \sum\limits_{i=0}^j \langle f(z) , \eta_i (z) \rangle\eta_i(z). $$
  Then, for all $ f \in zH^2 (\Dd, \Cc^m)$ and for almost all $z \in \Tt,$ 
 $$ (\eta_0 \telwe \eta_1 \telwe\cdots\telwe \eta_j \telwe f ) (z) = \eta_0 (z) \we \eta_1 (z) \we \cdots \we \eta_j (z) \we \left( f(z) - \sum\limits_{i=0}^j \langle f(z) , \eta_i (z) \rangle\eta_i(z)\right).$$
 Notice that $\left(f(z) - \sum\limits_{i=0}^j \langle f(z) , \eta_i (z) \rangle\eta_i(z)\right) \in V_j, $ thus
 $$ \eta_0 \telwe \eta_1 \telwe\cdots\telwe \eta_j \telwe  zH^2(\Dd, \Cc^m) \subset  \eta_0 \telwe \eta_1 \telwe\cdots\telwe \eta_j \telwe V_j .$$ The reverse inclusion holds by the definition of $V_j,$ 
 hence $$\eta_0 \telwe \eta_1 \telwe\cdots\telwe \eta_j \telwe  zH^2(\Dd, \Cc^m) = \eta_0 \telwe \eta_1 \telwe\cdots\telwe \eta_j \telwe V_j .$$
 Consequently, in order to prove the proposition it suffices to show that $\eta_0 \telwe \eta_1 \telwe\cdots\telwe \eta_j \telwe V_j $ is a closed subspace of $zH^2(\Dd, \we^{j+2} \Cc^m).$
 By Corollary \ref{vclosed}, $V_j$ is a closed subspace of $zH^2(\Dd, \Cc^m),$ being an  intersection of closed subspaces. 
 For any $f \in V_j,$ we have  
$$ \begin{array}{clllll}
 
 &\|\eta_0 \telwe \eta_1 \telwe\cdots\telwe \eta_j \telwe f\|_{L^2 ( \Tt, \we^{j+2} \Cc^m)}^2 \vspace{2ex} \\
 &=\displaystyle\frac{1}{2\pi} 
 \int_0^{2\pi} \det \begin{pmatrix}
                                 \|\eta_0 (e^{i\theta})\|_{\Cc^m}^2 & \cdots  & \langle \eta_0 (e^{i\theta}) , f(e^{i\theta})\rangle_{\Cc^m}\\
                                 \langle \eta_1 (e^{i\theta}) , \eta_0 (e^{i\theta}) \rangle_{\Cc^m}&  \|\eta_1 (e^{i\theta})\|_{\Cc^m}^2 & \cdots\\
                                 \vdots  & & \hspace{-20ex} \ddots \\
                                 \langle f(e^{i\theta}) , \eta_0 (e^{i\theta})\rangle_{\Cc^m} & \cdots & \|f(e^{i\theta})\|_{\Cc^m}^2   \end{pmatrix}
d\theta.\end{array}$$ Note that $f$ and $\eta_i$ are pointwise orthogonal almost everywhere on $\Tt$ and, by Proposition \ref{onxi}, $\{\eta_0(z), \dots, \eta_j(z)\}$ is an orthonormal set for almost every $z \in \Tt.$
Hence   
$$\begin{array}{clll}
&\|\eta_0 \telwe \eta_1 \telwe\cdots\telwe \eta_j \telwe f\|_{L^2 ( \Tt, \we^{j+2} \Cc^m)}^2  \vspace{2ex}\\
&=  \displaystyle\frac{1}{2\pi} 
\int_0^{2\pi} \det \begin{pmatrix}
  1  &    0  & \cdots    & 0\\
  0  &    1  &  \cdots   & 0\\
\vdots & \vdots&  \ddots  &  \vdots\\
  0  &   0  &  \cdots & \|f(e^{i\theta})\|_{\Cc^m}^2 \end{pmatrix}d\theta \vspace{2ex}\\
&= \|f\|_{L^2 ( \Tt, \Cc^m)}^2.\end{array}$$
Thus $$\eta_0 \telwe \eta_1 \telwe\cdots\telwe \eta_j \telwe \cdot \colon  V_j \to  \eta_0 \telwe \eta_1 \telwe\cdots\telwe \eta_j \telwe V_j $$ 
 is an isometry. 
Furthermore  $$(\eta_0 \telwe \eta_1 \telwe\cdots\telwe \eta_j \telwe \cdot)\colon V_j \to \eta_0 \telwe \eta_1 \telwe\cdots\telwe \eta_j \telwe V_j $$ is a surjective mapping, thus $V_j$ and $\eta_0 \telwe \cdots \telwe \eta_j \telwe V_j$ are isometrically isomorphic. Therefore, since $V_j$ is a closed subspace of $zH^2(\Dd,\Cc^m),$ the space $\eta_0 \telwe \cdots \telwe \eta_j \telwe V_j $ is a closed subspace of $zH^2(\Dd,\we^{j+2}\Cc^m).$ Hence 
$$\bar{ \eta}_0\telwe \cdots \telwe \bar{ \eta}_j \telwe H^2(\Dd,\Cc^m)^\perp $$ is a closed subspace of $H^2(\Dd,\we^{j+2}\Cc^m)^\perp.$
\end{proof}

\begin{corollary}\label{projwellgen}
	Let $0\leq j \leq m -2.$
The orthogonal projection 
$$P_{Y_j}\colon  L^2(\Tt, \we^{j+2}\Cc^m)\to Y_j$$ is well-defined. \end{corollary}

\begin{proof}
By Proposition \ref{H2subsetL2}, $H^2(\Dd,\we^{j+2}\Cc^m)$ can be identified with a closed subspace of $L^2(\Tt,\we^{j+2}\Cc^m),$ thus we have $$ H^2(\Dd,\we^{j+2}\Cc^m)^\perp= L^2(\Tt,\we^{j+2}\Cc^m)\ominus H^2(\Dd,\we^{j+2}\Cc^m).$$ Now the assertion follows immediately from Proposition \ref{clwegen}.
\end{proof}

%% file: algorithm_proofs.tex
\section{$T_j$ is a well-defined operator }\label{Tj-well-def}

 \begin{proposition}\label{Twell} 
Let $G\in H^\infty(\Dd,\CCmn)+C(\Tt,\CCmn)$ and let $0\leq j \leq\min(m,n)-2.$ Let the functions $\xi_{i}, \eta_{i}$ be defined by equations \eqref{xij+1etaj+1}, that is,
 \begin{equation}\label{xizetaiz}\xi_i = \displaystyle\frac{x_i}{h_i} , \quad \eta_i=\displaystyle\frac{\bar{z} \bar{ \eta_i}}{h_i} \end{equation} for $i=0,\cdots,j$ and let  
 
 $$X_i = \xi_0 \telwe \xi_1 \telwe \cdots \telwe \xi_{i-1} \dot{\we} H^2(\Dd, \Cc^n) \subset H^2(\Dd, \we^{i+1} \Cc^n),\quad i=0,\cdots, j,$$
 $$Y_i =  \bar{\eta}_{0} \telwe \bar{\eta}_{1} \telwe \cdots \telwe  \bar{\eta}_{i-1} \dot{\we} H^2 (\Dd, \Cc^m )^\perp \subset H^2 (\Dd, \we^{i+1} \Cc^m)^\perp,\quad i=0,\cdots, j.$$ Let $Q_i \in H^\infty(\Dd,\CCmn)$ satisfy
 
 \begin{equation}\label{G-qii}(G-Q_i)x_{k} = t_{k} y_{k}, \quad (G-Q_i)^*y_{k} =t_{k} x_{k}\end{equation} for $k=0,\dots,i-1.$
 
\noindent  Then, the operators $T_i\colon  X_i \to Y_i,$ $i=0,\cdots,j,$ given by 
 \be\label{defTi}
T_i(\xi_0 \dot{\we} \xi_1 \dot{\we} \cdots \dot{\we} \xi_{i-1} \dot{\we} x)= P_{Y_i} \left( \bar{\eta}_0 \dot{\we} \bar{\eta}_1 \dot{\we} \cdots \dot{\we} \bar{\eta}_{i-1}\dot{\we} (G-Q_i)x \right)
\ee
\noindent are well-defined and are independent of the choice of $Q_i\in H^\infty(\Dd, \CCmn)$ satisfying equations \eqref{G-qii}.
    \end{proposition}
\begin{proof}
By Corollary \ref{projwellgen}, the projections $P_{Y_i}$ are well-defined for all $i=0, \cdots, j.$
Hence it suffices to show that, for all $i=0,1,\cdots,j,$ $T_i$ maps a zero from its domain to a zero in its range and that $T_i$ does not depend on the choice of $Q_i,$ which satisfies equations (\ref{G-qii}). 

 For $i=0,$ the operator $T_0$ is the Hankel operator $H_G.$ If $f_0 \equiv 0,$ then $H_G f_0 =0$ and, moreover, $H_G$ is independent of the choice of any $Q\in H^\infty(\Dd,\CCmn)$ as $H_{G-Q} = H_G.$ Thus, $T_0$ is well-defined. 

 For $i=1,$ let $(x_0,y_0)$ be a Schmidt pair for the compact operator $H_G$ corresponding to $t_0=\|H_G\|,$ where $x_0 \in H^2(\Dd,\Cc^n)$ and $y_0 \in H^2(\Dd,\Cc^m)^\perp.$
By Lemma \ref{2.2}, $x_0, \bar{z}\bar{y}_0$ admit the inner-outer factorisations $x_0= \xi_0 h_0, \quad \bar{z}\bar{y}_0= \eta_0 h_0,$ where $\xi_0 \in H^\infty(\Dd,\Cc^n),$ $\eta_0 \in H^\infty(\Dd,\Cc^m)$ are inner vector-valued functions and $h_0 \in H^2(\Dd,\Cc)$ is scalar outer. The spaces $X_1$ and $Y_1$ are given by the formulas $$X_1 = \xi_0 \telwe H^2(\Dd,\Cc^n),\quad Y_1=\bar{\eta}_0 \telwe H^2(\Dd,\Cc^m)^\perp .$$ The operator $T_1\colon X_1\to Y_1$ is given by $$T_1(\xi_0\telwe x) = P_{Y_1} (\bar{\eta}_0 \telwe (G-Q_1)x)$$ for all $x \in H^2(\Dd,\Cc^n),$ where $Q_1 \in H^\infty(\Dd,\CCmn)$ satisfies equations (\ref{G-qii}).

\begin{lemma}\label{0to0} Let $\xi_0 \telwe u =\xi_0 \telwe v$ for some $u,v\in H^2(\Dd,\Cc^n)$. Then 
 $$ \bar{\eta}_0 \telwe (G-Q_1)u = \bar{\eta}_0 \telwe (G-Q_1)v. $$
\end{lemma}
\begin{proof}
 Suppose that $\xi_0 \telwe u =\xi_0 \telwe v$ for some $u,v\in H^2(\Dd,\Cc^n)$.  Let $x=u-v$, then $\xi_0 \telwe x =0$, and so $x$ and $\xi_0$ are pointwise linearly dependent in $\Cc^n$ on $\Dd.$ 
Therefore there exist maps $\beta,\lambda \colon  \Dd \to \Cc$, having no common zero in $\Dd$,
such that
   
 \begin{equation}\label{lindep}\beta(z) \xi_0(z) = \lambda(z) x(z) \;\; \text{in}\; \Cc^n ,\end{equation} for all $z\in \Dd.$ By assumption, $Q_1 \in H^\infty(\Dd, \CCmn)$ satisfies equations (\ref{G-qii}). Thus, for all $z \in \Dd,$
\begin{equation}\label{t0y0} t_0 y_0(z) = (G-Q_1)(z) x_0 (z).\end{equation}

\noindent By equations (\ref{xizetaiz}) and (\ref{lindep}),  \begin{equation}\label{x0lindep} \beta(z) x_0(z)=  \beta(z) h_0(z)\xi_0(z)= h_0(z) \lambda (z) x(z)\end{equation} 
for all $z \in \Dd.$ By equations (\ref{t0y0}) and (\ref{x0lindep}),  for all $z \in \Dd,$

$$\begin{array}{cllll}
 t_0 y_0(z) &= (G-Q_1)(z) x_0 (z) ,\\
\beta(z) t_0 z\displaystyle\frac{y_0(z)}{\bar{h}_0(z)} &= 
\displaystyle  (G-Q_1)(z)h_0 (z) \lambda (z) x(z)\frac{z}{\bar{h}_0(z)}. \end{array}  $$

\noindent Therefore, by equations (\ref{xizetaiz}), for all $z\in \Dd,$
$$t_0 \beta(z) \bar{\eta}_0 (z) = (G-Q_1)(z) x(z) \mu(z)\; \; \text{in} \; \Cc^m,$$
where 
$$\mu(z)= \frac{zh_0(z) \lambda(z)}{\bar{h}_0(z)}, \; \text{for all }\; z\in \Dd.$$ 
Hence, by Definition \ref{pointwiseld}, $\bar{\eta}_0$ and $(G-Q_1)x$ are pointwise linearly dependent in $\Cc^m $ on $\Dd,$ and so 
$$ \bar{\eta}_0 \dot{\we} (G-Q_1)x =0.
$$ 
Consequently, 
\[
 \bar{\eta}_0 \dot{\we}(G-Q_1)u= \bar{\eta}_0 \dot{\we}(G-Q_1)v.
\]
\end{proof}
Therefore the formula \eqref{defTi} (with $i=1$) does uniquely define $T_1u \in Y_1$.
Next, we show that the operator $T_1$ is independent of the choice of $Q_1 \in H^\infty(\Dd,\CCmn)$  satisfying equations (\ref{G-qii}). 

Suppose $Q_1 , Q_2 \in H^\infty(\Dd, \CCmn)$ satisfy 
\begin{equation}\label{q_1}
(G-Q_1 ) x_0 = t_0 y_0 \;, \; y_0^* (G-Q_1 )  = t_0 x_0^* \end{equation}
and
 \begin{equation}\label{q_2}
(G-Q_2 ) x_0 = t_0 y_0 \;, \; y_0^* (G-Q_2 )  = t_0 x_0^* .\end{equation}

\noindent Then, we claim that, for all $x \in H^2(\Dd, \Cc^n),$  
$$P_{Y_1} ( \bar{\eta}_0 \telwe (G-Q_1 ) x ) = P_{Y_1}( \bar{\eta}_0 \telwe (G-Q_2 ) x ) ,
$$
that is, 
$$ 
P_{Y_1} ( \bar{\eta}_0 \telwe (Q_1 - Q_2 ) x )=0. 
$$
The latter equation is equivalent to the statement that, for all $x \in H^2(\Dd, \Cc^n)$,
 $\bar{\eta}_0 \telwe (Q_2 - Q_1)x$ is orthogonal to $Y_1$, that is, to $ \bar{\eta}_0 \telwe \varrho$  for all
 $ \varrho \in H^2(\Dd, \Cc^m)^\perp.$ As a matter of convenience, set 
 $$ 
A x = (Q_2 - Q_1)x, \quad x \in H^2(\Dd,\Cc^n).
$$
We have to prove that 
$$\langle \bar{\eta}_0 \telwe Ax , \bar{\eta}_0 \telwe \varrho \rangle_{L^2(\Tt, \we^2\Cc^m)}=0$$ 
for all $x \in H^2 (\Dd, \Cc^n)$ and all $\varrho \in H^2(\Cc^m)^\perp.$
Note that
 $$\begin{array}{clllllll}
 \langle \bar{\eta}_0 \telwe Ax , \bar{\eta}_0 \telwe \varrho\rangle_{L^2(\Tt, \we^2\Cc^m)} =\displaystyle\frac{1}{2\pi}\int\limits_0^{2\pi}\langle \bar{\eta}_0 (e^{i\theta}) \telwe A(e^{i\theta}) x(e^{i\theta}) , 
 \bar{\eta}_0(e^{i\theta}) \telwe \varrho(e^{i\theta})\rangle_{\we^2\Cc^m} \; d\theta,\end{array}$$which by Proposition \ref{we} yields
 
 $$\begin{array}{llll}&\displaystyle\frac{1}{2\pi}\int\limits_0^{2\pi} \det \begin{pmatrix}
                                                       \langle \bar{\eta}_0 (e^{i\theta}) ,\bar{\eta}_0  (e^{i\theta}) \rangle_ {\Cc^m}& \langle \bar{\eta}_0  (e^{i\theta}) , \varrho(e^{i\theta})\rangle_{\Cc^m}\\
                                                       \langle A(e^{i\theta})x(e^{i\theta}), \bar{\eta}_0 (e^{i\theta})\rangle_{\Cc^m}
                                                       &  \langle A(e^{i\theta}) x(e^{i\theta}), \varrho(e^{i\theta})\rangle_ {\Cc^m}\\
                                                       
                                                      \end{pmatrix}d\theta\vspace{2ex}\\                                                     
&=\displaystyle\frac{1}{2\pi}\int\limits_0^{2\pi} \|\bar{\eta}_0  (e^{i\theta})\|_{\Cc^m}^2  \langle A(e^{i\theta})x(e^{i\theta}), \varrho (e^{i\theta})\rangle_{\Cc^m}\;d\theta \vspace{2ex} \\
&\hspace{5ex}- \displaystyle\frac{1}{2\pi}\int\limits_0^{2\pi} \langle A(e^{i\theta})x(e^{i\theta}), \bar{\eta}_0 (e^{i\theta})\rangle_{\Cc^m} \langle \bar{\eta}_0  (e^{i\theta}) , \varrho(e^{i\theta})\rangle_{\Cc^m} \; d\theta .
 \end{array}$$
By Proposition \ref{onxi}, $\|\bar{\eta}_0 (e^{i\theta})\|_{\Cc^m} =1$ for almost every $\eiu \in \Tt.$ Since $Ax \in H^2 (\Dd, \Cc^m)$ and $\varrho \in  H^2 (\Dd, \Cc^m)^\perp,$ 

$$\displaystyle\frac{1}{2\pi}\int\limits_0^{2\pi}  \langle A(e^{i\theta})x(e^{i\theta}), \varrho (e^{i\theta}\rangle_{\Cc^m} d\theta = \langle Ax, \varrho\rangle_{L^2 (\Tt, \Cc^m)}=0.$$
 Thus
$$
\begin{array}{cllllll}
\langle \bar{\eta}_0 \telwe Ax , \bar{\eta}_0 \telwe \varrho\rangle_{L^2(\Tt, \we^2\Cc^m)}&= \displaystyle\frac{1}{2\pi}\int\limits_0^{2\pi}
\langle A(e^{i\theta})x(e^{i\theta}), \bar{\eta}_0  (e^{i\theta})\rangle_{\Cc^m} \langle \bar{\eta}_0  (e^{i\theta}) , \varrho(e^{i\theta})\rangle_{\Cc^m} d\theta \vspace{2ex} \\

&= \displaystyle\frac{1}{2\pi}\int\limits_0^{2\pi} \bar{\eta}_0^* (e^{i\theta})A(e^{i\theta})x(e^{i\theta}) \langle \bar{\eta}_0  (e^{i\theta}) , \varrho(e^{i\theta})\rangle_{\Cc^m} d\theta. \end{array}$$ 

\noindent Recall that by equation  (\ref{xi0eta0}), $\bar{\eta}_0(z) =\displaystyle \frac{zy_0 (z)}{\bar{h}_0(z)}, \; z \in \Tt,$ so that
 $$\bar{\eta}_0^*  (e^{i\theta})= \left(\frac{e^{i\theta} y_0(e^{i\theta}) }{\bar{h}_0 (e^{i\theta}) }\right)^* = \frac{e^{-i\theta} y_0^*(e^{i\theta}) }{h_0 (e^{i\theta}) }.$$

\noindent Therefore 
$$\langle \bar{\eta}_0 \telwe Ax , \bar{\eta}_0 \telwe \varrho\rangle_{L^2(\Tt, \we^2\Cc^m)}=\displaystyle\frac{1}{2\pi}\int\limits_0^{2\pi} \frac{e^{-i\theta} y_0^*(e^{i\theta}) }{h_0 (e^{i\theta}) } 
 A(e^{i\theta})x(e^{i\theta}) \langle \bar{\eta}_0 (e^{i\theta}) , \varrho(e^{i\theta})\rangle_{\Cc^m} d\theta .$$

\noindent Recall our 
initial assumption was that $Q_1 , Q_2$ satisfy equations (\ref{q_1}) and (\ref{q_2}), consequently,
$$ y_0^* (G- Q_i ) =t_0 x_0 ^* , \; \text{for} \; i=1,2.  $$
Hence, for $z \in \Tt,$ 

$$\begin{array}{clll}y_0^* (z) A(z)x(z) &= y_0^* (z) (G-Q_1)(z) x(z)- y_0^* (z) (G-Q_2) (z) x(z) \vspace{2ex} \\
&= (t_0 x_0^* x - t_0 x_0^* x) (z) \vspace{2ex} \\&=0. \end{array}$$
 \noindent We deduce that $$ \displaystyle\frac{1}{2\pi}\int\limits_0^{2\pi} \bar{\eta}_0^* (e^{i\theta})A(e^{i\theta})x(e^{i\theta}) \langle \bar{\eta}_0 (e^{i\theta}) , \varrho(e^{i\theta})\rangle_{\Cc^m} d\theta =0. $$

\noindent To conclude, we have proved that, if $Q_1 , Q_2 \in H^\infty( \Dd, \CCmn)$ satisfy equations (\ref{q_1}) and (\ref{q_2}), then 
 $$P_{Y_1}(\bar{\eta}_0\telwe (G-Q_1)x)= P_{Y_1}(\bar{\eta}_0\telwe (G-Q_2)x),$$ that is, $T_1$ is independent of the choice of $Q_1$ subject to equations \eqref{q_1},\eqref{q_2}.

 Recursive step: suppose that functions $x_{i-1} \in L^2(\Tt,\Cc^n),$  $y_{i-1} \in L^2(\Tt,\Cc^m),$ outer functions $h_{i-1} \in H^2(\Dd,\Cc),$ positive numbers $t_i,$ matrix-valued functions $Q_i \in H^\infty(\Dd,\CCmn),$ spaces $X_i, Y_i$ and compact operators $T_i \colon  X_i \to Y_i$ are constructed inductively by the algorithm for $i=0,\dots,j.$

  Let us prove that $T_j\colon X_j \to Y_j,$ given by equation (\ref{T_j}), is well-defined for $0\leq j \leq r$. Note, by Corollary \ref{projwellgen}, the projection $P_{Y_{j}}$ is well-defined. 


We must show that if an element of $X_j$  has two different expressions as an element of $\xi_0\telwe \dots\telwe \xi_{j-1}\telwe \xi_j \telwe H^2(\Dd,\Cc^n),$ say 

\begin{equation}\label{tworeps}
u= \xi_0\telwe \dots\telwe \xi_{j-1} \telwe x= 
\xi_0\telwe \dots\telwe \xi_{j-1} \telwe \tilde{x}
\end{equation}
for some $x,\tilde{x} \in H^2(\Dd,\Cc^n),$ then the two corresponding formulae for $T_j u $ given by the defining equation \eqref{T_j} agree, that is,
$$P_{Y_j}(\bita_0 \telwe \bita_1\telwe  \dots \telwe \bita_{j-1}\telwe (G-Q_j)x)=P_{Y_j}(\bita_0 \telwe \bita_1\telwe  \dots \telwe \bita_{j-1}\telwe (G-Q_j)\tilde{x}),$$ 
or equivalently, 
$$ P_{Y_j}\bigg(\bita_0 \telwe \bita_1\telwe  \dots \telwe \bita_{j-1}\telwe (G-Q_j)(x-\tilde{x})\bigg)=0,$$
which is to say that we need to show that
\begin{equation}\label{weneed}
\bita_0 \telwe \bita_1 \telwe \dots \telwe \bita_{j-1}\telwe (G-Q_j)(x-\tilde{x}) \in Y_j^\perp.
\end{equation}
If $x,\tilde{x}$ satisfy equation \eqref{tworeps}, then 
$$\xi_0\telwe \dots\telwe \xi_{j-1} \telwe (x-\tilde{x})=0, $$ and so, by Corollary \ref{lin-depend}, $\xi_0,\xi_1,\dots, \xi_{j-1},x-\tilde{x}$ are pointwise linearly dependent almost everywhere on $\Tt.$ 

\noindent It follows immediately that, for almost all $z\in \Tt,$ the vectors $$ (G-Q_j)\xi_0(z),\dots,  (G-Q_j)\xi_{j-1}(z), (G-Q_j)\xi_j(z),  (G-Q_j)(x-\tilde{x})(z)$$ are linearly dependent in $\Cc^m.$

\noindent Since $y_j = \bar{z}\bar{h}_j \bar{\eta}_j,$ by equations \eqref{xij+1etaj+11}, equations \eqref{G-qii} imply, for $i=0,\dots,j-1$ and almost all $z\in \Tt,$
$$\displaystyle (G-Q_j)\xi_i(z) = (G-Q_j)\frac{x_i}{h_i}(z)= t_i \frac{y_i}{h_i}(z) =\frac{t_i}{h_i} \bar{z}\bar{h}_i\bar{\eta}_i(z).$$
Thus, for almost all $z\in \Tt,$ the vectors 
$$\frac{t_0}{h_0} \bar{z}\bar{h}_0\bar{\eta}_0(z),\;\dots,\; \frac{t_{j-1}}{h_{j-1}} \bar{z}\bar{h}_{j-1}\bar{\eta}_{j-1}(z), \;(G-Q_j)(x-\tilde{x})(z) $$
are linearly dependent in $\Cc^m.$ Since $t_0 \geq t_1 \geq \dots \geq t_j >0$,
it follows that
$$\bita_0(z),\; \dots,\; \bita_{j-1}(z), \;(G-Q_j)(x-\tilde{x})(z) $$ are linearly dependent for almost all $z \in \Tt$ and so, by Corollary \ref{lin-depend}, 
$$\bita_0 \telwe \dots \telwe \bita_{j-1} \telwe (G-Q_j)(x-\tilde{x}) =0,$$
which certainly implies the desired relation \eqref{weneed}. 
Thus $T_j\colon X_j \to Y_j$ is well-defined.

 For the operator $T_j$ to be uniquely defined in the algorithm, it remains to prove 
$T_j$ is independent of the choice of $Q_j \in H^\infty ( \Dd, \CCmn)$ subject to equations (\ref{G-qii}). Let $Q_1, Q_2 \in H^\infty(\Dd,\CCmn)$ satisfy

\begin{equation}\label{qixiyiti}(G-Q_1)x_i = t_i y_i,\quad(G-Q_2)x_i = t_i y_i,\quad  
 y_i^* (G-Q_1) = t_i x_i^*,\quad y_i^* (G-Q_2) = t_i x_i^* 
\end{equation}
 for  $i=0, \cdots, j-1.$
We shall  prove that, for all $x \in H^2(\Dd, \Cc^n),$  
$$P_{Y_j} ( \bar{\eta}_0\telwe \cdots \telwe \bar{\eta}_{j-1} \telwe (G-Q_1)x ) = P_{Y_j}( \bar{\eta}_0\telwe \cdots \telwe \bar{\eta}_{j-1} \telwe (G-Q_2)x )  .$$
The latter equality holds if and only if, for all $x \in H^2(\Dd, \Cc^n),$  
$$P_{Y_j} ( \bar{\eta}_0\telwe \cdots \telwe \bar{\eta}_{j-1} \telwe (Q_2 -Q_1)x ) =0 $$which is equivalent to the assertion that 
 $\bar{\eta}_0 \telwe \cdots \telwe \bar{\eta}_{j-1} \telwe (Q_2 -Q_1)x$ is orthogonal to $ \bar{\eta}_0\telwe \cdots \telwe \bar{\eta}_{j-1}  \telwe q$ for all $x \in H^2(\Dd, \Cc^n)$ and for all
 $q \in H^2(\Dd, \Cc^m)^\perp.$

\noindent Equivalently
$$\langle \bar{\eta}_0 \telwe \cdots \telwe \bar{\eta}_{j-1} \telwe (Q_2 -Q_1)x ,\bar{\eta}_0\telwe \cdots \telwe \bar{\eta}_{j-1}  \telwe q \rangle_{L^2(\Tt, \we^{j+1}\Cc^m)}=0$$
for all $x \in H^2 (\Dd, \Cc^n)$ and for all $q \in H^2(\Cc^m)^\perp.$
Set $Ax= (Q_2-Q_1)x,$ $x\in H^2(\Dd,\Cc^n).$

\noindent By Proposition \ref{we},
 $$
 \langle \bar{\eta}_0\telwe \cdots \telwe \bar{\eta}_{j-1} \telwe (Q_2 -Q_1 )x ,\bar{\eta}_0\telwe \cdots \telwe \bar{\eta}_{j-1}  \telwe q\rangle_{L^2(\Tt, \we^{j+1}\Cc^m)}$$ 
 \noindent is equal to
 $$\frac{1}{2\pi}\int\limits_0^{2\pi}
\footnotesize{
\det \begin{pmatrix}
                                                   \langle \bar{\eta}_0 (e^{i\theta}) , \bar{\eta}_0(e^{i\theta})\rangle_{\Cc^m} & \cdots &  
                                                   \langle \bar{\eta}_0(e^{i\theta}), \bar{\eta}_{j-1}(e^{i\theta})\rangle_{\Cc^m} & \langle \bar{\eta}_0(e^{i\theta}), q(e^{i\theta}) \rangle_{\Cc^m} \\
                                                   
        \vdots & \ddots& \vdots & \vdots\\
                                                   
                                                  \langle \bar{\eta}_{j-1}(e^{i\theta}), \bar{\eta}_0(e^{i\theta})\rangle_{\Cc^m} & \cdots
                                                  & \langle \bar{\eta}_{j-1}(e^{i\theta}), \bar{\eta}_{j-1}(e^{i\theta})\rangle_{\Cc^m} & \langle \bar{\eta}_{j-1}(e^{i\theta}), q (e^{i\theta})\rangle_{\Cc^m} \\
                                                   \langle A(e^{i\theta})x(e^{i\theta}) , \bar{\eta}_0(e^{i\theta})\rangle_{\Cc^m} & \cdots 
                                                   & \langle A(e^{i\theta})x(e^{i\theta}) , \bar{\eta}_{j-1}(\eiu)\rangle_{\Cc^m} &\langle A(e^{i\theta})x(e^{i\theta}) , q(e^{i\theta})\rangle_{\Cc^m}\end{pmatrix} 
                                                   d\theta
}.$$
Notice that $Ax$ and $q$ are orthogonal in $L^2(\Tt, \Cc^m) $ and, by Proposition \ref{onxi}, $\{\bar{\eta}_i(z)\}_{i=0}^{j-1}$ is an orthonormal sequence in  $\Cc^m$ almost everywhere on $\Tt.$ Also, for  $i=0, \cdots, j-1,$ by equations (\ref{qixiyiti}),
  $$
 \begin{array}{cllll}
 \langle A(e^{i\theta})x(e^{i\theta}) , \overline{ \eta}_i(e^{i\theta})\rangle\rangle_{\Cc^m} 
  &=\displaystyle  \eta_i^T(e^{i\theta}) A(e^{i\theta})x(e^{i\theta}) \\&= 
  \displaystyle  \frac{e^{-i\theta} y_i^*(e^{i\theta}) }{h_i (e^{i\theta}) } A(e^{i\theta})x(e^{i\theta})\\
   &=\displaystyle \frac{e^{-i\theta}}{h_i (e^{i\theta}) }\left( y_i^*(e^{i\theta}) (G-Q_1)(z)x(z)-y_i^*(e^{i\theta})(G-Q_2)(z)x(z)\right)\\
   &=\displaystyle \frac{e^{-i\theta}}{h_i (e^{i\theta}) } (t_i x_i^*x-t_ix_i^*x)=0.  \end{array}$$
Thus 
$$\begin{array}{clllllll}
&\langle \bar{\eta}_0\telwe \cdots \telwe \bar{\eta}_j \telwe (Q_2 -Q_1 )x ,\bar{\eta}_0\telwe \cdots \telwe \bar{\eta}_j  \telwe q\rangle_{L^2(\Tt, \we^{j+1}\Cc^m)} \vspace{2ex} \\
&= \displaystyle\frac{1}{2\pi}\int\limits_0^{2\pi} \det \begin{pmatrix} 1 & 0 & \cdots & \langle \bar{\eta}_0(e^{i\theta}), q(e^{i\theta}) \rangle_{\Cc^m} \\ 
0 & 1 & \cdots&  \langle \bar{\eta}_2(e^{i\theta}), q (e^{i\theta})\rangle_{\Cc^m}\\
\vdots & & \ddots & \vdots\\
0 & 0 &\cdots & \langle A(e^{i\theta})x(e^{i\theta}) , q(e^{i\theta})\rangle_{\Cc^m}  \end{pmatrix}\; d\theta \vspace{2ex} \\
&= \langle Ax , q\rangle_{L^2 ( \Tt, \Cc^m)} \vspace{2ex} =0.
\end{array}$$

\noindent Consequently $$P_{Y_j} ( \bar{\eta}_0\telwe \cdots \telwe \bar{\eta}_{j-1} \telwe (G-Q_1)x ) = P_{Y_j}( \bar{\eta}_0\telwe \cdots \telwe \bar{\eta}_{j-1} \telwe (G-Q_2)x )  ,$$ and so 
$T_j$ is independent of the choice of $Q_j$ subject to equations \eqref{G-qii}. 
\end{proof}


%% file: T_j-compact.tex
\section{Compactness of the operators $T_1$ and $T_2$}\label{T_1-compact}

Here we use notations from the algorithm of Subsection \ref{Alg_statement} to prove the compactness of the operator $T_j$ given by equation (\ref{T_j}) for $j=1, 2$. 
 The proof requires several steps. Let us first prove that the operator $T_1$ is compact.

 Recall that since $G\in H^\infty(\Dd,\CCmn)+C(\Tt,\CCmn),$ by Hartman's theorem, the operator $T_0 = H_G$ is compact and hence there exist $x_0 \in H^2(\Dd,\Cc^n)$ and $y_0 \in H^2(\Dd,\Cc^m)^\perp$ such that $(x_0,y_0)$ is a Schmidt pair for $H_G$ corresponding to the singular value $\|H_G\|=t_0.$ 
 
 By Lemma \ref{2.2}, $x_0, \bar{z} \bar{y}_0$ admit the inner-outer factorizations \begin{equation}\label{x00y00} x_0 = \xi_0 h_0, \quad \bar{z} \bar{y}_0=\eta_0 h_0,\end{equation} where $\xi_0 \in H^\infty(\Dd,\Cc^n)$, $\eta_0 \in H^\infty(\Dd,\Cc^m)$ are vector-valued inner functions and $h_0\in H^2(\Dd,\Cc)$ is a scalar outer function. Moreover there exist unitary-valued functions of types $n\times n, m\times m$ respectively, of the form  \begin{equation}\label{v00w00} V_0= \begin{pmatrix} \xi_0 & \bar{\alpha}_0 \end{pmatrix},\quad W_0= \begin{pmatrix}\eta_0 & \bar{\beta}_0  \end{pmatrix}^T, \end{equation} where $\alpha_0,\beta_0$ are inner, co-outer, quasi-continuous functions of types $n\times (n-1)$, $m\times (m-1)$ respectively and all minors on the first columns of $V_0,W_0^T$ are in $H^\infty$. Furthermore every $Q_1 \in H^\infty(\Dd,\CCmn)$ which is at minimal distance from $G$ satisfies
 
 $$W_0 (G-Q_1) V_0 = \begin{pmatrix}
 t_0 u_0 & 0\\
 0 & F_1
 \end{pmatrix} $$ for some $$F_1\in H^\infty(\Dd,\Cc^{(m-1)\times(n-1)})+C(\Tt,\Cc^{(m-1)\times(n-1)}) $$ and some quasi-continuous function $u_0$ with $|u_0(z)|=1$ almost everywhere on $\Tt.$
 
 \noindent Recall that  $$X_1 = \xi_0 \telwe H^2(\Dd,\Cc^n),\quad Y_1 = \bar{\eta}_0\telwe H^2(\Dd,\Cc^m)^\perp $$ and $T_1 \colon X_1 \to Y_1 $ is given by  $$
 T_1(\xi_0 \telwe x) = P_{Y_1}[ \bar{\eta}_0 \telwe(G-Q_1)x] \quad \text{for all}\; x \in H^2(\Dd,\Cc^n).
 $$ 
 Our first endeavour in this section is to prove the following theorem. 
 
\begin{theorem}\label{T0compact}
Let 	
 	\begin{equation}\label{calk1l1}
 	\mathcal{K}_1 \stackrel{\emph{def}}{=} V_0 \begin{pmatrix}0 \\ H^2(\Dd,\Cc^{n-1}) \end{pmatrix}, \quad \mathcal{L}_1 \stackrel{\emph{def}}{=} W_0^* \begin{pmatrix}0 \\ H^2(\Dd,\Cc^{m-1})^\perp \end{pmatrix},\end{equation}and let the maps $$U_1\colon  H^2(\Dd,\Cc^{n-1}) \to \mathcal{K}_1,$$ $$U_2\colon   H^2(\Dd,\Cc^{m-1})^\perp \to \mathcal{L}_1$$ be given by
 $$ U_1 x = V_0 \begin{pmatrix}0 \\ x \end{pmatrix}, \quad U_2 y = W_0^* \begin{pmatrix}0 \\ y \end{pmatrix}  $$for all $x\in H^2(\Dd,\Cc^{n-1}),\; y \in H^2(\Dd,\Cc^{m-1})^\perp.$ Consider the operator $\Gamma_1 = P_{\mathcal{L}_1} M_{G-Q_1} |_{\mathcal{K}_1}.$ 
 	\noindent Then
 	\begin{enumerate}
 		\item[\emph{(i)}] The maps $U_1,U_2$ are unitaries.
 		\item[\emph{(ii)}] The maps $(\xi_0\telwe \cdot )\colon\mathcal{K}_1\to H^2(\Dd,\we^2\Cc^n)$ and $(\bar{ \eta}_0 \telwe \cdot)\colon\mathcal{L}_1 \to H^2(\Dd,\Cc^m)^\perp$ are unitaries.
 	\item[\emph{(iii)}]  The following 
	diagram is commutative: 
	
	\begin{equation}\label{commdiagr}
\begin{array}{clllll}	
H^2(\Dd,\Cc^{n-1}) &\xrightarrow{U_1} & \mathcal{K}_1 &\xrightarrow{\xi_0 \telwe \cdot}& \xi_0 \telwe H^2 (\Dd, \Cc^n)=X_1\\
	\Big\downarrow\rlap{$\scriptstyle H_{F_1}  $} & ~  &\Big\downarrow\rlap{$\scriptstyle \Gamma_1$}  &~&\hspace{3ex}\Big\downarrow\rlap{$\scriptstyle T_1$} \\
H^2(\Dd,\Cc^{m-1})^\perp &\xrightarrow{U_2}&	\mathcal{L}_1 &\xrightarrow{\bar{\eta}_0 \telwe \cdot } & \bar{\eta}_0 \telwe H^2 (\Dd, \Cc^m)^\perp =Y_1
	\end{array}.\end{equation}
	
\item[\emph{(iv)}] $T_1$ is a compact operator. 
\item[\emph{(v)}]$ \|T_1\| = \|\Gamma_1\| = t_1$.
	 	\end{enumerate}
	
\end{theorem}
\begin{proof}
Statement {\rm (i)} follows from Lemma \ref{3.2constr}. 
Statement {\rm (ii)} follows from Propositions \ref{xi0wek0subset} and \ref{eta0telweh2} below, which are consequences of the following lemmas.

\begin{lemma}\label{maxvect}
	In the notation of Theorem \ref{T0compact}, the Hankel operator $H_G$ has a maximizing vector $x_0$ of unit norm such that $\xi_0,$ which is defined by $\xi_0 = \frac{x_0}{h_0},$ is a co-outer function. 
\end{lemma}

\begin{proof}
	Choose any maximizing vector $x_0.$ By Lemma \ref{2.2}, $x_0$ has the inner-outer factorization $x_0 = \xi_0 h_0,$ where $h_0$ is a scalar outer factor. Then, the closure of $ \xi_0^T H^2(\Dd,\Cc^n),$ denoted by $\clos(\xi_0^T H^2(\Dd,\Cc^n)),$ is a closed shift-invariant subspace of $H^2(\Dd,\Cc),$ so, by Beurling's theorem, $$ \clos(\xi_0^T H^2(\Dd,\Cc^n)) = \phi H^2(\Dd,\Cc)$$ for some scalar inner function $\phi.$ Hence 
	$$\bar{\phi} \xi_0^T H^2(\Dd,\Cc^n) \subset H^2(\Dd,\Cc). $$ 
	Thus, if $\xi_0^T = (\xi_{01}, \cdots,\xi_{0n} ),$ we have $\bar{\phi}\xi_{0j} \in H^\infty(\Dd,\Cc)$ for $j=1,\cdots, n,$ and so, $$ \overline{ \phi } \xi_0 \in H^\infty(\Dd,\Cc^n).$$ 
	Hence $$ \overline{ \phi}x_0 = \overline{ \phi} \xi_0 h_0 \in H^2(\Dd,\Cc^n).$$ 
	
\noindent 	Let $Q$ be a best $H^\infty$ approximation to $G.$ Since $x_0$ is a maximizing vector for $H_G$, by Theorem \ref{1.7},
	$$ (G-Q)x_0 \in H^2(\Dd,\Cc^m)^\perp $$ and $$\|(G-Q)(z)x_0(z)\|_{\Cc^m} = \|H_G\| \|x_0(z)\|_{\Cc^n}$$ for almost all $z \in \Tt.$ Thus 
	$$ (G-Q)\overline{\phi} x_0 \in H^2(\Dd,\Cc^m)^\perp $$ and $$\|(G-Q)\overline{\phi}x_0(z) \|_{\Cc^m} = \|H_G\| \| \overline{ \phi}x_0(z)\|_{\Cc^n} $$ for almost all $z \in \Tt.$ 
	
	\noindent Hence $\overline{ \phi}x_0 \in H^2(\Dd,\Cc^n)$ is a maximizing vector for $H_G,$ and $\overline{\phi} x_0$ is co-outer. Then $\frac{\bar{ \phi}x_0}{\|x_0\|}$ is a co-outer maximizing vector of unit norm for $H_G.$ \end{proof}

\begin{lemma}\label{utxi=1}
Let $x_0$ be a co-outer maximizing vector of unit norm for $H_G,$ and let \newline $x_0=\xi_0 h_0$ be the inner-outer factorisation of $x_0.$ Then

{\em (i)} $\xi_0$ is a quasi-continuous function and 

{\em (ii)} there exists a function $A\in H^\infty(\Dd,\Cc^n)$ such that $$A^T \xi_0 =1. $$ 
\end{lemma}
\begin{proof}
	Let us first show that $$\xi_0 \in (H^\infty(\Dd,\Cc^n) + C(\Tt, \Cc^n)) \cap      \overline{H^\infty(\Dd,\Cc^n) + C(\Tt, \Cc^n)}.$$
	
	\noindent Let $Q$ be a best $H^\infty$ approximation to $G.$ Then, by Theorem \ref{1.7}, the function $Q$ satisfies the equation 
	$$(G-Q)^* y_0 = t_0 x_0. $$ Taking complex conjugates in equations \eqref{x00y00}, we have
	$$(G-Q)^T \bar{y}_0 = t_0 \overline{x_0}.  $$ Hence, for $z \in \Tt,$
	$$(G-Q)^T z h_0 \eta_{0}  = t_0 \overline{h_0} \overline{ \xi_0} ,$$ and therefore
	$$\displaystyle \frac{(G-Q)^T z h_0 \eta_0 }{t_0 \overline{h_0}} = \overline{\xi_0}.$$  
	
	\noindent	Recall, by equation \eqref{eq223} (with $\phi=1$), $u_0 =\frac{\bar{z}\bar{h}_0}{h_0}.$ By Lemma \ref{2.2}, $u_0 \in QC ,$ hence $\overline{u_0} \in H^\infty+C.$ Note $ \overline{u_0} = \frac{zh_0}{\overline{h_0}}, $ and hence $$\overline{\xi_0} = \displaystyle \frac{(G-Q)^T \overline{u_0} \eta_0 }{t_0 }. $$  
	
	\noindent	Since $H^\infty +C$ is an algebra and $(G-Q)^T, \; \eta_0 \in H^\infty +C,$ it follows that $\overline{\xi_0} \in H^\infty+C,$ thus $$ \xi_0 \in (H^\infty(\Dd,\Cc^n) + C(\Tt, \Cc^n)) \cap     \overline{H^\infty(\Dd,\Cc^n) + C(\Tt, \Cc^n)}.$$ 
	
	\noindent	The conclusion that there exists a function $A\in H^\infty(\Dd,\Cc^n)$ such that $A^T \xi_0 =1 $ now follows directly from Lemma \ref{L6.2}.
\end{proof}

\begin{lemma}\label{a0h2}
In the notation of Theorem \ref{T0compact}, let $\xi_0 \in H^\infty(\Dd,\Cc^n)$ be a vector-valued inner, co-outer, quasi-continuous function and let $$V_0 = \begin{pmatrix}
\xi_0 & \bar{\alpha}_0 
\end{pmatrix}$$ be a thematic completion of $\xi_0$ as described in Lemma \ref{2.2}, where $\alpha_0$ is an inner, co-outer, quasi-continuous function of order $n\times (n-1) $ and all minors on the first column of $V_0$ are analytic. Then,
$$ \alpha_0^T H^2(\Dd,\Cc^n) = H^2(\Dd,\Cc^{n-1}). $$

\end{lemma}
\begin{proof} By Lemma \ref{L6.2}, for the given $\alpha_0,$ there exists $A_0\in H^\infty(\Dd, \Cc^{(n-1)\times n })$ such that $A_0\alpha_0 = I_{n-1}.$ Equivalently, $\alpha_0^T A_0^T = I_{n-1}.$ \medskip

\noindent Let $g\in H^2(\Dd,\Cc^{n-1}).$ Then $g = (\alpha_0^T A_0^T) g \in \alpha_0^T A_0^T H^2(\Dd,\Cc^{n-1}),$
which implies that $ g \in \alpha_0^T H^2(\Dd,\Cc^n).$ Hence $H^2(\Dd,\Cc^{n-1}) \subseteq \alpha_0^T H^2(\Dd,\Cc^n).$\medskip

\noindent For the reverse inclusion, note that since $\alpha_0$ is 
in $H^\infty(\Dd, \Cc^{n\times (n-1)})$, we have
$\alpha_0^T H^2(\Dd,\Cc^n) \subseteq  H^2(\Dd,\Cc^{n-1}).$ Thus $$ \alpha_0^T H^2(\Dd,\Cc^n) = H^2(\Dd,\Cc^{n-1}). $$
\end{proof}

\begin{proposition}\label{v*poc}
Let $\xi_0, \alpha_0$ and $V_0$ be as in Lemma \ref{a0h2}. Then 
$$V_0^* \Poc(\{\xi_0\}, L^2(\Tt,\Cc^n)) = \begin{pmatrix} 0 \\ L^2(\Tt,\Cc^{n-1})
\end{pmatrix}. $$
\end{proposition}

\begin{proof}
Let $g \in V_0^* \Poc(\{\xi_0\}, L^2(\Tt,\Cc^n)).$ Equivalently, $g$ can be written as 
$g= V_0^* f$ for some $f\in L^2(\Tt,\Cc^n)$ such that $f(z) \perp \xi_0(z)$ for almost all $z\in \Tt.$ 
This in turn is equivalent to the assertion that 
$g=V_0^*f$ for some $f\in L^2(\Tt,\Cc^n)$ such that $(V_0^*f)(z) \perp (V_0^*\xi_0)(z)$ for almost all $z \in \Tt,$ since $V_0(z)$ is unitary for almost all $z \in \Tt.$  

\noindent Note that, by the fact that $V_0$ is unitary-valued almost everywhere on $\Tt$, we have \begin{align}I_n &= V_0^*(z)V_0(z)\nonumber \vspace{2ex}\\  &= \begin{pmatrix}
\xi_0^*(z) \\ \alpha_0^T(z) 
\end{pmatrix} \begin{pmatrix}
\xi_0(z) & \bar{ \alpha}_0(z)
\end{pmatrix}\nonumber \vspace{2ex}\\ 
&= \begin{pmatrix}
\xi_0^*(z) \xi_0(z) & \xi_0^*(z) \bar{ \alpha}_0(z) \\
\alpha_0^T(z) \xi_0(z) & \alpha_0^T(z) \bar{\alpha}_0(z)
\end{pmatrix}\quad\text{almost everywhere on}\;\Tt\label{vounitary},\end{align} and so

 $$V_0^*\xi_0 = \begin{pmatrix}\xi_0^* \\ \alpha_0^T \end{pmatrix} \xi_0 = \begin{pmatrix} 1 \\ 0_{(n-1)\times 1} \end{pmatrix},$$ where 
$0_{(n-1)\times 1}$ denotes the zero vector in $\Cc^{n-1}.$

\noindent Hence $g=V_0^*f$ with $(V_0^*f)(z)$ orthogonal to $(V_0^*\xi_0)(z)$ for almost every $z\in \Tt,$ is equivalent to the statement 
$g \in L^2(\Tt,\Cc^n)$ and $$g(z) \perp  \begin{pmatrix} 1 \\ 0_{(n-1)\times 1} \end{pmatrix}$$ for almost all $z\in \Tt,$ or equivalently, $g\in \begin{pmatrix} 0\\ L^2(\Tt,\Cc^{n-1})\end{pmatrix}.$ \end{proof}  \index{ $0_{(n-1)\times 1}$}

\begin{proposition}\label{xi0wek0subset}
	Under the assumptions of Theorem \ref{T0compact}, where $x_0$ is a co-outer maximizing vector of unit norm for $H_G,$
 $\xi_0\in H^\infty(\Dd,\Cc^n)$ is a vector-valued inner function given by $\xi_0= \frac{x_0}{h_0} ,$ $V_0=\begin{pmatrix}
\xi_0 & \bar{\alpha}_0
\end{pmatrix}$ is a thematic completion of $\xi_0$ and $\mathcal{K}_1$ is defined by $$\mathcal{K}_1 = V_0 \begin{pmatrix}
0 \\ H^2(\Dd,\Cc^{n-1})
\end{pmatrix}\subseteq L^2(\Tt,\Cc^n),$$ we have 
	$$\xi_0 \telwe \mathcal{K}_1 = \xi_0 \telwe H^2(\Dd,\Cc^n) $$ and the operator
	 $$(\xi_0 \telwe \cdot) \colon  \mathcal{K}_1 \to \xi_0 \telwe H^2(\Dd,\Cc^n)$$ is unitary. 
\end{proposition}

\begin{proof}
	\noindent Let us first prove $\xi_0 \telwe H^2(\Dd,\Cc^n) \subset \xi_0 \telwe \mathcal{K}_1.$ Let $\phi\in H^2(\Dd,\Cc^n).$ Since $V_0$ is unitary-valued, 
	$$
\xi_0\xi_0^* + \bar{\alpha}_0\alpha_0^T = I_n. 
$$ 
Thus 
	$$
 \begin{array}{cllll} \xi_0 \telwe \phi &= \xi_0 \telwe (\xi_0\xi_0^*\phi +\bar{\alpha}_0\alpha_0^T\phi) \vspace{2ex} \\ 
	&= \xi_0 \telwe \xi_0 (\xi_0^*\phi) + \xi_0 \telwe \bar{\alpha}_0 (\alpha_0^T \phi) \vspace{2ex} \\ 
	&= 0 + 	\xi_0 \telwe \bar{\alpha}_0 (\alpha_0^T \phi) \end{array}
$$
on account of the pointwise linear dependence of $\xi_0$ and $\xi_0\xi_0^*\phi$ on $\Dd.$ Recall that, by Lemma \ref{a0h2}, $\alpha_0^T \phi \in H^2(\Dd,\Cc^{n-1})$ and, by the definition of $\mathcal{K}_1,$ 
$$\mathcal{K}_1 = \bar{\alpha}_0 H^2(\Dd,\Cc^{n-1}).$$ 
Hence, for $\phi\in H^2(\Dd,\Cc^n),$ $$\xi_0 \telwe \phi = \xi_0 \telwe \bar{\alpha}_0\alpha_0^T \phi \in \xi_0 \telwe \bar{ \alpha}_0 H^2(\Dd,\Cc^{n-1}),$$ and thus 
	\begin{equation}\label{xi0telwek0}\xi_0 \telwe H^2(\Dd,\Cc^n) \subseteq \xi_0 \telwe \mathcal{K}_1. \end{equation}
	
 Let us now show that $\xi_0 \telwe \mathcal{K}_1 \subseteq \xi_0 \telwe H^2(\Dd,\Cc^n).$ Since $\mathcal{K}_1 = \bar{\alpha}_0H^2(\Dd,\Cc^{n-1}),$ an arbitrary element $u \in \xi_0 \telwe \mathcal{K}_1$ is of the form 
$$ u = \xi_0 \telwe \bar{\alpha}_0g,$$ for some $g\in H^2(\Dd,\Cc^{n-1}).$
Note that, by Lemma \ref{a0h2}, there exists a function $f \in H^2(\Dd,\Cc^n)$ such that $g = \alpha_{0}^T f.$ Hence $u=\xi_0 \telwe \bar{\alpha}_0 \alpha_0^T f.$ By equation \eqref{vounitary}, $\xi_0 \xi_0^* + \bar{\alpha}_0 \alpha_0^T = I_n.$ Thus 
$$u= \xi_0 \telwe (I_{n}-\xi_0 \xi_0^*)f = \xi_0 \telwe f - \xi_0 \telwe \xi_0 \xi_0^*f = \xi_0 \telwe f \in \xi_0 \telwe H^2(\Dd,\Cc^n),  $$ and so, 
$ \xi_0 \telwe \mathcal{K}_1 \subseteq \xi_0 \telwe H^2(\Dd,\Cc^n).$
 Combining the latter inclusion with the relation \eqref{xi0telwek0}, we have
 $$
 \xi_0 \telwe \mathcal{K}_1 = \xi_0 \telwe H^2(\Dd,\Cc^n) .
$$  
	 Now, let us show that the operator $(\xi_0 \telwe \cdot) \colon  \mathcal{K}_1 \to \xi_0 \telwe H^2(\Dd,\Cc^n)$ is unitary. As we have shown above, the operator is surjective. We will show it is also an isometry.   Let $f \in \mathcal{K}_1.$ Then, 
	$$\begin{array}{cllllllll}
	\| \xi_0 \telwe f\|_{L^2(\Tt,\we^2\Cc^n)}^2 
	&= \langle \xi_0 \telwe f, \xi_0 \telwe f \rangle_{L^2(\Tt,\we^2\Cc^n)} \vspace{2ex} \\
		&= \displaystyle \frac{1}{2\pi} \int_0^{2\pi} \langle \xi_0(\eiu) \telwe f(\eiu), \xi_0(\eiu) \telwe f(\eiu) \rangle_{\we^2\Cc^n} d\theta\end{array}.$$ By Proposition \ref{we}, the latter integral is equal to
	
			$$\begin{array}{llll}
		
			&= \displaystyle \frac{1}{2\pi} \int_0^{2\pi} \det \begin{pmatrix}
			\langle \xi_0 (\eiu), \xi_0 (\eiu) \rangle_{\Cc^n} & \langle \xi_0(\eiu) , f(\eiu) \rangle_{\Cc^n}  \\
			\langle f(\eiu) , \xi_0 (\eiu) \rangle_{\Cc^n} & \langle f(\eiu) , f(\eiu) \rangle_{\Cc^n} 
			\end{pmatrix}\; d\theta \vspace{2ex} \\
			&=  \displaystyle \frac{1}{2\pi} \int_0^{2\pi} \|\xi_0 (\eiu)\|_{\Cc^n}^2\langle f(\eiu) , f(\eiu) \rangle_{\Cc^n}  - |\langle \xi_0(\eiu) , f(\eiu) \rangle_{\Cc^n}  |^2 \; d\theta. 
	\end{array} $$
\noindent Note that, by Proposition \ref{onxi}, $  \|\xi_0 (\eiu)\|_{\Cc^n} =1 $ for almost all $\eiu$ on $\Tt.$ Moreover, since $$\mathcal{K}_1= \bar{\alpha}_0H^2(\Dd,\Cc^{n-1}) ,$$  $f=\bar{\alpha}_0 g$ for some $g \in  H^2(\Dd,\Cc^{n-1}).$ 
Hence $$\langle \xi_0 (\eiu), f(\eiu) \rangle_{\Cc^n} = \langle \xi_0(\eiu), \bar{\alpha}_0(\eiu) g(\eiu)\rangle_{\Cc^n}  =\langle  \alpha_{0}^T(\eiu)\xi_0(\eiu),  g(\eiu)\rangle_{\Cc^{n-1}}= 0 $$ 
almost everywhere on $\Tt,$ since $ V_0 = \begin{pmatrix} \xi_0 & \bar{\alpha}_0 \end{pmatrix}$ 
is unitary-valued. Thus 

$$ \| \xi_0 \telwe f\|_{L^2(\Tt,\we^2\Cc^n)}^2 = \|f\|_{L^2(\Tt,\Cc^n)}^2,$$that is, the operator $(\xi_0 \telwe \cdot) \colon  \mathcal{K}_1 \to \xi_0 \telwe H^2(\Dd,\Cc^n)$ is an isometry. Therefore, the surjective operator $(\xi_0 \telwe \cdot)$ is unitary.  \end{proof}

\begin{lemma}\label{inner0iff}
	Let $u\in L^2(\Tt,\Cc^m)$ and let $\eta_0 \in H^\infty(\Dd,\Cc^m)$  be a vector-valued inner function. Then 
\begin{equation}\label{cond9.7}
\langle \bar{\eta}_0 \telwe u, \bar{\eta}_0 \telwe \bar{z} \bar{f} \rangle_{L^2(\Tt,\we^2\Cc^m)} = 0 \quad \text{for all}\quad f\in H^2(\Dd,\Cc^m)
\end{equation}
	if and only if the function 
	$$ z \mapsto u(z) - \langle u(z), \bar{ \eta}_{0}(z)\rangle_{\Cc^m} \bar{ \eta}_{0}(z)$$ belongs to $H^2(\Dd,\Cc^m).$
\end{lemma}

\begin{proof}
	\noindent The statement that $\bar{\eta}_0 \telwe u$ is orthogonal to $\bar{\eta}_0 \telwe \bar{z} \bar{f}$ in $L^2(\Tt,\we^2\Cc^m)$ is equivalent to the equation $I =0$, where
	$$I=\frac{1}{2\pi}\int_{0}^{2\pi} \langle \bar{\eta}_0(\eiu) \telwe u(\eiu) , \bar{\eta}_0(\eiu) \telwe e^{-i\theta} \bar{f}(\eiu)\rangle_{\we^2\Cc^m} \; d\theta.$$
By Proposition \ref{we}, 
	$$I=\displaystyle \frac{1}{2\pi}\int_{0}^{2\pi} \det \begin{pmatrix}
	\langle \bar{\eta}_0(\eiu), \bar{\eta}_0(\eiu)\rangle_{\Cc^m} & \langle \bar{\eta}_0(\eiu) , e^{-i\theta} \bar{f}(\eiu)\rangle_{\Cc^m} \\
	\langle u(\eiu) , \overline{ \eta_{0}}(\eiu) \rangle_{\Cc^m} & \langle u(\eiu), e^{-i\theta} \bar{f}(\eiu)\rangle_{\Cc^m}
	\end{pmatrix}d\theta. $$
Notice that, since $\eta_0$ is an inner function,  $\|\bar{\eta}_0(\eiu)\|_{\Cc^m}=1$ almost everywhere on $\Tt,$ and hence
	$$\begin{array}{cllllll}
	I&=\displaystyle \frac{1}{2\pi}\int_{0}^{2\pi} \det \begin{pmatrix}
	1 & \langle \bar{\eta}_0(\eiu) , e^{-i\theta} \bar{f}(\eiu)\rangle_{\Cc^m} \\
	\langle u(\eiu) , \overline{ \eta_{0}}(\eiu) \rangle_{\Cc^m} & \langle u(\eiu), e^{-i\theta} \bar{f}(\eiu)\rangle_{\Cc^m}
	\end{pmatrix}d\theta \vspace{2ex}\\	
	&=\displaystyle \frac{1}{2\pi} \displaystyle \int_{0}^{2\pi} \langle u(\eiu), e^{-i\theta}\bar{f}(\eiu)\rangle_{\Cc^m} \\ &\hspace{15ex}- \langle \bar{\eta}_0(\eiu) , e^{-i\theta} \bar{f}(\eiu)\rangle_{\Cc^m} \langle u(\eiu) , \bar{\eta}_0(\eiu)\rangle_{\Cc^m} d\theta \vspace{2ex} \\
	&=\displaystyle \frac{1}{2\pi}\displaystyle \int_{0}^{2\pi} \langle u(\eiu), e^{-i\theta}\bar{f}(\eiu)\rangle_{\Cc^m} \\ &\hspace{15ex}- \left\langle \langle u(\eiu), \bar{\eta}_0(\eiu) \rangle_{\Cc^m} \bar{\eta}_0(\eiu) , e^{-i\theta} \bar{f}(\eiu)  \right\rangle_{\Cc^m} d\theta \vspace{2ex} \\
	&=\displaystyle \frac{1}{2\pi} \displaystyle \int_{0}^{2\pi} \left\langle u(\eiu) -\langle u(\eiu), \bar{\eta}_0(\eiu)\rangle_{\Cc^m}\bar{\eta}_0(\eiu) , e^{-i\theta} \bar{f}(\eiu)\right\rangle_{\Cc^m} d\theta .
	\end{array}$$ 
	
\noindent 	Thus, the condition \eqref{cond9.7} holds  if and only if 
	$$\displaystyle\frac{1}{2\pi}\int_{0}^{2\pi} \langle \bar{\eta}_0(\eiu) \telwe u(\eiu) , \bar{\eta}_0(\eiu) \telwe e^{-i\theta} \bar{f}(\eiu)\rangle_{\we^2\Cc^m} \; d\theta = 0  \quad \text{for all}\quad f\in H^2(\Dd,\Cc^m) $$ if and only if 
	$$ \displaystyle\frac{1}{2\pi} \int_{0}^{2\pi} \left\langle u(\eiu) - \langle u(\eiu), \bar{\eta}_0(\eiu)\rangle_{\Cc^m}\bar{\eta}_0(\eiu) , e^{-i\theta} \bar{f}(\eiu)\right\rangle_{\Cc^m} d\theta  =0   $$for all $f\in H^2(\Dd,\Cc^m),$ and the latter equation holds if and only if 
	$$u(\eiu) - \langle u(\eiu), \bar{\eta}_0(\eiu)\rangle_{\Cc^m} \bar{\eta}_0(\eiu) $$ belongs to $H^2(\Dd,\Cc^m).$ \end{proof}

\begin{lemma}\label{l0perp}
In the notation of Theorem \ref{T0compact},
	$$ \mathcal{L}_1^\perp = \{f \in L^2(\Tt,\Cc^m) \; : \beta_0^*f \in H^2(\Dd,\Cc^{m-1}) \} .$$
\end{lemma}
\begin{proof}
It is easy to see that $\mathcal{L}_1 = \beta_0 H^2(\Dd,\Cc^{m-1})^\perp.$
The general element of $\beta_0 H^2(\Dd,\Cc^{m-1})^\perp$ is $\beta_0 \bar{z} \bar{g}$ with $g \in H^2(\Dd,\Cc^{m-1}).$
For $f \in L^2(\Tt,\Cc^m),$ 
$f \in \mathcal{L}_1^\perp$ if and only if 
$$ \langle f,\beta_0 \bar{z} \bar{g} \rangle_{L^2(\Tt,\Cc^m)} = 0 \quad \text{for all} \quad g \in H^2(\Dd,\Cc^{m-1}).$$ 
Equivalently, $f \in \mathcal{L}_1^\perp$ if and only if
$$ \displaystyle \frac{1}{2\pi} \int_{0}^{2\pi} \langle f(\eiu), \beta_0(\eiu) e^{-i\theta} \overline{g}(\eiu) \rangle_{\Cc^m} d\theta =0  \quad \text{for all} \quad g \in H^2(\Dd,\Cc^{m-1}) $$ if and only if 
$$\displaystyle \frac{1}{2\pi} \int_{0}^{2\pi} \langle \beta_0(\eiu)^* f(\eiu),  e^{-i\theta} \overline{g}(\eiu) \rangle_{\Cc^{m-1}} d\theta =0  \quad \text{for all} \quad g \in H^2(\Dd,\Cc^{m-1}) .$$ 
The latter statement is equivalent to the assertion that $ \beta_0^* f$ is orthogonal to $ H^2(\Dd,\Cc^{m-1})^\perp$ in $L^2(\Tt,\Cc^{m-1}),$ which holds if and only if $\beta_0^* f $ belongs to $H^2(\Dd,\Cc^{m-1}).$
	
\noindent Hence  
	$$ \mathcal{L}_1^\perp = \{f \in L^2(\Tt,\Cc^m) \; : \beta_0^*f \in H^2(\Dd,\Cc^{m-1}) \} $$ as required.	\end{proof}

\begin{proposition}\label{beta0*h2} 
Under the assumptions of Theorem \ref{T0compact}, let $\eta_0$ be defined by equation \eqref{x00y00} and let $W_0^T = \begin{pmatrix} \eta_0 & \bar{\beta}_0 \end{pmatrix}$ be a thematic completion of $\eta_0,$ where $\beta_0$ is an inner, co-outer, quasi-continuous function of type $m \times (m-1).$ 
Then, 
$$ \beta_0^* H^2(\Dd,\Cc^{m})^\perp = H^2(\Dd,\Cc^{m-1})^\perp.$$	
\end{proposition}

\begin{proof}
	By virtue of the fact that complex conjugation is a unitary operator on $L^2(\Tt,\Cc^m),$ an equivalent statement is that 
$\beta_0^T z H^2(\Dd,\Cc^{m}) = z H^2(\Dd,\Cc^{m-1}).$ 
By Lemma \ref{L6.2}, since $\beta_0$ is an inner, co-outer and quasi-continuous function, there exists a matrix-valued function 
$B_0 \in H^\infty( \Dd,\Cc^{(m-1)\times m})$ such that 
$$B_0 \beta_0 = I_{m-1}$$ 
or, equivalently, 
$$ \beta_0^T B_0^T = I_{m-1}.$$
	
	\noindent Let $g \in z H^2(\Dd,\Cc^{m-1}).$ 
Then, 
$$ g = (\beta_0^T B_0^T) g  \in  \beta_0^T B_0^T z H^2(\Dd,\Cc^{m-1}) \subseteq \beta_0^T z H^2(\Dd,\Cc^{m}).$$ 
Hence 
$$ z H^2(\Dd,\Cc^{m-1}) \subseteq \beta_0^T z H^2(\Dd,\Cc^{m}). $$
Note that, 
$\beta_0 \in  H^\infty(\Dd,\Cc^{m \times (m-1)}),$
and so, 
$$ zH^2(\Dd,\Cc^{m-1}) \subseteq \beta_0^T z H^2(\Dd,\Cc^{m}) 
\subseteq z H^2(\Dd,\Cc^{m-1}).$$ 
Thus 
$$
\beta_0^T z H^2(\Dd,\Cc^{m}) = z H^2(\Dd,\Cc^{m-1}).
$$
\end{proof}

\begin{proposition}\label{eta0telweh2}
In the notation of Theorem \ref{T0compact}, let $\eta_0 \in H^\infty(\Dd,\Cc^m)$ be a vector-valued inner function given by equation \eqref{x00y00}, let $W_0^T = \begin{pmatrix}
	\eta_0 & \bar{ \beta}_0
	\end{pmatrix}$ be a thematic completion of $\eta_0$ given by equation \eqref{v00w00}, and let 
$$
\mathcal{L}_1 = W_0^* \begin{pmatrix}
	0 \\ H^2(\Dd,\Cc^{m-1})^\perp\end{pmatrix}.
$$ 
	 Then 
$$
\bar{\eta}_0 \telwe \mathcal{L}_1 = \bar{\eta}_0 \telwe H^2(\Dd,\Cc^m)^\perp
$$
and the operator 
$$
(\bar{\eta}_0 \telwe \cdot)\colon \mathcal{L}_1 \to 
\bar{\eta}_0\telwe H^2(\Dd,\Cc^m)^\perp 
$$ 
is unitary. 
\end{proposition}

\begin{proof}
Let us first prove that 
$\bar{\eta}_0 \telwe H^2(\Dd,\Cc^m)^\perp \subseteq \bar{\eta}_0 \telwe \mathcal{L}_1.$ Consider an element  $f \in H^2(\Dd,\Cc^m)^\perp.$ Note that, since $W_0^T$ is unitary valued, we have 
\begin{equation}\label{wounit} 
\bar{\eta}_0\eta_0^T  +\beta_0\beta_0^*=I_m.
\end{equation} 
Thus
$$
\begin{array}{cllll}
\bar{\eta}_0 \telwe f &= \bar{\eta}_0 \telwe (\bar{\eta}_0\eta_0^T +\beta_0\beta_0^*)f \vspace{2ex} \\ 
&= \bar{\eta}_0 \telwe \bar{\eta}_0\eta_0^T f + \bar{\eta}_0 \telwe \beta_0\beta_0^*f \vspace{2ex}\\
&=  0+\bar{\eta}_0 \telwe \beta_0\beta_0^*f,
\end{array} 
$$
the last equality following by the pointwise linear dependence of $\bar{ \eta}_0$ and $\bar{ \eta}_0 (\eta_0^T f)$ on $\Dd.$  By Proposition \ref{beta0*h2}, $$ \beta_0^* H^2(\Dd,\Cc^{m})^\perp =  H^2(\Dd,\Cc^{m-1})^\perp,
$$
and,  by the definition of $\mathcal{L}_1$, we have
 $$
\mathcal{L}_1 = \beta_0 H^2(\Dd,\Cc^{m-1})^\perp. 
$$ 
Hence, for $f \in H^2(\Dd,\Cc^m)^\perp,$
$$ 
\bar{\eta}_0 \telwe f = \bar{\eta}_0 \telwe \beta_0\beta_0^*f \in \bar{\eta}_0 \telwe \beta_0 H^2(\Dd,\Cc^{m-1})^\perp, 
$$
and thus 
$$
\bar{\eta}_0 \telwe H^2(\Dd,\Cc^m)^\perp \subseteq \bar{\eta}_0 \telwe \mathcal{L}_1. 
$$ 

\noindent Let us show 
$$
\bar{\eta}_0 \telwe \mathcal{L}_1 \subseteq \bar{\eta}_0 \telwe H^2(\Dd,\Cc^m)^\perp.
$$ 
A typical element of $\bar{\eta}_0 \telwe \mathcal{L}_1$ is of the form $\bar{\eta}_0 \telwe \beta_0 g ,$ for some $g \in H^2(\Dd,\Cc^{m-1})^\perp.$ By Proposition \ref{beta0*h2}, there exists a $\phi \in H^2(\Dd,\Cc^m)^\perp$ such that $\beta_0^*\phi= g.$ Then
$$ 
\bar{\eta}_0 \telwe \beta_0 g = \bar{\eta}_0 \telwe \beta_0 \beta_0^*\phi.
$$
 By equation \eqref{wounit}, we have
$$ \bar{\eta}_0 \telwe \beta_0 g = \bar{\eta}_0 \telwe (I_{m}- \bar{\eta}_0\eta_0^T)\phi =  \bar{\eta}_0 \telwe \phi,  $$ 
the last equality following by pointwise linear dependence of $\bar{ \eta}_0$ and $\bar{ \eta}_0(\eta_0^T\phi) $ on $\Dd$. Thus $$ \bar{\eta}_0 \telwe \beta_0 g \in \bar{\eta}_0 \telwe H^2(\Dd,\Cc^m)^\perp ,$$ and so  $  \bar{\eta}_0 \telwe \mathcal{L}_1 \subseteq \bar{\eta}_0 \telwe H^2(\Dd,\Cc^m)^\perp.$
 Consequently 
 $$ \bar{\eta}_0 \telwe \mathcal{L}_1 = \bar{\eta}_0 \telwe H^2(\Dd,\Cc^m)^\perp.$$
 
 To prove that the operator 
$$
(\bar{\eta}_0 \telwe \cdot)\colon \mathcal{L}_1 \to 
\bar{\eta}_0\telwe H^2(\Dd,\Cc^m)^\perp 
$$ 
is unitary, it suffices to show that it is an isometry, since the preceding discussion asserts that it is surjective. To this end, let $s \in \mathcal{L}_1.$ Then,
$$\begin{array}{lll} \|\bar{\eta}_0 \telwe s\|_{L^2(\Tt,\we^2\Cc^m)}^2 &=  \langle \bar{\eta}_0 \telwe s , \bar{\eta}_0 \telwe s \rangle_{L^2(\Tt,\we^2\Cc^m)} \vspace{2ex} \\  
&= \displaystyle\frac{1}{2\pi} \int_0^{2\pi} \langle \bar{\eta}_0 (\eiu) \telwe s(\eiu) , \bar{\eta}_0(\eiu) \telwe s(\eiu) \rangle_{\we^2\Cc^m}\;d\theta \vspace{2ex} \\ 
&= \displaystyle \frac{1}{2\pi}\int_{0}^{2\pi} \det \begin{pmatrix}
\langle \bar{\eta}_0 (\eiu) , \bar{\eta}_0 (\eiu) \rangle_{\Cc^m}	& \langle \bar{\eta}_0 (\eiu) , s(\eiu)\rangle_{\Cc^m} \\
	\langle s(\eiu) , \bar{\eta}_0 (\eiu) \rangle_{\Cc^m} & \langle s(\eiu), s(\eiu)\rangle_{\Cc^m}
	\end{pmatrix}d\theta.
 \end{array} $$ By Proposition \ref{onxi}, $\|\bar{\eta}_0(z)\|_{\Cc^m} =1$ almost everywhere on $\Tt.$ Moreover, since $s \in \mathcal{L}_1,$ there exists a function $\psi \in H^2(\Dd,\Cc^{m-1})^\perp$ such that $s=\beta_0 \psi.$ Then
 $$\langle \bar{\eta}_0(\eiu), s(\eiu)\rangle_{\Cc^m} = \langle \bar{\eta}_0(\eiu), \beta_0(\eiu) \psi(\eiu)\rangle_{\Cc^m}=\langle \beta_0^*(\eiu)\bar{\eta}_0(\eiu),  \psi(\eiu)\rangle_{\Cc^m}=0$$almost everywhere on $\Tt,$ which follows by the fact that $W_0$ is unitary-valued, and so $$(W_0W_0^*)(z)=\begin{pmatrix}\eta_0^T(z) \\ \beta_0^*(z) \end{pmatrix}\begin{pmatrix}
 \bar{ \eta}_0(z) & \beta_0(z) 
 \end{pmatrix}=\begin{pmatrix}
 \eta_0^T(z)\bar{ \eta}_0(z) & \eta_0^T(z) \beta_0^T(z)\\ \beta_0^*(z) \bar{ \eta}_0(z) & \beta_0^*(z) \beta_0(z)
 \end{pmatrix} =I_m$$ almost everywhere on $\Tt.$

\noindent Thus, for all $s \in \mathcal{L}_1$,
 $$\|\bar{\eta}_0 \telwe s\|_{L^2(\Tt,\we^2\Cc^m)}^2 = \|s\|_{L^2(\Tt,\Cc^m)}^2 , $$ 
which shows that the operator 
$$
(\bar{\eta}_0 \telwe \cdot)\colon \mathcal{L}_1 \to 
\bar{\eta}_0\telwe H^2(\Dd,\Cc^m)^\perp 
$$ 
is an isometry. We have proved it is also surjective, hence the operator $(\bar{\eta}_0\telwe \cdot)$ is unitary.\end{proof}

{\it Continuation of the proof of Theorem \ref{T0compact}.}\\
{\rm (iii)}	We have to prove that diagram \eqref{commdiagr} commutes. 	
Recall,	by Lemma \ref{3.2constr}, the left hand square commutes, so it suffices to show that that the right hand square, namely
\begin{equation}\label{commdiagrr}
\begin{array}{clllll}

\mathcal{K}_1 &\xrightarrow{\xi_0 \telwe \cdot}& \xi_0 \telwe H^2 (\Dd, \Cc^n)=X_1\\
\Big\downarrow\rlap{$\scriptstyle \Gamma_1$}  &~&\hspace{3ex}\Big\downarrow\rlap{$\scriptstyle T_1$} \\
	\mathcal{L}_1 &\xrightarrow{\bar{\eta}_0 \telwe \cdot } & \bar{\eta}_0 \telwe H^2 (\Dd, \Cc^m)^\perp =Y_1
\end{array},\end{equation}also commutes. That is, we wish to prove that, for all $x \in \mathcal{K}_1,$ 
$$	T_1(\xi_0 \telwe x)= \bar{ \eta}_0 \telwe \Gamma_1(x), $$where
	$\Gamma_1(x) = P_{\mathcal{L}_1}((G-Q_1)x)$ for any function $Q_1 \in H^\infty(\Dd,\CCmn)$ that satisfies the following equations
	$$(G-Q_1)x_0 = t_0 y_0, \quad y_0^*(G-Q_1) = t_0 x_0^*. $$ By Proposition \ref{xi0wek0subset}, 
	$$\xi_0 \telwe \mathcal{K}_1 = \xi_0 \telwe H^2(\Dd,\Cc^n), $$ and so, for every $x \in \mathcal{K}_1, $ there exists $\tilde{x} \in H^2(\Dd,\Cc^n)$ such that 
	$$\xi_0 \telwe x = \xi_0 \telwe \tilde{x}. $$ Thus, for $ x \in \mathcal{K}_1,$
	
	$$T_1(\xi_0 \telwe x)= T_1(\xi_0\telwe \tilde{x}) = P_{Y_1}(\bar{ \eta}_0 \telwe (G-Q_1)\tilde{x}), $$ and 
	$$ \bar{ \eta}_0 \telwe \Gamma_1(x) = \bar{ \eta}_0 \telwe P_{\mathcal{L}_1}(G-Q_1)x.$$ 
Hence to prove the commutativity of diagram (\ref{commdiagrr}), it suffices to show that, for all $x\in \mathcal{K}_1,$ 
$$P_{Y_1}[ \bar{\eta}_0 \telwe(G-Q_1)\tilde{x})] = \bar{\eta}_0 \telwe P_{\mathcal{L}_1}(G-Q_1)x$$
in $Y_1,$ where $\xi_0\telwe (x-\tilde{x})=0.$ By Proposition \ref{eta0telweh2},  
 $$\bar{\eta}_0 \telwe \mathcal{L}_1 = \bar{\eta}_0 \telwe H^2(\Dd,\Cc^m)^\perp=Y_1,$$ and so, for all $x \in \mathcal{K}_1,$	$\bar{ \eta}_0 \telwe P_{\mathcal{L}_1}(G-Q_1)x\in Y_1.$ 
	\noindent  Let us show that, for $x \in \mathcal{K}_1,$
	$$\bar{\eta}_0 \telwe  (G-Q_1)\tilde{x} - \bar{\eta}_0\telwe P_{\mathcal{L}_1} (G-Q_1)x  $$ is orthogonal to $Y_1$ in $L^2(\Tt,\we^2 \Cc^m),$ or equivalently, that for every $f \in H^2(\Dd,\Cc^m),$ 
	\begin{equation}\label{7.1}	\left\langle \bar{\eta}_0 \telwe [ (G-Q_1)\tilde{x} - P_{\mathcal{L}_1}(G-Q_1)x ] , \bar{\eta}_0 \telwe \bar{z}\bar{f} \right\rangle_{L^2(\Tt, \we^2\Cc^m)} =0  \end{equation}  for $x \in \mathcal{K}_1$ and for any $\tilde{x} \in H^2(\Dd,\Cc^n)$ such that $\xi_0 \telwe \tilde{x} = \xi_0 \telwe x.$ By Lemma \ref{0to0}, $$\bar{ \eta}_0 \telwe (G-Q_1)x=\bar{ \eta}_0 \telwe (G-Q_1)\tilde{x}.$$ Then equation (\ref{7.1}) is equivalent to the equation 
	\begin{equation}\label{7.11} \langle \bar{\eta}_0 \telwe P_{\mathcal{L}_1^\perp} (G-Q_1)x, \bar{\eta}_0 \telwe \bar{z} \bar{f} \rangle_{L^2(\Tt, \we^2\Cc^m)} =0  \end{equation} for any $x \in \mathcal{K}_1.$
 By Lemma \ref{inner0iff}, equation (\ref{7.11}) holds if and only if the function 
\begin{equation}\label{zmapsto1} z \mapsto [P_{\mathcal{L}_1^\perp}(G-Q_1)x](z) - \langle [P_{\mathcal{L}_1^\perp}(G-Q_1)x](z) ,\bar{\eta}_0(z)\rangle_{\Cc^m} \bar{\eta}_0(z)\end{equation}belongs to $H^2(\Dd,\Cc^m).$
By Lemma \ref{l0perp}, there exists a function $\psi \in L^2(\Tt,\Cc^m)$ such that 
\begin{equation}\label{PLperp}
P_{\mathcal{L}_1^\perp} (G-Q_1)x=\psi
\end{equation}
and
$$\beta_0^*\psi \in H^2(\Dd,\Cc^{m-1}).$$
 Equation \eqref{PLperp} implies
 $$(G-Q_1)x -\psi \in \mathcal{L}_1 = \beta_0 H^2(\Dd,\Cc^{m-1})^\perp.$$
	
	\noindent Hence, to prove that the function defined by equation (\ref{zmapsto1}) belongs to $H^2(\Dd,\Cc^m)$, we have to show that 
	$$\psi - (\eta_0^T \psi) \bar{\eta}_0 \in H^2(\Dd,\Cc^m).$$	 
Since $W_0 =\begin{pmatrix}
		\eta_0 & \overline{ \beta_{0}}
	\end{pmatrix}^T $ is a unitary-valued function,
	$$\bar{ \eta}_0(z) \eta_0^T(z)+ \beta_0(z)\beta_0^*(z)=I_m$$
almost everywhere on $\Tt.$
	 \noindent Since $\eta_0^T \psi$ is a scalar-valued function,
$$ \psi - (\eta_0^T \psi) \overline{ \eta_0 } = 
(I_m - \eta_0^T  \overline{ \eta_0 })\psi = \beta_0 \beta_0^*\psi  \in H^2(\Dd,\Cc^m) .$$ 
Thus diagram (\ref{commdiagrr}) commutes. 
	
{\rm (iv)} By Lemma \ref{2.2}, $$F_1 \in H^\infty(\Dd,\Cc^{(m-1)\times (n-1)})+C(\Tt,\Cc^{(m-1)\times (n-1)}).$$ Then, by Hartman's Theorem \ref{2.04}, the Hankel operator $H_{F_1}$ is compact, and by (iii), $$(\bar{\eta}_0 \telwe \cdot) \circ (U_2 H_{F_1} U_1^*) \circ(\xi_0\telwe \cdot )^*= T_1 .$$ By (i) and (ii), the operators $U_1, U_2,$ $(\xi_0 \telwe \cdot)$ and $(\bar{ \eta}_0 \telwe \cdot)$ are unitary. Hence $T_1$ is a compact operator.
  
{\rm (v)} Since diagram \eqref{commdiagr} is commutative and $U_1,U_2,(\xi_0\telwe\cdot)$ and $(\eta_0\telwe \cdot)$ are unitaries,
$$\|T_1\| = \|\Gamma_1\| = \|H_{F_1}\|. $$ 	
\end{proof}

\noindent In what follows, we will prove an analogous statement to Theorem \ref{T0compact} for $T_2.$ To this end, we need the following results. 
\begin{lemma}\label{connofschpairs1}
	In the notation of Theorem \ref{T0compact}, $v_1 \in H^2(\Dd,\Cc^n)$ and $w_1 \in H^2(\Dd,\Cc^m)^\perp$ are such that $(\xi_0 \telwe v_1, \bar{\eta}_0 \telwe w_1)$ is a Schmidt pair for the operator $T_1$ corresponding to $\|T_1\|.$ Then

{\em (i)} there exist $x_1 \in \mathcal{K}_1$ and $y_1\in \mathcal{L}_1$ such that $(x_1,y_1) $ is a Schmidt pair for the operator $\Gamma_1$;

{\em (ii)} for any $x_1 \in \mathcal{K}_1$ and $y_1\in \mathcal{L}_1$ such that $$
	\xi_0 \telwe x_1 = \xi_0 \telwe v_1,\quad 	
 \bar{ \eta}_0 \telwe y_1 = \bar{ \eta}_0 \telwe w_1,$$ the pair $(x_1,y_1)$ is a Schmidt pair for $\Gamma_1$ corresponding to $\|\Gamma_1\|.$ 
\end{lemma}

\begin{proof}
{\rm (i)} By Theorem \ref{T0compact}, the diagram (\ref{commdiagr}) commutes, $(\xi_0 \telwe \cdot)$ is unitary from $\mathcal{K}_1$ to $X_1,$ and $(\bar{ \eta}_0\telwe \cdot)$ is unitary from $\mathcal{L}_1$ to $Y_1.$ 
Thus $\|\Gamma_1\| =\|T_1\|=t_1.$ 
Moreover, by Lemma \ref{3.2constr}, the operator $\Gamma_1\colon \mathcal{K}_1\to \mathcal{L}_1$ is compact, hence there exist $x_1 \in \mathcal{K}_1,$ $y_1 \in \mathcal{L}_1$ such that $(x_1,y_1)$ is a Schmidt pair for $\Gamma_1$ corresponding to $\|\Gamma_1\|=t_1.$ 

{\rm (ii)} Suppose that $x_1\in \mathcal{K}_1,y_1\in \mathcal{L}_1$ satisfy \begin{equation}\label{xitel1}
\xi_0 \telwe x_1 = \xi_0 \telwe v_1,\end{equation}	
\begin{equation}\label{eta0tel1} \bar{ \eta}_0 \telwe y_1 = \bar{ \eta}_0 \telwe w_1. \end{equation} Let us show that $(x_1,y_1)$ is a Schmidt pair for $\Gamma_1$ corresponding to $t_1$, that is, 
$$\Gamma_1 x_1 = t_1y_1,\quad \Gamma_1^*y_1=t_1x_1. $$ Since diagram (\ref{commdiagrr}) commutes,
\begin{equation}\label{commt1gamma1}T_1 \circ (\xi_0\telwe \cdot  )=(\bar{ \eta}_0 \telwe \cdot)\circ\Gamma_1 , \quad (\xi_0\telwe \cdot  )^*\circ T_1^* = \Gamma_1^* \circ (\bar{\eta}_0 \telwe \cdot)^*. \end{equation} 
By hypothesis,
\begin{equation}\label{hypt1} T_1 (\xi_0 \telwe v_1)= t_1 (\bar{ \eta}_0 \telwe w_1), \quad T_1^*(\bar{ \eta}_0 \telwe w_1)= t_1 (\xi_0 \telwe v_1). \end{equation} Thus, by equations \eqref{eta0tel1}, \eqref{commt1gamma1} and \eqref{hypt1},
$$\begin{array}{clllll} \Gamma_1  x_1&= (\bar{\eta}_0 \telwe \cdot)^*  T_1  (\xi_0\telwe v_1 ) \vspace{2ex} \\
&=  (\bar{\eta}_0 \telwe \cdot)^*  t_1 (\bar{ \eta}_0 \telwe w_1) \vspace{2ex} \\ &= t_1 (\bar{\eta}_0 \telwe \cdot)^* (\bar{ \eta}_0 \telwe y_1).\end{array}$$  Hence  $$\Gamma_1  x_1= t_1 (\bar{\eta}_0 \telwe \cdot)^* (\bar{ \eta}_0 \telwe \cdot)y_1= t_1 y_1.$$ 

\noindent By equation (\ref{xitel1}), 
$$x_1 = (\xi_0\telwe \cdot  )^* (\xi_0\telwe v_1  ),$$ 
and, by equation (\ref{eta0tel1}), 
$$(\bar{\eta}_0 \telwe \cdot)^*(\bar{ \eta}_0 \telwe w_1)=y_1. $$Thus $$\begin{array}{clll}\Gamma_1^* y_1 &= \Gamma_1^*(\bar{ \eta}_0 \telwe \cdot)^*(\bar{\eta}_0 \telwe w_1)\vspace{2ex} \\ 
&= (\xi_0 \telwe \cdot )^*  T_1^* (\bar{ \eta}_0 \telwe w_1),\end{array}$$ the last equality following by the second equation of (\ref{commt1gamma1}). 
By equations \eqref{xitel1} and (\ref{hypt1}), we have 
$$ T_1^* (\bar{ \eta}_0 \telwe w_1) = t_1 (\xi_0\telwe v_1)= t_1(\xi_0\telwe x_1),$$ and so, $$ \Gamma_1^* y_1 = t_1 x_1.$$ Therefore $(x_1,y_1)$ is a Schmidt pair for 
$\Gamma_1$ corresponding to $\|\Gamma_1\|=\|T_1\|=t_1.$ \end{proof}

\begin{lemma}\label{schfohf1}
	Suppose $(\xi_0 \telwe v_1, \bita_0 \telwe w_1)$ is a Schmidt pair for $T_1$ corresponding to $t_1.$ Let 
	$$x_1 = (I_{n} - \xi_0 \xi_0^*)v_1,\quad y_1= (I_{m} - \bita_0\eta_0^T)w_1,$$ and let 
	$$\hx_1 = \alpha_0^T x_1,\quad \hy_1=\beta_0^*y_1. $$
Then 

{\em (i)} \begin{equation}\label{x1=alpha0alphaTx1}
x_1= \bar{\alpha}_0 \alpha_0^Tx_1,\quad y_1=\beta_0\beta_0^*y_1;
\end{equation}

{\em (ii)} the pair $(\hx_1,\hy_1)$ is a Schmidt pair for $H_{F_1}$ corresponding to $\|H_{F_1}\|=t_1.$
	
\end{lemma}

\begin{proof} {\rm (i)} 
Since $V_0=\big(\begin{matrix}\xi_0 & \bar{\alpha}_0 \end{matrix}\big)$ is unitary-valued,  $I_{n} - \xi_0 \xi_0^* = \balpha_0 \alpha_0^T,$ and so
\begin{align} \label{x1=alpha0alphaTx1-pr}
\bar{\alpha}_0 \alpha_0^Tx_1& =(I_{n} - \xi_0 \xi_0^*) (I_{n} - \xi_0 \xi_0^*)v_1 \nn\\
& = (I_{n} - 2\xi_0 \xi_0^* + \xi_0 \xi_0^*\xi_0 \xi_0^*)v_1 \nn\\
& = (I_{n} - \xi_0 \xi_0^*)v_1 = x_1.
\end{align}
Similarly, since $W_0^T=\big( \begin{matrix}\eta_0 & \bar{\beta}_0 \end{matrix}\big)$ 
is unitary valued,  $I_{m} - \bita_0 \eta_0^T = \beta_0 \beta_0^*,$
and so 
\begin{align} \label{y1=beta0beta*y1-pr}
\beta_0\beta_0^*y_1 &= (I_{m} - \bita_0\eta_0^T)(I_{m} - \bita_0\eta_0^T)w_1 \nn\\
& = (I_{m} - 2\bita_0\eta_0^T +\bita_0\eta_0^T\bita_0\eta_0^T)w_1 \nn\\
& = (I_{m} - \bita_0\eta_0^T)w_1= y_1.
\end{align}
{\rm (ii)} Recall that, by Lemma \ref{3.2constr}, the maps $$U_1 \colon  H^2(\Dd,\Cc^{n-1}) \to \mathcal{K}_1,\quad U_2 \colon  H^2(\Dd,\Cc^{m-1})^\perp \to \mathcal{L}_1,$$ defined by 
$$U_1 \chi = V_0 \begin{pmatrix}
0 \\ \chi 
\end{pmatrix}=\bar{\alpha}_0\chi, \quad U_2 \psi = W_0^* \begin{pmatrix}
0 \\ \psi
\end{pmatrix}=\beta_0\psi \quad $$  for all $\chi \in H^2(\Dd,\Cc^{n-1})$ and all $\psi \in H^2(\Dd,\Cc^{m-1})^\perp,$ are unitaries. By the commutativity of the diagram \eqref{commdiagr},
 \begin{equation}\label{hf1g1} H_{F_1} = U_2^* \Gamma_1 U_1.\end{equation}

By Part (i), $x_1 \in \mathcal{K}_1$ and $y_1 \in \mathcal{L}_1$ and
by Proposition \ref{onxi}, 
$$	\xi_0 \telwe x_1 = \xi_0 \telwe v_1,\quad 	
 \bar{ \eta}_0 \telwe y_1 = \bar{ \eta}_0 \telwe w_1.$$ 
Thus, by Lemma \ref{connofschpairs1}, $(x_1,y_1)$ is a Schmidt pair for the operator $\Gamma_1$ corresponding to $t_1=\|\Gamma_1\|,$ that is,
\begin{equation}\label{schmg1}
\Gamma_1 x_1 =t_1y_1,\quad \Gamma_1^*y_1=t_1x_1.
\end{equation}

To prove that  the pair $(\hx_1,\hy_1)$ is a Schmidt pair for $H_{F_1}$ corresponding to $\|H_{F_1}\|=t_1$, we need to show that 
$$H_{F_1}\hx_1 = t_1 \hy_1,\; \text{and} \; H_{F_1}^*\hy_1 = t_1\hx_1. $$
By equations \eqref{hf1g1} and \eqref{x1=alpha0alphaTx1}, we have
\begin{align}\label{u2star} 
H_{F_1}\hat{x}_1&= H_{F_1}\alpha_0^T\hx_1\nonumber\vspace{2ex}\\&= U_2^*\Gamma_1 U_1 \alpha_0^T x_1 = U_2^* \Gamma_1 \bar{\alpha}_0\alpha_0^T x_1 \nonumber\vspace{2ex}\\
&= U_2^*\Gamma_1 x_1 = t_1 \beta_0^* y_1 = t_1 \hy_1.
\end{align}

Let us show that $H_{F_1}^*\hy_1 = t_1\hx_1.$  By equations \eqref{hf1g1} and
\eqref{x1=alpha0alphaTx1}, we have
\begin{align}\label{g1*}
H_{F_1}^*\hy_1 &=H_{F_1}^* \beta_0^* y_1\nonumber\vspace{2ex}\\
&= U_1^* \Gamma_1^* U_2 \beta_0^* y_1 = U_1^*\Gamma_1^*\beta_0\beta_0^* y_1\nonumber \vspace{2ex}\\ 
&= U_1^*\Gamma_1^* y_1= t_1 U_1^* x_1 = t_1 \alpha_0^T x_1 = t_1 \hx_1.
\end{align}

Therefore
$(\hx_1,\hy_1)$ is a Schmidt pair for $H_{F_1}$ corresponding to 
$\|H_{F_1}\|=t_1.$ \end{proof}

\begin{proposition}\label{x0wev1eta1wew1}
	Let $(\xi_0\telwe v_1,\bita_0\telwe w_1)$ be a Schmidt pair for $T_1$ corresponding to $t_1$ for some $v_1\in H^2(\Dd,\Cc^n),w_1\in H^2(\Dd,\Cc^m)^\perp,$ let $h_1 \in H^2(\Dd,\Cc)$ be the scalar outer factor of $\xi_0 \telwe v_1,$ let 
$$x_1 = (I_{n}- \xi_0 \xi_0^*)v_1,\quad y_1=(I_{m} - \bar{\eta}_0 \eta_0^T)w_1,$$
 and let 
$$\hx_1 = \alpha_0^T x_1,\quad \hy_1=\beta_0^*y_1. $$
Then
$$\|\hx_1 (z) \|_{\Cc^{n-1}} =  \|\hy_1(z)\|_{\Cc^{m-1}} = |h_1(z)|,$$
		$$\| x_1(z) \|_{\Cc^n} =  \|y_1(z)\|_{\Cc^m} = |h_1(z)|$$
and	
	$$\| \xi_0 (z) \we v_1(z) \|_{\we^2\Cc^n} = \| \bar{\eta}_0(z) \we w_1(z)\|_{\we^2\Cc^m} = |h_1(z)|$$
almost everywhere on $\Tt.$	
\end{proposition}

\begin{proof}
By Lemma \ref{schfohf1}, $(\hx_1,\hy_1)$ is a Schmidt pair for $H_{F_1}$ corresponding to $\|H_{F_1}\|=t_1$. Hence
$$H_{F_1}\hx_1 = t_1\hy_1 \quad \text{and}\quad H_{F_1}^* \hy_i = t_1 \hx_1. $$
By Theorem \ref{1.7},  for the Hankel operator $H_{F_1}$ and the Schmidt pair
$(\hx_1,\hy_1)$, we have 
 \begin{equation}\label{hatseq} \|\hy_1(z)\|_{\Cc^{m-1}}= \|\hx_1(z)\|_{\Cc^{n-1}} \end{equation} almost everywhere on $\Tt.$ 

By equations \eqref{x1=alpha0alphaTx1},
$$x_1= \bar{\alpha}_0 \alpha_0^Tx_1 =\bar{\alpha}_0 \hx_1,\quad y_1= \beta_0\beta_0^*y_1=\beta_0 \hy_1.$$
Since $\bar{\alpha}_0(z)$ and $\beta_0(z)$ are isometric for almost every $z\in \Tt$,
$$\|x_1(z)\|_{\Cc^n}=\|\hx_1(z)\|_{\Cc^{n-1}} \;\; \text{and} \;\; 
\|y_1(z)\|_{\Cc^m}=\|\hy_1(z)\|_{\Cc^{m-1}} $$
almost everywhere on $\Tt$. By equations \eqref{hatseq}, we deduce
\begin{equation}\label{x1isy1}
\|x_1(z)\|_{\Cc^n}=\|y_1(z)\|_{\Cc^m} 
\end{equation} almost everywhere on $\Tt.$

By Theorem \ref{T0compact}, $(\xi_0 \telwe \cdot)$ is an isometry from $\mathcal{K}_1$ to $X_1,$ and $(\bar{ \eta}_0\telwe \cdot)$ is an isometry from $\mathcal{L}_1$ to $Y_1.$ By Proposition \ref{onxi}, 
$$	\xi_0 \telwe x_1 = \xi_0 \telwe v_1,\quad 	
 \bar{ \eta}_0 \telwe y_1 = \bar{ \eta}_0 \telwe w_1.$$ 
Hence 
$$\begin{array}{lll}\|\xi_0(z)\we v_1(z)\|_{\we^2\Cc^n}&= \|\xi_0(z)\we x_1(z)\|_{\we^2\Cc^n}\vspace{2ex}=\|x_1(z)\|_{\Cc^n} \end{array}$$
almost everywhere on $\Tt.$  Also
$$\begin{array}{llll}
\| \bita_0(z) \we w_1(z)\|_{\we^2\Cc^m}&=\| \bita_0(z) \we y_1(z) \|_{\we^2\Cc^m}\vspace{2ex}=\|y_1(z)\|_{\Cc^m} 
\end{array} $$
almost everywhere on $\Tt.$  Thus, by equation \eqref{x1isy1},
  $$\| \xi_0 (z) \we v_1(z) \|_{\we^2\Cc^n} = \| \bar{\eta}_0(z) \we w_1(z)\|_{\we^2\Cc^m}$$almost everywhere on $\Tt.$	
Recall that $h_1$ is the scalar outer factor of $\xi_0 \telwe v_1$.
Hence
$$\| \xi_0 (z) \we v_1(z) \|_{\we^2\Cc^n} = \| \bar{\eta}_0(z) \we w_1(z)\|_{\we^2\Cc^m} = |h_1(z)|,$$
$$\| x_1(z) \|_{\Cc^n} =  \|y_1(z)\|_{\Cc^m} = |h_1(z)|$$ and 
$$\|\hx_1 (z) \|_{\Cc^{n-1}} =  \|\hy_1(z)\|_{\Cc^{m-1}} = |h_1(z)|$$ almost everywhere on $\Tt.$	\end{proof}


\index{level $j-$superoptimal error function} \index{$\mathcal{E}_j$}
\begin{definition}\label{epsilonj}
	Given $G\in H^\infty(\Dd,\CCmn)+C(\Tt,\CCmn)$ and $0\leq j \leq\min(m,n),$ define $\Omega_j$ to be \emph{the set of level $j$ superoptimal analytic approximants to $G$}, that is, the set of $Q\in H^\infty(\Dd,\CCmn)$ which minimize the tuple
	$$\big(s_0^\infty(G-Q),s_1^\infty(G-Q), \dots, s_j^\infty(G-Q)\big)  $$with respect to the lexicographic ordering over $Q \in H^\infty(\Dd,\CCmn).$  For $Q\in \Omega_j$ we call $G-Q$ a \emph{level $j$ superoptimal error function}, and we denote by $\mathcal{E}_j$ the \emph{set of all level $j$ superoptimal error functions}, that is 
	$$\mathcal{E}_j =\{G-Q \; : \; Q\in \Omega_j\}. $$
\end{definition}

\begin{proposition}\label{tildew1v1}
	Let $m,n$ be positive integers such that $\min(m,n)\geq2.$ Let \linebreak $G\in H^\infty(\Dd,\CCmn)+C(\Tt,\CCmn).$ In line with the algorithm from Subsection \ref{Alg_statement},
	let $Q_1 \in H^\infty(\Dd,\CCmn)$ satisfy	
	$$(G-Q_1)x_0 = t_0 y_0,\quad (G-Q_1)^*y_0=t_0x_0. $$ Let the spaces $X_1 , Y_1$ be given by
	$$ X_1 = \xi_0 \telwe H^2(\Dd,\Cc^n) \subset H^2(\Dd,\we^2\Cc^n), \quad Y_1 = \bar{\eta}_0 \telwe H^2(\Dd,\Cc^m)^\perp \subset H^2(\Dd,\we^2\Cc^m)^\perp,$$ 
	and consider the compact operator $T_1\colon X_1 \to Y_1$ given by 
$$T_1(\xi_0 \telwe x) = P_{Y_1} (\bar{\eta}_0 \telwe (G-Q_1)x)$$ 
for all $x \in H^2(\Dd,\Cc^n).$ 
Let  $(\xi_0 \telwe v_1,\bar{\eta}_0\telwe  w_1)$ be a Schmidt pair for the  operator $T_1$ corresponding to $t_1 = \|T_1\|,$  let $h_1 \in H^2(\Dd,\Cc)$ be the scalar outer factor of $\xi_0 \telwe v_1,$ let 
$$x_1 = (I_{n} - \xi_0 \xi_0^*)v_1, \quad y_1=(I_{m}-\bar{\eta}_0\eta_0^T)w_1 $$ and let $$\xi_1 =\frac{{x}_1}{h_1}, \quad \eta_1 =\frac{\bar{z}\bar{y}_1}{h_1}. $$
	
	Then, there exist unitary-valued functions $\tilde{V}_1, \tilde{W}_1$ of types $(n-1)\times(n-1),$\linebreak $(m-1)\times (m-1)$ respectively of the form \begin{equation}\label{V1}	
\tilde{V}_1 \stackrel{\emph{def}}{=} \begin{pmatrix}\alpha_0^T \xi_1 & \overline{\alpha}_1 \end{pmatrix}	\end{equation}and \begin{equation}\label{W1} \tilde{W}_1^T \stackrel{\emph{def}}{=} \begin{pmatrix} \beta_0^T \eta_1 & \overline{\beta}_1 \end{pmatrix},\end{equation}where $\alpha_1,\beta_1$ are inner, co-outer, quasi-continuous functions of types $(n-1)\times (n-2),$ \linebreak $(m-1)\times (m-2)$ respectively, and all minors on the first columns of $\tilde{V}_1,\tilde{W}_1^T$ are in $H^\infty.$ 

Furthermore, the set of all level $1$ superoptimal functions $\mathcal{E}_1$ satisfies  
	 \begin{equation}\label{g-qv0v1w0w1}
	\mathcal{E}_1 = W_0^* \begin{pmatrix} 1 & 0 \\ 0& \tilde{W}_1^*\end{pmatrix} \begin{pmatrix} t_0 u_0 & 0&0 \\ 0& t_1 u_1 &0\\ 0&0 & \left(F_2+H^\infty(\Dd,\Cc^{(m-2)\times(n-2)})\right)\cap B(t_1) \end{pmatrix}\begin{pmatrix} 1 & 0 \\ 0 & \tilde{V}_1^*  \end{pmatrix} V_0^* , \end{equation}
	where $F_2 \in   H^\infty(\Dd,\Cc^{(m-2)\times (n-2)})+C(\Tt,\Cc^{(m-2)\times (n-2)}),$ $u_1 = \frac{\bar{z} \bar{h}_1}{h_1}$ is a quasi-continuous unimodular function and $V_0,W_0^T$ are as in Theorem \ref{T0compact}, and $B(t_1)$ is the closed ball of radius $t_1$ in $L^\infty(\Tt,\Cc^{(m-2) \times (n-2)})$.
\end{proposition} 

\begin{proof}
By Theorem \ref{T0compact}, the following diagram commutes 
\begin{equation}\label{commdiagrhf1}
\begin{array}{clllll}
H^2(\Dd,\Cc^{n-1}) &\xrightarrow{U_1} & \mathcal{K}_1 &\xrightarrow{\xi_0 \telwe \cdot}& \xi_0 \telwe H^2 (\Dd, \Cc^n)=X_1\\
\Big\downarrow\rlap{$\scriptstyle H_{F_1}  $} & ~  &\Big\downarrow\rlap{$\scriptstyle \Gamma_1$}  &~&\hspace{3ex}\Big\downarrow\rlap{$\scriptstyle T_1$} \\
H^2(\Dd,\Cc^{m-1})^\perp &\xrightarrow{U_2}&	\mathcal{L}_1 &\xrightarrow{\bar{\eta}_0 \telwe \cdot } & \bar{\eta}_0 \telwe H^2 (\Dd, \Cc^m)^\perp =Y_1.
\end{array}\end{equation}
Let $\hx_1=\alpha_0^Tx_1,\; \hy_1 = \beta_0^* y_1.$
By Lemma \ref{schfohf1}, $(\hx_1,\hy_1)$ is a Schmidt pair for $H_{F_1}$ corresponding to $t_1$. 
By equations \eqref{x1=alpha0alphaTx1},
$$x_1= \bar{\alpha}_0 \alpha_0^Tx_1 =\bar{\alpha}_0 \hx_1\; \text{and} \; y_1= \beta_0\beta_0^*y_1=\beta_0 \hy_1.$$

 We want to apply Lemma \ref{2.2} to $H_{F_1}$ and the Schmidt pair $(\hx_1,\hy_1)$ to find unitary-valued functions $\tilde{V}_1,\tilde{W}_1$ such that, for any function $\tilde{Q}_1\in H^\infty(\Dd,\Cc^{(m-1)\times (n-1)})$ which is at minimal distance from $F_1,$ the following equation holds $$F_1-\tilde{Q}_1 = \tilde{W}_1^* \begin{pmatrix}
 t_1 u_1 & 0 \\ 0 & F_2 
 \end{pmatrix}\tilde{V}_1^*, $$
for some $F_2 \in   H^\infty(\Dd,\Cc^{(m-2)\times (n-2)})+C(\Tt,\Cc^{(m-2)\times (n-2)}).$
For this purpose we find the inner-outer factorisations $\hat{x}_1$ and 
$\bar{z}\bar{\hy}_1$.  By Proposition \ref{x0wev1eta1wew1},
\begin{equation}\label{h1common}
\begin{aligned}&\|\hx_1(z)\|_{\Cc^{n-1}}=\|x_1(z)\|_{\Cc^n}=\| \xi_0(z)\telwe v_1(z)\|_{\we^2{\Cc^n}} = |h_1(z)|\;\\
\text{and}\\
 &\| \hy_1(z)\|_{\Cc^{m-1}}= \|y_1(z)\|_{\Cc^m}=\| \bar{\eta}_0(z) \telwe w_1(z)\|_{\we^2\Cc^m} =|h_1(z)|\end{aligned}\end{equation} almost everywhere on $\Tt.$ 
Equations \eqref{h1common} imply that $h_1\in H^2(\Dd,\Cc)$ is the scalar outer factor of both $\hat{x}_1$ and $\bar{z}\bar{\hat{y}}_1.$ By Lemma \ref{2.2},  $\hat{x}_1,\bar{z}\bar{\hat{y}}_1$ admit the inner-outer factorisations 
$$\hat{x}_1 = \hat{\xi}_1 h_1, \quad \bar{z}\bar{\hat{y}}_1=\hat{\eta}_1 h_1 ,$$ for some inner vector-valued $\hat{\xi}_1\in H^\infty(\Dd,\Cc^{n-1})$ and $\hat{\eta}_1 \in H^\infty(\Dd,\Cc^{m-1}). $ 
Recall that  
$$\hat{x}_1 =  \alpha_0^T x_1 =  \alpha_0^T \xi_1 h_1,\quad \bar{z}\bar{\hat{y}}_1=\bar{z}\beta_0^T \bar{y}_1 = \beta_0^T \eta_1 h_1, $$ 
which imply
$$ \hat{\xi}_1 = \alpha_0^T\xi_1 \quad \text{and}\quad \hat{\eta}_1 = \beta_0^T \eta_1  .$$ 

 Let us  show that $\alpha_0^T \xi_1,\;\beta_0^T \eta_1 $ are inner in order to apply Lemma \ref{2.2}. 

\noindent Recall that, since $V_0,W_0^T$ are unitary-valued, we have 
$$I_n -\xi_0 \xi_0^* =\bar{\alpha}_0 \alpha_0^T, \quad I_m - \bar{\eta}_0 \eta_0^T = \beta_0 \beta_0^* .$$
Therefore
$$x_1 = (I_{n} - \xi_0 \xi_0^*)v_1=\bar{\alpha}_0 \alpha_0^Tv_1,\quad y_1= (I_{m} - \bita_0\eta_0^T)w_1=\beta_0\beta_0^*w_1 .$$ 
Then,
\begin{equation}\label{aox1aov1}\alpha_0^Tx_1 = \alpha_0^T v_1, \quad \beta_0^T \bar{y}_1 =\beta_0^T\bar{w}_1  ,\end{equation} and since
$$\xi_1 =\frac{x_1}{h_1},\quad \eta_1 = \frac{\bar{z}\bar{y}_1}{h_1}, $$ the functions
$$\alpha_0^T\xi_1=\frac{\alpha_0^Tv_1}{h_1},\quad  \beta_0^T\eta_1 = \frac{\beta_0^T\bar{z}\bar{w}_1}{h_1}$$ are analytic. 
Furthermore, by Proposition \ref{x0wev1eta1wew1},  
$$\|x_1(z)\|_{\Cc^n}= \|y_1(z)\|_{\Cc^m}=|h_1(z)|= \|\hx_1(z)\|_{\Cc^{n-1}}=\| \hy_1(z)\|_{\Cc^{m-1}}$$
almost everywhere on $\Tt$. 
Thus
$$\|\alpha_0^T(z)x_1(z)\|_{\Cc^{n-1}}=\|\alpha_0^T(z)v_1(z)\|_{\Cc^{n-1}}=|h_1(z)|$$
and
$$ \|\beta_0^T(z)\bar{z} \bar{y}_1(z)\|_{\Cc^{m-1}}=\|\beta_0^T (z)\bar{z}\bar{w}_1(z)\|_{\Cc^{m-1}}=|h_1(z)| $$
almost everywhere on $\Tt.$ Hence
$$\|\alpha_0^T(z)\xi_1(z)\|_{\Cc^{n-1}}=1,\quad \|\beta_0^T(z)\eta_1(z)\|_{\Cc^{m-1}}=1 $$
almost everywhere on $\Tt.$ Therefore $\alpha_0^T\xi_1,\; \beta_0^T\eta_1$ are inner functions. By Lemma \ref{2.2}, there exist inner, co-outer, quasi-continuous functions $\alpha_1,\beta_1$ of types $(n-1)\times (n-2)$ and $(m-1)\times (m-2)$ respectively such that 
$$\tilde{V}_1 =\begin{pmatrix}\alpha_0^T \xi_1 & \overline{\alpha}_1 \end{pmatrix}	,\quad \tilde{W}_1^T = \begin{pmatrix} \beta_0^T \eta_1 & \overline{\beta}_1 \end{pmatrix}$$ are unitary-valued and all minors on the first columns are in $H^\infty.$ 
Furthermore, by Lemma \ref{2.2}, every $\hat{Q}_1\in H^\infty(\Dd,\Cc^{(m-1)\times(n-1)})$ which is at minimal distance from $F_1$ satisfies 
$$F_1-\hat{Q}_1 = \tilde{W}_1^* \begin{pmatrix}
t_1 u_1 & 0 \\
0 & F_2
\end{pmatrix}\tilde{V}_1^*, $$where $F_2 \in   H^\infty(\Dd,\Cc^{(m-2)\times (n-2)})+C(\Tt,\Cc^{(m-2)\times (n-2)})$ and $u_1$ is a quasi-continuous unimodular function given by $u_1 = \frac{\bar{z} \bar{h}_1}{h_1}.$ 

By Lemma \ref{f+hinfty}, the set $$\tilde{\mathcal{E}}_{0} =\{F_{1} - \hat{Q} : \hat{Q} \in H^\infty(\Dd,\Cc^{(m-1)\times (n-1)}), \| F_{1} - \hat{Q}\|_{L^\infty}=t_{1}  \}$$ satisfies  $$\tilde{\mathcal{E}}_{0} = \tilde{W}_{1}^* 
\begin{pmatrix}
t_{1}u_{1} & 0 \\
0 & \left(F_2+H^\infty(\Dd,\Cc^{(m-2)\times(n-2)})\right) \cap B(t_1)
\end{pmatrix}V_{1}^*,
$$ 
for some 
$F_2 \in H^\infty(\Dd,\Cc^{(m-2)\times(n-2)}) + C(\Tt, \Cc^{(m-2)\times(n-2)})$ and for the closed ball $B(t_1)$ of radius $t_1$ in $L^\infty(\Tt,\Cc^{(m-2)\times (n-2)}).$ Thus, by Lemma \ref{f+hinfty}, $\mathcal{E}_1$ admits the factorisation \eqref{g-qv0v1w0w1} as claimed.
\end{proof} 

\begin{proposition}\label{g-q1y1t1x1}
	
Suppose the function $Q_2\in H^\infty(\Dd,\CCmn)$ minimises $$(s_0^\infty(G-Q),s_1^\infty(G-Q)).$$ Then $Q_2$ satisfies
$$(G-Q_2)x_0 = t_0y_0,\quad (G-Q_2)^*y_0 = t_0 x_0  $$and 
$$(G-Q_2)x_1 = t_1y_1,\quad (G-Q_2)^*y_1=t_1x_1, $$where $x_0,x_1,y_0,y_1,t_0,t_1$ are as in Theorem \ref{T0compact}.	
	
\end{proposition}
\begin{proof}
	Let $(x_0,y_0)$ be a Schmidt pair for the Hankel operator $H_G$ corresponding to \linebreak$\|H_G\|=t_0.$ Then, by Theorem \ref{1.7}, every $Q_2 \in H^\infty(\Dd,\CCmn)$ which is at minimal distance from $G$ satisfies 
	$$(G-Q_2)x_0 = t_0y_0,\quad (G-Q_2)^*y_0 = t_0 x_0 , $$ and, by Lemma \ref{2.2}, 
	$$W_0 (G-Q_2)V_0 = \begin{pmatrix}
	t_0u_0 & 0 \\
	0 & F_1
	\end{pmatrix} ,$$where $F_1 \in H^\infty(\Dd,\Cc^{(m-1)\times (n-1)})+C(\Tt,\Cc^{(m-1)\times (n-1)}).$

Moreover, by Lemma \ref{f+hinfty}, the set  $\mathcal{E}_0 = \{ G-Q : Q \in \Omega_0\}$ of all level $0$ superoptimal error functions satisfies  
	\begin{equation}\label{wevv}W_0 \mathcal{E}_0 V_0 = \begin{pmatrix}
	t_0 u_0  & 0\\ 0 & F_1 +H^\infty(\Dd, \Cc^{m-1\times n-1})
	\end{pmatrix}\cap B(t_0).\end{equation}
Suppose $Q_2\in \Omega_0$. Then  
$$ W_0 (G-Q_2) V_0 =\begin{pmatrix}
\eta_0^T \\\beta_0^* 
\end{pmatrix}(G-Q_2) \begin{pmatrix}
\xi_0 & \balpha_0
\end{pmatrix}= \begin{pmatrix}
\eta_0^T (G-Q_2)\xi_0  & \eta_0^T (G-Q_2)\balpha_0\\
\beta_0^* (G-Q_2)\balpha_0 & \beta_0^*(G-Q_2)\balpha_0
\end{pmatrix}.$$ By equation \eqref{wevv}, for $\tilde{Q}_1\in H^\infty(\Dd,\Cc^{(m-1)\times(n-1)})$ at minimal distance from $F_1,$ 

\begin{equation}\label{wog-qo}\begin{pmatrix}
\eta_0^T (G-Q_2)\xi_0  & \eta_0^T (G-Q_2)\balpha_0\\
\beta_0^* (G-Q_2)\balpha_0 & \beta_0^*(G-Q_2)\balpha_0
\end{pmatrix}=\begin{pmatrix}
	t_0 u_0 &0  \\ 0 & F_1 -\tilde{Q}_1
\end{pmatrix} \end{equation}

 \noindent Note that, by Theorem \ref{nehtmatr}, $$\| F_1 - \tilde{Q}_1\|_\infty=\|H_{F_1}\|, $$and, by Theorem \ref{T0compact} (part (v)), $\|H_{F_1}\|=t_1.$ 
	
	Consideration of the $(2,2)$ entries of equation \eqref{wog-qo} yields
	\begin{equation}\label{fq-q1} F_1-\tilde{Q}_1 =\beta_0^* (G-Q_2) \bar{\alpha}_0 .\end{equation}   Note that, if $(\hat{x}_1,\hat{y}_1)$ is a Schmidt pair for $H_{F_1}$ corresponding to $t_1=\|H_{F_1}\|,$ then, by Theorem \ref{1.7}, 
	
	$$(F_1-\tilde{Q}_1)\hat{x}_1 = t_1 \hat{y}_1,\quad  (F-\hat{Q}_1)^*\hat{y}_1=t_1\hat{x}_1 .$$In view of equation \eqref{fq-q1}, the latter equations imply
	\begin{equation}\label{bg-q1a}\beta_0^* (G-Q_2) \bar{\alpha}_0 \hat{x}_1 = t_1 \hat{y}_1,\end{equation}and \begin{equation}\label{bg-q1ab} \alpha_0^T (G-Q_2)^*\beta_0\hy_1 = t_1\hat{x}_1. \end{equation}

\noindent	By Lemma \ref{schfohf1}, we may choose the Schmidt pair for $H_{F_1}$ corresponding to $\|H_{F_1}\|$ to be \begin{equation}\label{schmhf1} \hat{x}_1 =\alpha_0^T x_1,\quad \hat{y}_1  = \beta_0^* y_1 .\end{equation} 
Recall that, by equations \eqref{x1=alpha0alphaTx1},
\begin{equation}\label{ex1} x_1 =\bar{\alpha}_0 \alpha_0^Tx_1\end{equation}
and 
\begin{equation}\label{yai1} y_1 =\beta_0 \beta_0^* y_1 .\end{equation}

	In view of equations \eqref{bg-q1a}  and \eqref{schmhf1}, we obtain 
	$$ \beta_0^* (G-Q_2) \bar{\alpha}_0 \alpha_0^Tx_1= t_1 \beta_0^*y_1.$$ Multiplying both sides of the latter equation by $\beta_0,$ we have
	$$ \beta_0 \beta_0^*(G-Q_2)\balpha_0 \alpha_0^Tx_1  = t_1\beta_0 \beta_0^* y_1,$$ which, by equation \eqref{ex1}, implies 
	$$ \beta_0 \beta_0^*(G-Q_2)x_1 = t_1\beta_0\beta_0^*y_1 ,$$ or equivalently,
	$$\beta_0\beta_0^* \bigg( (G-Q_2)x_1-t_1y_1\bigg)=0. $$ 
Since, by Theorem \ref{T0compact}, $U_2^* =M_{\beta_0\beta_0^*}$ is unitary, the latter equation yields
	$$ (G-Q_2)x_1 = t_1y_1.$$
	
	Moreover, by equations \eqref{bg-q1ab} and \eqref{schmhf1}, we obtain $$\alpha_0^T(G-Q_2)^*\beta_0\beta_0^*y_1= t_1 \alpha_0^Tx_1.$$ Multiplying both sides of the latter equation by $\balpha_0,$ we have 
	$$ \balpha_0\alpha_0^T(G-Q_2)^*\beta_0\beta_0^*y_1= t_1 \balpha_0\alpha_0^Tx_1.$$ 
In view of equation \eqref{yai1}, the latter expression is equivalent to the equation
	$$ \balpha_0\alpha_0^T(G-Q_2)^*y_1= t_1 \balpha_0\alpha_0^Tx_1, $$
or equivalently,
	$$\balpha_0\alpha_0^T \bigg( (G-Q_2)^*y_1-t_1x_1 \bigg) =0.$$
Since, by Theorem \ref{T0compact},  $U_1^* = M_{\balpha_0\alpha_0^T}$ is unitary, the latter equation yields
	$$ (G-Q_2)^*y_1=t_1x_1.$$ Therefore $Q_2$ satisfies the required equations.
\end{proof}

The next few propositions are in preparation for Theorem \ref{T2compactt} on the compactness of $T_2$.

\begin{proposition}\label{beta1h2}
For a thematic completion of the inner matrix-valued function $\beta_0^T\eta_1$ of the form
	 $\tilde{W}_1^T =\big( \begin{matrix} \beta_0^T \eta_1 & \bar{ \beta}_{1}\end{matrix}\big)$, where $\beta_1$ is an inner, co-outer, quasi-continuous function of type $(m-1) \times (m-2),$ the following equation holds 
$$\beta_1^* H^2(\Dd, \Cc^{m-1})^\perp=  H^2(\Dd, \Cc^{m-2})^\perp.$$	
	
\end{proposition}

\begin{proof}
	By virtue of the fact that complex conjugation is a unitary operator on $L^2(\Tt,\Cc^m),$ an equivalent statement is that $\beta_1^T z H^2(\Dd,\Cc^{m-1}) = z H^2(\Dd,\Cc^{m-2}).$ By Lemma \ref{L6.2}, there exists a matrix-valued function $B_1 \in H^\infty( \Dd,\Cc^{(m-2)\times (m-1)})$ such that $$B_1 \beta_1 = I_{m-2}$$ or, equivalently, $$ \beta_1^T B_1^T = I_{m-2}.$$
	
	\noindent Let $f \in z H^2(\Dd,\Cc^{m-2}).$ Then, $$ f = (\beta_1^T B_1^T) f  \in \beta_1^T B_1^T z H^2(\Dd,\Cc^{m-2}) \subseteq \beta_1^T z H^2(\Dd,\Cc^{m-1}).$$ Hence $$ z H^2(\Dd,\Cc^{m-2}) \subseteq \beta_1^T z H^2(\Dd,\Cc^{m-1}). $$
	
	\noindent Note that, since $\beta_1 \in H^\infty(\Dd, \Cc^{(m-1) \times (m-2)})$, we have $$ \beta_1^T z H^2(\Dd,\Cc^{m-1}) \subseteq z H^2(\Dd,\Cc^{m-2}).$$ Thus 
	$$\beta_1^T z H^2(\Dd,\Cc^{m-1}) = z H^2(\Dd,\Cc^{m-2}).$$
\end{proof}

\begin{lemma}\label{a1h2}
For a thematic completion of the inner matrix-valued function  $\alpha_0^T\xi_1$ of the form $ \tilde{V}_1 =\big(\begin{matrix}
	\alpha_0^T \xi_1 & \bar{ \alpha}_{1}
	\end{matrix}\big), $ where $\alpha_1$ is an inner, co-outer, quasi-continuous function of type $(n-1) \times (n-2),$ the following equation holds
	$$ \alpha_1^T H^2(\Dd,\Cc^{n-1}) = H^2(\Dd,\Cc^{n-2}). $$

\end{lemma}

\begin{proof}
	\noindent By Lemma \ref{L6.2}, for the given $\alpha_1,$ there exists $A_1\in H^\infty(\Dd, \Cc^{(n-2)\times (n-1) })$ such that $A_1\alpha_1 = I_{n-2}.$ Equivalently, $\alpha_1^T A_1^T = I_{n-2}.$ \medskip

	\noindent Let $g\in H^2(\Dd,\Cc^{n-2}).$ Then $g = (\alpha_1^T A_1^T) g \in \alpha_1^T A_1^T H^2(\Dd,\Cc^{n-2}),$
	which implies that \linebreak$ g \in \alpha_1^T H^2(\Dd,\Cc^{n-1}).$ Hence $H^2(\Dd,\Cc^{n-2}) \subseteq \alpha_1^T H^2(\Dd,\Cc^{n-1}).$\medskip

	\noindent For the reverse inclusion, note that, since 
$\alpha_1 \in H^\infty(\Dd, \Cc^{(n-1) \times (n-2)})$, we have
 $$\alpha_1^T H^2(\Dd,\Cc^{n-1}) \subseteq  H^2(\Dd,\Cc^{n-2}).$$ Thus $$ \alpha_1^T H^2(\Dd,\Cc^{n-1}) = H^2(\Dd,\Cc^{n-2}). $$
\end{proof}

\begin{remark}\label{V1V1*unit}
Let $V_0$ and $\tilde{V}_1$ be given by equations \eqref{V0W0} and \eqref{V1} respectively and let 
$V_1 = \begin{pmatrix}
1 & 0 \\ 0 & \tilde{V}_1
\end{pmatrix}.$ 
Since $V_0,$ $\tilde{V}_1$ and $V_1$ are unitary-valued, we have 
	\begin{equation}\label{V0V0*}
	I_n = V_0 V_0^*= \xi_0 \xi_0^* + \bar{\alpha}_0 \alpha_0^T,
	 \end{equation} 	
	\begin{equation}\label{V1V1*} 
I_{n-1} = \tilde{V}_1 \tilde{V}_1^* = \alpha_0^T \xi_1\xi_1^*
\bar{\alpha}_0 + \bar{\alpha}_1 \alpha_1^T.
  \end{equation} 
	\end{remark}
		
\begin{lemma}\label{xi=A1A1*}
Let $V_0$ and $\tilde{V}_1$ be given by equations \eqref{V0W0} and \eqref{V1} respectively. Let $A_1= \alpha_0\alpha_1$.
Then
\begin{equation} \label{1-xi=A1A1*}
I_n - \xi_0 \xi_0^* -\xi_1 \xi_1^*=
\bar{\alpha}_0 \bar{\alpha}_1 \alpha_1^T\alpha_0^T = \bar{A_1}A_1^T.
\end{equation}
 almost everywhere on $\Tt$.
\end{lemma}

\begin{proof} By equation \eqref{V1V1*}
$$ \bar{\alpha}_1 \alpha_1^T= I_{n-1} - \alpha_0^T \xi_1\xi_1^*
\bar{\alpha}_0,$$
thus 
$$
\bar{\alpha}_0 \bar{\alpha}_1\alpha_1^T\alpha_0^T =\bar{\alpha}_0( I_{n-1} - \alpha_0^T \xi_1\xi_1^* \bar{\alpha}_0)\alpha_0^T .$$
By equation \eqref{V0V0*},
$$ \bar{\alpha}_0 \alpha_0^T = I_n - \xi_0 \xi_0^* .$$
Hence 
$$
\bar{\alpha}_0 \bar{\alpha}_1 \alpha_1^T \alpha_0^T =(I_n - \xi_0 \xi_0^* )
- (I_n - \xi_0 \xi_0^* )\xi_1\xi_1^* (I_n - \xi_0 \xi_0^* ).$$
Since, by Proposition \ref{onxi}, the set $\{\xi_0(z), \xi_1(z) \}$ is orthonormal in $\Cc^m$ for almost every $z \in \Tt$,
$$
\bar{\alpha}_0 \bar{\alpha}_1 \alpha_1^T\alpha_0^T =
I_n - \xi_0 \xi_0^* -\xi_1 \xi_1^*
$$
almost everywhere on $\Tt$. \end{proof}

Let us state certain identities that are useful for the next statements.

\begin{remark}\label{w1w1*unit}
Let $W_0$ and $\tilde{W}_1$ be given by equations \eqref{V0W0} and \eqref{W1} respectively and let $W_1 = \begin{pmatrix}
1 & 0 \\ 0 & \tilde{W}_1
\end{pmatrix}.$ Then
	\begin{equation}\label{w0w0*}
	I_m = W_0^* W_0= \bar{\eta}_0 \eta_0^T + \beta_0 \beta_0^*,
	 \end{equation} 
	
	\begin{equation}\label{w1w1*} 
I_{m-1} = \tilde{W}_1^* \tilde{W}_1 = \beta_0^* \bar{\eta}_1 \eta_1^T \beta_0 + \beta_1 \beta_1^*. \end{equation} 

\begin{equation} \label{W1W0*}
\begin{array}{cllllll}
	 W_0^*  \begin{pmatrix}
	1 & 0 \\ 0 & \tilde{W}_1^* 
	\end{pmatrix}=
\begin{pmatrix}
	\bar{\eta}_0 & \beta_0 
	\end{pmatrix} \begin{pmatrix} 1 & 0 \\ 0 & \begin{pmatrix}
	\beta_0^* \bar{\eta}_1 & \beta_1 
	\end{pmatrix} \end{pmatrix}= \begin{pmatrix}
		\bar{\eta}_0 & \beta_0 \beta_0^* \bar{\eta}_1 & \beta_0 \beta_1 
		 \end{pmatrix}.
\end{array}
\end{equation}
  
	$$\begin{array}{cllllll}
	 W_0^*  \begin{pmatrix}
	1 & 0 \\ 0 & \tilde{W}_1^* 
	\end{pmatrix} 
        \begin{pmatrix}
	1 & 0 \\ 0 & \tilde{W}_1
	\end{pmatrix} W_0  \vspace{2ex} &= 
        \begin{pmatrix} \bar{\eta}_0 & \beta_0 
	\end{pmatrix} \begin{pmatrix} 1 & 0 \\ 0 & \begin{pmatrix}
	\beta_0^* \bar{\eta}_1 & \beta_1 
	\end{pmatrix} \end{pmatrix}   \begin{pmatrix} 1 & 0 \\0 & \begin{pmatrix}
	 \eta_1^T  \beta_0 \\ \beta_1^* 
	\end{pmatrix} \end{pmatrix}  \begin{pmatrix}
	\eta_{0}^T \\ \beta_0^* 
		\end{pmatrix}\vspace{2ex} \\\end{array}$$ 

$$\begin{array}{llll}&= \begin{pmatrix}
		\bar{\eta}_0 & \beta_0 \beta_0^* \bar{\eta}_1 & \beta_0 \beta_1 
		 \end{pmatrix} \begin{pmatrix}
		 \eta_0^T  \\ \eta_1^T \beta_0 \beta_0^* \\ \beta_1^* \beta_0^*
		 \end{pmatrix} \vspace{2ex} \\&= \bar{\eta}_0 \eta_0^T + \beta_0 \beta_0^* \bar{\eta}_1  \eta_1^T \beta_0 \beta_0^* + \beta_0 \beta_1 \beta_1^* \beta_{0}^* .\end{array}$$   Furthermore, \begin{equation}\label{b0b1conn}\bar{\eta}_0 \eta_0^T + \beta_0 \beta_0^* \bar{\eta}_1  \eta_1^T \beta_0 \beta_0^* + \beta_0 \beta_1 \beta_1^* \beta_{0}^* =\bar{\eta}_0 \eta_0^T + \beta_0 (I_{m-1} - \beta_1 \beta_1^* + \beta_1 \beta_1^* ) \beta_0^* \vspace{2ex} = I_m .\end{equation}  
	\end{remark}
Equations \eqref{w0w0*} and \eqref{w1w1*} follow from the facts that $W_0, \tilde{W}_1$ and $W_1$ are unitary-valued on $\Tt$.
Equations \eqref{b0b1conn} follow from equations (\ref{w0w0*}) and (\ref{w1w1*}).
		
\begin{lemma}\label{eta=B1B1*}
Let $W_0$ and $\tilde{W}_1$ be given by equations \eqref{V0W0} and \eqref{W1} respectively. Let $B_1= \beta_0 \beta_1$.
Then
\begin{equation} \label{1-eta=B1B1*}
I_m - \bar{\eta}_0 \eta_0^T -\bar{\eta}_1 \eta_1^T=
\beta_0 \beta_1 \beta_1^*\beta_0^*= B_1 B_1^*.
\end{equation}
 almost everywhere on $\Tt$.
\end{lemma}

\begin{proof} By equation \eqref{w1w1*}
$$\beta_1 \beta_1^*= I_{m-1} -\beta_0^* \bar{\eta}_1 \eta_1^T \beta_0 ,$$
thus 
$$
\beta_0 \beta_1 \beta_1^*\beta_0^*= \beta_0(I_{m-1} -\beta_0^* \bar{\eta}_1 \eta_1^T \beta_0 )\beta_0^*.$$
By equation \eqref{w0w0*},
$$ \beta_0 \beta_0^*= I_m - \bar{\eta}_0 \eta_0^T.$$
Hence 
$$
\beta_0 \beta_1 \beta_1^*\beta_0^*= (I_m - \bar{\eta}_0 \eta_0^T) -
(I_m - \bar{\eta}_0 \eta_0^T)\bar{\eta}_1 \eta_1^T  (I_m - \bar{\eta}_0 \eta_0^T).$$
Since, by Proposition \ref{onxi}, the set $\{\bar{\eta}_0(z), \bar{\eta}_1(z) \}$ is orthonormal in $\Cc^m$ for almost every $z \in \Tt$,
$$
\beta_0 \beta_1 \beta_1^*\beta_0^* =
I_m - \bar{\eta}_0 \eta_0^T -\bar{\eta}_1 \eta_1^T
$$
almost everywhere on $\Tt$. \end{proof}

\begin{proposition}\label{xi12telweunit}
With the notation of Proposition \ref{tildew1v1}, let unitary completions of $\xi_0$ and $\alpha_0^T\xi_1$ be given by 
	$$V_0 = \begin{pmatrix}
	\xi_0 & \bar{\alpha_{0}}
	\end{pmatrix}, \quad \tilde{V}_1 =  \begin{pmatrix}\alpha_0^T\xi_1 & \bar{\alpha_{1}}
	\end{pmatrix} ,$$ 
where $\alpha_0, \alpha_1$ are inner, co-outer, quasi-continuous matrix-valued functions  of types $n\times (n-1)$ and $(n-1)\times (n-2)$ respectively. Let $$V_1= \begin{pmatrix}
		1 & 0 \\ 0 & \tilde{V}_1\end{pmatrix}$$ and let 
\begin{equation}\label{K2}
\mathcal{K}_2 = V_0 V_1 \begin{pmatrix}
	0_{2\times 1} \\ H^2(\Dd,\Cc^{n-2})
	\end{pmatrix}.
\end{equation}
Then $$\xi_0 \telwe \xi_1 \telwe   H^2(\Dd,\Cc^n) = \xi_0 \telwe \xi_1 \telwe  \mathcal{K}_2$$ and the operator $(\xi_0 \telwe \xi_1 \telwe\cdot)\colon \mathcal{K}_2 \to \xi_0 \telwe \xi_1 \telwe H^2(\Dd,\Cc^n)$ is unitary. 
\end{proposition}

\begin{proof} First let us prove that 
$$\xi_0 \telwe \xi_1 \telwe H^2(\Dd,\Cc^n) = \xi_0 \telwe \xi_1 \telwe \mathcal{K}_2. $$ 
Recall that $A_1 = \alpha_0\alpha_1$.
Observe that, by definition, 
\begin{equation} \label{K2A}
\mathcal{K}_2 = \balpha_0 \balpha_1 H^2(\Dd,\Cc^{n-2})= \bar{A}_1 H^2(\Dd,\Cc^{n-2}).
\end{equation} 
By Lemmas \ref{a0h2} and \ref{a1h2}, 
\begin{equation}\label{alphaH2-K2}
\begin{array}{lll}
H^2(\Dd,\Cc^{n-2})&= \alpha_1^T H^2(\Dd,\Cc^{n-1})\\
&=  \alpha_1^T \alpha_0^T H^2(\Dd,\Cc^{n})\\
&= A_1^T H^2(\Dd,\Cc^{n}).
\end{array}
\end{equation}
By equations \eqref{K2A} and  \eqref{alphaH2-K2},
\begin{equation}\label{K2-2}
\mathcal{K}_{2} = \balpha_0 \balpha_1 H^2(\Dd,\Cc^{n-2})= \bar{A}_{1} A_1^T H^2(\Dd,\Cc^{n}).
\end{equation}
By Lemma \ref{xi=A1A1*},
\begin{equation}\label{AA1-K} 
\bar{A}_1 A_1^T = I_n - \sum\limits_{k=0}^1 \xi_k \xi_k^*.
\end{equation} 
By Proposition \ref{onxi}, $\{\xi_i(z)\}_{i=0}^1$ is an orthonormal set in $\Cc^n$ for almost every $z \in \Tt.$
Therefore, by equations \eqref{K2-2} and \eqref{AA1-K},
\begin{equation}\label{K2H2}
\begin{array}{lll}
\xi_0 \telwe \xi_1 \telwe \mathcal{K}_{2} &= \xi_0 \telwe \xi_1 \telwe \bar{A}_{1} A_1^T H^2(\Dd,\Cc^{n})\\
&= \xi_0 \telwe  \xi_1 \telwe (I_n - \sum\limits_{k=0}^1 \xi_k \xi_k^*)H^2(\Dd,\Cc^n)\\
&= \xi_0 \telwe \xi_1 \telwe H^2(\Dd,\Cc^n). 
\end{array}
\end{equation}

 Let us show that the operator $(\xi_0 \telwe \xi_1 \telwe\cdot)\colon \mathcal{K}_2 \to \xi_0 \telwe \xi_1 \telwe H^2(\Dd,\Cc^n)$ is unitary. The foregoing paragraph asserts that the operator is surjective. It remains to be shown that it is an isometry.  To this end, let $f \in \mathcal{K}_2.$ Then
	
$$\begin{array}{cllllllll}
\| \xi_0 \telwe \xi_1\telwe f\|_{L^2(\Tt,\we^3\Cc^n)}^2 
&= \langle \xi_0 \telwe \xi_1\telwe f, \xi_0 \telwe \xi_1\telwe f \rangle_{L^2(\Tt,\we^3\Cc^n)} \vspace{2ex} \\
&= \displaystyle \frac{1}{2\pi} \int_0^{2\pi} \langle \xi_0(\eiu) \telwe \xi_1(\eiu) \telwe f(\eiu), \xi_0(\eiu) \telwe \xi_1(\eiu) \telwe f(\eiu) \rangle_{\we^3\Cc^n}\; d\theta.\end{array}$$ By Proposition \ref{we}, the latter integral is equal to

$$ \displaystyle \frac{1}{2\pi} \int_0^{2\pi} \det \begin{pmatrix}
\langle \xi_0 (\eiu), \xi_0 (\eiu) \rangle_{\Cc^n} & \langle \xi_0(\eiu) , \xi_1 (\eiu) \rangle_{\Cc^n} &\langle \xi_0(\eiu) , f(\eiu) \rangle_{\Cc^n}  \\
\langle \xi_1(\eiu) , \xi_0 (\eiu) \rangle_{\Cc^n} & \langle \xi_1(\eiu) , \xi_1 (\eiu) \rangle_{\Cc^n} & \langle \xi_1(\eiu) , f (\eiu) \rangle_{\Cc^n} \\
\langle f(\eiu) , \xi_0 (\eiu) \rangle_{\Cc^n} & \langle f(\eiu) , \xi_1(\eiu) \rangle_{\Cc^n} & \langle f(\eiu) , f(\eiu) \rangle_{\Cc^n}
\end{pmatrix}\; d\theta. 
 $$
\noindent Note that, by Proposition \ref{onxi}, $  \{\xi_0(\eiu), \xi_1(\eiu)\} $ is an orthonormal set for almost all $\eiu$ on $\Tt.$ Moreover, since $\mathcal{K}_2= \bar{ \alpha}_{0} \bar{\alpha}_1 H^2(\Dd,\Cc^{n-2}) ,$ then $f=\bar{\alpha}_0 \bar{\alpha}_1 \varphi$ for some $\varphi \in  H^2(\Dd,\Cc^{n-2}).$ Hence $$\begin{array}{lll}\langle \xi_0 (\eiu), f(\eiu) \rangle_{\Cc^n} &= \langle \xi_0(\eiu), \bar{ \alpha}_{0}(\eiu)\bar{\alpha}_1(\eiu) \varphi(\eiu)\rangle_{\Cc^n}\vspace{2ex}\\  &= \langle \alpha_0^T(\eiu) \xi_0(\eiu) , \bar{\alpha}_1(\eiu)\varphi(\eiu) \rangle_{\Cc^{n-1}}\vspace{2ex}\\&=0\end{array} $$almost everywhere on $\Tt,$ since $ V_0$ is unitary-valued. Similarly, since $\tilde{V}_1$ is unitary valued, we deduce that 
$$\langle \xi_1(\eiu), f(\eiu) \rangle_{\Cc^n}= \langle \alpha_1^T(\eiu)\alpha_0^T(\eiu) \xi_1(\eiu), \varphi(\eiu)\rangle_{\Cc^{n-2}}= 0 $$almost everywhere on $\Tt.$
 Therefore 

$$ \| \xi_0 \telwe \xi_1 \telwe f\|_{L^2(\Tt,\we^3\Cc^n)}^2 =\frac{1}{2\pi} \int_0^{2\pi} \det \begin{pmatrix}
1 & 0 &0  \\
0 & 1 & 0\\
0 & 0 &  \|f(\eiu)\|_{\Cc^n}^2
\end{pmatrix}\; d\theta= \|f\|_{L^2(\Tt,\Cc^n)}^2,$$that is, $(\xi_0 \telwe \xi_1 \telwe\cdot)\colon \mathcal{K}_2 \to \xi_0 \telwe \xi_1 \telwe H^2(\Dd,\Cc^n)$ is an isometric operator. Thus, the operator $(\xi_0 \telwe \xi_1 \telwe \cdot)\colon \mathcal{K}_2 \to \xi_0 \telwe \xi_1 \telwe H^2(\Dd,\Cc^n)$ is unitary. 
\end{proof}


\begin{proposition}\label{eta0eta1b0b1}
	Let $\eta_0, \eta_1$ be defined by equations \eqref{xi0eta0} and \eqref{311} respectively, and let $\beta_0, \beta_1$ be inner, co-outer, quasi-continuous functions of types $m \times (m-1)$ and $(m-1) \times (m-2)$ respectively, such that the functions $$W_0^T= \begin{pmatrix} \eta_0 & \bar{\beta}_0 \end{pmatrix}, \quad  \tilde{W}_1^T = \begin{pmatrix}
	\beta_0^T \eta_1 & \bar{\beta}_1
	\end{pmatrix}  $$ are unitary-valued. Let
	$$W_1^T = \begin{pmatrix}
	1 & 0\\0& \tilde{W}_1^T
	\end{pmatrix} $$
and let 
\begin{equation}\label{L2}
\mathcal{L}_2 = W_0^* W_1^* \begin{pmatrix}
	0_{2\times 1}\\H^2(\Dd,\Cc^{m-2})^\perp
	\end{pmatrix}.
\end{equation}
Then  
\begin{equation} \label{L2=H2}
\bar{\eta}_0 \telwe \bar{\eta}_1 \telwe \mathcal{L}_2 = \bar{\eta}_0 \telwe \bar{\eta}_1 \telwe H^2(\Dd,\Cc^{m})^\perp, 
\end{equation}
and the operator $(\bar{\eta}_0 \telwe \bar{\eta}_1 \telwe \cdot)\colon  \mathcal{L}_2 \to \bar{\eta}_0 \telwe \bar{\eta}_1 \telwe H^2(\Dd,\Cc^m)^\perp$ is unitary. 
\end{proposition}

\begin{proof}
First let us prove that 
$$\bar{ \eta}_0 \telwe \bar{ \eta}_1\telwe \mathcal{L}_{2} = \bar{ \eta}_0 \telwe \bar{ \eta}_1  \telwe H^2(\Dd,\Cc^m)^\perp.$$ 
Let $B_1 = \beta_0 \beta_1$.
By equations \eqref{L2} and \eqref{W1W0*},
\begin{equation}\label{L2B}
\mathcal{L}_{2} = B_1  H^2(\Dd,\Cc^{m-2})^\perp.
\end{equation}
By Lemmas \ref{beta0*h2} and \ref{beta1h2},
\begin{equation}\label{betaH2-L2}
\begin{array}{lll}
H^2(\Dd,\Cc^{m-2})^\perp &= \beta_1^*(H^2(\Dd,\Cc^{m-1})^\perp\\
& = \beta_1^*\beta_0^*(H^2(\Dd,\Cc^{m})^\perp\\
& = B_1^*H^2(\Dd,\Cc^m)^\perp.
\end{array}
\end{equation}
By  equations \eqref{L2B} and  \eqref{betaH2-L2},
\begin{equation}\label{L2-2}
\mathcal{L}_{2} = B_1  H^2(\Dd,\Cc^{m-2})^\perp = B_1 B_1^*H^2(\Dd,\Cc^m)^\perp.
\end{equation}
By Lemma \ref{eta=B1B1*},
\begin{equation}\label{BB_1*-L} 
B_1 B_1^* = I_m - \sum\limits_{i=0}^1 \bita_i \eta_i^T .
\end{equation} 
Thus 
\begin{equation}\label{L2-eta}
\mathcal{L}_{2} = (I_m - \sum\limits_{i=0}^1 \bita_i \eta_i^T)H^2(\Dd,\Cc^m)^\perp. 
\end{equation}
By Proposition \ref{onxi}, $\{\bar{\eta}_i(z)\}_{i=0}^1$ is an orthonormal set in $\Cc^m$ for almost every $z \in \Tt.$
Therefore, by equations \eqref{L2-2} and \eqref{L2-eta},
\begin{equation}\label{L2H2}
\begin{array}{lll}
\bar{ \eta}_0 \telwe \bar{ \eta}_1 \telwe \mathcal{L}_{2} &=
\bar{ \eta}_0  \telwe \bar{ \eta}_1 \telwe (I_m - \sum\limits_{i=0}^1 \bita_i \eta_i^T)H^2(\Dd,\Cc^m)^\perp\\
&=  \bar{ \eta}_0  \telwe \bar{ \eta}_1 \telwe H^2(\Dd,\Cc^m)^\perp.
\end{array}
\end{equation}


To complete the proof, let us show that the operator $$(\bita_0 \telwe \bita_1 \telwe \cdot)\colon \mathcal{L}_2 \to \bita_0 \telwe \bita_1 \telwe H^2(\Dd,\Cc^{m})^\perp$$ is unitary. Observe that the 
foregoing paragraph asserts the operator is surjective. Hence it suffices to prove that it is an isometry. To this end, let $\upsilon \in \mathcal{L}_2.$ Then
$$\| \bita_0 \telwe \bita_1 \telwe \upsilon \|_{L^2(\Tt,\we^3\Cc^m)}^2 = \langle \bita_0 \telwe \bita_1 \telwe \upsilon, \bita_0 \telwe \bita_1 \telwe \upsilon \rangle_{L^2(\Tt,\we^3\Cc^m)}   ,$$ and, by Proposition \ref{we}, $$ \begin{array}{ll}&\langle \bita_0 \telwe \bita_1 \telwe \upsilon, \bita_0 \telwe \bita_1 \telwe \upsilon \rangle_{L^2(\Tt,\we^3\Cc^m)}\vspace{2ex}\\  &= \displaystyle \frac{1}{2\pi} \int_0^{2\pi} \det \begin{pmatrix} \langle \bita_0(\eiu) , \bita_0(\eiu)\rangle_{\Cc^m} &\langle\bita_0(\eiu),\bita_1 (\eiu)\rangle_{\Cc^m}&\langle \bita_0(\eiu), \upsilon(\eiu) \rangle_{\Cc^m}\\
\langle\bita_1(\eiu), \bita_0(\eiu)\rangle_{\Cc^m} &\langle \bita_1(\eiu) , \bita_1(\eiu)\rangle_{\Cc^m} &\langle\bita_1(\eiu) , \upsilon(\eiu)\rangle_{\Cc^m}\\
\langle\upsilon(\eiu) , \bita_0(\eiu)\rangle_{\Cc^m} & \langle\upsilon(\eiu), \bita_1(\eiu)\rangle_{\Cc^m} & \langle \upsilon (\eiu) , \upsilon(\eiu)\rangle_{\Cc^m} \end{pmatrix}\;d\theta.\end{array}$$
Notice that, by Proposition \ref{onxi}, $\{\bita_0(\eiu),\bita_1(\eiu) \}$ is an orthonormal set almost everywhere on $\Tt.$ Further, since $\mathcal{L}_2 = \beta_0 \beta_1 H^2(\Dd,\Cc^{m-2})^\perp,$  $\upsilon = \beta_0 \beta_1 \varphi$ for some $\varphi\in H^2(\Dd,\Cc^{m-2})^\perp.$ Hence 
$$\begin{array}{lll}\langle \bita_0 (\eiu),v(\eiu)\rangle_{\Cc^m} &= \langle \bita_0(\eiu), \beta_0 (\eiu) \beta_1(\eiu)\varphi(\eiu)\rangle_{\Cc^m}\vspace{2ex}\\
&= \langle \beta_0^* (\eiu)\bita_0(\eiu) ,  \beta_1(\eiu)\varphi(\eiu)\rangle_{\Cc^{m-1}} = 0,  \end{array}$$since $W_0^T$ is unitary-valued almost everywhere on $\Tt.$ Similarly, since, by Proposition \ref{tildew1v1}, $\tilde{W}_1^T$ is unitary-valued almost everywhere on $\Tt,$ we obtain
$$ \langle \bita_1(\eiu), \upsilon(\eiu)\rangle_{\Cc^m} = \langle \beta_1^*(\eiu) \beta_0^*(\eiu) \bita_1(\eiu), \varphi(\eiu)\rangle_{\Cc^{m-2}}=0.$$
Therefore 
$$ \| \bita_0 \telwe \bita_1 \telwe \upsilon \|_{L^2(\Tt,\we^3\Cc^m)}^2 = \displaystyle \frac{1}{2\pi} \int_0^{2\pi} \det \begin{pmatrix} 1 &0&0\\
0 &1 &0\\
0& 0 & \| \upsilon (\eiu) \|_{\Cc^m}^2 \end{pmatrix}\;d\theta = \|\upsilon\|^2_{L^2(\Tt,\Cc^m)},$$
that is, the operator $(\bita_0 \telwe \bita_1 \telwe \cdot)\colon \mathcal{L}_2 \to \bita_0 \telwe \bita_1 \telwe H^2(\Dd,\Cc^{m})^\perp$ is an isometry. Thus the operator is unitary.\end{proof}

\begin{proposition}\label{l2perp}
		Let $\eta_0, \eta_1$ be defined by equations \eqref{xi0eta0} and \eqref{311} respectively and let $\beta_0, \beta_1$ be inner, co-outer, quasi-continuous functions of types $m \times (m-1)$ and \\ $(m-1) \times (m-2)$ respectively, such that the functions $$W_0^T= \begin{pmatrix} \eta_0 & \bar{\beta}_0 \end{pmatrix}, \quad  \tilde{W}_1^T = \begin{pmatrix}
	\beta_0^T \eta_1 & \bar{\beta}_1
	\end{pmatrix}  $$ 
are unitary-valued. Let $$\mathcal{L}_2 = W_0^* \begin{pmatrix}
	1 & 0 \\0& \tilde{W}_1^*
	\end{pmatrix} \begin{pmatrix}
	0_{2\times 1} \\ H^2(\Dd,\Cc^{m-2})^\perp
	\end{pmatrix}.$$ Then
	$$ \mathcal{L}_2^\perp = \{ f \in L^2(\Tt,\Cc^m) : \beta_1^* \beta_0^* f \in H^2(\Dd,\Cc^{m-2})\}.$$
\end{proposition}

\begin{proof}
Clearly $\mathcal{L}_2 = \beta_0 \beta_1 H^2(\Dd,\Cc^{m-2})^\perp.$ 
The general element of $\beta_0 \beta_1 H^2(\Dd,\Cc^{m-2})^\perp$ is 
$\beta_0 \beta_1 \bar{z} \bar{g}$ with $ g \in H^2 (\Dd,\Cc^{m-2})$. 
A function $f \in L^2(\Tt,\Cc^m)$ belongs to  $\mathcal{L}_2^\perp$ if and only if 

$$\langle f, \beta_0 \beta_1 \bar{z} \bar{g} \rangle_{L^2(\Tt,\Cc^m)} =0 \quad \text{for all} \quad  g\in H^2(\Dd,\Cc^{m-2}) $$if and only if

$$\displaystyle \frac{1}{2\pi} \int_0^{2\pi} \langle f(\eiu), \beta_0(\eiu) \beta_1(\eiu) e^{-i\theta} \bar{g}(\eiu) \rangle_{\Cc^m}d\theta =0  \quad \text{for all} \quad  g\in H^2(\Dd,\Cc^{m-2})  $$ if and only if
$$ \displaystyle \frac{1}{2\pi} \int_0^{2\pi} \langle  \beta_1^*(\eiu)\beta_0^*(\eiu)f(\eiu),  e^{-i\theta} \bar{g}(\eiu) \rangle_{\Cc^{m-2}}d\theta =0  \quad \text{for all} \quad  g\in H^2(\Dd,\Cc^{m-2}),$$which in turn is equivalent to the assertion that 
$\beta_1^* \beta_0^* f $ is orthogonal to $H^2(\Dd,\Cc^{m-2})^\perp$ in $L^2(\Tt,\Cc^{m-2}),$ which holds if and only if $\beta_1^* \beta_0^* f  $ belongs to $H^2(\Dd,\Cc^{m-2}).$ Thus 

$$ \mathcal{L}_2^\perp = \{ f \in L^2(\Tt,\Cc^m) : \beta_1^* \beta_0^* f \in H^2(\Dd,\Cc^{m-2})\}$$as required. \end{proof}

\index{$\tilde{V}_2$} \index{$\tilde{W}_2$} \index{$\mathcal{E}_2$}  

\begin{theorem}\label{T2compactt}
Let $m,n$ be positive integers such that $\min(m,n)\geq2.$ Let  $G$ be in $H^\infty(\Dd,\CCmn)+C(\Tt,\CCmn).$ Let  $(\xi_0 \telwe v_1,\bar{\eta}_0\telwe  w_1)$ be a Schmidt pair for the  operator $T_1,$ as given in equation \eqref{T_0},  corresponding to $t_1 = \|T_1\|\neq 0,$  let $h_1 \in H^2(\Dd,\Cc)$ be the scalar outer factor of $\xi_0 \telwe v_1,$ let
$$x_1 = (I_{n} - \xi_0 \xi_0^*)v_1, \quad y_1 = (I_m - \bar{\eta}_0 \eta_0^T)w_1, $$ and let $$\xi_1 = \frac{x_1}{h_1} , \quad \bar{\eta}_1 = \frac{zy_1}{\bar{h}_1}. $$

\noindent Let $$V_0=\big(\begin{matrix}\xi_0 & \balpha_0\end{matrix}\big),\quad W_0^T=\big(\begin{matrix} \eta_0 & \bar{\beta}_0 \end{matrix} \big)$$ be given by equations \eqref{V0W0}, and let  $$\tilde{V}_1=\big(\begin{matrix} \alpha_{0}^T\xi_1 & \balpha_1 \end{matrix}\big), \quad \tilde{W}_1^T=\big( \begin{matrix}
\beta_0^T\eta_1 & \bar{\beta}_1
\end{matrix}\big)$$ be given by equations \eqref{V1} and \eqref{W1} respectively. Let 
$$X_2 = \xi_0 \telwe \xi_1 \telwe H^2(\Dd,\Cc^n), \quad Y_2 = \bar{\eta}_0 \telwe \bar{\eta}_1 \telwe H^2(\Dd,\Cc^m)^\perp ,$$
  let
\begin{equation}\label{k2l2}
\mathcal{K}_2 = V_0 \begin{pmatrix} 1 & 0 \\ 0& \tilde{V}_1\end{pmatrix} \begin{pmatrix}
0_{2\times 1} \\ H^2(\Dd,\Cc^{n-2})\end{pmatrix}, \quad 
\mathcal{L}_2 = W_0^* \begin{pmatrix} 1 & 0 \\ 0& \tilde{W}_1^*\end{pmatrix} \begin{pmatrix}
0_{2\times 1} \\ H^2(\Dd,\Cc^{m-2})^\perp \end{pmatrix}.
\end{equation}
Consider the operator $T_2 \colon  X_2 \to Y_2$ given by 
\begin{equation}\label{TT2} T_2 (\xi_0 \telwe \xi_1 \telwe x ) = P_{Y_2} (\bar{\eta}_0 \telwe \bar{\eta}_1 \telwe (G-Q_2) x ),\end{equation}
where $Q_2 \in H^\infty(\Dd,\CCmn)$ satisfies 
\begin{equation} \label{G_Q2xy}
(G-Q_2) x_i = t_i y_i,\quad (G-Q_2)^*y_i = t_ix_i \quad \text{for}\; i=0,1.
\end{equation}
Let the operator $\Gamma_2 \colon  \mathcal{K}_2 \to \mathcal{L}_2$ be given by $\Gamma_2 = P_{\mathcal{L}_2} M_{G-Q_2}|_{\mathcal{K}_2}.$  Then  \begin{itemize}	
 	\item[{\rm (i)}] The maps 
$$ M_{\bar{\alpha}_0\bar{\alpha}_1} \colon H^2(\Dd,\Cc^{n-2}) \to  
\mathcal{K}_2  \colon x \mapsto \bar{\alpha_0}\bar{\alpha}_1 x,
$$
and
$$M_{\beta_0\beta_1}: 
H^2(\Dd,\Cc^{m-2})^\perp \to \mathcal{L}_2: y \mapsto  
\beta_0 \beta_1 y $$ 
are unitaries.
 	\item[{\rm (ii)}] The maps $(\xi_0\telwe \xi_1\telwe \cdot)\colon\mathcal{K}_2\to X_2,$ $(\bita_0\telwe \bita_1\telwe\cdot)\colon\mathcal{L}_2\to Y_2 $ are unitaries.
 	
 	\item[{\rm (iii)}] The following diagram commutes 
 \begin{equation}\label{commdiagrt2}	
\begin{array}{clllll}
H^2(\Dd,\Cc^{n-2}) &\xrightarrow{M_{\bar{\alpha_0}\bar{\alpha}_1}} & \mathcal{K}_2 &\xrightarrow{\xi_0 \telwe \xi_1 \telwe \cdot}& \xi_0 \telwe \xi_1 \telwe H^2 (\Dd, \Cc^n)=X_2\\
\Big\downarrow\rlap{$\scriptstyle H_{F_2}  $} & ~  &\Big\downarrow\rlap{$\scriptstyle \Gamma_2$}  &~&\hspace{3ex}\Big\downarrow\rlap{$\scriptstyle T_2$} \\
H^2(\Dd,\Cc^{m-2})^\perp &\xrightarrow{M_{\beta_0 \beta_1} }&	\mathcal{L}_2 &\xrightarrow{\bar{\eta}_0 \telwe  \bar{\eta}_1\telwe  \cdot } & \bar{\eta}_0 \telwe  \bar{\eta}_1 \telwe  H^2 (\Dd, \Cc^m)^\perp =Y_2,
\end{array}\end{equation}
where $F_2 \in H^\infty(\Dd,\Cc^{(m-2)\times(n-2)})+ C(\Tt,\Cc^{(m-2)\times(n-2)})$  is the function defined in Proposition \ref{tildew1v1}.
\item [{\rm (iv)}] $T_2$ is a compact operator.
\item[{\rm (v)}] $\|T_2\|=\|\Gamma_2\|=\|H_{F_2}\|=t_2.$

 \end{itemize} 
\end{theorem}

\begin{proof}
{\rm (i)} follows from Lemma \ref{3.1constr}.
{\rm (ii)} follows from Propositions \ref{xi12telweunit} and \ref{eta0eta1b0b1}.	
	
{\rm (iii)} By Proposition \ref{Twell}, $T_2$ is well-defined and is independent of the choice of $Q_2 \in H^\infty(\Dd,\CCmn)$ satisfying equations
\eqref{G_Q2xy}. We can choose $Q_2$ which  minimises 
\[
(s_0^\infty(G-Q),s_1^\infty(G-Q)),
\]
 and therefore  satisfies  equations \eqref{G_Q2xy}.
	By Lemma \ref{3.2constr} and Theorem \ref{1.7}, the left hand side of diagram \eqref{commdiagrt2} commutes. Let us show the right hand side also commutes. A typical element of $\mathcal{K}_2$ is of the form $\bar{\alpha}_0 \bar{\alpha}_1x$ where $x \in  H^2(\Dd,\Cc^{n-2}).$ Then, by equation (\ref{TT2}),
$$\begin{array}{cllll}	T_2 (\xi_0 \telwe \xi_1 \telwe \bar{\alpha}_0 \bar{\alpha}_1x)= P_{Y_2} \left( \bar{ \eta}_{0}    \telwe \bar{ \eta}_{1}  \telwe (G-Q_2) \bar{\alpha}_0 \bar{\alpha}_1 x \right) .
	  \end{array}$$ 
By Proposition \ref{tildew1v1}, every $Q_2\in H^\infty(\Dd,\CCmn)$ which minimises $(s_0^\infty(G-Q),s_1^\infty(G-Q))$ satisfies the following equation (see equation \eqref{g-qv0v1w0w1}),
	  \begin{equation}\label{g-Q2}(G-Q_2) V_0 \begin{pmatrix}
	  1 & 0 \\ 0 & \tilde{V}_1
	  \end{pmatrix} \begin{pmatrix}
	  0 \\ 0 \\ x
	  \end{pmatrix} = W_0^* \begin{pmatrix}
	  1 & 0 \\ 0 & \tilde{W}_1^* 
	  \end{pmatrix}\begin{pmatrix}
	  0 \\ 0 \\ F_2 x
	  \end{pmatrix} ,
\end{equation}
for some $F_2 \in H^\infty(\Dd,\Cc^{(m-2)\times(n-2)})+ C(\Dd,\Tt^{(m-2)\times(n-2)})$. This implies that
\begin{equation}\label{comm2-0}
(G-Q_2)\bar{\alpha_0}\bar{\alpha}_1 x =  \beta_0 \beta_1 F_2 x,
\end{equation}
for $x \in H^2(\Dd,\Cc^{n-2}).$ Hence
\begin{equation}\label{comm2-1}
	T_2 (\xi_0 \telwe \xi_1 \telwe \bar{\alpha}_0 \bar{\alpha}_1 x) =  P_{Y_2} ( \bar{\eta}_0 \telwe \bar{ \eta}_1 \telwe \beta_0 \beta_1 F_2 x ).
\end{equation}
Furthermore, 
 $$\begin{array}{clll}
 ( \bar{\eta}_0 \telwe \bar{ \eta}_1 \telwe \cdot ) \Gamma_2 (\bar{\alpha}_0 \bar{\alpha}_1 x) &= \bar{\eta}_0 \telwe \bar{ \eta}_1 \telwe P_{\mathcal{L}_2} [ (G-Q_2) \bar{\alpha}_0 \bar{\alpha}_1 x].
 \end{array}  $$ 
Hence, by equation \eqref{comm2-0},
\begin{equation}\label{comm2-2}
 ( \bar{\eta}_0 \telwe \bar{ \eta}_1 \telwe \cdot ) \Gamma_2 \left(\bar{\alpha}_0 \bar{\alpha}_1 x \right) =  \bar{\eta}_0 \telwe \bar{ \eta}_1 \telwe P_{\mathcal{L}_2}( \beta_0 \beta_1 F_2 x ).
\end{equation}
To show  commutativity of the right hand square in the diagram (\ref{commdiagrt2}), we need to prove that, for every $x \in H^2(\Dd,\Cc^{n-2})$,
\begin{equation}\label{comm2-3}
	T_2 (\xi_0 \telwe \xi_1 \telwe \bar{\alpha}_0 \bar{\alpha}_1 x) = ( \bar{\eta}_0 \telwe \bar{ \eta}_1 \telwe \cdot ) \Gamma_2 (\bar{\alpha}_0 \bar{\alpha}_1 x).
\end{equation}
By equations \eqref{comm2-1} and \eqref{comm2-2}, it is equivalent to
show that
\begin{equation}\label{comm2-4}
	 P_{Y_2} ( \bar{\eta}_0 \telwe \bar{ \eta}_1 \telwe \beta_0 \beta_1 F_2 x ) = \bar{\eta}_0 \telwe \bar{ \eta}_1 \telwe P_{\mathcal{L}_2}( \beta_0 \beta_1 F_2 x ).
\end{equation}
Therefore, we need to show that
  $$ \bar{ \eta}_{0} \telwe \bar{ \eta}_{1} \telwe P_{\mathcal{L}_2}(\beta_0 \beta_1 F_2 x)\in Y_2 $$ and that  
  $$ \bar{\eta}_0 \telwe \bar{ \eta}_1 \telwe \beta_0 \beta_1 F_2 x - \bar{ \eta}_{0} \telwe \bar{ \eta}_{1} \telwe P_{\mathcal{L}_2}(\beta_0 \beta_1 F_2 x) $$ is orthogonal to $Y_2.$ 

By Proposition \ref{eta0eta1b0b1}, $\bar{ \eta}_{0} \telwe \bar{ \eta}_{1} \telwe P_{\mathcal{L}_2}(\beta_0 \beta_1 F_2 x)$ is indeed an element of $Y_2.$ 
Furthermore, 
$$ \begin{array}{lll}\bar{\eta}_0 \telwe \bar{ \eta}_1 \telwe \beta_0 \beta_1 F_2 x - \bar{ \eta}_{0} \telwe \bar{ \eta}_{1} \telwe P_{\mathcal{L}_2}(\beta_0 \beta_1 F_2 x) &=\bar{\eta}_0 \telwe \bar{ \eta}_1 \telwe [\beta_0 \beta_1 F_2 x - P_{\mathcal{L}_2}(\beta_0 \beta_1 F_2 x)] \vspace{2ex} \\ &= \bar{\eta}_0 \telwe \bar{ \eta}_1 \telwe P_{\mathcal{L}_2^\perp}(\beta_0 \beta_1 F_2 x) \end{array}.$$ 
Let us show that $ \bar{\eta}_0 \telwe \bar{ \eta}_1 \telwe P_{\mathcal{L}_2^\perp}(\beta_0 \beta_1 F_2 x)$ is  orthogonal to $Y_2$.

It is so  if and only if
  \begin{equation}\label{orthtoy2}\left\langle  \bar{\eta}_0 \telwe \bar{ \eta}_1 \telwe P_{\mathcal{L}_2^\perp}(\beta_0 \beta_1 F_2 x) , \bar{ \eta}_0 \telwe \bar{ \eta}_{1} \telwe g \right\rangle_{L^2(\Tt,\we^3\Cc^m)} = 0\quad \text{for every}\; g \in H^2(\Dd,\Cc^m)^\perp.\end{equation}  
Let $\Phi = P_{\mathcal{L}_2^\perp}(\beta_0 \beta_1 F_2 x) \in L^2(\Tt,\Cc^m)$.
By Proposition \ref{l2perp}, 
\begin{equation}\label{phil2} \beta_1^* \beta_0^* \Phi \in H^2(\Dd,\Cc^{m-2}). 
\end{equation} 
Then, by Proposition \ref{we}, assertion (\ref{orthtoy2}) is equivalent to the following assertion 
  $$\displaystyle\frac{1}{2\pi} \int_0^{2\pi} \det \begin{pmatrix}
  \langle \bar{ \eta}_{0} (\eiu), \bar{ \eta}_{0} (\eiu) \rangle_{\Cc^m} & \langle \bar{ \eta}_{0} (\eiu), \bar{ \eta}_{1} (\eiu)\rangle_{\Cc^m} & \langle \bar{ \eta}_{0} (\eiu) , g(\eiu) \rangle_{\Cc^m}\\
  \bar{ \eta}_{1} (\eiu) , \bar{ \eta}_{0} (\eiu) \rangle_{\Cc^m} & \langle \bar{ \eta}_{1} (\eiu), \bar{ \eta}_{1} (\eiu)\rangle_{\Cc^m} & \langle \bar{ \eta}_{1} (\eiu), g(\eiu) \rangle_{\Cc^m}\\
  \langle \Phi(\eiu) , \bar{ \eta}_{0} (\eiu) \rangle_{\Cc^m}  &\langle  \Phi(\eiu), \bar{ \eta}_{1} (\eiu) \rangle_{\Cc^m}& \langle \Phi(\eiu), g (\eiu) \rangle_{\Cc^m}
  \end{pmatrix} \; d\theta = 0  $$ 
for every $g \in H^2(\Dd,\Cc^m)^\perp,$ which in turn, by Proposition \ref{onxi}, is equivalent to the assertion  
$$ \displaystyle\frac{1}{2\pi} \int_0^{2\pi} \det \begin{pmatrix}
1 & 0 & \langle \bar{ \eta}_{0} (\eiu) , g(\eiu) \rangle_{\Cc^m}\\
0 & 1 & \langle \bar{ \eta}_{1} (\eiu), g(\eiu) \rangle_{\Cc^m}\\
\langle \Phi(\eiu) , \bar{ \eta}_{0} (\eiu) \rangle_{\Cc^m}&\langle  \Phi(\eiu), \bar{ \eta}_{1} (\eiu) \rangle_{\Cc^m}& \langle \Phi(\eiu), g (\eiu) \rangle_{\Cc^m}
\end{pmatrix} \; d\theta = 0$$ 
for every $g \in H^2(\Dd,\Cc^m)^\perp.$ The latter statement is equivalent to the assertion 
$$\begin{array}{clll}\displaystyle\frac{1}{2\pi} \int_0^{2\pi} \langle \Phi(\eiu), g (\eiu) \rangle_{\Cc^m} &- \langle  \Phi(\eiu), \bar{ \eta}_{1} (\eiu) \rangle_{\Cc^m}  \langle \bar{ \eta}_{1} (\eiu), g(\eiu) \rangle_{\Cc^m} \\&- \langle \Phi(\eiu) , \bar{ \eta}_{0} (\eiu) \rangle_{\Cc^m}   \langle \bar{ \eta}_{0} (\eiu) , g(\eiu) \rangle_{\Cc^m} \; d\theta =0 \end{array}$$ 
for every $g \in H^2(\Dd,\Cc^m)^\perp,$ which in turn is equivalent to the statement that 
$$\begin{array}{lll}\displaystyle\frac{1}{2\pi} \int_0^{2\pi} g^*(\eiu)\Phi(\eiu) &- g^*(\eiu) \bar{ \eta}_{0}  \eta_0^T(\eiu) \Phi(\eiu)  (\eiu)\\ &-g^*(\eiu) \bar{ \eta}_1(\eiu) \eta_1^T(\eiu) \Phi(\eiu)   \; d\theta =0 \end{array}$$ for every $g \in H^2(\Dd,\Cc^m)^\perp.$ Equivalently $$\displaystyle\frac{1}{2\pi} \int_0^{2\pi} g^*(\eiu) \bigg(I_{m} - \bar{ \eta}_0  (\eiu)\eta_0^T(\eiu) -\bar{\eta}_1 (\eiu)\eta_1^T(\eiu) \bigg)\Phi(\eiu)\;d\theta  =0 $$ for every $g \in H^2(\Dd,\Cc^m)^\perp$ if and only if
$$ \bigg(I_{m} - \bar{ \eta}_0  (\eiu)\eta_0^T(\eiu) -\bar{\eta}_1 (\eiu)\eta_1^T(\eiu) \bigg)\Phi(\eiu)$$ \text{is orthogonal to} $H^2(\Dd,\Cc^m)^\perp,$ which occurs if and only if 
 $$\left(I_{m} - \bar{ \eta}_0  \eta_0^T - \bar{ \eta}_1 \eta_1^T\right)\Phi \in H^2(\Dd,\Cc^m).$$ 
By Lemma \ref{eta=B1B1*}, 
 $$\left(I_{m} - \bar{ \eta}_0  \eta_0^T - \bar{ \eta}_1 \eta_1^T\right)\Phi = \beta_0 \beta_1 \beta_1^* \beta_0^*\Phi. $$
Recall that, by assertions (\ref{phil2}), $\beta_1^* \beta_0^*\Phi \in H^2(\Dd,\Cc^{m-2}),$ and so $$\beta_0 \beta_1 \beta_1^* \beta_0^*\Phi \in H^2(\Dd,\Cc^m).$$ 
Thus the right hand square in the diagram (\ref{commdiagrt2}) commutes, and so the diagram (\ref{commdiagrt2}) commutes.

{\rm(iv)} By Proposition \ref{tildew1v1},
$$F_2 \in H^\infty(\Dd,\Cc^{(m-2)\times(n-2)})+ 
C(\Tt, \Cc^{(m-2)\times(n-2)}). $$
Thus, by Hartman's Theorem, the Hankel operator $H_{F_2}$ is compact. By (iii), 
$$(\bita_0 \telwe \bita_1 \telwe \cdot)\circ (M_{\beta_0 \beta_1} H_{F_2} M_{\alpha_0^T\alpha_1^T} ) \circ (\xi_0\telwe \xi_1 \telwe \cdot)^*=T_2 .$$ By (i) and (ii), the operators $ M_{\balpha_0\balpha_1},\;M_{\beta_0 \beta_1},$ $(\xi_0\telwe \xi_1 \telwe \cdot)$ and $(\bita_0 \telwe \bita_1 \telwe \cdot)$ are unitaries. Hence $T_2$ is a compact operator.

{\rm (v)} Since diagram \eqref{commdiagrt2} commutes and 
the operators $ M_{\balpha_0\balpha_1},\;M_{\beta_0 \beta_1},$ $(\xi_0\telwe \xi_1 \telwe \cdot)$ and\linebreak $(\bita_0 \telwe \bita_1 \telwe \cdot)$ are unitaries,
$$\|T_2\|=\|\Gamma_2\|=\|H_{F_2}\|=t_2. $$
\end{proof}

\begin{lemma}\label{coofschpairs2}
	In the notation of Theorem \ref{T2compactt}, let $v_2 \in H^2(\Dd,\Cc^n)$ and $w_2 \in H^2(\Dd,\Cc^m)^\perp$ be such that $(\xi_0 \telwe\xi_1\telwe v_2, \bar{\eta}_0 \telwe\bita_1\telwe  w_2)$ is a Schmidt pair for the operator $T_2$ corresponding to $\|T_2\|.$ Then

{\em (i)} there exist $x_2 \in \mathcal{K}_2$ and $y_2\in \mathcal{L}_2$ such that $(x_2,y_2) $ is a Schmidt pair for the operator $\Gamma_2$;

	{\em (ii)} for any $x_2 \in \mathcal{K}_2$ and $y_2\in \mathcal{L}_2$ such that $$
		\xi_0 \telwe \xi_1 \telwe x_2= \xi_0 \telwe \xi_1 \telwe  v_2,\quad 	
		\bar{ \eta}_0 \telwe \bita_1 \telwe y_2 = \bar{ \eta}_0 \telwe  \bita_1\telwe w_2,$$ the pair $(x_2,y_2)$ is a Schmidt pair for $\Gamma_2$ corresponding to $\|\Gamma_2\|.$ 
\end{lemma}

\begin{proof}
{\rm (i)} By Theorem \ref{T2compactt}, the diagram (\ref{commdiagrt2}) commutes, $(\xi_0 \telwe \xi_1\telwe \cdot)$ is unitary from $\mathcal{K}_2$ to $X_2,$ and $(\bar{ \eta}_0\telwe \bita_1 \telwe \cdot)$ is unitary from $\mathcal{L}_2$ to $Y_2$ and $\|\Gamma_2\| =\|T_2\|=t_2.$ Moreover, by the commutativity of diagram (\ref{commdiagrt2}), the operator $\Gamma_2\colon \mathcal{K}_2\to \mathcal{L}_2$ is compact, hence there exist $x_2 \in \mathcal{K}_2,$ $y_2 \in \mathcal{L}_2$ such that $(x_2,y_2)$ is a Schmidt pair for $\Gamma_2$ corresponding to $\|\Gamma_2\|=t_2.$ 

{\rm (ii)} Suppose that $x_2\in \mathcal{K}_2,y_2\in \mathcal{L}_2$ satisfy \begin{equation}\label{xitel2}
		\xi_0 \telwe \xi_1 \telwe x_2= \xi_0 \telwe \xi_1 \telwe v_2\quad \text{and}\end{equation}	
		\begin{equation}\label{eta0tel2}
 \bar{ \eta}_0 \telwe \bita_1 \telwe y_2 = 
\bar{ \eta}_0\telwe \bita_1  \telwe w_2. 
\end{equation} 
Let us show that $(x_2,y_2)$ is a Schmidt pair for $\Gamma_2,$ that is, 
		$$\Gamma_2 x_2 = t_2y_2,\quad \Gamma_2^*y_2=t_2x_2. $$ Since diagram (\ref{commdiagrt2}) commutes,
		\begin{equation}\label{commt2gamma2}
T_2 \circ (\xi_0\telwe \xi_1 \telwe \cdot  )=(\bar{ \eta}_0 \telwe \bita_1 \telwe  \cdot)\circ\Gamma_2 , \quad (\xi_0\telwe\xi_1\telwe  \cdot  )^*\circ T_2^* = \Gamma_2^* \circ (\bar{\eta}_0 \telwe\bita_1\telwe \cdot)^*. \end{equation} 
By hypothesis,
		\begin{equation}\label{hypt2} T_2 (\xi_0 \telwe  \xi_1 \telwe v_2)= t_2 (\bar{ \eta}_0 \telwe \bita_1 \telwe w_2), \quad T_2^*(\bar{ \eta}_0 \telwe \bita_1\telwe w_2)= t_2 (\xi_0 \telwe\xi_1\telwe  v_2). \end{equation} 
Thus, by equations \eqref{xitel2}, \eqref{eta0tel2} and \eqref{hypt2},

$$\begin{array}{clllll} \Gamma_2  x_2&= (\bar{\eta}_0 \telwe \bita_1 \telwe \cdot)^*  T_2  (\xi_0\telwe \xi_1\telwe v_2) \vspace{2ex} \\
		&=  (\bar{\eta}_0 \telwe \bita_1 \telwe \cdot)^*  t_2 (\bar{ \eta}_0 \telwe \bita_1  \telwe w_2) \vspace{2ex} \\ &= t_2 (\bar{\eta}_0 \telwe \bita_1  \telwe \cdot)^* (\bar{ \eta}_0 \telwe \bita_1  \telwe y_2).\end{array}$$  
Hence  $$\Gamma_2 x_2= t_2 (\bar{\eta}_0 \telwe \bita_1  \telwe \cdot)^* (\bar{ \eta}_0 \telwe \bita_1  \telwe \cdot)y_2= t_2 y_2.$$ 
		
		\noindent By equation (\ref{xitel2}), $$x_2 = (\xi_0\telwe \xi_1 \telwe \cdot  )^* (\xi_0\telwe \xi_1 \telwe  v_2  ),$$ and, by equation (\ref{eta0tel2}), 
		$$(\bar{\eta}_0 \telwe \bita_1  \telwe \cdot)^*(\bar{ \eta}_0 \telwe \bita_1  \telwe w_2)=y_2. $$
Thus 
$$\begin{array}{clll}\Gamma_2^* y_2	 &= \Gamma_2^*(\bar{ \eta}_0 \telwe\bita_1  \telwe \cdot)^*(\bar{\eta}_0 \telwe \bita_1  \telwe w_2)\vspace{2ex} \\ 
		&= (\xi_0 \telwe \xi_1 \telwe \cdot )^*  T_2^* (\bar{ \eta}_0 \telwe \bita_1  \telwe w_2),\end{array}$$ the last equality following by the second equation of (\ref{commt2gamma2}). By equations \eqref{xitel2} and (\ref{hypt2}), we have 
		$$ T_2^* (\bar{ \eta}_0 \telwe \bita_1  \telwe w_2) = t_2 (\xi_0\telwe \xi_1 \telwe v_2)= t_2(\xi_0\telwe \xi_1\telwe  x_2),$$ and so, $$ \Gamma_2^* y_2 = t_2 x_2.$$Therefore $(x_2,y_2)$ is a Schmidt pair for $\Gamma_2$ corresponding to $t_2$. \end{proof}

\begin{lemma}\label{schfohf2}
	Suppose $(\xi_0 \telwe \xi_1\telwe v_2, \bita_0 \telwe \bita_1  \telwe w_2)$ is a Schmidt pair for $T_2$ corresponding to $t_2.$ Let 
	$$x_2 = (I_{n} - \xi_0 \xi_0^*-\xi_1\xi_1^*)v_2,\quad y_2= (I_{m} - \bita_0\eta_0^T-\bita_1\eta_1^T)w_2,$$ and let 
	$$\hx_2 = \alpha_1^T\alpha_0^T x_2,\quad \hy_2=\beta_1^*\beta_0^*y_2. $$Then the pair $(\hx_2,\hy_2)$ is a Schmidt pair for $H_{F_2}$ corresponding to $\|H_{F_2}\|=t_2.$
	
\end{lemma}

\begin{proof}
Let us first show that  $\hat{x}_2\in H^2(\Dd,\Cc^{n-2})$ and $x_2\in \mathcal{K}_2.$ Recall that $V_0=\big( \begin{matrix} \xi_0 & \balpha_0 \end{matrix}\big)$ and $\tilde{V}_1=\big(\begin{matrix}
\alpha_0^T\xi_1& \balpha_1
\end{matrix} \big)$ are unitary-valued, that is, $\alpha_0^T\xi_0=0,$ $\alpha_1^T\alpha_0^T\xi_1=0,$ 
 \begin{equation}\label{vounit}I_{n}-\xi_0\xi_0^* =\balpha_0 \alpha_0^T, \end{equation}and \begin{equation}\label{tildv1unit}I_{n-1}-\alpha_0^T\xi_1\xi_1^* \balpha_0 = \balpha_1\alpha_1^T .\end{equation} Then
\begin{align}\hx_2&=\alpha_1^T\alpha_0^T x_2\nonumber\vspace{2ex}\\&= \alpha_1^T\alpha_0^T (I_{n}-\xi_0\xi_0^*-\xi_1\xi_1^*)v_2\nonumber\vspace{2ex}\\
&=\alpha_1^T\alpha_0^T v_2 - \alpha_1^T\alpha_0^T\xi_0\xi_0^*v_2 - \alpha_1^T\alpha_0^T\xi_1\xi_1^*v_2 \nonumber\vspace{2ex}\\
&= \alpha_1^T\alpha_0^T v_2,\label{hx2} \end{align}
which, by Lemmas \ref{a0h2} and \ref{a1h2}, implies that $\hx_2 \in H^2(\Dd,\Cc^{n-2}).$ Moreover, by Lemma \ref{xi=A1A1*},  we obtain
$$\begin{array}{ll}\balpha_0\balpha_1\hx_2&=\balpha_0\balpha_1\alpha_1^T\alpha_0^T v_2\vspace{2ex}\\
&=(I_{n}-\xi_0\xi_0^*-\xi_1\xi_1^* )v_2 =x_2. \end{array}$$
Hence \begin{equation}\label{x2}x_2 =\balpha_0\balpha_1\alpha_1^T\alpha_0^Tv_2=\balpha_0\balpha_1\hx_2,\end{equation}and thus $x_2\in \mathcal{K}_2.$ 

Next, we shall show that $\hy_2\in H^2(\Dd,\Cc^{n-2})^\perp$ and $y_2 \in \mathcal{L}_2.$ Notice that since $\tilde{W}_1^T=\big(\begin{matrix} \beta_0^T \eta_1 & \bar{\beta}_1\end{matrix} \big)$ and $W_0^T=\big( \begin{matrix}
\eta_0 & \bar{\beta}_0
\end{matrix}\big)$ are unitary-valued, 
 $\beta_0^*\bita_0 =0, $ $\beta_1^* \beta_0^*\bita_1=0,$ \begin{equation}\label{tw1}
 \big(I_{m-1} - \beta_0^* \bita_1 \eta_1^T \beta_0 \big)= \beta_1 \beta_1^*
\end{equation}and
\begin{equation}\label{wo}
\big(I_m-\bita_0 \eta_0^T\big)=\beta_0 \beta_0^*.
\end{equation}
We have 
\begin{align}\hy_2 &= \beta_1^*\beta_0^*y_2 \nonumber\vspace{2ex}\\
&=\beta_1^*\beta_0^*(I_{m} - \bita_0\eta_0^T-\bita_1\eta_1^T)w_2\nonumber\vspace{2ex}\\
&=\beta_1^*\beta_0^* w_2 - \beta_1^*\beta_0^*\bita_0\eta_0^Tw_2- \beta_1^*\beta_0^*\bita_1\eta_1^Tw_2\nonumber\vspace{2ex}\\
&= \beta_1^*\beta_0^* w_2\label{betaw2}, \end{align}
which, by Propositions \ref{beta0*h2} and \ref{beta1h2}, implies that $\hy_2 \in H^2(\Dd,\Cc^{m-2})^\perp.$ 
By Lemma \ref{eta=B1B1*}, we have 
$$\begin{array}{lllllll}\beta_0 \beta_1 \hy_2&= 
\beta_0 \beta_1 \beta_1^*\beta_0^* w_2 \vspace{2ex}\\ 
&= (I_m -\bita_0 \eta_0^T- \bita_1 \eta_1^T)w_2\vspace{2ex}= y_2. \end{array}$$
 Hence
\begin{equation}\label{y2}y_2 = \beta_0 \beta_1 \beta_1^*\beta_0^* w_2 = \beta_0 \beta_1 \hy_2,\end{equation}and therefore $y_2\in \mathcal{L}_2.$

By Theorem \ref{T2compactt}, the maps $$M_{\bar{\alpha}_0\balpha_1} \colon  H^2(\Dd,\Cc^{n-2}) \to \mathcal{K}_2,\quad M_{\beta_0\beta_1} \colon  H^2(\Dd,\Cc^{m-2})^\perp \to \mathcal{L}_2,$$ are unitaries and 
\begin{equation}\label{hf2g2} H_{F_2} = M^*_{\beta_0\beta_1} \Gamma_2 M_{ \balpha_0\balpha_1}.\end{equation} We need to show that 
	$$H_{F_2}\hx_2 = t_2 \hy_2,\quad H_{F_2}^*\hy_2 = t_2\hx_2. $$By equations \eqref{hx2} and \eqref{x2},
	\begin{equation}\label{x_2alphax_2} x_2= \balpha_0 \balpha_1 \alpha_1^T\alpha_0^Tx_2.\end{equation} Hence equation \eqref{hf2g2} yields
	\begin{align} 
	H_{F_2}\hat{x}_2&= \beta_1^*\beta_0^* \Gamma_2 \balpha_0 \balpha_1 \hx_2=\beta_1^*\beta_0^* \Gamma_2 x_2.\label{h_f2.1}\end{align}

By Proposition \ref{onxi}(ii), 
$$		\xi_0 \telwe \xi_1 \telwe x_2= \xi_0 \telwe \xi_1 \telwe  v_2,\quad 	
		\bar{ \eta}_0 \telwe \bita_1 \telwe y_2 = \bar{ \eta}_0 \telwe  \bita_1\telwe w_2.$$
Thus,  by Lemma \ref{coofschpairs2},
 $(x_2,y_2)$ is a Schmidt pair for the operator $\Gamma_2$ corresponding to $t_2= \|\Gamma_2\|$, that is,
	\begin{equation}\label{schmg2}
	\Gamma_2 x_2 =t_2y_2,\quad \Gamma_2^*y_2=t_2x_2.
	\end{equation} Thus equation \eqref{h_f2.1} yields 
	$$H_{F_2}\hx_2= \beta_1^*\beta_0^* \Gamma_2 x_2 = \beta_1^*\beta_0^* t_2 y_2 = t_2\hy_2 $$
as required.
	Let us show that $H_{F_2}^*\hy_2 = t_2\hx_2.$ By equations \eqref{betaw2} and \eqref{y2},  
	\begin{equation}\label{y2betay2}y_2=\beta_0 \beta_1 \beta_1^* \beta_0^* y_2.\end{equation} 
By equation \eqref{hf2g2},
$$H_{F_2}^* = M_{ \balpha_0\balpha_1}^*\circ \Gamma_2^* \circ M_{\beta_0\beta_1} .$$
Hence
	\begin{align}\label{h_f2.4}
	H_{F_2}^*\hy_2= \alpha_1^T \alpha_0^T \Gamma_2^* \beta_0\beta_1 \hy_2
=\alpha_1^T \alpha_0^T \Gamma_2^*y_2,\end{align} 
and, by equations \eqref{schmg2} and \eqref{h_f2.4}, 	
	$$H_{F_2}^*\hy_2 = \alpha_1^T \alpha_0^T \Gamma_2^*y_2 = \alpha_1^T \alpha_0^T t_2 x_2 =t_2\hx_2 .$$ 
Therefore $(\hx_2,\hy_2)$ is a Schmidt pair for $H_{F_2}$ corresponding to  $\|H_{F_2}\|=t_2.$
\end{proof}

\begin{proposition}\label{x0wev2eta2wew2}
	Let $(\xi_0\telwe \xi_1\telwe  v_2,\bita_0\telwe\bita_0 \telwe  w_2)$ be a Schmidt pair for $T_2$ corresponding to $t_2$ for some $v_2\in H^2(\Dd,\Cc^n), w_2\in H^2(\Dd,\Cc^m)^\perp,$ let $h_2 \in H^2(\Dd,\Cc)$ be the scalar outer factor of $\xi_0\telwe \xi_1\telwe  v_2$,  let 
$$x_2 = (I_{n}- \xi_0 \xi_0^*-\xi_1\xi_1^*)v_2,\quad y_2=(I_{m} - \bar{\eta}_0 \eta_0^T-\bita_1 \eta_1^T)w_2,$$
and let \begin{equation}\label{hatx2}\hx_2=\alpha_1^T \alpha_0^T x_2\quad \text{and}\quad \hy_2=\beta_1^*\beta_0^*y_2.\end{equation}
Then
$$\|\hx_2 (z) \|_{\Cc^{n-2}} =  \|\hy_2(z)\|_{\Cc^{m-2}} = |h_2(z)|,$$
$$\| x_2(z) \|_{\Cc^n} =  \|y_2(z)\|_{\Cc^m} =|h_2(z)| $$
and	
	$$\| \xi_0 (z) \we \xi_1(z) \we v_2(z) \|_{\we^3\Cc^n} = \| \bar{\eta}_0(z) \we \bita_1(z) \we w_2(z)\|_{\we^3\Cc^m} =|h_2(z)| $$
almost everywhere on $\Tt.$	
\end{proposition}

\begin{proof} By Lemma \ref{schfohf2}, $(\hx_2,\hy_2)$ is a Schmidt pair for $H_{F_2}$ corresponding to $\|H_{F_2}\|=t_2$ (see Theorem \ref{T2compactt} (v)). Hence
	$$H_{F_2}\hx_2 = t_2\hy_2 \quad \text{and}\quad H_{F_2}^* \hy_2 = t_2 \hx_2. $$ Then, by Theorem \ref{1.7}, 
	$$\|\hy_2(z)\|_{\Cc^{m-2}}=\|\hx_2(z)\|_{\Cc^{n-2}} $$ almost everywhere on $\Tt$.
	Notice that, by equations \eqref{hatx2},  $$x_2 =\bar{\alpha}_0\balpha_1\hx_2,	$$
and since $\bar{\alpha}_0(z),\balpha_1(z)$ are isometric for almost every $z\in \Tt,$ we obtain
	$$\|x_2(z)\|_{\Cc^n}=\|\hx_2(z)\|_{\Cc^{n-2}}. $$
	
	\noindent Furthermore, by equations \eqref{hatx2}, 
	$$y_2 =\beta_0 \beta_1 \hy_2, $$ and since $\beta_0(z),\beta_1(z)$ are isometries almost everywhere on $\Tt,$ we have  
	$$\|y_2(z)\|_{\Cc^m}=\|\hy_2(z)\|_{\Cc^{m-2}} $$
almost everywhere on $\Tt.$ By equations \eqref{hatseq}, we deduce
	\begin{equation}\label{x2isy2}
\|x_2(z)\|_{\Cc^n}=\|\hx_2(z)\|_{\Cc^{n-2}}=\|\hy_2(z)\|_{\Cc^{m-2}}=\|y_2(z)\|_{\Cc^m} \end{equation} 
almost everywhere on $\Tt.$  By Proposition \ref{onxi}, 
$$\xi_0\we  \xi_1 \telwe v_2 = \xi_0\we  \xi_1 \telwe x_2, \quad 	\bita_0 \we \bita_1 \we w_2= \bita_0 \we \bita_1 \we y_2.$$ 
Hence, by Proposition  \ref{weon},
$$\begin{array}{lll}
&\|\xi_0(z)\we  \xi_1(z) \telwe v_2(z)\|_{\we^3\Cc^n}=\|\xi_0(z)\we \xi_1(z)\we x_2(z)\|_{\we^3\Cc^n}\\
& = \| x_2(z) - \displaystyle\sum\limits_{i=0}^1 \langle x_2(z), \xi_i(z)\rangle \xi_i(z)\|_{\Cc^n}=\|x_2(z)\|_{\Cc^n} 
\end{array}$$
almost everywhere on $\Tt$.  Furthermore 
	$$\begin{array}{llll}
	&\| \bita_0(z) \we \bita_1(z) \we w_2(z)\|_{\we^3\Cc^m}=\| \bita_0(z) \we \bita_1(z)\we y_2(z) \|_{\we^3\Cc^m}\\
&=\| y_2(z) - \displaystyle\sum\limits_{i=0}^1 \langle y_2(z), \bita_i(z)\rangle \bita_i(z)\|_{\Cc^m} =\|y_2(z)\|_{\Cc^m} 	
\end{array}$$
almost everywhere on $\Tt$.   Thus, by equation \eqref{x2isy2},
	$$\| \xi_0 (z) \we \xi_1(z)\we v_2(z) \|_{\we^3\Cc^n} = \| \bar{\eta}_0(z) \we\bita_1(z)\we w_2(z)\|_{\we^3\Cc^m}$$
almost everywhere on $\Tt.$

Recall that $h_2$ is the scalar outer factor  of $\xi_0\telwe \xi_1\telwe  v_2$. Hence
$$\|\hx_2 (z) \|_{\Cc^{n-2}} =  \|\hy_2(z)\|_{\Cc^{m-2}} = |h_2(z)|,$$
$$\| x_2(z) \|_{\Cc^n} =  \|y_2(z)\|_{\Cc^m} =|h_2(z)| $$
and	
	$$\| \xi_0 (z) \we \xi_1(z) \we v_2(z) \|_{\we^3\Cc^n} = \| \bar{\eta}_0(z) \we \bita_1(z) \we w_2(z)\|_{\we^3\Cc^m} =|h_2(z)| $$
almost everywhere on $\Tt.$ \end{proof}

\begin{proposition}\label{epsilon2}
Let $m,n$ be positive integers such that $\min(m,n)\geq2.$ Let \linebreak $G\in H^\infty(\Dd,\CCmn)+C(\Tt,\CCmn).$ In line with the algorithm from Subsection \ref{Alg_statement},
	let $Q_2 \in H^\infty(\Dd,\CCmn)$ satisfy
\begin{equation}\label{Q2cond}	
\begin{array}{llll}
(G-Q_2)x_0 &= t_0 y_0,\quad (G-Q_2)^*y_0&=t_0x_0, \\
(G-Q_2)x_1 &= t_1 y_1,\quad (G-Q_2)^*y_1&=t_1x_1.
\end{array}
\end{equation}
Let the spaces $X_2 , Y_2$ be given by
$$X_2 = \xi_0 \telwe \xi_1 \telwe H^2(\Dd,\Cc^n), \quad Y_2 = \bar{\eta}_0 \telwe \bar{\eta}_1 \telwe H^2(\Dd,\Cc^m)^\perp ,$$
and consider the compact operator $T_2\colon X_2 \to Y_2$ given by $$T_2(\xi_0 \telwe \xi_1 \telwe x) = P_{Y_2} (\bar{\eta}_0 \telwe \bar{\eta}_1 \telwe (G-Q_2)x)$$ 
for all $x \in H^2(\Dd,\Cc^n).$ 
Let  $(\xi_0 \telwe \xi_1\telwe v_2,\bar{\eta}_0\telwe \bar{\eta}_1 \telwe  w_2)$ be a Schmidt pair for the  operator $T_2$ corresponding to $t_2 = \|T_2\|,$  let $h_2 \in H^2(\Dd,\Cc)$ be the scalar outer factor of $
\xi_0 \telwe \xi_1\telwe v_2$, let 
$$x_2 = (I_{n} - \xi_0 \xi_0^*- \xi_1 \xi_1^*)v_2, \quad y_2=(I_{m}-\bar{\eta}_0\eta_0^T-\bar{\eta}_1\eta_1^T)w_2 $$ 
and let 
$$\xi_2 =\frac{{x}_2}{h_2}, \quad \eta_2 =\frac{\bar{z}\bar{y}_2}{h_2}. $$
Then there exist unitary-valued functions $\tilde{V}_2, \tilde{W}_2$ of types $(n-2)\times (n-2)$ and $(m-2) \times (m-2)$ respectively of the form 
	$$\tilde{V}_2 =\begin{pmatrix}
	\alpha_1^T \alpha_0^T \xi_2  &  \bar{\alpha}_2
	\end{pmatrix}, \quad \tilde{W}_2^T = \begin{pmatrix}
	{\beta}_1^T \beta_0^T \eta_2 & \bar{\beta}_2
	
\end{pmatrix} ,
$$
where $\alpha_2, \beta_2$ are inner, co-outer, quasi-continuous and all minors on the first columns of $\tilde{V}_2, \tilde{W}_2^T$ are in $H^\infty.$ Furthermore, the set $\mathcal{E}_2$ of all level $2$ superoptimal error functions for $G$ satisfies
$$\mathcal{E}_2 = W_0^* W_1^* \begin{pmatrix}
I_2 & 0 \\ 
0& \tilde{W}_2^*

\end{pmatrix} \begin{pmatrix}
t_0 u_0 & 0 & 0&0 \\
0 & t_1 u_1 &0 &0 \\
0& 0 & t_2 u_2&0 \\
0&0 & 0 & (F_3 +H^\infty)\cap B(t_2)

\end{pmatrix}\begin{pmatrix}
I_2 & 0 \\ 0 & \tilde{V}_2^*
\end{pmatrix} V_1^* V_0^* , $$
for some  $F_3 \in H^\infty(\Dd,\Cc^{(m-3)\times (n-3)}) +C(\Tt,\Cc^{(m-3)\times (n-3)}),$ where $u_3=\frac{\bar{z}\bar{h}_3}{h_3}$ is a quasi-continuous unimodular function and $B(t_2)$ is the closed ball of radius $t_2$ in $L^\infty(\Tt,\Cc^{(m-3)\times (n-3)}).$

\end{proposition}
\begin{proof} By Theorem \ref{T2compactt}, the following diagram commutes
 \begin{equation}\label{commdiagrt2-2}	\begin{array}{clllll}
H^2(\Dd,\Cc^{n-2}) &\xrightarrow{M_{\bar{\alpha}_0\bar{\alpha}_1}} & \mathcal{K}_2 &\xrightarrow{\xi_0 \telwe \xi_1 \telwe \cdot}& \xi_0 \telwe \xi_1 \telwe H^2 (\Dd, \Cc^n)=X_2\\
\Big\downarrow\rlap{$\scriptstyle H_{F_2}  $} & ~  &\Big\downarrow\rlap{$\scriptstyle \Gamma_2$}  &~&\hspace{3ex}\Big\downarrow\rlap{$\scriptstyle T_2$} \\
H^2(\Dd,\Cc^{m-2})^\perp &\xrightarrow{M_{\beta_0 \beta_1} }&	\mathcal{L}_2 &\xrightarrow{\bar{\eta}_0 \telwe  \bar{\eta}_1\telwe  \cdot } & \bar{\eta}_0 \telwe  \bar{\eta}_1 \telwe  H^2 (\Dd, \Cc^m)^\perp =Y_2.
\end{array}\end{equation}
Recall that the operators $ M_{\balpha_0\balpha_1},\;M_{\beta_0 \beta_1},$ $(\xi_0\telwe \xi_1 \telwe \cdot)$ and $(\bita_0 \telwe \bita_1 \telwe \cdot)$ are unitaries. 

By Proposition \ref{Twell}, $T_2$ is well defined and is independent of the choice of $Q_2 \in H^\infty(\Dd,\Cc^{m\times n}) $ satisfying conditions \eqref{Q2cond}. Hence we may choose $Q_2$ to minimise 
$(s_0^\infty(G-Q), s_1^\infty(G-Q))$, and then, by Proposition \ref{g-q1y1t1x1}, the conditions \eqref{Q2cond} hold.

By Lemma \ref{coofschpairs2}, $(x_2,y_2)$ defined above is a Schmidt pair for $\Gamma_2$ corresponding to $t_2$.  By Lemma \ref{schfohf2}, $(\hx_2,\hy_2)$ is a Schmidt pair for $H_{F_2}$ corresponding to $t_2,$ where 
$$\hx_2 = \alpha_1^T\alpha_0^T x_2,\quad \hy_2=\beta_1^*\beta_0^*y_2. $$	
We intend to apply Lemma \ref{2.2} to $H_{F_2}$ and the Schmidt pair 
$(\hx_2,\hy_2)$ to find unitary-valued functions $\tilde{V}_2,\tilde{W}_2$ such that, for every $\tilde{Q}_2\in H^\infty(\Dd,\Cc^{(m-2)\times(n-2)})$ which is at minimal distance from $F_2,$ a  factorisation of the form
$$F_2-\tilde{Q}_2 = \tilde{W}_2^* \begin{pmatrix}t_2 u_2 &0\\ 0 & F_3
\end{pmatrix}\tilde{V}_2^* $$is obtained, for some $F_3 \in H^\infty(\Dd,\Cc^{(m-2)\times(n-2)})+C(\Tt,\Cc^{(m-2)\times(n-2)}).$ For this purpose we find the inner-outer factorisations of $\hx_2$ and $\bar{z}\bar{\hy}_2.$ 

By Lemma \ref{x0wev2eta2wew2} 
	\begin{equation}\label{h2common}
\|\hx_2(z)\|_{\Cc^{n-2}} = |h_2(z)|\;\text{and}\; \|\hy_2(z)\|_{\Cc^{m-2}} =|h_2(z)|\end{equation} 
almost everywhere on $\Tt.$ Equations \eqref{h2common} imply that $h_2\in H^2(\Dd,\Cc)$ is the scalar outer factor of both $\hat{x}_2$ and $\bar{z}\bar{\hat{y}}_2.$
	
	 By Lemma \ref{2.2}, 
	$\hat{x}_2, \bar{z}\bar{\hat{y}}_2$ admit the inner outer factorisations 
	$$\hat{x}_2 = \hat{\xi}_2 h_2, \quad \bar{z} \bar{y}_2 = \hat{\eta}_2 h_2,$$for some inner $\hat{\xi}_2 \in H^\infty(\Dd,\Cc^{n-2}), \hat{\eta}_2\in  H^\infty(\Dd,\Cc^{m-2}).$ Then 
	
	$$ \hat{x}_2 = \hat{\xi}_2 h_2 =\alpha_1^T\alpha_0^T x_2,\quad  \bar{z}\bar{\hat{y}}_2 =\hat{\eta}_2 h_2=\bar{z}\beta_1^T \beta_0^T \bar{y}_2,$$from which we obtain 
	$$ \hat{\xi}_2 = \alpha_1^T\alpha_0^T \xi_2,\quad  \hat{\eta}_2 = \beta_1^T \beta_0^T \eta_2. $$ We show that $\alpha_1^T\alpha_0^T \xi_2,\;\beta_1^T\beta_0^T \eta_2 $ are inner in order to apply Lemma \ref{2.2} and obtain $\tilde{V}_2$ and $\tilde{W}_2.$  Recall that, by Lemma \ref{schfohf2},
	$$x_2 = (I_{n} - \xi_0 \xi_0^*-\xi_1\xi_1^*)v_2=\bar{\alpha}_0\balpha_1\alpha_1^T \alpha_0^Tv_2,\quad y_2= (I_{m} - \bita_0\eta_0^T-\bita_1\eta_1^T)w_2=\beta_0\beta_1\beta_1^*\beta_0^*w_2 .$$ 
Thus
	$$\alpha_1^T\alpha_0^Tx_2 = \alpha_1^T\alpha_0^T v_2, \quad \beta_1^T\beta_0^T \bar{y}_2 =\beta_1^T\beta_0^T\bar{w}_2  ,$$ 
and since
$$\xi_2 =\frac{x_2}{h_2},\quad \eta_2 = \frac{\bar{z}\bar{y}_2}{h_2}, $$
we deduce that the functions
	$$\alpha_1^T\alpha_0^T\xi_2=\frac{\alpha_1^T\alpha_0^Tv_2}{h_2},\quad  \beta_1^T\beta_0^T\eta_2 = \frac{\beta_1^T\beta_0^T\bar{z}\bar{w}_2}{h_2}$$ are analytic. Furthermore,  $\|\xi_2(z)\|_{\Cc^n}=1$ and $\|\eta_2(z)\|_{\Cc^m}=1$ almost everywhere on $\Tt,$ and, by equations \eqref{h2common},

	$$\|\alpha_1^T(z)\alpha_0^T(z)x_2(z)\|_{\Cc^{n-2}}=\|\alpha_1^T(z)\alpha_0^T(z)v_2(z)\|_{\Cc^{n-2}}=|h_2(z)|$$and$$ \|\beta_1^T(z)\beta_0^T(z)\bar{y}_2(z)\|_{\Cc^{m-2}}=\|\beta_1^T(z)\beta_0^T (z)\bar{w}_2(z)\|_{\Cc^{m-2}}=|h_2(z)| $$almost everywhere on $\Tt.$ Hence 
	$$\|\alpha_1^T(z)\alpha_0^T(z)\xi_2(z)\|_{\Cc^{n-2}}=1,\quad \|\beta_1^T(z)\beta_0^T(z)\eta_2(z)\|_{\Cc^{m-2}}=1 $$almost everywhere on $\Tt.$ Thus $\alpha_1^T\alpha_0^T\xi_2,\; \beta_1^T\beta_0^T\eta_2$ are inner functions.
	
	By Lemma \ref{2.2}, there exist inner, co-outer, quasi-continuous functions $\alpha_2, \beta_2$  of types $(n-2)\times (n-3), (m-2)\times (m-3)$ respectively such that the functions 
	$$\tilde{V}_2 =\begin{pmatrix}
	\alpha_1^T \alpha_0^T \xi_2  &  \bar{\alpha}_2
	\end{pmatrix}, \quad \tilde{W}_2^T = \begin{pmatrix}
	{\beta}_1^T \beta_0^T \eta_2 & \bar{\beta}_2
	
\end{pmatrix} $$
are unitary-valued with all minors on the first columns in $H^\infty.$

\noindent Furthermore, by Lemma \ref{2.2}, every $\hat{Q}_2\in H^\infty(\Dd,\Cc^{(m-2)\times(n-2)})$ which is at minimal distance from $F_2$ satisfies 
$$F_2-\hat{Q}_2 = \tilde{W}_2^* \begin{pmatrix}
t_2 u_2 & 0 \\
0 & F_3
\end{pmatrix}\tilde{V}_2^*, $$
for some  $F_3 \in   H^\infty(\Dd,\Cc^{(m-3)\times (n-3)})+C(\Tt,\Cc^{(m-3)\times (n-3)})$ and for the quasi-continuous unimodular function given by $u_2 = \frac{\bar{z} \bar{h}_2}{h_2}.$ 

By Lemma \ref{f+hinfty}, the set $$\tilde{\mathcal{E}}_{2} =\{F_{2} - \hat{Q} : \hat{Q} \in H^\infty(\Dd,\Cc^{(m-2)\times (n-2)}), \| F_{2} - \hat{Q}\|_{L^\infty}=t_{2}  \}$$ satisfies  $$\tilde{\mathcal{E}}_{2} = \tilde{W}_{2}^* 
\begin{pmatrix}
t_{2}u_{2} & 0 \\
0 & (F_3+H^\infty)\cap B(t_2)
\end{pmatrix}V_{2}^*,
$$ where $B(t_2)$ is the closed ball of radius $t_2$ in $L^\infty(\Tt,\Cc^{(m-3)\times (n-3)}).$ Thus, by Proposition \ref{tildew1v1}, $\mathcal{E}_2$ admits the factorisation claimed. \end{proof}

\begin{proposition}\label{g-q2}
	Every $Q_3 \in H^\infty(\Dd,\CCmn)$ which minimises $$(s_0^\infty(G-Q),s_1^\infty(G-Q),s_2^\infty(G-Q)) $$satisfies 
	$$(G-Q_3)x_i = t_i y_i,\quad (G-Q_3)^*y_i=t_ix_i\quad \text{for}\quad i=0,1,2. $$
\end{proposition}
\begin{proof}
	By Proposition \ref{g-q1y1t1x1}, every $Q_3 \in H^\infty(\Dd,\CCmn)$ that minimises $$(s_0^\infty(G-Q),s_1^\infty(G-Q))$$ satisfies
$$(G-Q_3)x_i = t_i y_i,\quad (G-Q_3)^*y_i=t_ix_i\quad \text{for}\quad i=0,1. $$Hence it suffices to show that $Q_3$ satisfies 
$$(G-Q_3)x_2 = t_2y_2, \quad (G-Q_3)^*y_2=t_2x_2. $$

By Theorem \ref{T2compactt}, the following diagram commutes
 $$	\begin{array}{clllll}
H^2(\Dd,\Cc^{n-2}) &\xrightarrow{M_{\bar{\alpha}_0\bar{\alpha}_1}} & \mathcal{K}_2 &\xrightarrow{\xi_0 \telwe \xi_1 \telwe \cdot}& \xi_0 \telwe \xi_1 \telwe H^2 (\Dd, \Cc^n)=X_2\\
\Big\downarrow\rlap{$\scriptstyle H_{F_2}  $} & ~  &\Big\downarrow\rlap{$\scriptstyle \Gamma_2$}  &~&\hspace{3ex}\Big\downarrow\rlap{$\scriptstyle T_2$} \\
H^2(\Dd,\Cc^{m-2})^\perp &\xrightarrow{M_{\beta_0 \beta_1} }&	\mathcal{L}_2 &\xrightarrow{\bar{\eta}_0 \telwe  \bar{\eta}_1\telwe  \cdot } & \bar{\eta}_0 \telwe  \bar{\eta}_1 \telwe  H^2 (\Dd, \Cc^m)^\perp =Y_2,
\end{array}$$ where the operator $\Gamma_2 \colon  \mathcal{K}_2 \to \mathcal{L}_2$ is given by $\Gamma_2 = P_{\mathcal{L}_2} M_{G-Q_2}|_{\mathcal{K}_2}$ and \linebreak$F_2\in H^\infty(\Dd,\Cc^{(m-2)\times(n-2)})+C(\Tt,\Cc^{(m-2)\times(n-2)})$ is constructed as follows. 

\noindent By Lemma \ref{2.2} and Proposition \ref{tildew1v1}, there exist unitary-valued functions 
$$\tilde{V}_1=\big( \begin{matrix}
\alpha_0^T \xi_1 & \balpha_1
\end{matrix}\big),\quad \tilde{W}_1^T =\big( \begin{matrix}
\beta_0^T \eta_1 & \bar{\beta}_1
\end{matrix}\big), $$where $\alpha_1, \beta_1$ are inner, co-outer, quasi-continuous functions of types $(n-1)\times(n-2)$ and $(m-1)\times (m-2)$ respectively, and all minors on the first columns of $\tilde{V}_1, \tilde{W}_1^T$ are in $H^\infty.$ Furthermore, the set of all level $1$ superoptimal functions $\mathcal{E}_1=\{G-Q:Q\in \Omega_1\}$ satisfies 
 
\begin{equation}\label{g-qv0v1w0w111}
\mathcal{E}_1 = W_0^* \begin{pmatrix} 1 & 0 \\ 0& \tilde{W}_1^*\end{pmatrix} \begin{pmatrix} t_0 u_0 & 0&0 \\ 0& t_1 u_1 &0\\ 0&0 & (F_2+H^\infty(\Dd,\Cc^{(m-2)\times(n-2)}))\cap B(t_1) \end{pmatrix}\begin{pmatrix} 1 & 0 \\ 0 & \tilde{V}_1^*  \end{pmatrix} V_0^*, \end{equation}
for some $F_2 \in   H^\infty(\Dd,\Cc^{(m-2)\times (n-2)})+C(\Tt,\Cc^{(m-2)\times (n-2)})$, for a quasi-continuous unimodular function $u_1 = \frac{\bar{z} \bar{h}_1}{h_1}$, and for the closed ball $B(t_1)$ of radius $t_1$ in $L^\infty(\Tt, \Cc^{(m-2)\times (n-2)})$. 

Consider some $Q_3\in \Omega_1,$ so that, according to equation \eqref{g-qv0v1w0w111},

$$ \begin{pmatrix} 1 & 0 \\ 0& \tilde{W}_1\end{pmatrix}W_0 (G-Q_3)V_0  \begin{pmatrix}
1 & 0 \\ 0 & \tilde{V}_1 \end{pmatrix} =\begin{pmatrix} t_0 u_0 & 0&0 \\ 0& t_1 u_1 &0\\ 0&0 & F_2-\tilde{Q}_2 \end{pmatrix},$$
for some $\tilde{Q}_2\in H^\infty(\Dd,\Cc^{(m-2)\times(n-2)}),$ that is,

\begin{equation}\label{g-q3}
\begin{pmatrix}
1 & 0 \\
0 & \begin{pmatrix}\eta_1^T\beta_0 \\\beta_1^*\end{pmatrix}
\end{pmatrix} \begin{pmatrix}\eta_0^T \\\beta_0^*\end{pmatrix}(G-Q_3) \big(\begin{matrix}
\xi_0&\balpha_0
\end{matrix}\big)\begin{pmatrix}
1 & 0 &0 \\
0 & 
\alpha_0^T \xi_1 & \balpha_1
\end{pmatrix}=\begin{pmatrix} t_0 u_0 & 0&0 \\ 0& t_1 u_1 &0\\ 0&0 & F_2-\tilde{Q}_2 \end{pmatrix}.\end{equation}Observe

$$ \begin{pmatrix}\eta_0^T \\\beta_0^*\end{pmatrix}(G-Q_3) \big(\begin{matrix}
\xi_0&\balpha_0
\end{matrix}\big)= \begin{pmatrix}
t_0 u_0 & 0 \\ 
0& \beta_0^* (G-Q_3)\balpha_0
\end{pmatrix},$$hence 

$$\begin{pmatrix}
	1 & 0 \\
	0 & \begin{pmatrix}\eta_1^T\beta_0 \\\beta_1^*\end{pmatrix}
\end{pmatrix} \begin{pmatrix}
t_0 u_0 & 0 \\ 
0& \beta_0^* (G-Q_3)\balpha_0
\end{pmatrix} \begin{pmatrix}
1 & 0 &0 \\
0 & 
\alpha_0^T \xi_1 & \balpha_1
\end{pmatrix}$$ is equal to

$$\begin{pmatrix}
t_0 u_0 & 0 & 0 \\
0& \eta_1^T \beta_0 \beta_0^* (G-Q_3)\balpha_0 \alpha_0^T\xi_1 & \eta_1^T \beta_0 \beta_0^* (G-Q_3)\balpha_0 \balpha_1 \\ 0& \beta_1^* \beta_0^* (G-Q_3)\balpha_0 \alpha_0^T \xi_1 & \beta_1^* \beta_0^* (G-Q_3)\balpha_0 \balpha_1
\end{pmatrix}, $$ and so
 equation \eqref{g-q3} yields 

$$\begin{pmatrix}
t_0 u_0 & 0 & 0 \\
0& \eta_1^T \beta_0 \beta_0^* (G-Q_3)\balpha_0 \alpha_0^T\xi_1 & \eta_1^T \beta_0 \beta_0^* (G-Q_3)\balpha_0 \balpha_1 \\ 0& \beta_1^* \beta_0^* (G-Q_3)\balpha_0 \alpha_0^T \xi_1 & \beta_1^* \beta_0^* (G-Q_3)\balpha_0 \balpha_1
\end{pmatrix} =\begin{pmatrix} t_0 u_0 & 0&0 \\ 0& t_1 u_1 &0\\ 0&0 & F_2-\tilde{Q}_2 \end{pmatrix},$$which is equivalent to the following equations

$$\eta_1^T \beta_0 \beta_0^* (G-Q_3)\balpha_0 \alpha_0^T\xi_1= t_1 u_1, $$

$$\eta_1^T \beta_0 \beta_0^* (G-Q_3)\balpha_0 \balpha_1 =0, $$

$$\beta_1^* \beta_0^* (G-Q_3)\balpha_0 \alpha_0^T \xi_1= 0, $$
and 

\begin{equation}\label{g-q_3}
\beta_1^* \beta_0^* (G-Q_3)\balpha_0 \balpha_1 = F_2 -\tilde{Q}_2.
\end{equation}

By Theorem \ref{1.7} applied to $H_{F_2},$ if $(\hx_2,\hy_2)$ is a Schmidt pair for $H_{F_2}$ corresponding to $t_2=\|H_{F_2}\|,$ then, for any $\tilde{Q}_2$ which is at minimal distance from $F_2,$ we have 
\begin{equation}\label{f2-q2}(F_2-\tilde{Q}_2)\hx_2 = t_2\hy_2,\quad (F_2-\tilde{Q}_2)^*\hy_2 =t_2\hx_2. \end{equation}

By equations \eqref{g-q_3} and \eqref{f2-q2},

\begin{equation}\label{g-q31}\beta_1^* \beta_0^* (G-Q_3)\balpha_0 \balpha_1\hx_2= t_2\hy_2\end{equation}and \begin{equation}\label{g-q32}\alpha_1^T \alpha_0^T (G-Q_3)^*\beta_0 \beta_1 \hy_2=t_2\hx_2. \end{equation}Recall that, by equations \eqref{x2} and \eqref{y2}, \begin{equation}\label{a01x1}
\balpha_0\balpha_1\hx_2 = x_2\quad \text{and}\quad \hy_2 = \beta_1^*\beta_0^*y_2.\end{equation} 
Hence, by equation \eqref{g-q31}, we obtain 
$$\beta_1^* \beta_0^* (G-Q_3)x_2=t_2\beta_1^* \beta_0^* y_2, $$or equivalently, 
$$ \beta_1^* \beta_0^* \bigg( (G-Q_3)x_2-t_2y_2\bigg)=0.$$ Since, by Theorem \ref{T2compactt}, $M_{\beta_0\beta_1}$ is unitary, the latter equation yields 
$$ (G-Q_3)x_2=t_2y_2.$$

Moreover, in view of equations \eqref{g-q_3}, \eqref{f2-q2} and \eqref{a01x1}, equation \eqref{g-q32} implies
$$ \alpha_1^T \alpha_0^T (G-Q_3)^*y_2=t_2\alpha_1^T\alpha_0^Tx_2,$$which in turn is equivalent to the equation 
$$ \alpha_1^T \alpha_0^T\bigg(  (G-Q_3)^*y_2-t_2x_2\bigg)=0.$$By Theorem \ref{T2compactt}, $M_{\balpha_0 \balpha_1}$ is unitary, hence the latter equation yields
$$(G-Q_3)^*y_2=t_2x_2 $$and therefore the assertion has been proved. \end{proof}

\section{Compactness of the operator $T_{j+1}$}\label{T_j-compact}

At this point, the reader has seen the proof of  the compactness of the operators $T_1$ and $T_2$. 
Suppose we have applied steps $0,\dots,j$ of the superoptimal analytic approximation algorithm from Subsection \ref{Alg_statement} to $G$, we have constructed 
$$\begin{array}{lll}	&t_0 \geq t_1 \geq \cdots \geq t_j > 0\\
&x_0, x_1, \cdots, x_j \in L^2 (\Tt, \Cc^n)\\
&y_0 , y_1 , \cdots , y_j \in L^2(\Tt, \Cc^m) \\
&h_0, h_1, \cdots, h_j \in H^2(\Dd,\Cc) \; \text{outer}\\
& \xi_0,\xi_1, \cdots, \xi_j \in L^\infty(\Tt,\Cc^n)\; \text{pointwise orthonormal on}\;\Tt\\
& \eta_0, \eta_1, \cdots , \eta_j \in L^\infty(\Tt,\Cc^m)\; \text{pointwise orthonormal on}\;\Tt\\
&X_0 = H^2(\Dd,\Cc^n),X_1, \cdots, X_j \\
&Y_0 = H^2(\Dd,\Cc^m)^\perp, Y_1, \cdots, Y_j\\
&T_0, T_1, \cdots, T_j \; \text{compact operators},\end{array} $$
and all the claimed properties hold.
We shall apply a similar method to show that the operator $T_{j+1}$ as given in equation \eqref{T_j} is compact.

\begin{proposition}\label{tildvjwj} \index{$\mathcal{E}_j$}
	Let $m,n$ be positive integers such that $\min(m,n)\geq2.$ Let \linebreak $G\in H^\infty(\Dd,\CCmn)+C(\Tt,\CCmn).$ In line with the algorithm from Subsection \ref{Alg_statement}, let $Q_j \in H^\infty(\Dd,\CCmn)$ satisfy
\begin{equation}\label{g-qj}
(G-Q_{j})x_i = t_i y_i,\quad (G-Q_{j})^*y_i=t_ix_i\quad \text{for}\quad i=0,1,\dots,j-1.
\end{equation}
Let the spaces $X_j , Y_j$ be given by
$$X_j = \xi_0 \telwe \xi_1 \telwe \dots \telwe \xi_{j-1} \telwe H^2(\Dd,\Cc^n), \quad Y_j = \bar{\eta}_0 \telwe \bar{\eta}_1 \telwe \dots \telwe \bar{\eta}_{j-1} \telwe H^2(\Dd,\Cc^m)^\perp ,$$
and consider the compact operator 
$T_j\colon X_j \to Y_j$ given by 
$$T_j(\xi_0 \telwe \xi_1 \telwe \dots \telwe \xi_{j-1} \telwe x) = P_{Y_j} (\bar{\eta}_0 \telwe \bar{\eta}_1 \telwe \dots \telwe\bar{\eta}_{j-1} \telwe (G-Q_j)x)$$ for all $x \in H^2(\Dd,\Cc^n).$ 
Let  $(\xi_0 \telwe \xi_1 \telwe \dots \telwe \xi_{j-1} \telwe v_j,
\bar{\eta}_0\telwe \bar{\eta}_1 \telwe \dots \bar{\eta}_{j-1} \telwe  w_j)$ be a Schmidt pair for the  operator $T_j$ corresponding to $t_j = \|T_j\|,$  let $h_j \in H^2(\Dd,\Cc)$ be the scalar outer factor of 
$\xi_0 \telwe \xi_1 \telwe \dots \telwe \xi_{j-1} \telwe v_j$, 
let 
$$x_{j} =(I_{n}-\xi_0 \xi_0^* \dots- \xi_{j-1}\xi_{j-1}^*)v_{j},\quad y_{j} = (I_{m} - \bar{\eta}_0 \eta_0^T - \dots - \bar{\eta}_{j-1}\eta_{j-1}^T)w_{j} $$
and let 
\begin{equation} \label{etajxij}\xi_{j} = \frac{x_{j}}{h_{j}},\quad \eta_{j} = \frac{\bar{z} \bar{y}_{j}}{h_{j}}.
\end{equation}
Let, for $i=0,1,\dots, j-1$,
\begin{equation} \label{tildevjwjj}
\tilde{V}_{i} = \begin{pmatrix}
\alpha_{i-1}^T \cdots \alpha_0^T \xi_{i} & \bar{ \alpha}_{i} 
\end{pmatrix}, \quad \tilde{W}_{i}^T = \begin{pmatrix}
\beta_{i-1}^T  \cdots \beta_0^T \eta_{i} & \bar{\beta}_{i}
\end{pmatrix}
\end{equation} be unitary-valued functions, 
as described in Lemma \ref{2.2} (see also Proposition \ref{epsilon2} for 
$\tilde{V}_{2}$ and $\tilde{W}_{2}^T$), 
  $u_i = \frac{\bar{z} \bar{h}_i}{h_i}$ are quasi-continuous unimodular functions, and $$V_{i} = \begin{pmatrix}
I_{i} & 0 \\
0 & \tilde{V}_{i} 
\end{pmatrix}, \quad W_{i} = \begin{pmatrix}
I_{i} & 0 \\
0 & \tilde{W}_{i}
\end{pmatrix} .$$

	 There exist unitary-valued functions $\tilde{V}_{j},\tilde{W}_{j}$ of the form \begin{equation}\label{tildevjwj} \tilde{V}_{j} = \begin{pmatrix}
	  \alpha_{j-1}^T \cdots \alpha_0^T \xi_{j} & \bar{ \alpha}_{j} 
	 \end{pmatrix}, \quad \tilde{W}_{j}^T = \begin{pmatrix}
	 \beta_{j-1}^T  \cdots \beta_0^T \eta_{j} & \bar{\beta}_{j}
	 \end{pmatrix} ,\end{equation}
	 where $\alpha_0, \dots, \alpha_{j-1}$ and $\beta_0,\dots,\beta_{j-1}$ are of types $n\times (n-1),\dots,(n-j-1)\times (n-j-2)$ and $m\times (m-1),\dots,(m-j-1)\times (m-j-2)$ respectively, and are inner, co-outer and quasi-continuous. 
	
	Furthermore, the set of all level $j$ superoptimal error functions $\mathcal{E}_j$ satisfies 
	\index{$\mathcal{E}_j$}
	\begin{equation}\label{epsilonnj}
\mathcal{E}_j = W_0^* W_1^* \cdots W_{j}^* \begin{pmatrix}
	t_0 u_0  &    0   &  \cdots &  0          &0_{1\times (n-j-1)}\\
	0        &t_1u_1  & \dots   &  0        &0_{1\times (n-j-1)} \\
	\vdots   &  \vdots &  \ddots & \vdots    &\vdots \\
	0        &   0    &   \cdots  &t_{j}u_{j} & 0\\ 
	0_{(m-j-1)\times 1} & 0_{(m-j-1)\times 1} & \dots &    \dots    &     (F_{j+1}+H^\infty)\cap B(t_{j})
	\end{pmatrix}V_{j}^* \cdots V_0^* , \end{equation} 
for some  $F_{j+1} \in H^\infty(\Dd,\Cc^{(m-j-1)\times (n-j-1)})+C(\Tt,\Cc^{(m-j-1)\times (n-j-1)}),$ for the quasi-continuous unimodular functions  $u_i = \frac{\bar{z} \bar{h}_i}{h_i}$, for all $i=0,\dots, j ,$ for the closed ball $B(t_{j})$ of radius $t_j$ in $L^\infty(\Tt,\Cc^{(m-j-1)\times (n-j-1)})$, and for the unitary valued functions $$V_{j} = \begin{pmatrix}
	I_{j} & 0 \\
	0 & \tilde{V}_{j} 
	\end{pmatrix}, \quad W_{j} = \begin{pmatrix}
	I_{j} & 0 \\
	0 & \tilde{W}_{j}
	\end{pmatrix}. $$
\end{proposition} \index{$\tilde{V}_j$} \index{$\tilde{W}_j$}


\begin{proof}
	 Suppose we have applied steps $0,\dots,j$ of the algorithm from Subsection \ref{Alg_statement} and that the following diagram commutes
	\begin{equation}\label{tj-1comm}
	\begin{array}{clllll}		
	H^2(\Dd,\Cc^{n-j}) &\xrightarrow{M_{\bar{\alpha}_0\cdots \bar{ \alpha}_{j-1}}} & \mathcal{K}_{j} &\xrightarrow{\xi_{(j-1)} \telwe \cdot}& \xi_{(j-1)} \telwe H^2 (\Dd, \Cc^n)=X_{j}\\
	\Big\downarrow\rlap{$\scriptstyle H_{F_{j}}  $} & ~  &\Big\downarrow\rlap{$\scriptstyle \Gamma_{j}$}  &~&\hspace{3ex}\Big\downarrow\rlap{$\scriptstyle T_{j}$} \\
	H^2(\Dd,\Cc^{m-j})^\perp &\xrightarrow{M_{\beta_0 \cdots \beta_{j-1}}}&\mathcal{L}_{j} &\xrightarrow{\bar{\eta}_{(j-1)} \telwe \cdot } & \bar{\eta}_{(j-1)} \telwe H^2 (\Dd, \Cc^m)^\perp =Y_{j},
	\end{array}\end{equation}
where the maps 
$$M_{\bar{\alpha}_0\cdots \bar{ \alpha}_{j-1}}\colon
H^2(\Dd,\Cc^{n-j}) \to \mathcal{K}_{j}\colon x \mapsto \bar{\alpha}_0\cdots \bar{ \alpha}_{j-1}x,$$ 
$$ M_{\beta_0 \cdots \beta_{j-1}}\colon
H^2(\Dd,\Cc^{m-j})^\perp \to \mathcal{L}_{j}\colon y \mapsto 
\beta_0 \cdots \beta_{j-1}y,$$ 
$$(\xi_{(j-1)} \telwe \cdot)\colon \mathcal{K}_{j}\to X_{j} \; \text{ and} \; (\bar{\eta}_{(j-1)} \telwe \cdot)\colon\mathcal{L}_{j}\to Y_{j} $$ 
are unitaries. 

Let $(\xi_{(j-1)} \telwe v_{j}, \bar{\eta}_{(j-1)}\telwe w_{j})$ be a Schmidt pair for the compact operator $T_{j}.$ Then 
\linebreak$x_{j} \in \mathcal{K}_{j},$ $y_{j} \in \mathcal{L}_{j}$ are such that $(x_{j},  y_{j})$ is a Schmidt pair for $\Gamma_{j}$ corresponding to $t_{j}=\|\Gamma_{j}\|,$ and $(\hat{x}_{j},\hat{y}_{j})$ is a Schmidt pair for $H_{F_{j}}$ corresponding to $t_{j}=\| H_{F_{j}} \|,$ where  
\begin{equation}\label{hatxj-1}\hat{x}_{j}= \alpha_{j-1}^T\cdots \alpha_0^T x_{j}, \quad  \hat{y}_{j}=\beta_{j-1}^* \cdots \beta_0^* y_{j}. \end{equation} 
We intend to apply Lemma \ref{2.2} to $H_{F_j}$ and the Schmidt pair $(\hx_j,\hy_j)$ to find unitary-valued functions $\tilde{V}_j,\tilde{W}_j$ such that, for every $\tilde{Q}_j\in H^\infty(\Dd,\Cc^{(m-j)\times(n-j)})$ which is at minimal distance from $F_j,$ a  factorisation of the form
$$F_j-\tilde{Q}_j = \tilde{W}_j^* \begin{pmatrix}t_j u_j &0\\ 0 & F_{j+1}
\end{pmatrix}\tilde{V}_j^* $$ is obtained, for some 
$F_{j+1} \in H^\infty(\Dd,\Cc^{(m-j-1)\times (n-j-1)})+C(\Tt,\Cc^{(m-j-1)\times (n-j-1)}).$ 
 For this purpose we find the inner-outer factorisations of $\hx_j$ and $\bar{z}\bar{\hy}_j.$

By the inductive hypothesis (see Lemma \ref{x0wev2eta2wew2} for $j=2$), we have
	\begin{equation}\label{hjcommon}
\begin{array}{lll}
|h_j(z)|&=	\| \xi_0 (z) \we \dots \we \xi_{j-1}(z) \we v_j(z) \|_{\we^{j+1}\Cc^n} = \| \bar{\eta}_0(z) \we \dots \we\bita_{j-1}(z) \we w_j(z)\|_{\we^{j+1}\Cc^m}, \\
\|\hx_j (z) \|_{\Cc^{n-j}} &=  \|\hy_j(z)\|_{\Cc^{m-j}} = |h_j(z)|,\;\text{and}\\
\| x_j(z) \|_{\Cc^n} &=  \|y_j(z)\|_{\Cc^m} =|h_j(z)|,  \;
\end{array}
\end{equation} 
almost everywhere on $\Tt.$
Equations \eqref{hjcommon} imply that $h_j\in H^2(\Dd,\Cc)$ is the scalar outer factor of both $\hat{x}_j$ and $\bar{z}\bar{\hat{y}}_j.$


By Lemma \ref{2.2}, $\hat{x}_{j}, \bar{z}\bar{\hat{y}}_{j}$ admit the inner-outer factorisations 
\begin{equation}\label{facthatj-1}\hat{x}_{j} = \hat{\xi}_{j} h_{j},\quad \bar{z}\bar{\hat{y}}_{j} = \hat{\eta}_{j} h_{j}, \end{equation}where $\hat{\xi}_{j} \in H^\infty(\Dd,\Cc^{n-j})$ and $\hat{\eta}_{j}
\in 	H^\infty(\Dd,\Cc^{m-j})$ are vector-valued inner functions. 

By equations \eqref{hatxj-1} and \eqref{facthatj-1}, we deduce that 
$$\hat{\xi}_{j} = \alpha_{j-1}^T\cdots \alpha_0^T \xi_{j},\quad \hat{\eta}_{j} = \beta_{j-1}^T \cdots \beta_0^T \eta_{j}.$$ 

We shall show that $\alpha_{j-1}^T\cdots\alpha_0^T \xi_j,\;\beta_{j-1}^T\cdots\beta_0^T \eta_j $ are inner in order to apply Lemma \ref{2.2} and obtain $\tilde{V}_j$ and $\tilde{W}_j$ as required.  We have
$$\begin{array}{lll}\hx_j&=\alpha_{j-1}^T\cdots \alpha_0^Tx_j \vspace{2ex}\\&= \alpha_{j-1}^T\cdots \alpha_0^T(I_{n} - \xi_0 \xi_0^*-\dots-\xi_{j-1}\xi_{j-1}^*)v_j\vspace{2ex}\\
&= \alpha_{j-1}^T\cdots \alpha_0^Tv_j- \alpha_{j-1}^T\cdots \alpha_0^T\xi_0\xi_0^*v_j-\cdots-\alpha_{j-1}^T\cdots \alpha_0^T\xi_{j-1}\xi_{j-1}^*v_j. \end{array}$$

Recall that, by the inductive hypothesis, for $i=0,\dots,j-1,$ each $$\tilde{V}_{i}= \begin{pmatrix}
\alpha_{i-1}^T \cdots \alpha_0^T \xi_{i} & \bar{ \alpha}_{i} 
\end{pmatrix}$$ is unitary-valued, and so $\alpha_{i}^T \alpha_{i-1}^T \cdots \alpha_0^T \xi_{i} =0$. Hence, if $0 \le i \le j-1$, we have 
$$\alpha_{j-1}^T\cdots\alpha_{i+1}^T \alpha_{i}^T\cdots\alpha_0^T\xi_i=0.$$ 
Thus
$$\hx_j=\alpha_{j-1}^T\cdots \alpha_0^Tx_j =\alpha_{j-1}^T\cdots \alpha_0^Tv_j,$$ that is, $\hx_j \in H^2(\Dd,\Cc^{n-j})$ and
$$ \alpha_{j-1}^T\cdots \alpha_0^T\xi_j= \frac{1}{h_j} \alpha_{j-1}^T\cdots \alpha_0^Tx_j=\frac{1}{h_j} \alpha_{j-1}^T\cdots \alpha_0^Tv_j$$is analytic. Moreover, by equations \eqref{hjcommon},
$$\|\alpha_{j-1}^T(z)\cdots\alpha_0^T(z)x_j(z)\|_{\Cc^{n-j}}=  \|\alpha_{j-1}^T(z)\cdots\alpha_0^T(z)v_j(z)\|_{\Cc^{n-j}}=|h_j(z)|$$
almost everywhere on $\Tt,$ and hence
$$\|\alpha_{j-1}^T(z)\cdots \alpha_0^T(z)\xi_j(z)\|_{\Cc^{n-j}} =1$$ almost everywhere on $\Tt.$ Therefore $\alpha_{j-1}^T\cdots \alpha_0^T\xi_j$ is inner.

Furthermore 
$$ \begin{array}{llll}\hy_j &=\beta_{j-1}^*\cdots \beta_0^*y_j\vspace{2ex}\\
&= \beta_{j-1}^*\cdots \beta_0^* (I_{m} - \bar{\eta}_0 \eta_0^T - \dots - \bar{\eta}_{j-1}\eta_{j-1}^T)w_{j}\vspace{2ex}\\
&= \beta_{j-1}^*\cdots \beta_0^* w_j - \beta_{j-1}^*\cdots \beta_0^* \bar{\eta}_0 \eta_0^T w_j-\cdots-\beta_{j-1}^*\cdots \beta_0^* \bar{\eta}_{j-1} \eta_{j-1}^Tw_j.
\end{array}$$ 
Notice that, by the inductive hypothesis, for $i=0,\dots,j-1,$ each 
$$\tilde{W}_{i}^T = \begin{pmatrix}
\beta_{i-1}^T  \cdots \beta_0^T \eta_{i} & \bar{\beta}_{i}
\end{pmatrix}
$$
 is unitary-valued, and so $\beta_{i}^*\cdots \beta_0^* \bar{\eta}_{i} =0.$
Hence, if $0 \le i \le j-1$, we have 
$$\beta_{j-1}^*\cdots \beta_{i+1}^*\beta_{i}^*\cdots \beta_0^* \bar{\eta}_{i} =0.$$ 
Thus 
$$\hy_j =\beta_{j-1}^*\cdots \beta_0^*y_j = \beta_{j-1}^*\cdots \beta_0^* w_j,   $$ that is,
$\hy_j \in H^2(\Dd,\Cc^{m-j})^\perp$ and 
$$\beta_{j-1}^T\cdots \beta_0^T  \eta_j= \frac{1}{h_j} \beta_{j-1}^T\cdots \beta_0^T \bar{z}\bar{y}_j =\frac{1}{h_j} \beta_{j-1}^T\cdots \beta_0^T \bar{z}\bar{w}_j$$ is analytic. Further, by equations \eqref{hjcommon},
$$\| \beta_{j-1}^T(z)\cdots \beta_0^T (z) \bar{z}\bar{y}_j(z) \|_{\Cc^{m-j}}=  \| \beta_{j-1}^T(z)\cdots \beta_0^T (z) \bar{z}\bar{w}_j(z) \|_{\Cc^{m-j}}=|h_j(z)| $$ almost everywhere on $\Tt,$ and therefore 
$$ \|\beta_{j-1}^T(z)\cdots \beta_0^T (z)\eta_j(z)\|_{\Cc^m}=1 $$almost everywhere on $\Tt,$ that is, $\beta_{j-1}^T\cdots \beta_0^T  \eta_j$ is inner.

We apply  Lemma \ref{2.2} to the Hankel operator $H_{F_j}$ and the Schmidt pair
$(\hx_j, \hy_j)$ to deduce that there exist inner, co-outer, quasi-continuous functions $\alpha_{j}, \beta_{j}$ of types $(n-j)\times (n-j-1),$ $ (m-j)\times (m-j-1)$ respectively such that $$\tilde{V}_{j} = \begin{pmatrix}
\alpha_{j-1}^T\cdots \alpha_0^T \xi_{j} & \bar{ \alpha}_{j}
\end{pmatrix}, \quad \tilde{W}_{j}^T = \begin{pmatrix}
\beta_{j-1}^T \cdots \beta_0^T \eta_{j} & \bar{ \beta}_{j}
\end{pmatrix} $$ are unitary-valued and all minors on the first columns of $\tilde{V}_{j},\tilde{W}_{j}$ are in $H^\infty.$ Moreover, every function $\hat{Q}_{j}\in H^\infty(\Dd,\Cc^{(m-j)\times (n-j)}),$ which is at minimal distance from $F_{j},$ satisfies
$$F_{j} - \hat{Q}_{j} = \tilde{W}_{j}^* \begin{pmatrix}
t_{j}u_{j} & 0 \\
0 & F_{j+1}
\end{pmatrix}\tilde{V}_{j}^* ,$$ 
for some $F_{j+1} \in H^\infty(\Dd,\Cc^{(m-j-1)\times (n-j-1)})+C(\Tt,\Cc^{(m-j-1)\times (n-j-1)})$ and for the quasi-continuous unimodular function $u_j = \frac{\bar{z} \bar{h}_j}{h_j}.$ 

By Lemma \ref{f+hinfty}, the set $$\tilde{\mathcal{E}}_{j} =\{F_{j} - \hat{Q} : \hat{Q} \in H^\infty(\Dd,\Cc^{(m-j)\times (n-j)}), \| F_{j} - \hat{Q}\|_{L^\infty}=t_{j}  \}$$ satisfies  $$\tilde{\mathcal{E}}_{j} = \tilde{W}_{j}^* 
\begin{pmatrix}
t_{j}u_{j} & 0 \\
0 & (F_{j+1}+H^\infty)\cap B(t_{j})
\end{pmatrix}\tilde{V}_{j}^*,
$$ 
where $B(t_{j})$ is the closed ball of radius $t_{j}$ in $L^\infty(\Tt,\Cc^{(m-j-1)\times (n-j-1)}).$ 

By the inductive hypothesis,  the set of all level $j$ superoptimal error functions $\mathcal{E}_j$ satisfies 
\index{$\mathcal{E}_j$}
\begin{equation}\label{epsilonnj-1}
\mathcal{E}_{j-1} = W_0^* W_1^* \cdots W_{j-1}^* \begin{pmatrix}
	t_0 u_0  &    0   &  \cdots &  0          &0_{1\times (n-j)}\\
	0        &t_1u_1  & \dots   &  0        &0_{1\times (n-j)} \\
	\vdots   &  \vdots &  \ddots & \vdots    &\vdots \\
	0        &   0    &   \cdots  &t_{j-1}u_{j-1} & 0\\ 
	0_{(m-j)\times 1} & 0_{(m-j)\times 1} & \dots &    \dots    &     (F_{j}+H^\infty)\cap B(t_{j-1})
\end{pmatrix}V_{j-1}^* \cdots V_0^* , \end{equation} 
for some  $F_{j} \in H^\infty(\Dd,\Cc^{(m-j)\times (n-j)})+C(\Tt,\Cc^{(m-j)\times (n-j)}),$  $u_i = \frac{\bar{z} \bar{h}_i}{h_i}$ are quasi-continuous unimodular functions for all $i=0,\dots, j-1$, and for the closed ball  $B(t_{j-1})$ of radius $t_{j-1}$ in $L^\infty(\Tt,\Cc^{(m-j)\times (n-j)}).$  

Thus, by equation \eqref{epsilonnj-1}, $\mathcal{E}_j$ admits the factorisation \eqref{epsilonnj} as claimed. \end{proof}


\begin{remark}\label{V0VjandW0Wj} Let, for $i=0,1,\dots, j$,
\begin{equation} \label{tildevjwjj-1}
\tilde{V}_{i} = \begin{pmatrix}
\alpha_{i-1}^T \cdots \alpha_0^T \xi_{i} & \bar{ \alpha}_{i} 
\end{pmatrix}, \quad \tilde{W}_{i}^T = \begin{pmatrix}
\beta_{i-1}^T  \cdots \beta_0^T \eta_{i} & \bar{\beta}_{i}
\end{pmatrix}
\end{equation} be unitary-valued functions, 
as described in Lemma \ref{2.2}. Let 
 $$V_{j} = \begin{pmatrix}
	I_{j} & 0 \\
	0 & \tilde{V}_{j} 
	\end{pmatrix}, \quad W_{j} = \begin{pmatrix}
	I_{j} & 0 \\
	0 & \tilde{W}_{j}
	\end{pmatrix}. $$
Let $A_j= \alpha_0\alpha_1 \dots \alpha_j$, $A_{-1} = I_n$,  $B_j= \beta_0 \beta_1\dots\beta_j$ and $B_{-1} = I_m$.
		
	Note $$W_1W_0 = \begin{pmatrix}
	1 & 0\\
	0& \eta_1^T\beta_0\\
	0& \beta_1^*
	\end{pmatrix}\begin{pmatrix}
	\eta_0^T\\
	\beta_0^*
	\end{pmatrix}= \begin{pmatrix}
	\eta_0^T \\ \eta_1^T B_0B_0^*\\B_1^*
	\end{pmatrix}$$
and 
	$$W_2W_1W_0 = \begin{pmatrix}
	I_2 & 0\\
	0& \eta_2^TB_1\\
	0& \beta_2^*
	\end{pmatrix}\begin{pmatrix}
	\eta_0^T \\ \eta_1^T B_0 B_0^*\\B_1^*
	\end{pmatrix} = \begin{pmatrix}
	\eta_0^T  \\ \eta_1^TB_0B_0^* \\ \eta_2^T B_1B_1^*\\ B_2^*
	\end{pmatrix}.$$ 
Similarly one obtains
\begin{equation}\label{WjW0}
W_jW_{j-1}\cdots W_0 = \begin{pmatrix}
	\eta_0^T \\ \eta_1^TB_0B_0^* \\\vdots\\ \eta_j^TB_{j-1}B_{j-1}^* \\ B_j^*
	\end{pmatrix}.
\end{equation}
Therefore 
\begin{equation}\label{WjW0*}
W_0^* W_1^*\cdots W_j^*= \begin{pmatrix}
\bar{\eta}_{0} & B_0B_0^*\bar{\eta}_{1} & \dots & B_{j-1} B_{j-1}^* \bar{\eta}_{j} & B_{j} \end{pmatrix}. 
\end{equation}
Thus 
\begin{equation} \label{W0*Wj*WjW0}
I_m= W_0^* W_1^* \cdots W_j^* W_j \cdots W_{1} W_0 =
 \sum\limits_{i=0}^{j} B_{i-1}B_{i-1}^* \bita_i \eta_i^T B_{i-1}B_{i-1}^* +
B_j B_j^*.
\end{equation}

	Furthermore 
	$$V_0 V_1 = \begin{pmatrix}
	\xi_0 & \balpha_0 
	\end{pmatrix} \begin{pmatrix}
	1 & 0 &0\\
	0& \alpha_0^T\xi_1 & \balpha_1
	
\end{pmatrix}= \begin{pmatrix}
	\xi_0  & A_0^T \xi_1 & \bar{A}_1
\end{pmatrix}$$ and 
$$V_0 V_1 V_2 = \begin{pmatrix}
	\xi_0  & \bar{A}_0 A_0^T \xi_1 & \bar{A}_1
\end{pmatrix} \begin{pmatrix}
	I_2 & 0 & 0 \\
	0 &A_1^T\xi_2 &\bar{\alpha}_2 \\
\end{pmatrix}=\begin{pmatrix}
	\xi_0 & \bar{A}_0 A_0^T \xi_1 &  \bar{A}_1 A_1^T \xi_2 & \bar{A}_2
\end{pmatrix}.  $$ 
One can easily show by induction that 
\begin{equation}\label{V0Vj}
 V_0\cdots V_j = \begin{pmatrix}
	\xi_0 & \bar{A}_0 A_0^T \xi_1 & \bar{A}_1A_1^T\xi_2& \dots & \bar{A}_{j-1}A_{j-1}^T \xi_{j}& \bar{A}_{j}
\end{pmatrix}.	
\end{equation}
Therefore,
\begin{equation} \label{V0VjVj*V0*}
I_n= V_0\cdots V_j V_j^*\cdots V_0^* = \xi_0 \xi_0^* + \bar{A}_0 A_0^T \xi_1 \xi_1^* \bar{A}_0 A_0^T + \dots \bar{A}_{j-1} A_{j-1}^T \xi_j \xi_j^* \bar{A}_{j-1} A_{j-1}^T +\bar{A}_j A_j^T.
\end{equation}

\end{remark}

\begin{lemma}\label{xi=AkAkT}
Let
\begin{equation}\label{tildevjwjj-2}
\tilde{V}_{i} = \begin{pmatrix}
\alpha_{i-1}^T \cdots \alpha_0^T \xi_{i} & \bar{\alpha}_{i} 
\end{pmatrix}
\end{equation} 
be unitary-valued functions, for $i=0,1,\dots, j$, 
as described in Lemma \ref{2.2}
 For $i = 0,1,\dots,j$, let $A_i= \alpha_0\alpha_1 \dots \alpha_i$ and $A_{-1} = I_n$.
 Then, for $i=0,1,\dots, j$,
\begin{equation}\label{AAk} 
\bar{A}_i A_i^T = I_n - \sum\limits_{k=0}^i \xi_k \xi_k^*
\end{equation} 
almost everywhere on $\Tt$.	
\end{lemma}
\begin{proof}
By equation \eqref{V0VjVj*V0*}, for $k= 0, \dots,j$,	
\begin{equation}\label{AkAkT} 
\bar{A}_k A_k^T 
 = I_n - \sum\limits_{i=0}^{k} \bar{A}_{i-1} A_{i-1}^T \xi_i \xi_i^* \bar{A}_{i-1} A_{i-1}^T.
\end{equation} 
Thus to prove condition \eqref{AAk} it suffices to show that, for  $k=0,\dots,j,$
$$\bar{A}_{k-1} A_{k-1}^T \xi_k \xi_k^* \bar{A}_{k-1} A_{k-1}^T =\xi_k \xi_k^*.$$
	
	For $k=0,$ 
$$\bar{A}_{-1} A_{-1}^T \xi_0 \xi_0^* \bar{A}_{-1} A_{-1}^T =\xi_0 \xi_0^*,$$ 
and so, equation \eqref{AkAkT} yields 
$$\bar{A}_0 A_0^T = I_n -\xi_0 \xi_0^*.$$
	
	For $k=1,$  
$$\bar{A}_{0} A_{0}^T \xi_1 \xi_1^* \bar{A}_{0} A_{0}^T=
(I_n -\xi_0 \xi_0^*)\xi_1 \xi_1^*(I_n -\xi_0 \xi_0^*)
$$ 
By Proposition \ref{onxi}, $\xi_1$ and $\xi_0$ are pointwise orthogonal almost everywhere on $\Tt,$ hence 
$$\bar{A}_{0} A_{0}^T \xi_1 \xi_1^* \bar{A}_{0} A_{0}^T=\xi_1 \xi_1^*,$$ and in view of equation \eqref{AkAkT}, we get
	$$\bar{A}_1 A_1^T 
 = I_n - \xi_0 \xi_0^* -\xi_1 \xi_1^*.$$ 
	
	Suppose 
\begin{equation}\label{AAkT}
\bar{A}_{\ell-1} A_{\ell-1}^T \xi_\ell \xi_\ell^* \bar{A}_{\ell-1} A_{\ell-1}^T =\xi_\ell \xi_\ell^*
\end{equation}
 holds for every $\ell \le k$, where $0\le k \le j$, almost everywhere on $\Tt$. By equations \eqref{AkAkT} and \eqref{AAkT}, this implies 
$$ \bar{A}_{k} A_{k}^T = I_n - \sum\limits_{i=0}^{k} \xi_i \xi_i^* .$$
Let us show that 
$$\bar{A}_{k} A_{k}^T \xi_{k+1} \xi_{k+1}^* \bar{A}_{k} A_{k}^T =\xi_{k+1} \xi_{k+1}^*.$$ 
Note that
	$$\bar{A}_{k} A_{k}^T \xi_{k+1} \xi_{k+1}^* \bar{A}_{k} A_{k}^T
= (I_n - \sum\limits_{i=0}^{k} \xi_i \xi_i^*)\xi_{k+1} \xi_{k+1}^* 
(I_n - \sum\limits_{i=0}^{k} \xi_i \xi_i^*).
$$
By Proposition \ref{onxi}, the set $\{ \xi_i(z)\}_{i=0}^{k+1}$ is pointwise orthogonal almost everywhere on $\Tt,$ and therefore
	$$\bar{A}_{k} A_{k}^T \xi_{k+1} \xi_{k+1}^* \bar{A}_{k} A_{k}^T
= \xi_{k+1} \xi_{k+1}^*.$$
Thus, by equation \eqref{AkAkT},
	$$\bar{A}_{k+1} A_{k+1}^T = I_n - \sum\limits_{i=0}^{k+1} \xi_i \xi_i^*$$
almost everywhere on $\Tt$, and the assertion has been proved. \end{proof}

\begin{lemma}\label{eta=BkBk*}
Let
\begin{equation} \tilde{W}_{i}^T = \begin{pmatrix}
\beta_{i-1}^T  \cdots \beta_0^T \eta_{i} & \bar{\beta}_{i}
\end{pmatrix}
\end{equation} be unitary-valued functions, for $i=0,1, \dots, j$,
as described in Lemma \ref{2.2}.
For $i=0,1,\dots, j$, let $B_i= \beta_0 \beta_1\dots\beta_i$ and $B_{-1} = I_m$. Then, for $k=0,1,\dots, j$,
\begin{equation}\label{BB_k*} 
B_k B_k^* = I_m - \sum\limits_{i=0}^k \bita_i \eta_i^T 
\end{equation} 	
almost everywhere on $\Tt$.
\end{lemma}

\begin{proof} 
By equation \eqref{W0*Wj*WjW0}, for $k= 0, \dots,j$,	
\begin{equation}\label{BB_k**} 
B_k B_k^* = I_m - \sum\limits_{i=0}^{k} B_{i-1}B_{i-1}^*\bita_i\eta_i^T B_{i-1}B_{i-1}^*.\end{equation} 
Thus to prove condition \eqref{BB_k*} it suffices to show that, for  $k=0,\dots,j,$  
	$$B_{k-1}B_{k-1}^*\bita_k\eta_k^T B_{k-1}B_{k-1}^* =\bita_k\eta_k^T.$$
	
	For $k=0,$ 
$$B_{-1} B_{-1}^* \bita_0 \eta_0^T B_{-1} B_{-1}^*=I_m \bita_0 \eta_0^T I_m =\bita_0 \eta_0^T ,$$ 
and so, equation \eqref{BB_k**} yields 
$$B_0 B_0^* = I_m - \bita_0 \eta_0^T. $$
	
	For $k=1,$  
$$B_{0}B_{0}^*\bita_1\eta_1^T B_{0}B_{0}^* =  (I_m - \bita_0 \eta_0^T) \bita_1 \eta_1^T (I_m - \bita_0 \eta_0^T).$$ 
By Proposition \ref{onxi}, $\eta_1$ and $\eta_0$ are pointwise orthogonal almost everywhere on $\Tt,$ hence 
	$$B_{0}B_{0}^*\bita_1\eta_1^T B_{0}B_{0}^*= \bita_1 \eta_1^T,$$ and in view of equation \eqref{BB_k**}, we get
	$$B_1B_1^* = I_m - \bita_0 \eta_0^T - \bita_1\eta_1^T. $$ 
	
	Suppose 
\begin{equation}\label{BBk*}
B_{\ell-1}B_{\ell-1}^*\bita_\ell\eta_\ell^T B_{\ell-1}B_{\ell-1}^* =\bita_\ell \eta_\ell^T
\end{equation}
holds for every $\ell \le k$, where $0\le k \le j$, almost everywhere on $\Tt$. By equations \eqref{BB_k**} and 
\eqref{BBk*}, this implies 
$$ B_k B_k^* = I_m - \sum\limits_{i=0}^k \bita_i \eta_i^T .$$
Let us show that 
$$B_{k}B_{k}^*\bita_{k+1}\eta_{k+1}^T B_{k}B_{k}^* =\bita_{k+1}\eta_{k+1}^T .$$ 
Note that
	$$B_{k}B_{k}^*\bita_{k+1}\eta_{k+1}^T B_{k}B_{k}^* = (I_m - \sum\limits_{i=0}^k \bita_i \eta_i^T) \bita_{k+1}\eta_{k+1}^T (I_m - \sum\limits_{i=0}^k \bita_i \eta_i^T).$$ By Proposition \ref{onxi}, the set $\{ \bita_i(z)\}_{i=0}^{k+1}$ is pointwise orthogonal almost everywhere on $\Tt,$ and therefore
	$$B_{k}B_{k}^*\bita_{k+1}\eta_{k+1}^T B_{k}B_{k}^* = \bita_{k+1}\eta_{k+1}^T .$$ Thus, by equation \eqref{BB_k**},
	$$B_{k+1} B_{k+1}^* = I_m - \sum\limits_{i=0}^{k+1} \bita_i \eta_i^T$$ 
almost everywhere on $\Tt$, and the assertion has been proved. \end{proof}


The following statement asserts that any function $Q_{j+1} \in \Omega_{j}$ necessarily satisfies equations \eqref{g-qi}.  
\begin{proposition}\label{g-qjj}
	Every $Q_{j+1} \in H^\infty(\Dd,\CCmn)$ which minimises $$(s_0^\infty(G-Q),s_1^\infty(G-Q),\dots,s_j^\infty(G-Q)) $$ satisfies 
	$$(G-Q_{j+1})x_i = t_i y_i, (G-Q_{j+1})^*y_i=t_ix_i,\quad \text{for}\; i=0,1,\dots,j. $$
\end{proposition}
\begin{proof}
	By the recursive step of the algorithm from Subsection \ref{Alg_statement}, every $Q_{j+1} \in H^\infty(\Dd,\CCmn)$ that minimises $$(s_0^\infty(G-Q),\dots, s_{j-1}^\infty(G-Q))$$ satisfies
	$$(G-Q_{j+1})x_i = t_i y_i,\quad (G-Q_{j+1})^*y_i=t_ix_i\quad \text{for}\quad i=0,1,\dots,j-1. $$
	Hence it suffices to show that $Q_{j+1}$ satisfies 
	$$ (G-Q_{j+1})x_j = t_j y_j, (G-Q_{j+1})^*y_j=t_jx_j .$$ Notice that, by the inductive step, the following diagram commutes
	\begin{equation}\label{tj-1commj}
	\begin{array}{clllll}
	H^2(\Dd,\Cc^{n-j}) &\xrightarrow{M_{\bar{\alpha}_0\cdots \bar{ \alpha}_{j-1}}} & \mathcal{K}_{j} &\xrightarrow{\xi_{(j-1)} \telwe \cdot}& \xi_{(j-1)} \telwe H^2 (\Dd, \Cc^n)=X_{j}\\
	\Big\downarrow\rlap{$\scriptstyle H_{F_{j}}  $} & ~  &\Big\downarrow\rlap{$\scriptstyle \Gamma_{j}$}  &~&\hspace{3ex}\Big\downarrow\rlap{$\scriptstyle T_{j}$} \\
	H^2(\Dd,\Cc^{m-j})^\perp &\xrightarrow{M_{\beta_0 \cdots \beta_{j-1}}}&\mathcal{L}_{j} &\xrightarrow{\bar{\eta}_{(j-1)} \telwe \cdot } & \bar{\eta}_{(j-1)} \telwe H^2 (\Dd, \Cc^m)^\perp =Y_{j},
	\end{array}\end{equation}
where the maps $M_{\bar{\alpha}_0\cdots \bar{ \alpha}_{j-1}},$ $ M_{\beta_0 \cdots \beta_{j-1}},$ $(\xi_{(j-1)} \telwe \cdot)\colon \mathcal{K}_{j}\to X_{j} $ and $(\bar{\eta}_{(j-1)} \telwe \cdot)\colon\mathcal{L}_{j}\to Y_{j} $ are unitaries, and $F_{j} \in H^\infty(\Dd,\Cc^{(m-j)\times(n-j)})+C(\Tt,\Cc^{(m-j)\times(n-j)}).$

	By equation \eqref{epsilonnj-1}, the set of all level $j-1$ superoptimal error functions $$\mathcal{E}_{j-1} =\{G-Q:Q\in \Omega_{j-1} \}$$ satisfies \begin{equation}\label{epsilonnj-11}
\mathcal{E}_{j-1} = W_0^* W_1^* \cdots W_{j-1}^* \begin{pmatrix}
	t_0 u_0  &    0   &  \cdots &  0          &0_{1\times (n-j)}\\
	0        &t_1u_1  & \dots   &  0        &0_{1\times (n-j)} \\
	\vdots   &  \vdots &  \ddots & \vdots    &\vdots \\
	0        &   0    &   \cdots  &t_{j-1}u_{j-1} & 0\\ 
	0_{(m-j)\times 1} & 0_{(m-j)\times 1} & \dots &\dots & (F_{j}+H^\infty)\cap B(t_{j-1})
\end{pmatrix}V_{j-1}^* \cdots V_0^* , \end{equation} 
for some $F_{j} \in H^\infty(\Dd,\Cc^{(m-j)\times (n-j)})+C(\Tt,\Cc^{(m-j)\times (n-j)}),$ for quasi-continuous unimodular functions  $u_i = \frac{\bar{z} \bar{h}_i}{h_i}$, $i=0,\dots, j-1,$ and the closed ball $B(t_{j-1})$ of radius $t_{j-1}$ in $L^\infty(\Tt,\Cc^{(m-j)\times (n-j)}).$
	Consider some  $Q_{j+1}\in \Omega_{j-1},$ so that, according to equation \eqref{epsilonnj-11},  
	\begin{equation}\label{wow1vovj} \footnotesize{\begin{pmatrix} I_{j-1} & 0 \\ 
	0& \tilde{W}_{j-1} \end{pmatrix}\cdots W_0 (G-Q_{j+1}) V_0\cdots \begin{pmatrix} I_{n-j-1} & 0 \\ 
	0& \tilde{V}_{j-1} \end{pmatrix} = \scriptsize{\begin{pmatrix}
	t_0 u_0 & 0 &\dots &0 \\
	0&t_1u_1 &\dots &0 \\
	\vdots & \hspace{10ex}\ddots &  &\vdots \\
	0  &\cdots &t_{j-1}u_{j-1} & 0\\ 
	0 & \dots & \dots & F_{j}-\tilde{Q}_{j}
	\end{pmatrix}}}, \end{equation}
where $\tilde{Q}_{j} \in H^\infty(\Dd,\Cc^{(m-j)\times(n-j)})$ is at minimal distance from $F_j.$ Let $B_j = \beta_0 \cdots \beta_j $ and let $A_j = \alpha_0 \cdots \alpha_j$.
\index{$B_j$} \index{$A_j$} 
By equations \eqref{tildevjwjj}, we have	
	
$$\begin{array}{lll}	&\begin{pmatrix} I_{j-1} & 0 \\ 
		0& \tilde{W}_{j-1} \end{pmatrix}\cdots W_0 (G-Q_{j+1}) V_0\cdots \begin{pmatrix} I_{n-j-1} & 0 \\ 
		0& \tilde{V}_{j-1} \end{pmatrix}\vspace{2ex}\\ &=\footnotesize{\begin{pmatrix}
		t_0 u_0 & 0 &\dots &0 \\
		\vdots &\ddots  & \dots &\vdots \\	
		0 & \dots  &\eta_{j-1}^T B_{j-2} B_{j-2}^*  (G-Q_{j+1})\bar{A}_{j-2}A_{j-2}^T\xi_{j-1} & \eta_{j-1}^T B_{j-2} B_{j-2}^*  (G-Q_{j+1}) \bar{A}_{j-1} \\ 
		0& \dots & B_{j-1}^*(G-Q_{j+1})\bar{A}_{j-2}A_{j-2}^T\xi_{j-1} & B_{j-1}^* (G-Q_{j+1})\bar{A}_{j-1}
		\end{pmatrix}},\end{array}$$
which, combined with equation \eqref{wow1vovj}, yields 
		
		\begin{equation}\label{f-qj} B_{j-1}^* (G-Q_{j+1})\bar{A}_{j-1}= F_j -\tilde{Q}_j. \end{equation}
	
	 Since $\tilde{Q}_{j}$ is at minimal distance from $F_{j},$ $$\|F_{j}-\tilde{Q}_{j}\|_\infty=\|H_{F_{j}}\|=t_{j}.$$ Note that, if $(\hat{x}_{j},\hat{y}_{j})$ is a Schmidt pair for $H_{F_{j}}$ corresponding to $t_{j},$ then, by Theorem \ref{1.7}, 
	$$(F_{j}-\tilde{Q}_{j})\hat{x}_{j} = t_{j} \hat{y}_{j}, \quad  (F_{j}-\tilde{Q}_{j})^*\hat{y}_{j} = t_{j} \hat{x}_{j}. $$In view of equation \eqref{f-qj}, the latter equations imply 
	$$ B_{j-1}^* (G-Q_{j+1}) \bar{A}_{j-1} \hat{x}_{j} = t_{j} \hat{y}_{j},\quad A_{j-1}^T (G-Q_{j+1})^* B_{j-1}\hat{y}_{j}  = t_{j}\hat{x}_{j}.$$
By equation \eqref{hatxj-1},
	$$\hat{x}_{j} = A_{j-1}^T x_{j}, \quad \hat{y}_{j} =B_{j-1}^* y_{j}$$
Thus 
	$$B_{j-1}^* (G-Q_{j+1}) \bar{A}_{j-1} \hat{x}_{j} = B_{j-1}^* (G-Q_{j+1}) x_{j} = t_{j} B_{j-1}^*y_j, $$
or equivalently, $$B_{j-1}^* \big( (G-Q_{j+1})x_j - t_jy_j\big)=0, $$
and since, by the inductive hypothesis, $M_{B_{j-1}}$ is a unitary map, we have
	$$ (G-Q_{j+1})x_j =t_jy_j.$$
	Furthermore
$$A_{j-1}^T (G-Q_{j+1})^* B_{j-1}\hat{y}_{j}  = A_{j-1}^T (G-Q_{j+1})^*y_{j} =t_j A_{j-1}^Tx_j,   $$
or equivalently,
$$A_{j-1}^T \big((G-Q_{j+1})^*y_{j} - t_jx_j \big)=0. $$ By the inductive hypothesis, $M_{\bar{A}_{j-1}}$ is a unitary map, hence $$(G-Q_{j+1})^*y_{j} = t_jx_j,$$
and therefore $Q_{j+1}$ satisfies the required equations.  \end{proof}

\begin{lemma}\label{corona} Let
\begin{equation} \label{tildevjwjj-c}
\tilde{V}_{i} = \begin{pmatrix}
\alpha_{i-1}^T \cdots \alpha_0^T \xi_{i} & \bar{ \alpha}_{i} 
\end{pmatrix}, \quad \tilde{W}_{i}^T = \begin{pmatrix}
\beta_{i-1}^T  \cdots \beta_0^T \eta_{i} & \bar{\beta}_{i}
\end{pmatrix}, \;\;
i=0, 1, \dots, j,
\end{equation} be unitary-valued functions, 
as described in Lemma \ref{2.2}.  Then 
$$\alpha_l^T H^2(\Dd,\Cc^{n-l}) = H^2(\Dd,\Cc^{n-l-1})  $$ and 
$$ \beta_l^*(H^2(\Dd,\Cc^{m-l})^\perp = H^2(\Dd,\Cc^{m-l-1})^\perp, $$	
for all $l=0,\dots,j.$
\end{lemma}

\begin{proof}
Recall that, by Lemma \ref{L6.2}, for all $l=0,\dots,j,$ the inner, co-outer, quasi-continuous functions $\alpha_l, \beta_l$ of types $(n-l)\times (n-l-1)$ and $(m-l)\times (m-l-1)$ respectively, are left invertible. The rest of the proof is similar to
Lemmas \ref{a0h2} and \ref{beta0*h2}. \end{proof}

As a preparation for proof of the main inductive step we prove several propositions.

\begin{proposition}\label{xitelakxik0}
Let
\begin{equation}\label{tildevjwjj-Kj}
\tilde{V}_{i} = \begin{pmatrix}
\alpha_{i-1}^T \cdots \alpha_0^T \xi_{i} & \bar{\alpha}_{i} 
\end{pmatrix}
\end{equation} 
be unitary-valued functions, for $i=0,1,\dots, j$, 
as described in Lemma \ref{2.2}
Let $A_i= \alpha_0\alpha_1 \dots \alpha_i$, for $i = 0,1,\dots,j$,  and $A_{-1} = I_n$. Let  
$$V_i = \begin{pmatrix}
I_i & 0 \\ 0 & \tilde{V}_i
\end{pmatrix}, \; \text{for} \; i = 0,1,\dots,j,  $$ 
and let \begin{equation}\label{calkj+1}\mathcal{K}_{j+1}= V_0 \cdots V_{j} \begin{pmatrix}
0_{(j+1)\times 1} \\ H^2(\Dd,\Cc^{n-j-1})
\end{pmatrix}.\end{equation} Let $\xi_{(j)} = \xi_0 \telwe \dots \telwe \xi_j.$ Then, 
$$\xi_{(j)} \telwe \mathcal{K}_{j+1} = \xi_{(j)} \telwe H^2(\Dd,\Cc^n) $$ and the operator $(\xi_{(j)} \telwe \cdot) \colon  \mathcal{K}_{j+1} \to \xi_{(j)} \telwe H^2(\Dd,\Cc^n)$ is unitary.	
\index{$\xi_{(i)}$}
\end{proposition}

\begin{proof}
First let us prove that 
$$\xi_0 \telwe \dots \telwe \xi_j \telwe \mathcal{K}_{j+1} = \xi_0 \telwe \dots \telwe \xi_j \telwe H^2(\Dd,\Cc^n). $$
By equations \eqref{calkj+1} and \eqref{V0Vj}
\begin{equation}\label{K(j+1)}
\mathcal{K}_{j+1} = \bar{A}_{j}H^2(\Dd,\Cc^{n-j-1}).
\end{equation}
By Lemma \ref{corona},
$$\alpha_l^T H^2(\Dd,\Cc^{n-l}) = H^2(\Dd,\Cc^{n-l-1})  $$ 
for all $l=0,\dots,j.$ Thus
\begin{equation}\label{alphaH2}
\begin{array}{lll}
H^2(\Dd,\Cc^{n-j-1})&= \alpha_j^T H^2(\Dd,\Cc^{n-j})\\
&= \alpha_j^T \alpha_{j-1}^T H^2(\Dd,\Cc^{n-j+1})\\
&= \alpha_j^T \alpha_{j-1}^T \alpha_{j-2}^T H^2(\Dd,\Cc^{n-j+2})\\
&= \dots\\
&=  \alpha_j^T \alpha_{j-1}^T \dots \alpha_1^T H^2(\Dd,\Cc^{n-1})\\
&= \alpha_j^T \alpha_{j-1}^T \dots \alpha_1^T \alpha_0^T H^2(\Dd,\Cc^{n})\\
&= A_j^T H^2(\Dd,\Cc^{n}).
\end{array}
\end{equation}
By equations \eqref{K(j+1)} and  \eqref{alphaH2},
\begin{equation}\label{K(j+1)-2}
\mathcal{K}_{j+1} = \bar{A}_{j}H^2(\Dd,\Cc^{n-j-1})= \bar{A}_{j} A_j^T H^2(\Dd,\Cc^{n}).
\end{equation}
By Lemma \ref{xi=AkAkT},
\begin{equation}\label{AAk-K} 
\bar{A}_j A_j^T = I_n - \sum\limits_{k=0}^j \xi_k \xi_k^*.
\end{equation} 
By Proposition \ref{onxi}, $\{\xi_i(z)\}_{i=0}^j$ is an orthonormal set in $\Cc^n$ for almost every $z \in \Tt.$
Therefore, by equations \eqref{K(j+1)-2} and \eqref{AAk-K},
\begin{equation}\label{Kj+1H2}
\begin{array}{lll}
\xi_0 \telwe \dots \telwe \xi_j \telwe \mathcal{K}_{j+1} &= \xi_0 \telwe \dots \telwe \xi_j \telwe \bar{A}_{j} A_j^T H^2(\Dd,\Cc^{n})\\
&= \xi_0 \telwe \dots \telwe \xi_j \telwe (I_n - \sum\limits_{k=0}^j \xi_k \xi_k^*)H^2(\Dd,\Cc^n)\\
&= \xi_0 \telwe \dots \telwe \xi_j \telwe H^2(\Dd,\Cc^n). 
\end{array}
\end{equation}

 To show that the operator $(\xi_{(j)} \telwe \cdot)\colon\mathcal{K}_{j+1} \to \xi_{(j)} \telwe H^2(\Dd,\Cc^n)$ is unitary, it suffices to prove that, for every $\vartheta \in \mathcal{K}_{j+1},$ 
$$\| \xi_{(j)} \telwe \vartheta\|_{L^2(\Tt,\we^{j+2}\Cc^n)} = \| \vartheta \|_{L^2(\Tt,\Cc^n)} .$$

\noindent Let $\vartheta \in \mathcal{K}_{j+1}.$ Then, by Proposition \ref{we}, we have 
$$\begin{array}{llll} \| \xi_{(j)} \telwe \vartheta\|_{L^2(\Tt,\we^{j+2}\Cc^n)}^2 &= \langle \xi_{(j)} \telwe \vartheta , \xi_{(j)} \telwe \vartheta \rangle_{L^2(\Tt,\we^{j+2}\Cc^n)} \vspace{2ex} \\
&= \displaystyle\frac{1}{2\pi} \int_0^{2\pi}\langle \xi_{(j)}(\eiu) \telwe \vartheta(\eiu) , \xi_{(j)}(\eiu) \telwe \vartheta(\eiu) \rangle_{\we^{j+2}\Cc^n} d\theta \vspace{2ex} \\
&= \displaystyle\frac{1}{2\pi} \int_0^{2\pi} \det \begin{pmatrix}
\langle \xi_0(\eiu), \xi_0 (\eiu)\rangle_{\Cc^n} &\dots & \langle\xi_0(\eiu),\vartheta (\eiu)\rangle_{\Cc^n}\\
\langle\xi_1(\eiu) , \xi_0(\eiu) \rangle_{\Cc^n} & \dots &\langle \xi_1(\eiu) , \vartheta(\eiu) \rangle_{\Cc^n}\\
\vdots & \ddots & \vdots \\
\langle\vartheta(\eiu) , \xi_0(\eiu) \rangle_{\Cc^n}& \dots &\langle \vartheta(\eiu) , \vartheta(\eiu) \rangle_{\Cc^n}
\end{pmatrix}d\theta .\end{array}$$
By Proposition \ref{onxi}, $\{\xi_i(z)\}_{i=0}^j$ is an orthonormal set in $\Cc^n$ for almost every $z \in \Tt.$ 
Thus the latter integral is equal to
$$ \displaystyle\frac{1}{2\pi} \int_0^{2\pi} \det \begin{pmatrix}
1 & 0& \dots & \langle\xi_0(\eiu),\vartheta (\eiu)\rangle_{\Cc^n}\\
0 &1& \dots &\langle \xi_1(\eiu) , \vartheta(\eiu) \rangle_{\Cc^n}\\
\vdots & ~&\ddots & \vdots \\
\langle\vartheta(\eiu) , \xi_0(\eiu) \rangle_{\Cc^n}& ~&\dots &\langle \vartheta(\eiu) , \vartheta(\eiu) \rangle_{\Cc^n}
\end{pmatrix}d\theta.  $$ 
Note that since $\vartheta \in \mathcal{K}_{j+1},$ 
$$
\vartheta = \bar{A}_l A_j^T \psi= (I_n - \sum\limits_{i=0}^j \xi_i \xi_i^*)\psi$$ for some $\psi \in H^2(\Dd,\Cc^{n}).$ Then, for almost every $\eiu \in \Tt,$ 
$$\begin{array}{lll}
 \langle \xi_k (\eiu), \vartheta(\eiu) \rangle_{\Cc^n} 
&= \langle \xi_k (\eiu) , (I_n - \sum\limits_{i=0}^j \xi_i \xi_i^*)(\eiu) \psi(\eiu) \rangle_{ \Cc^n}\\
& = \langle \xi_k (\eiu) ,\psi(\eiu) \rangle_{ \Cc^n} -  \langle \xi_k, \xi_k \rangle_{\Cc^n} \langle \xi_k, \psi(\eiu) \rangle_{ \Cc^n} =0.
\end{array}$$ 
Hence 
$$ \| \xi_{(j)} \telwe \vartheta\|_{L^2(\Tt,\we^{j+2}\Cc^n)}^2=
\displaystyle\frac{1}{2\pi} \int_0^{2\pi} \det \begin{pmatrix}
1 & ~&0 &\dots & 0\\
0 & ~& 1&  \dots &0\\
\vdots & ~& ~& \ddots & \vdots \\
0&~ &~ & \dots &\langle \vartheta(\eiu) , \vartheta(\eiu) \rangle_{\Cc^n}
\end{pmatrix}d\theta ,$$
which yields
$$\begin{array}{clll}\displaystyle\frac{1}{2\pi} \int_0^{2\pi} \| \vartheta(\eiu)\|_{\Cc^n}^2 d\theta = \| \vartheta \|_{L^2(\Tt,\Cc^n)}^2. \end{array} $$  Therefore, the operator $ (\xi_{(j)} \telwe \cdot)\colon \mathcal{K}_{j+1} \to \xi_{(j)} \telwe H^2(\Dd,\Cc^n)$ is unitary. \end{proof}

\begin{proposition}\label{etatelwelj}
Let
\begin{equation} \tilde{W}_{i}^T = \begin{pmatrix}
\beta_{i-1}^T  \cdots \beta_0^T \eta_{i} & \bar{\beta}_{i}
\end{pmatrix}
\end{equation} be unitary-valued functions, for $i=0,1, \dots, j$,
as described in Lemma \ref{2.2}.
Let $B_i= \beta_0 \beta_1\dots\beta_i$, for $i=0,1,\dots, j$, and $B_{-1} = I_m$. 
Let  $W_{i}^T = \begin{pmatrix}
	I_i & 0 \\ 0 & \tilde{W}_i^T 
	\end{pmatrix},\; \text{for} \; i=0,1,\dots, j,
$ and let 
\begin{equation}\label{callj+1}
\mathcal{L}_{j+1} = W_0^* \cdots W_j^* \begin{pmatrix}
	0_{(j+1)\times 1} \\ H^2(\Dd,\Cc^{m-j-1})^\perp
	\end{pmatrix},\end{equation}
Let $\bar{ \eta}_{(j)} = \bar{ \eta}_0 \telwe \dots \telwe \bar{ \eta}_j.$ 
Then, 
	$$\bar{ \eta}_{(j)} \telwe \mathcal{L}_{j+1} = \bar{ \eta}_{(j)} \telwe H^2(\Dd,\Cc^m)^\perp$$ 
and the operator $(\bar{\eta}_{(j)} \telwe \cdot)\colon  \mathcal{L}_{j+1} \to \bar{\eta}_{(j)} \telwe   H^2(\Dd,\Cc^m)^\perp$ is unitary.
\index{$\bar{\eta}_{(j)}$}
\end{proposition}
\begin{proof}
First let us prove that 
$$\bar{ \eta}_0 \telwe \dots \telwe \bar{ \eta}_j\telwe \mathcal{L}_{j+1} = \bar{ \eta}_0 \telwe \dots \telwe \bar{ \eta}_j  \telwe H^2(\Dd,\Cc^m)^\perp.$$ 
By equations \eqref{callj+1} and \eqref{WjW0*},
\begin{equation}\label{L(j+1)}
\mathcal{L}_{j+1} = B_j  H^2(\Dd,\Cc^{m-j-1})^\perp.
\end{equation}
By Lemma \ref{corona},
$$ \beta_\ell^*(H^2(\Dd,\Cc^{m-\ell})^\perp = H^2(\Dd,\Cc^{m-\ell-1})^\perp, $$	
for all $\ell=0,\dots,j.$ Thus
\begin{equation}\label{betaH2}
\begin{array}{lll}
H^2(\Dd,\Cc^{m-j-1})^\perp &= \beta_j^*(H^2(\Dd,\Cc^{m-j})^\perp\\
& = \beta_j^*\beta_{j-1}^*(H^2(\Dd,\Cc^{m-j+1})^\perp\\
& = \beta_j^*\beta_{j-1}^*\beta_{j-2}^*(H^2(\Dd,\Cc^{m-j+2})^\perp\\
&= \dots\\
& = \beta_j^*\beta_{j-1}^* \dots \beta_1^*(H^2(\Dd,\Cc^{m-1})^\perp\\
& = \beta_j^*\beta_{j-1}^* \dots \beta_1^*\beta_0^*(H^2(\Dd,\Cc^{m})^\perp\\
& = B_j^*H^2(\Dd,\Cc^m)^\perp.
\end{array}
\end{equation}
By  equations \eqref{L(j+1)} and  \eqref{betaH2},
\begin{equation}\label{L(j+1)-2}
\mathcal{L}_{j+1} = B_j  H^2(\Dd,\Cc^{m-j-1})^\perp = B_j B_j^*H^2(\Dd,\Cc^m)^\perp.
\end{equation}
By Lemma \ref{eta=BkBk*},
\begin{equation}\label{BB_k*-L} 
B_j B_j^* = I_m - \sum\limits_{i=0}^j \bita_i \eta_i^T .
\end{equation} 
Thus 
\begin{equation}\label{L(j+1)-eta}
\mathcal{L}_{j+1} = (I_m - \sum\limits_{i=0}^j \bita_i \eta_i^T)H^2(\Dd,\Cc^m)^\perp. 
\end{equation}
By Proposition \ref{onxi}, $\{\bar{\eta}_i(z)\}_{i=0}^j$ is an orthonormal set in $\Cc^m$ for almost every $z \in \Tt.$
Therefore, by equations \eqref{L(j+1)-2} and \eqref{L(j+1)-eta},
\begin{equation}\label{Lj+1=H2}
\begin{array}{lll}
\bar{ \eta}_0 \telwe \dots \telwe \bar{ \eta}_j\telwe \mathcal{L}_{j+1} &=
\bar{ \eta}_0 \telwe \dots \telwe \bar{ \eta}_j\telwe (I_m - \sum\limits_{i=0}^j \bita_i \eta_i^T)H^2(\Dd,\Cc^m)^\perp\\
&=  \bar{ \eta}_0 \telwe \dots \telwe \bar{ \eta}_j  \telwe H^2(\Dd,\Cc^m)^\perp.
\end{array}
\end{equation}

\noindent To show that the operator $(\bar{\eta}_{(j)} \telwe \cdot)\colon  \mathcal{L}_{j+1} \to \bar{\eta}_{(j)} \telwe   H^2(\Dd,\Cc^m)^\perp$ is unitary, it suffices to prove that, for every $\varphi \in \mathcal{L}_{j+1},$ 
$$\|\bar{\eta}_{(j)} \telwe \varphi \|_{L^2(\Tt,\we^{j+2}\Cc^m)} = \| \varphi \|_{L^2(\Tt,\Cc^m)}. $$
		
\noindent By Proposition \ref{we}, we have 
$$ \begin{array}{clllll}\|\bar{\eta}_{(j)} \telwe \varphi \|_{L^2(\Tt,\we^{j+2}\Cc^m)}^2 &= \langle \bar{\eta}_{(j)} \telwe \varphi , \bar{\eta}_{(j)} \telwe \varphi \rangle_{ L^2(\Tt,\we^{j+2}\Cc^m)} \vspace{2ex} \\

&= \displaystyle\frac{1}{2\pi} \int_0^{2\pi} \det \begin{pmatrix} \langle \bar{\eta}_0(\eiu) , \bar{\eta}_0(\eiu) \rangle_{\Cc^m} & \dots & \langle\bar{\eta}_0(\eiu) , \varphi (\eiu) \rangle_{\Cc^m} \\
 \langle\bar{\eta}_1 (\eiu), \bar{\eta}_0(\eiu) \rangle_{\Cc^m}& \dots & \langle\bar{\eta}_1 (\eiu),\varphi(\eiu) \rangle_{\Cc^m} \\
 \vdots & \ddots & \vdots \\
\langle \varphi(\eiu), \bar{\eta}_0(\eiu) \rangle_{\Cc^m}& \dots &\langle  \varphi(\eiu), \varphi(\eiu) \rangle_{\Cc^m} 
  \end{pmatrix} d\theta .\end{array} $$
By Proposition \ref{onxi}, the set $\{ \bar{\eta}_i(z)\}_{i=0}^j$ is orthonormal in $\Cc^n$ for almost every $z \in \Tt.$ 
Then the latter integral is equal to 
  $$\displaystyle\frac{1}{2\pi} \int_0^{2\pi} \det \begin{pmatrix}1 & 0& \dots & \langle\bar{\eta}_0(\eiu) , \varphi (\eiu) \rangle_{\Cc^m} \\
  0& 1&\dots & \langle\bar{\eta}_1 (\eiu),\varphi(\eiu) \rangle_{\Cc^m} \\
  \vdots & ~ &\ddots & \vdots \\
  \langle \varphi(\eiu), \bar{\eta}_0(\eiu) \rangle_{\Cc^m}& \dots& \dots &\langle  \varphi(\eiu), \varphi(\eiu) \rangle_{\Cc^m} 
  \end{pmatrix} d\theta .$$ 
Note that since $\varphi \in \mathcal{L}_{j+1},$
$$
\varphi = (I_m - \sum\limits_{i=0}^j \bita_i \eta_i^T) \psi,$$
for some $\psi \in H^2(\Dd,\Cc^m)^\perp. $ Then, for almost every 
$\eiu \in \Tt,$ 
$$\begin{array}{lll}
   \langle\bar{\eta}_k (\eiu),\varphi(\eiu) \rangle_{\Cc^m} 
&= \langle\bar{\eta}_k (\eiu) , (I_m - \sum\limits_{i=0}^j \bita_i \eta_i^T)(\eiu) \psi(\eiu) \rangle_{\Cc^m} \vspace{2ex} \\
   &= \langle\bar{\eta}_k (\eiu) ,\psi(\eiu) \rangle_{\Cc^m} - 
\langle\bar{\eta}_k (\eiu) ,\bar{\eta}_k (\eiu)\rangle_{\Cc^m}
\langle\bar{\eta}_k (\eiu) ,\psi(\eiu) \rangle_{\Cc^m} =0. \vspace{2ex} 
   \end{array} $$ 
Thus   
  $$
\begin{array}{clll}
 \|\bar{\eta}_{(j)} \telwe \varphi \|_{L^2(\Tt,\we^{j+2}\Cc^m)}^2 = &\displaystyle\frac{1}{2\pi} \int_0^{2\pi} \det \begin{pmatrix}1 & 0& \dots & 0\\
  0& 1&\dots &0\\
  \vdots & ~ &\ddots & \vdots \\
  0& \dots& \dots &\langle  \varphi(\eiu), \varphi(\eiu) \rangle_{\Cc^m} 
  \end{pmatrix} d\theta \vspace{2ex} \\
  &=  \displaystyle\frac{1}{2\pi} \int_0^{2\pi} \langle  \varphi(\eiu), \varphi(\eiu) \rangle_{\Cc^m}  d\theta = \| \varphi \|_{L^2(\Tt,\Cc^m)}^2. \end{array}$$ 
Hence the operator $(\bar{\eta}_{(j)} \telwe \cdot)\colon  \mathcal{L}_{j+1} \to \bar{\eta}_{(j)} \telwe   H^2(\Dd,\Cc^m)^\perp$ is unitary.\end{proof}

\begin{proposition}\label{callj+1perp}
	With the notation of Proposition \ref{etatelwelj}
	$$\mathcal{L}_{j+1}^\perp = \{ f \in L^2(\Tt,\Cc^m): \beta_{j}^* \cdots \beta_{0}^*f \in H^2(\Dd,\Cc^{m-j-1})\}. $$
\end{proposition}
\begin{proof}
Clearly $\mathcal{L}_{j+1} = \beta_0 \cdots \beta_j H^2(\Dd,\Cc^{m-j-1})^\perp.$ 
The general element of $\beta_0 \cdots \beta_j H^2(\Dd,\Cc^{m-j-1})^\perp$ is $\beta_0 \cdots \beta_j \bar{z} \bar{g}$ with $ g \in H^2 (\Dd,\Cc^{m-j-1})$. A function $f\in L^2(\Tt,\Cc^m)$ belongs to $\mathcal{L}_{j+1}^\perp$ 
if and only if 
$$\langle f, \beta_0 \cdots \beta_j \bar{z} \bar{g}\rangle_{L^2(\Tt,\Cc^m)}=0 \quad \text{for all} \quad g \in H^2(\Dd,\Cc^{m-j-1}) $$
 if and only if 
$$\displaystyle \frac{1}{2\pi} \int_0^{2\pi} \langle f(\eiu), \beta_0(\eiu)\cdots \beta_j(\eiu) e^{-i\theta} \bar{g}(\eiu) \rangle_{\Cc^m}d\theta =0  \quad \text{for all} \quad  g \in H^2(\Dd,\Cc^{m-j-1})  $$ if and only if
$$ \displaystyle \frac{1}{2\pi} \int_0^{2\pi} \langle  \beta_j^*(\eiu)\cdots\beta_0^*(\eiu)f(\eiu),  e^{-i\theta} \bar{g}(\eiu) \rangle_{\Cc^{m-2}}d\theta =0  \quad \text{for all} \quad g \in H^2(\Dd,\Cc^{m-j-1}).$$

\noindent The latter statement is equivalent to the assertion that 
$\beta_j^*\cdots \beta_0^* f $ is orthogonal to \linebreak$H^2(\Dd,\Cc^{m-j-1})^\perp$ in $L^2(\Tt,\Cc^{m-j-1}),$ which holds if and only if $$\beta_j^*\cdots \beta_0^* f \in H^2(\Dd,\Cc^{m-j-1}).$$ Thus 

	$$\mathcal{L}_{j+1}^\perp = \{ f \in L^2(\Tt,\Cc^m): \beta_{j}^* \cdots \beta_{0}^*f \in H^2(\Dd,\Cc^{m-j-1})\} $$as required. \end{proof}

Let us proceed to the main theorem of this section. 

\begin{theorem}\label{Tkismultipleofhankel}
Let $m,n$ be positive integers such that $\min(m,n)\geq2.$ Let  $G$ be in $H^\infty(\Dd,\CCmn)+C(\Tt,\CCmn)$. In the notation of the algorithm \ref{Alg_statement}, let 
$$(\xi_0 \telwe \dots \telwe \xi_{j-1}\telwe v_{j}, \bar{\eta}_0\telwe \dots \telwe \bar{\eta}_{j-1}  \telwe w_{j} ) $$ 
be a Schmidt pair for $T_{j}$ corresponding to $t_{j}= \|T_{j}\| \neq 0.$ Let 
 $h_{j} \in H^2(\Dd,\Cc)$ be the scalar outer factor of $$\xi_0 \telwe \dots \telwe \xi_{j-1}\telwe v_{j}.$$
 Let 
$$x_{j} = (I_{n} - \xi_0 \xi_0^* - \dots - \xi_{j-1}\xi_{j-1}^* ) v_{j} ,$$
$$y_{j} = (I_{m} - \bar{\eta}_0 \eta_0^T - \dots - \bar{\eta}_{j-1} \eta_{j-1}^T )w_{j}$$ and
$$\xi_{j} = \frac{x_{j}}{h_{j}} , \quad \bar{\eta}_{j} = \frac{z y_{j}}{\bar{h}_{j}} .$$
For $i=0,1,\dots, j$, let
\begin{equation} \label{tildevjwjj+1}
\tilde{V}_{i} = \begin{pmatrix}
\alpha_{i-1}^T \cdots \alpha_0^T \xi_{i} & \bar{ \alpha}_{i} 
\end{pmatrix}, \quad \tilde{W}_{i}^T = \begin{pmatrix}
\beta_{i-1}^T  \cdots \beta_0^T \eta_{i} & \bar{\beta}_{i}
\end{pmatrix}
\end{equation} be unitary-valued functions, 
as described in Lemma \ref{2.2}. Let 
 $$V_{j} = \begin{pmatrix}
	I_{j} & 0 \\
	0 & \tilde{V}_{j} 
	\end{pmatrix}, \quad W_{j} = \begin{pmatrix}
	I_{j} & 0 \\
	0 & \tilde{W}_{j}
	\end{pmatrix}. $$
Let $A_j= \alpha_0\alpha_1 \dots \alpha_j$, $A_{-1} = I_n$,  $B_j= \beta_0 \beta_1\dots\beta_j$ and $B_{-1} = I_m$.
Let
$$X_{j+1} = \xi_0 \telwe \dots \telwe \xi_{j} \telwe H^2(\Dd,\Cc^n) \subset H^2(\Dd,\we^{j+2}{\Cc^n}), $$ 
and let
$$Y_{j+1} = \bar{ \eta}_0 \telwe \dots \telwe \bar{ \eta}_{j} \telwe H^2(\Dd,\Cc^m)^\perp \subset H^2(\Dd,\we^{j+2}\Cc^m)^\perp. $$  
Let 
$$T_{j+1} (\xi_0 \telwe \dots \telwe \xi_j \telwe x) = P_{Y_{j+1}} (\bar{ \eta}_0 \telwe \dots \telwe \bar{ \eta}_{j} \telwe (G-Q_{j+1})x )$$ 
for all $x \in H^2(\Dd,\Cc^n)$, where $Q_{j+1}$ satisfies 
\begin{equation} \label{G_Q(j+1)xy}
	(G-Q_{j+1})x_i = t_i y_i, \; \text{and}\; \;(G-Q_{j+1})^*y_i=t_ix_i,\quad \text{for}\; i=0,1,\dots,j.
\end{equation}
Let
\begin{equation}\label{k(j+1)l{j+1}}
\mathcal{K}_{j+1} =  V_0 \cdots V_{j} 
\begin{pmatrix}
0_{(j+1)\times 1} \\ H^2(\Dd,\Cc^{n-j-1})
\end{pmatrix}, \quad 
\mathcal{L}_{j+1} = W_0^* \cdots W_j^* \begin{pmatrix}
	0_{(j+1)\times 1} \\ H^2(\Dd,\Cc^{m-j-1})^\perp
	\end{pmatrix}.
\end{equation}
Let the operator $\Gamma_{j+1} \colon  \mathcal{K}_{j+1} \to \mathcal{L}_{j+1}$ be given by $$\Gamma_{j+1} = P_{\mathcal{L}_{j+1}} M_{G-Q_{j+1}}|_{\mathcal{K}_{j+1}}.$$  Then  \begin{itemize}	
 	\item[{\rm (i)}] The maps 
$$ M_{\bar{A}_j} \colon H^2(\Dd,\Cc^{n-j-1}) \to  
\mathcal{K}_{j+1}  \colon x \mapsto \bar{A}_j x,\;\; \text{and} \;\;
M_{B_j}\colon H^2(\Dd,\Cc^{m-j-1})^\perp \to \mathcal{L}_{j+1} \colon y \mapsto  {B_j} y $$ 
are unitaries.
 	\item[{\rm (ii)}] The maps $(\xi_0\telwe \dots \telwe \xi_j \telwe \cdot)\colon\mathcal{K}_{j+1}\to X_{j+1},$ $(\bita_0\telwe \dots \telwe \bita_j \telwe\cdot)\colon\mathcal{L}_{j+1}\to Y_{j+1} $ are unitaries.
 	\item[{\rm (iii)}] the following diagram commutes 
	\begin{equation}\label{t(j+1)comm}
	\begin{array}{clllll}		
	H^2(\Dd,\Cc^{n-j}) &\xrightarrow{M_{\bar{\alpha}_0\cdots \bar{ \alpha}_{j}}} & \mathcal{K}_{j+1} &\xrightarrow{\xi_{(j)} \telwe \cdot}& \xi_{(j)} \telwe H^2 (\Dd, \Cc^n)=X_{j+1}\\
	\Big\downarrow\rlap{$\scriptstyle H_{F_{j+1}}  $} & ~  &\Big\downarrow\rlap{$\scriptstyle \Gamma_{j+1}$}  &~&\hspace{3ex}\Big\downarrow\rlap{$\scriptstyle T_{j+1}$} \\
	H^2(\Dd,\Cc^{m-j})^\perp &\xrightarrow{M_{\beta_0 \cdots \beta_{j}}}&\mathcal{L}_{j+1} &\xrightarrow{\bar{\eta}_{(j)} \telwe \cdot } & \bar{\eta}_{(j)} \telwe H^2 (\Dd, \Cc^m)^\perp =Y_{j+1},
	\end{array}\end{equation}
where $F_{j+1} \in H^\infty(\Dd,\Cc^{(m-j-1)\times(n-j-1)})+ C(\Tt,\Cc^{(m-j-1)\times(n-j-1)})$  is the function defined in Proposition \ref{tildvjwj};
\item [{\rm (iv)}] ${\Gamma}_{j+1}$ and $T_{j+1}$ are compact operators;
\item[{\rm (v)}] $\|T_{j+1}\|=\|\Gamma_{j+1}\|=\|H_{F_{j+1}}\|=t_{j+1}.$
 \end{itemize} 
\end{theorem}

\begin{proof}
{\rm (i)} It follows from Lemma \ref{3.1constr}.
{\rm (ii)} follows from Propositions \ref{xitelakxik0} and \ref{etatelwelj}.		
{\rm (iii)}
 By Theorem \ref{1.8}, there exists a function $Q_{j+1} \in H^\infty(\Dd, \CCmn)$ such that
the sequence 
$$\left(s_0^\infty(G-Q_{j+1}),s_1^\infty(G-Q_{j+1}), \dots , s_{j+1}^\infty(G-Q_{j+1}) \right)$$
 is lexicographically minimized.  
By Proposition \ref{g-qjj}, any such  $Q_{j+1}$ satisfies 
\begin{equation} \label{G_Qjxy}
	(G-Q_{j+1})x_i = t_i y_i, (G-Q_{j+1})^*y_i=t_ix_i,\quad \text{for}\; i=0,1,\dots,j.
\end{equation}
By Proposition \ref{Twell}, $T_{j+1}$ is well-defined and is independent of the choice of $Q_{j+1}\in H^\infty(\Dd,\CCmn)$ satisfying equations
\eqref{G_Qjxy}. We can choose $Q_{j+1}$ which  minimises $$\left(s_0^\infty(G-Q_{j+1}),s_1^\infty(G-Q_{j+1}), \dots , s_{j+1}^\infty(G-Q_{j+1}) \right),$$ and therefore   satisfies  equations
\eqref{G_Qjxy}.	 
Consider the following diagram. 	 	 
	 \begin{equation}\label{diagr1122}	\begin{array}{clllll}	 
	 
	 &\mathcal{K}_{j+1} &\xrightarrow{\xi_0 \telwe\cdots \telwe \xi_{j} \telwe \cdot}& \xi_0 \telwe \cdots \telwe \xi_{j} \telwe  H^2 (\Dd, \Cc^n)=X_{j+1}\\
	 &\Big\downarrow\rlap{$\scriptstyle \Gamma_{j+1}$}  &~ &\hspace{13ex}\Big\downarrow\rlap{$\scriptstyle T_{j+1}$} \\
	 &	\mathcal{L}_{j+1} &\xrightarrow{\bar{\eta}_0 \telwe \cdots \telwe \bar{{\eta}}_{j} \telwe \cdot } & \bar{\eta}_0 \telwe \cdots \telwe \bar{\eta}_{j} \telwe H^2 (\Dd, \Cc^m)^\perp =Y_{j+1}.
	 \end{array}\end{equation} 
Let us prove first that diagram \eqref{diagr1122} commutes. 
By Proposition \ref{tildvjwj}, every  $Q_{j+1}\in H^\infty(\Dd,\CCmn)$, which  minimises $$\left(s_0^\infty(G-Q_{j+1}),s_1^\infty(G-Q_{j+1}), \dots , s_{j+1}^\infty(G-Q_{j+1}) \right),$$  satisfies the following equation (see equation \eqref{epsilonnj}).
\begin{equation}\label{G-Q(j+1)}
G - Q_{j+1} =  W_0^* W_1^* \cdots W_{j}^* \begin{pmatrix}
	t_0 u_0  &    0   &  \cdots &  0          &0_{1\times (n-j-1)}\\
	0        &t_1u_1  & \dots   &  0        &0_{1\times (n-j-1)} \\
	\vdots   &  \vdots &  \ddots & \vdots    &\vdots \\
	0        &   0    &   \cdots  &t_{j}u_{j} & 0\\ 
	0_{(m-j-1)\times 1} & 0_{(m-j-1)\times 1} & \dots &    \dots    &     (F_{j+1}+H^\infty)\cap B(t_{j})
	\end{pmatrix}V_{j}^* \cdots V_0^* , \end{equation} 
Thus, for every $\chi \in H^2(\Dd,\Cc^{n-j-1})$,
\begin{equation}\label{G-Q(j+1)-2}
(G - Q_{j+1}) V_0 \cdots V_{j}
\begin{pmatrix}
	0_{(j+1)\times 1}\\H^2(\Dd,\Cc^{m-j-1})^\perp 
\end{pmatrix} = \hspace{6cm}
\end{equation} 
$$ W_0^* W_1^* \cdots W_{j}^*
\begin{pmatrix}
	t_0 u_0  &    0   &  \cdots &  0          &0_{1\times (n-j-1)}\\
	0        &t_1u_1  & \dots   &  0        &0_{1\times (n-j-1)} \\
	\vdots   &  \vdots &  \ddots & \vdots    &\vdots \\
	0        &   0    &   \cdots  &t_{j}u_{j} & 0\\ 
	0_{(m-j-1)\times 1} & 0_{(m-j-1)\times 1} & \dots &    \dots    &     (F_{j+1}+H^\infty)\cap B(t_{j})
	\end{pmatrix}
  \begin{pmatrix}
	0_{(j+1)\times 1}\\H^2(\Dd,\Cc^{m-j-1})^\perp 
\end{pmatrix},
$$ 
for some  $F_{j+1} \in H^\infty(\Dd,\Cc^{(m-j-1)\times (n-j-1)})+C(\Tt,\Cc^{(m-j-1)\times (n-j-1)}),$ for the quasi-continuous unimodular functions  $u_i = \frac{\bar{z} \bar{h}_i}{h_i}$, for all $i=0,\dots, j ,$ for the closed ball $B(t_{j})$ of radius $t_j$ in $L^\infty(\Tt,\Cc^{(m-j-1)\times (n-j-1)})$.
By equation \eqref{WjW0*}, 
$$W_0^* W_1^* \cdots W_{j}^*= \begin{pmatrix}
\bar{\eta}_{0} & B_0B_0^*\bar{\eta}_{1} & \dots & B_{j-1} B_{j-1}^* \bar{\eta}_{j} & B_{j} \end{pmatrix}. $$
By equation \eqref{V0Vj},
\begin{equation}\label{V0Vj-2}
 V_0\cdots V_j = \begin{pmatrix}
	\xi_0 & \bar{A}_0 A_0^T \xi_1 & \bar{A}_1A_1^T\xi_2& \dots & \bar{A}_{j-1}A_{j-1}^T \xi_{j-1}& \bar{A}_{j}
\end{pmatrix}.	
\end{equation}
Therefore, by equation \eqref{G-Q(j+1)-2}, for every $\chi \in H^2(\Dd,\Cc^{n-j-1})$,
\begin{equation}\label{(G-Q)Ax=BFX}
(G-Q_{j+1}) \bar{A}_j \chi = B_j F_{j+1} \chi.
\end{equation}

\noindent A typical element $x \in \mathcal{K}_{j+1}$ is of the form $ x = \bar{A}_j \chi,$ for some $\chi \in H^2(\Dd,\Cc^{n-j-1}).$ Then, by Proposition \ref{xitelakxik0},
	 $$ (\xi_0 \telwe\dots \telwe \xi_{j} \telwe \cdot)\bar{A}_j \chi = \xi_0 \telwe\dots \telwe \xi_{j} \telwe \bar{A}_j \chi \in X_{j+1}.$$ 
Therefore, by the definition of $T_{j+1}$ and by equation \eqref{(G-Q)Ax=BFX},
 $$\begin{array}{cllll}
T_{j+1} ( \xi_0 \telwe\dots \telwe \xi_{j} \telwe \bar{A}_j \chi) &= P_{Y_{j+1}} (\bar{\eta}_0 \telwe \cdots \telwe \bar{\eta}_{j} \telwe (G-Q_{j+1}) \bar{A}_j \chi)\vspace{2ex} \\
	 &=P_{Y_{j+1}} (\bar{\eta}_0 \telwe \cdots \telwe \bar{\eta}_{j} \telwe B_j F_{j+1} \chi).\end{array}$$
Furthermore, by the definition of $\Gamma_{j+1}$ and by equation \eqref{(G-Q)Ax=BFX},
	 $$\begin{array}{cllll}
(\bar{\eta}_0 \telwe \cdots \telwe \bar{\eta}_{j} \telwe \cdot)\Gamma_{j+1}(\bar{A}_j\chi) &= 
\bar{\eta}_0 \telwe \cdots \telwe \bar{\eta}_{j} \telwe  P_{\mathcal{L}_{j+1}} (G-Q_{j+1})(\bar{A}_j\chi)\\
	 &=
\bar{\eta}_0 \telwe \cdots \telwe \bar{\eta}_{j} \telwe P_{\mathcal{L}_{j+1}} B_j F_{j+1} \chi. 
\end{array}	 $$ 
In order to prove the commutativity of diagram \eqref{diagr1122}, we need to show that, for every  $\chi \in H^2(\Dd,\Cc^{n-j-1})$,
$$T_{j+1} ( \xi_0 \telwe\dots \telwe \xi_{j} \telwe \bar{A}_j \chi)= 
(\bar{\eta}_0 \telwe \cdots \telwe \bar{\eta}_{j} \telwe \cdot)\Gamma_{j+1}(\bar{A}_j\chi) .$$
Hence we must prove that,  for every  $\chi \in H^2(\Dd,\Cc^{n-j-1})$,
	 $$ \bar{\eta}_0 \telwe \cdots \telwe \bar{\eta}_{j} \telwe P_{\mathcal{L}_{j+1}} B_j F_{j+1} \chi \in Y_{j+1}$$
	 and that 
$$ \bar{\eta}_0 \telwe \cdots \telwe \bar{\eta}_{j} \telwe \left(B_j F_{j+1} \chi- P_{\mathcal{L}_{j+1}} B_j F_{j+1} \chi \right), \; \text{which is equal to } \; \bar{\eta}_0 \telwe \cdots \telwe \bar{\eta}_{j} \telwe  P_{\mathcal{L}_{j+1}^\perp} B_j F_{j+1} \chi, $$ 
is orthogonal to $Y_{j+1}.$
Observe that, by Proposition \ref{etatelwelj}, for any $ \chi \in H^2(\Dd,\Cc^{n-j-1})$, $ \bar{\eta}_0 \telwe \cdots \telwe \bar{\eta}_{j} \telwe P_{\mathcal{L}_{j+1}} B_j F_{j+1} \chi$ is indeed an element of $Y_{j+1}.$ 
To prove that $$ \bar{\eta}_0 \telwe \cdots \telwe \bar{\eta}_{j} \telwe  P_{\mathcal{L}_{j+1}^\perp} B_j F_{j+1} \chi $$ is orthogonal to $Y_{j+1},$
 it suffices to prove that 
	 $$\langle \bar{\eta}_{(j)} \telwe \Phi, \bar{\eta}_{(j)} \telwe \psi \rangle_{L^2(\Tt,\we^{j+2}\Cc^m)}=0$$ 
for $\Phi = P_{\mathcal{L}_{j+1}^\perp} B_j F_{j+1} \chi$, for all $ \chi \in H^2(\Dd,\Cc^{n-j-1})$ and
for all $\psi \in H^2(\Dd,\Cc^{m})^\perp.$ By Proposition \ref{we},
	 $$\begin{array}{llll}
	 &\langle \bar{\eta}_{(j)} \telwe \Phi, \bar{\eta}_{(j)} \telwe \psi \rangle_{L^2(\Tt,\we^{j+2}\Cc^m)} \vspace{3ex}\\	 
	 &\vspace{1ex}= \displaystyle\frac{1}{2\pi} \int_0^{2\pi} \det \begin{pmatrix} \langle \bar{\eta}_0(\eiu) , \bar{\eta}_0(\eiu) \rangle_{\Cc^m} & \dots & \langle\bar{\eta}_0(\eiu) , \psi (\eiu) \rangle_{\Cc^m}  \\
	 \langle\bar{\eta}_1 (\eiu), \bar{\eta}_0(\eiu) \rangle_{\Cc^m}& \dots & \langle\bar{\eta}_1 (\eiu),\psi(\eiu) \rangle_{\Cc^m} \\
	 \vdots & \ddots & \vdots \\
	 \langle \Phi(\eiu), \bar{\eta}_0(\eiu) \rangle_{\Cc^m}& \dots &\langle  \Phi(\eiu), \psi(\eiu) \rangle_{\Cc^m} 
	 \end{pmatrix} d\theta
	 \end{array} $$ for all $ \psi \in H^2(\Dd,\Cc^{m})^\perp.$ Recall, by Proposition \ref{onxi}, the set $\{\eta_i\}_{i=0}^j$ is an orthonormal set in $\Cc^m$ almost everywhere on $\Tt.$ Hence
	 
	 $$\begin{array}{llll}
	 &\langle \bar{\eta}_{(j)} \telwe \Phi, \bar{\eta}_{(j)} \telwe \psi \rangle_{L^2(\Tt,\we^{j+2}\Cc^m)} \vspace{3ex}\\	 
	 & \vspace{1ex} = \displaystyle\frac{1}{2\pi} \int_0^{2\pi} \det \begin{pmatrix} 1 & 0&\dots & \langle\bar{\eta}_0(\eiu) , \psi (\eiu) \rangle_{\Cc^m}  \\
	 0& 1& \dots & \langle\bar{\eta}_1 (\eiu),\psi(\eiu) \rangle_{\Cc^m} \\
	 \vdots & ~&\ddots & \vdots \\
	 \langle \Phi(\eiu), \bar{\eta}_0(\eiu) \rangle_{\Cc^m}& ~&\dots &\langle  \Phi(\eiu), \psi(\eiu) \rangle_{\Cc^m} 	 
	 \end{pmatrix} d\theta. 
	 \end{array}$$ 
Multiplying the $k$-th column by $\langle \bar{\eta}_k(\eiu), \psi (\eiu) \rangle_{\Cc^m}  $ and subtracting it from the last column of the determinant above, we obtain
	 $$\displaystyle\frac{1}{2\pi} \int_0^{2\pi} \scriptsize{\det \begin{pmatrix} 1 & 0&\dots & 0  \\
	 	0& 1& \dots & 0 \\
	 	\vdots & ~&\ddots & \vdots \\
	 	\langle \Phi(\eiu), \bar{\eta}_0(\eiu) \rangle_{\Cc^m}& ~&\dots &\langle  \Phi(\eiu), \psi(\eiu) \rangle_{\Cc^m} \\ &~ &~  &-\sum_{i=0}^{j} \langle \Phi(\eiu), \bar{\eta}_i (\eiu) \rangle_{\Cc^m} \langle \bar{\eta}_i (\eiu), \psi(\eiu)\rangle_{\Cc^m} \end{pmatrix}}d\theta, $$which is equal to
	 	$$ \begin{array}{llll}
	  \displaystyle\frac{1}{2\pi} \int_0^{2\pi} \psi^*(\eiu) \Phi(\eiu) &- \sum_{i=0}^{j}\psi^*(\eiu) \bar{\eta}_i(\eiu) \eta_i^T(\eiu)\Phi(\eiu)d\theta  \vspace{2ex}  \\
	 &= \displaystyle\frac{1}{2\pi} \int_0^{2\pi} \psi^*(\eiu) \bigg(I_m - \sum_{i=0}^{j} \bar{\eta}_i(\eiu) \eta_i^T(\eiu)\bigg)\Phi(\eiu)d\theta.
	 \end{array}$$
	 
	 \noindent Then 
	 $$\langle \bar{\eta}_{(j)} \telwe \Phi, \bar{\eta}_{(j)} \telwe \psi \rangle_{L^2(\Tt,\we^{j+2}\Cc^m)}=0$$ for all $\psi \in H^2(\Dd,\Cc^{m})^\perp$ if and only if 
	 $$\displaystyle\frac{1}{2\pi} \int_0^{2\pi}\left\langle \bigg(I_m - \sum_{i=0}^{j} \bar{\eta}_i(\eiu) \eta_i^T(\eiu)\bigg)\Phi(\eiu), \psi(\eiu) \right\rangle_{\Cc^m} =0 $$ for all $\psi \in H^2(\Dd,\Cc^{m})^\perp,$
	 which holds if and only if 
$$\bigg(I_m - \sum_{i=0}^{j} \bar{\eta}_i \eta_i^T\bigg)\Phi \in H^2(\Dd,\Cc^m). $$
By Lemma \ref{eta=BkBk*},
\begin{equation}\label{BB_j*} 
B_j B_j^* = I_m - \sum\limits_{i=0}^j \bita_i \eta_i^T. 
\end{equation} 	
Thus	 
	 $$\langle \bar{\eta}_{(j)} \telwe \Phi, \bar{\eta}_{(j)} \telwe \psi \rangle_{L^2(\Tt,\we^{j+2}\Cc^m)}=0$$ for all $\psi \in H^2(\Dd,\Cc^{m})^\perp$ if and only if $$\displaystyle\frac{1}{2\pi} \int_0^{2\pi} \langle B_{j-1}(\eiu)B_{j-1}^*(\eiu)\Phi(\eiu), \psi(\eiu)\rangle_{\Cc^m}=0, $$ which holds if and only if  $B_j B_j^* \Phi \in H^2(\Dd,\Cc^m). $
By Proposition \ref{callj+1perp}, $\Phi = P_{\mathcal{L}_{j+1}^\perp} B_j F_{j+1} \chi $ satisfies the following property $B_j^* \Phi \in H^2(\Dd,\Cc^{m-j-1}).$ 
 Hence diagram \eqref{diagr1122} commutes. 

Recall that, by Lemma \ref{3.2constr}, the following diagram also commutes
	 
\begin{equation}\label{hfj+1}	\begin{array}{clllll}
&H^2(\Dd,\Cc^{n-j-1})  &\xrightarrow{M_{\bar{A}_j}}& \mathcal{K}_{j+1} \\
&\Big\downarrow\rlap{$\scriptstyle H_{F_{j+1}}$}  &~ &\Big\downarrow\rlap{$\scriptstyle {\Gamma}_{j+1}$} \\
&	H^2(\Dd,\Cc^{m-j-1})^\perp &\xrightarrow{M_{B_j}} & \mathcal{L}_{j+1}.
\end{array}\end{equation}

{\rm(iv)} 
 Since $F_{j+1} \in  H^\infty(\Dd,\Cc^{(n-j-1)\times (m-j-1)}) + C(\Tt, \Cc^{(n-j-1)\times (m-j-1)}),$ by Hartman's Theorem, the Hankel operator $H_{F_{j+1}}$ is compact.
Since diagram \eqref{hfj+1} commutes and 
the operators $ M_{\bar{A}_{j}}$ and $\;M_{B_j}$  are unitaries, ${\Gamma}_{j+1}$
is compact.
By (iii), 
$$ (\bita_0 \telwe \dots \telwe \bita_j\telwe \cdot) \circ (M_{B_j} \circ H_{F_{j+1}} \circ M_{\bar{A}_{j}}^*) \circ (\xi_0\telwe \dots  \telwe \xi_j \telwe \cdot)^*=T_{j+1} .$$ 
By (i) and (ii), the operators
$ M_{\bar{A}_{j}},\;M_{B_j},$ $(\xi_0\telwe  \dots \telwe \xi_j \telwe \cdot)$ and $(\bita_0 \telwe \dots \telwe \bita_j \telwe \cdot)$ are unitaries,
Hence $T_{j+1}$ is a compact operator.


{\rm (v)} Since diagram \eqref{t(j+1)comm} commutes and 
the operators $ M_{\bar{A}_{j}},\;M_{B_j},$ $(\xi_0\telwe  \dots \telwe \xi_j \telwe \cdot)$ and $(\bita_0 \telwe \dots \telwe \bita_j \telwe \cdot)$ are unitaries,
$$\|T_{j+1}\|=\|\Gamma_{j+1}\|=\|H_{F_{j+1}}\|=t_{j+1}. $$ \end{proof}

\begin{lemma}\label{coofschpairsj}
Let $v_{j+1} \in H^2(\Dd,\Cc^n)$ and $w_{j+1} \in H^2(\Dd,\Cc^m)^\perp$ be such that $$(\xi_0 \telwe\cdots \telwe \xi_{j}\telwe v_{j+1}, \bar{\eta}_0 \telwe \cdots \telwe \bita_j\telwe  w_{j+1})$$ is a Schmidt pair for the operator $T_{j+1}$ corresponding to $\|T_{j+1}\|.$ Then\begin{enumerate}
		\item[\emph{(i)}] there exist $x_{j+1} \in \mathcal{K}_{j+1}$ and $y_{j+1}\in \mathcal{L}_{j+1}$ such that $(x_{j+1},y_{j+1}) $ is a Schmidt pair for the operator $\Gamma_{j+1}$;
		\item[\emph{(ii)}] for any $x_{j+1} \in \mathcal{K}_{j+1}$ and $y_{j+1}\in \mathcal{L}_{j+1}$ such that $$
		\xi_0 \telwe \cdots \telwe \xi_j\telwe  x_{j+1}= \xi_0 \telwe \cdots  \telwe \xi_j \telwe v_{j+1},\quad 	
		\bar{ \eta}_0 \telwe \cdots \telwe \bita_j \telwe  y_{j+1} = \bar{ \eta}_0 \telwe  \cdots\telwe \bita_j\telwe w_{j+1},$$ the pair $(x_{j+1},y_{j+1})$ is a Schmidt pair for $\Gamma_{j+1}$ corresponding to $\|\Gamma_{j+1}\|.$ \end{enumerate} 
\end{lemma}

\begin{proof}
		{\rm (i)} By Theorem \ref{Tkismultipleofhankel}, the operator $\Gamma_{j+1}\colon \mathcal{K}_{j+1}\to \mathcal{L}_{j+1}$ is compact and 
$$\|\Gamma_{j+1}\| =\|T_{j+1}\|=t_{j+1}.$$
Hence there exist $x_{j+1} \in \mathcal{K}_{j+1},$ $y_{j+1} \in \mathcal{L}_{j+1}$ such that $(x_{j+1},y_{j+1})$ is a Schmidt pair for $\Gamma_{j+1}$ corresponding to $\|\Gamma_{j+1}\|=t_{j+1}.$ 

{\rm (ii)} Suppose that $x_{j+1}\in \mathcal{K}_{j+1},y_{j+1}\in \mathcal{L}_{j+1}$ satisfy \begin{equation}\label{xitelj}
		\xi_0 \telwe \cdots \telwe \xi_j\telwe  x_{j+1}= \xi_0 \telwe \cdots  \telwe \xi_j \telwe v_{j+1},\end{equation}	
		\begin{equation}\label{eta0telj} \bar{ \eta}_0 \telwe \cdots \telwe \bita_j \telwe  y_{j+1} = \bar{ \eta}_0 \telwe  \cdots\telwe \bita_j\telwe w_{j+1}. \end{equation} Let us show that $(x_{j+1},y_{j+1})$ is a Schmidt pair for $\Gamma_{j+1},$ that is, 
		$$\Gamma_{j+1} x_{j+1} = t_{j+1}y_{j+1},\quad \Gamma_{j+1}^*y_{j+1}=t_{j+1}x_{j+1}. $$ Since diagram \eqref{diagr1122} commutes,
		\begin{equation}\begin{aligned}\label{commt2gammaj}
&T_{j+1} \circ (\xi_0\telwe \cdots \telwe \xi_j \telwe \cdot  )=(\bar{ \eta}_0 \telwe \cdots \telwe  \bita_j \telwe \cdot)\circ\Gamma_{j+1} \;\\
\text{and} \\
 &(\xi_0\telwe \cdots \telwe \xi_j\telwe \cdot  )^*\circ T_{j+1}^* = \Gamma_{j+1}^* \circ (\bar{\eta}_0 \telwe\cdots \telwe \bita_j \telwe  \cdot)^*. \end{aligned}\end{equation} 
By hypothesis,
		\begin{equation}\begin{aligned}\label{hyptj} 
&T_{j+1} (\xi_0 \telwe  \cdots \telwe \xi_{j}\telwe v_{j+1})= t_{j+1} (\bar{ \eta}_0 \telwe \cdots \telwe \bita_j\telwe  w_{j+1}) \;\\
\text{and}\\ 
& T_{j+1}^*(\bar{ \eta}_0 \telwe \cdots \telwe \bita_j\telwe  w_{j+1})= t_{j+1} (\xi_0 \telwe\cdots \telwe \xi_{j}\telwe   v_{j+1}).\end{aligned} \end{equation}
Thus, by equations \eqref{xitelj}, \eqref{eta0telj} and \eqref{hyptj},
		$$\begin{array}{clllll} \Gamma_{j+1}  x_{j+1}&= (\bar{\eta}_0 \telwe \cdots \telwe \bita_j \telwe  \cdot)^*  T_{j+1}  (\xi_0\telwe \cdots \telwe \xi_j \telwe  v_{j+1}) \vspace{2ex} \\
		&=  (\bar{\eta}_0 \telwe \cdots \telwe \bita_j \telwe \cdot)^*  t_{j+1} (\bar{ \eta}_0 \telwe \cdots \telwe \bita_j \telwe w_{j+1}) \vspace{2ex} \\ &= t_{j+1} (\bar{\eta}_0 \telwe \cdots \telwe \bita_j \telwe \cdot)^* (\bar{ \eta}_0 \telwe \cdots \telwe \bita_j \telwe y_{j+1}).\end{array}$$  
Hence  
$$\Gamma_{j+1} x_{j+1}= t_{j+1} (\bar{\eta}_0 \telwe \cdots \telwe \bita_j \telwe \cdot)^* (\bar{ \eta}_0 \telwe\cdots \telwe \bita_j \telwe \cdot)y_{j+1}= t_{j+1} y_{j+1}.$$ 
		
		\noindent By equation \eqref{xitelj}, 
$$x_{j+1} = (\xi_0\telwe  \cdots \telwe \xi_{j}\telwe \cdot)^* (\xi_0 \telwe  \cdots \telwe \xi_{j}\telwe v_{j+1}  ),$$ 
and, by equation \eqref{eta0telj}, 
		$$(\bar{\eta}_0 \telwe\cdots \telwe \bita_j \telwe \cdot)^*(\bar{ \eta}_0 \telwe\cdots \telwe \bita_j \telwe w_{j+1})=y_{j+1}. $$
Thus $$\begin{array}{clll}\Gamma_{j+1}^* y_{j+1}	 &= \Gamma_{j+1}^*(\bar{ \eta}_0 \telwe\cdots \telwe \bita_j \telwe \cdot)^*(\bar{\eta}_0 \telwe\cdots \telwe \bita_j \telwe w_{j+1})\vspace{2ex} \\ 
		&= (\xi_0 \telwe \cdots  \telwe \xi_j \cdot )^*  T_{j+1}^* (\bar{ \eta}_0 \telwe\cdots \telwe \bita_j \telwe w_{j+1}),\end{array}$$ 
the last equality following by the second equation of \eqref{commt2gammaj}. By equations \eqref{xitelj} and \eqref{hyptj}, we have 
		$$ T_{j+1}^* (\bar{ \eta}_0 \telwe\cdots \telwe \bita_j \telwe w_{j+1}) = t_{j+1} (\xi_0\telwe \cdots \telwe \xi_j\telwe v_{j+1})= t_{j+1}(\xi_0\telwe \cdots \telwe \xi_j\telwe   x_{j+1}),$$ and so, $$ \Gamma_{j+1}^* y_{j+1} = t_{j+1} x_{j+1}.$$Therefore $(x_{j+1},y_{j+1})$ is a Schmidt pair for $\Gamma_{j+1}$ corresponding to $\|\Gamma_{j+1}\|=t_{j+1}.$ \end{proof}

\begin{lemma}\label{schfohfj}
	Suppose that $$(\xi_0 \telwe\cdots \telwe \xi_{j}\telwe v_{j+1}, \bar{\eta}_0 \telwe \cdots \telwe \bita_j\telwe  w_{j+1})$$ is a Schmidt pair for the operator $T_{j+1}$ corresponding to $\|T_{j+1}\|=t_{j+1}.$ Let 
	$$x_{j+1} = (I_{n} - \xi_0 \xi_0^*-\cdots-\xi_j\xi_j^*)v_{j+1},$$ 
$$ y_{j+1}= (I_{m} - \bita_0\eta_0^T- \cdots-  \bita_j\eta_j^T)w_{j+1},$$ and let 
	$$\hx_{j+1} = A_{j}^T x_{j+1},\quad \hy_{j+1}=B_j^*y_{j+1}. $$
Then 

{\rm (i)} the pair $(x_{j+1},y_{j+1})$ is a Schmidt pair for the operator $\Gamma_{j+1}$ corresponding to $t_{j+1}$;

{\rm (ii)} the pair $(\hx_{j+1},\hy_{j+1})$ is a Schmidt pair for $H_{F_{j+1}}$ corresponding to $\|H_{F_{j+1}}\|=t_{j+1}.$
	
\end{lemma}

\begin{proof}
By Lemmas \ref{xi=AkAkT}  and \ref{eta=BkBk*},
\begin{equation}\label{xjAA}
x_{j+1} = (I_{n} - \xi_0 \xi_0^*-\cdots-\xi_j\xi_j^*)v_{j+1}=\bar{A}_{j}A_{j}^Tv_{j+1}
\end{equation}
and 
\begin{equation}\label{yjBB}
 y_{j+1}= (I_{m} - \bita_0\eta_0^T-\cdots -\bita_j\eta_j^T)w_{j+1}=B_{j}B_{j}^*w_{j+1} .
\end{equation}
Hence
\begin{equation}\label{hxj+1}
\hx_{j+1} = A_{j}^T x_{j+1}= A_{j}^T \bar{A}_{j} A_{j}^Tv_{j+1}=  A_{j}^Tv_{j+1}
\end{equation}
and
\begin{equation}\label{bjwj}
\hy_{j+1}=B_j^*y_{j+1}= B_j^*B_{j}B_{j}^*w_{j+1} = B_{j}^*w_{j+1}.
\end{equation}
These  imply that $\hat{x}_{j+1}\in H^2(\Dd,\Cc^{n-j-1})$, $x_{j+1}\in \mathcal{K}_{j+1}$, $\hy_{j+1}\in H^2(\Dd,\Cc^{n-j-1})^\perp$ and $y_{j+1} \in \mathcal{L}_{j+1}.$
By Proposition \ref{onxi},
$$
		\xi_0 \telwe \cdots \telwe \xi_j\telwe  x_{j+1}= \xi_0 \telwe \cdots  \telwe \xi_j \telwe v_{j+1}\quad \text{and} \;	
		\bar{ \eta}_0 \telwe \cdots \telwe \bita_j \telwe  y_{j+1} = \bar{ \eta}_0 \telwe  \cdots\telwe \bita_j\telwe w_{j+1}.$$ 
Thus, by Proposition \ref{coofschpairsj},	
the pair $(x_{j+1},y_{j+1})$ is a Schmidt pair for $\Gamma_{j+1}$ corresponding to $\|\Gamma_{j+1}\|.$ Therefore,
\begin{equation}\label{gammaJ+1}
\Gamma_{j+1}x_{j+1}= t_{j+1}y_{j+1},\quad \Gamma_{j+1}^*y_{j+1}=t_{j+1}x_{j+1}. \end{equation}
To show that the pair $(\hx_{j+1},\hy_{j+1})$ is a Schmidt pair for $H_{F_{j+1}}$ corresponding to $t_{j+1}$, we need to prove that 
	$$H_{F_{j+1}}\hx_{j+1} = t_{j+1} \hy_{j+1},\quad H_{F_{j+1}}^*\hy_{j+1} = t_{j+1}\hx_{j+1}. $$ 

By Theorem \ref{Tkismultipleofhankel}, 
\begin{equation}\label{hf2gj} 
H_{F_{j+1}} = (M_{B_j})^* \circ \Gamma_{j+1} \circ M_{\bar{A}_j},\end{equation} 
and
\begin{equation}\label{hf2*gj} 
H_{F_{j+1}}^* = M_{\bar{A}_j}^* \circ \Gamma_{j+1}^* \circ M_{B_j}.\end{equation} 
By equation \eqref{hf2gj}, we have
	\begin{align} \label{h_fj.1}
	H_{F_{j+1}}\hat{x}_{j+1}&= H_{F_{j+1}}A_{j}^T x_{j+1}\nonumber\vspace{2ex}\\
	&= B_j^* \Gamma_{j+1} \bar{A}_{j}A_{j}^Tx_{j+1}. \end{align} 	
Notice that, by equations \eqref{xjAA} and \eqref{hxj+1}, 
\begin{equation}\label{x_j+1}
x_{j+1}=\bar{A}_jA_j^Tx_{j+1}.\end{equation} 
Hence, by equations \eqref{gammaJ+1} and \eqref{h_fj.1}, we obtain	
	$$H_{F_{j+1}}\hx_{j+1} = B_{j}^* \Gamma_{j+1} x_{j+1}=t_{j+1}B_j^*y_{j+1}=t_{j+1}\hy_{j+1}.$$
	
	Let us show that $H_{F_{j+1}}^*\hy_{j+1} = t_{j+1}\hx_{j+1}$. By equations \eqref{hf2*gj} and  \eqref{gammaJ+1}, we have
	\begin{align}
	H_{F_{j+1}}^*\hy_{j+1}&= H_{F_{j+1}}^* B_{j}^* y_{j+1}\nonumber\vspace{2ex}\\ &= A_{j}^T\Gamma_{j+1}^* B_j B_j^* y_{j+1} \label{h_fj.4}\end{align} 
Observe that, in view of equations \eqref{yjBB} and \eqref{bjwj}, 
 \begin{equation}\label{y_j+1}
	y_{j+1}=B_jB_j^*y_{j+1}.\end{equation} 
Hence, by equations \eqref{gammaJ+1} and \eqref{h_fj.4}, we obtain 
		 	$$H_{F_{j+1}}^*\hy_{j+1}=  A_{j}^T\Gamma_{j+1}^* y_{j+1} = t_{j+1} A_{j}^Tx_{j+1} = t_{j+1} \hx_{j+1}.$$Therefore $(\hx_{j+1},\hy_{j+1})$ is a Schmidt pair for the Hankel operator $H_{F_{j+1}}$ corresponding to $\| H_{F_{j+1}}\| = t_{j+1}.$ \end{proof}

\begin{proposition}\label{xjwevjetajwewj}
	Let $$(\xi_0\telwe \cdots \telwe \xi_j  \telwe v_{j+1},\bita_0\telwe \cdots \telwe \bita_j \telwe  w_{j+1})$$ be a Schmidt pair for $T_{j+1}$ corresponding to $t_{j+1}$ for some $v_{j+1}\in H^2(\Dd,\Cc^n),$ $w_{j+1}\in H^2(\Dd,\Cc^m)^\perp.$ Let $$x_{j+1} = (I_{n}- \xi_0 \xi_0^*-\cdots-\xi_j\xi_j^*)v_{j+1},\quad y_{j+1}=(I_{m} - \bar{\eta}_0 \eta_0^T-\cdots-\bita_j \eta_j^T)w_{j+1},$$
and let 
\begin{equation}\label{hatxj}\hx_{j+1}=A_j^T x_{j+1}\quad \text{and}\quad \hy_{j+1}=B_j^*y_{j+1}.\end{equation} 
Then
	\begin{equation}\label{hj+1common}
\begin{array}{lll}
\| \xi_0 (z) \we \dots \we \xi_{j}(z) \we v_{j+1}(z) \|_{\we^{j+2}\Cc^n} &=
 \| \bar{\eta}_0(z) \we \dots \we\bita_{j}(z) \we w_{j+1}(z)\|_{\we^{j+2}\Cc^m}= |h_{j+i}(z)|, \\
\|\hx_{j+1} (z) \|_{\Cc^{n-j-1}} &=  \|\hy_{j+1}(z)\|_{\Cc^{m-j-1}} = |h_{j+1}(z)|,\;\text{and}\\
\| x_{j+1}(z) \|_{\Cc^n} &=  \|y_{j+1}(z)\|_{\Cc^m} =|h_{j+1}(z)|,  \;
\end{array}
\end{equation} 
almost everywhere on $\Tt.$
\end{proposition}

\begin{proof}

By Lemma \ref{schfohfj}, $(\hx_{j+1},\hy_{j+1})$ is a Schmidt pair for $H_{F_{j+1}}$ corresponding to $\|H_{F_{j+1}}\|=t_{j+1}$. Hence
	$$H_{F_{j+1}}\hx_{j+1} = t_{j+1}\hy_{j+1} \quad \text{and}\quad H_{F_{j+1}}^* \hy_{j+1} = t_{j+1} \hx_{j+1}. $$
By Theorem \ref{1.7}, 
	$$t_{j+1} \|\hy_{j+1}(z)\|_{\Cc^{m-j-1}}=\|H_{F_{j+1}}\| \|\hx_{j+1}(z)\|_{\Cc^{n-j-1}} $$almost everywhere on $\Tt.$ 
Thus 
\begin{equation}\label{hatseqj} \|\hy_{j+1}(z)\|_{\Cc^{m-j-1}}= \|\hx_{j+1}(z)\|_{\Cc^{n-j-1}} \end{equation} almost everywhere on $\Tt.$ 
	
	Notice that $\bar{A}_j(z)$ is isometric for almost every $z\in \Tt,$ and therefore, by equations \eqref{hatxj}, we obtain
	$$\|x_{j+1}(z)\|_{\Cc^{n}}=\|\hx_{j+1}(z)\|_{\Cc^{n-j-1}}. $$
	
	Moreover, since $B_j(z)$ are isometries almost everywhere on $\Tt,$ by equations \eqref{hatxj}, we have  
	$$\|y_{j+1}(z)\|_{\Cc^m}=\|\hy_{j+1}(z)\|_{\Cc^{m-j-1}} $$
almost everywhere on $\Tt.$ By equations \eqref{hatseqj}, we deduce
	\begin{equation}\label{x2isyj}\|x_{j+1}(z)\|_{\Cc^n}=\|y_{j+1}(z)\|_{\Cc^m} \end{equation} almost everywhere on $\Tt.$

By Proposition \ref{onxi},
\begin{equation} \label{xitelvj}	
 \xi_0 \telwe \cdots \telwe \xi_j \telwe x_{j+1}=\xi_0\telwe \cdots \telwe \xi_j  \telwe v_{j+1}
\end{equation}
and
\begin{equation} \label{etatelwj}
\bar{\eta}_0 \telwe \cdots \telwe \bita_j \telwe y_{j+1} =
\bar{\eta}_0 \telwe \cdots \telwe \bita_j \telwe w_{j+1}.
\end{equation}
Hence, by Proposition  \ref{weon},
$$\begin{array}{lll}
&\|\xi_0(z)\we \cdots \we \xi_j(z)\we v_{j+1}(z)\|_{\we^{j+2}\Cc^n}\\
&= \|\xi_0(z)\we \cdots \we \xi_j(z)\we x_{j+1}(z)\|_{\we^{j+2}\Cc^n}, \\
& = \| x_{j+1}(z) - \displaystyle\sum\limits_{i=0}^j \langle x_{j+1}(z), \xi_i(z)\rangle \xi_i(z)\|_{\Cc^n}=\|x_{j+1}(z)\|_{\Cc^n},
\end{array}$$
almost everywhere on $\Tt$.  Furthermore 
$$\begin{array}{llll}
&\| \bita_0(z) \we \cdots \we \bita_j(z) \we w_{j+1}(z)\|_{\we^{j+2}\Cc^m}\\
&=\| \bita_0(z) \we \cdots \we \bita_j(z) \we y_{j+1}(z)\|_{\we^{j+2}\Cc^m}\\
&=\| y_{j+1}(z) - \displaystyle\sum\limits_{i=0}^j \langle y_{j+1}(z), \bita_i(z)\rangle \bita_i(z)\|_{\Cc^m} =\|y_{j+1}(z)\|_{\Cc^m}, 	
\end{array}$$
almost everywhere on $\Tt$.

Thus, by equation \eqref{x2isyj},	
$$\| \bita_0(z) \we \cdots \we \bita_j(z) \we w_{j+1}(z)\|_{\we^{j+2}\Cc^m} = \|\xi_0(z)\we \cdots \we \xi_j(z)\we v_{j+1}(z)\|_{\we^{j+2}\Cc^n}$$
almost everywhere on $\Tt.$
	
Recall that $h_{j+1}$ is the scalar outer factor  of $\xi_0\we \cdots \we \xi_j\we v_{j+1}$. Hence
$$\|\hx_{j+1} (z) \|_{\Cc^{n-j-1}} =  \|\hy_{j+1}(z)\|_{\Cc^{m-j-1}} = |h_{j+1}(z)|,$$
$$\| x_{j+1}(z) \|_{\Cc^n} =  \|y_{j+1}(z)\|_{\Cc^m} =|h_{j+1}(z)|, $$
and	
$$\|\xi_0(z)\we \cdots \we \xi_j(z)\we v_{j+1}(z)\|_{\we^{j+2}\Cc^n}
 = \| \bita_0(z) \we \cdots \we \bita_j(z) \we w_{j+1}(z)\|_{\we^{j+2}\Cc^m}\|
=|h_{j+1}(z)|, $$
almost everywhere on $\Tt.$	\end{proof}

\begin{proposition}\label{epsilonjj}
	In the notation of Theorem \ref{Tkismultipleofhankel}, there exist unitary-valued functions $\tilde{V}_{j+1}, \tilde{W}_{j+1}$ of types $(n-j-1)\times (n-j-2)$ and $(m-j-1) \times (m-j-2)$ respectively of the form 
	$$\tilde{V}_{j+1} =\begin{pmatrix}
	A_{j} \xi_{j+1}  &  \bar{\alpha}_{j+1}
	\end{pmatrix}, \quad \tilde{W}_{j+1}^T = \begin{pmatrix}
B_{j}^T \eta_{j+1} & \bar{\beta}_{j+1}
	
\end{pmatrix} ,$$ 
where $\alpha_{j+1}, \beta_{j+1}$ are inner, co-outer, quasi-continuous and all minors on the first columns of $\tilde{V}_{j+1}, \tilde{W}_{j+1}^T$ are in $H^\infty.$ Furthermore, the set $\mathcal{E}_{j+1}$ of all level $j+1$ superoptimal error functions for $G$ is equal to the following set
$$\scriptsize{ W_0^*\cdots \begin{pmatrix}
	I_{j+1} & 0 \\ 
	0& \tilde{W}_{j+1}^*
	
\end{pmatrix} \begin{pmatrix}
	t_0 u_0 & 0 & 0&0 \\
	0 & t_1 u_1 &0 &0 \\
		\vdots  &\hspace{8ex} \ddots&~&\vdots \\
	0& 0 & t_{j+1} u_{j+1}&0 \\
	0&0 & 0 & (F_{j+2} +H^\infty)\cap B(t_{j+1})
	
\end{pmatrix}\begin{pmatrix}
	I_{j+1} & 0 \\ 0 & \tilde{V}_{j+1}^*
\end{pmatrix}\cdots V_0^* , }$$
where $F_{j+2} \in H^\infty(\Dd,\Cc^{(m-j-2)\times (n-j-2)}) +C(\Tt,\Cc^{(m-j-2)\times (n-j-2)}),$ $u_{j+1}=\frac{\bar{z}\bar{h}_{j+1}}{h_{j+1}}$ is a quasi-continuous unimodular function and $B(t_{j+1})$ is the closed ball of radius $t_{j+1}$ in $L^\infty(\Tt,\Cc^{(m-j-2)\times (n-j-2)}).$

\end{proposition}
\begin{proof}
Recall that, in diagrams \eqref{diagr1122} and \eqref{hfj+1}, the operators 
$M_{\bar{A}_j},\;M_{B_j},$ $(\xi_0\telwe \cdots \telwe \xi_j \cdot)$ and $(\bita_0 \telwe \cdots\telwe  \bita_{j} \telwe \cdot)$ are unitaries. Since both diagrams commute and $(x_{j+1},y_{j+1})$ defined above is a Schmidt pair for $\Gamma_{j+1}$ corresponding to $t_{j+1},$ by Lemma \ref{schfohfj}, $(\hx_{j+1},\hy_{j+1})$ is a Schmidt pair for $H_{F_{j+1}}$ corresponding to $t_{j+1},$ where 
$$\hx_{j+1} = A_{j} x_{j+1},\quad \hy_{j+1}=B_{j}^*y_{j+1}. $$	
We intend to apply Lemma \ref{2.2} to $H_{F_{j+1}}$ and the Schmidt pair $(\hx_{j+1},\hy_{j+1})$ to find unitary-valued functions $\tilde{V}_{j+1},\tilde{W}_{j+1}$ such that, for every $\tilde{Q}_{j+1}\in H^\infty(\Dd,\Cc^{(m-j-1)\times(n-j-1)})$ which is at minimal distance from $F_{j+1},$ we obtain a  factorisation of the form
$$
F_{j+1}-\tilde{Q}_{j+1} = \tilde{W}_{j+1}^* \begin{pmatrix}t_{j+1} u_{j+1} &0\\ 0 & F_{j+2}
\end{pmatrix}\tilde{V}_{j+1}^*, 
$$
for some $F_{j+2} \in H^\infty(\Dd,\Cc^{(m-j-2)\times(n-j-2)})+C(\Tt,\Cc^{(m-j-2)\times(n-j-2)}).$ For this purpose we find the inner-outer factorisations of $\hx_{j+1}$ and $\bar{z}\bar{\hy}_{j+1}.$

By Proposition \ref{xjwevjetajwewj},
\begin{equation}\label{hj+1common1}
\|\hx_{j+1}(z)\|_{\Cc^{n-j-1}}= |h_{j+1}(z)|
\end{equation}
and
\begin{equation}\label{hj+1common2}
\|\hy_{j+1}(z)\|_{\Cc^{m-j-1}}=|h_{j+1}(z)|,
\end{equation} almost everywhere on $\Tt.$ 
Equations \eqref{hj+1common1} and \eqref{hj+1common2} imply that $h_{j+1}\in H^2(\Dd,\Cc)$ is the scalar outer factor of both $\hat{x}_{j+1}$ and $\bar{z}\bar{\hat{y}}_{j+1}.$ 
Hence, by Lemma \ref{2.2}, 
$\hat{x}_{j+1}, \bar{z}\bar{\hat{y}}_{j+1}$ admit the inner outer factorisations 
$$\hat{x}_{j+1} = \hat{\xi}_{j+1} h_{j+1}, \quad \bar{z} \bar{y}_{j+1} = \hat{\eta}_{j+1} h_{j+1},$$for some inner $\hat{\xi}_{j+1} \in H^\infty(\Dd,\Cc^{n-j-1}), \hat{\eta}_{j+1}\in  H^\infty(\Dd,\Cc^{m-j-1}).$ Then 

$$ \hat{x}_{j+1} = \hat{\xi}_{j+1} h_{j+1} =A_j^T x_{j+1},\quad  \bar{z}\bar{\hat{y}}_{j+1} =\hat{\eta}_{j+1} h_{j+1}=\bar{z}B_j^T \bar{y}_{j+1},$$from which we obtain 
$$ \hat{\xi}_{j+1} = A_j^T \xi_{j+1},\quad  \hat{\eta}_{j+1} =B_j^T \eta_{j+1}. $$

 We wish to show that $A_j^T \xi_{j+1},\; B_j^T \eta_{j+1}$ are inner functions in order to apply Lemma \ref{2.2}.  Observe that, by equations \eqref{xjAA} and \eqref{yjBB}, 
$$x_{j+1}=\bar{A}_jA_j^Tv_{j+1},\quad y_{j+1}= B_j B_j^*w_{j+1} .$$ Then,
$$A_j^T x_{j+1}= A_j^Tv_{j+1}, \quad B_j^T\bar{y}_{j+1} =B_j^T\bar{w}_{j+1}  ,$$ and since
$$\xi_{j+1} =\frac{x_{j+1}}{h_{j+1}},\quad \eta_{j+1} = \frac{\bar{z}\bar{y}_{j+1}}{h_{j+1}}, $$ 
the functions
$$A_j^T \xi_{j+1}=\frac{A_j^Tv_{j+1}}{h_{j+1}},\quad B_j^T \eta_{j+1}= \frac{B_j^*w_{j+1}}{h_2}$$ are analytic. Furthermore, by Proposition \ref{onxi} $\|\xi_{j+1}(z)\|_{\Cc^n}=1$ and $\|\eta_{j+1}(z)\|_{\Cc^m}=1$ almost everywhere on $\Tt,$ and, by equations \eqref{hj+1common}, 
 
$$\|A_j^T(z) \xi_{j+1}(z)\|_{\Cc^{n-j-1}}=1,\quad \|B_j^T(z) \eta_{j+1}(z)\|_{\Cc^{m-j-1}}=1 $$almost everywhere on $\Tt.$ Thus $A_j^T \xi_{j+1},\; B_j^T \eta_{j+1}$ are inner functions.

By Lemma \ref{2.2}, there exist inner, co-outer, quasi-continuous functions $\alpha_{j+1}, \beta_{j+1}$  of types $(n-j-1)\times (n-j-2), (m-j-1)\times (m-j-2)$ respectively such that the functions 
$$\tilde{V}_{j+1} =\begin{pmatrix}
A_j^T \xi_{j+1}  &  \bar{\alpha}_{j+1}
\end{pmatrix}, \quad \tilde{W}_{j+1}^T = \begin{pmatrix}
B_j^T \eta_{j+1} & \bar{\beta}_{j+1}

\end{pmatrix} $$are unitary-valued with all minors on the first columns in $H^\infty.$   Furthermore, by Lemma \ref{2.2}, every $\hat{Q}_{j+1}\in H^\infty(\Dd,\Cc^{(m-j-1)\times(n-j-1)})$ which is at minimal distance from $F_{j+1}$ satisfies 
$$F_{j+1}-\hat{Q}_{j+1} = \tilde{W}_{j+1}^* \begin{pmatrix}
t_{j+1} u_{j+1} & 0 \\
0 & F_{j+2}
\end{pmatrix}\tilde{V}_{j+1}^*, $$
for some $F_{j+2} \in   H^\infty(\Dd,\Cc^{(m-j-2)\times (n-j-2)})+C(\Tt,\Cc^{(m-j-2)\times (n-j-2)})$, where $u_{j+1}$ is a quasi-continuous unimodular function given by $u_{j+1} = \frac{\bar{z} \bar{h}_{j+1}}{h_{j+1}}.$ 

By Lemma \ref{f+hinfty}, the set $$\tilde{\mathcal{E}}_{j+1} =\{F_{j+1} - \hat{Q} : \hat{Q} \in H^\infty(\Dd,\Cc^{(m-j-1)\times (n-j-1)}), \| F_{j+1} - \hat{Q}\|_{L^\infty}=t_{j+1}  \}$$ satisfies  $$\tilde{\mathcal{E}}_{j+1} = \tilde{W}_{j+1}^* 
\begin{pmatrix}
t_{_{j+1}}u_{_{j+1}} & 0 \\
0 & (F_{j+2}+H^\infty)\cap B(t_{j+1})
\end{pmatrix}V_{j+1}^*,
$$
where $B(t_{j+1})$ is the closed ball of radius $t_{j+1}$ in $L^\infty(\Tt,\Cc^{(m-j-2)\times (n-j-2)}).$ Thus, by Proposition \ref{tildvjwj}, $\mathcal{E}_{j+1}$ admits the factorisation claimed.\end{proof}

\begin{theorem}\label{mathcalAG}

Let $G\in H^\infty(\Dd, \Cc^{m\times n})+C(\Tt, \Cc^{m\times n}).$ Let $T_i, x_i, y_i, h_i$, for $i\ge 0$, be defined by the algorithm from
Subsection \ref{Alg_statement}. Let $r$ be the least index $j \ge 0$ such that $T_j=0$. Then $r\leq \min(m,n)$ and the superoptimal approximant $\mathcal{A}G$ is given by the formula
	 	$$ G-\mathcal{A}G= \displaystyle \sum\limits_0^{r-1} \frac{t_i y_i x_i^*}{|h_i|^2} .$$
\end{theorem}
\begin{proof}
	First observe that, if $T_0=H_G=0,$ then this implies $G\in H^\infty(\Dd,\CCmn),$ and so $$\mathcal{A}G=G.$$ 
Otherwise, let $t_0=\|H_G\|>0.$ 
If  $T_1=0$, by Theorem \ref{T0compact},
$H_{F_{1}} =0,$ that is, $$F_1 \in H^\infty(\Dd, \Cc^{(m-1)\times (n-1)}).$$ Then, by Lemma \ref{f+hinfty}, we have 
	
	$$W_0 (G-\mathcal{A}G) V_0 = \begin{pmatrix}
	t_0 u_0 & 0 \\
	0 & 0
	\end{pmatrix}.
	$$ 
	
	\noindent Equivalently
	$$
	\begin{array}{cllll}
	G-\mathcal{A}G  &=W_0^* \begin{pmatrix}
	t_0 u_0 & 0 \\
	0 & 0
	\end{pmatrix}V_0^* \vspace{2ex} \\
	&= \begin{pmatrix}
	\bar{\eta}_0& \beta_0 
	\end{pmatrix} \begin{pmatrix}  t_0 u_0 & 0 \\
	0 & 0 \end{pmatrix} \begin{pmatrix} \xi_0^*\\\alpha_0^T \end{pmatrix}\vspace{2ex} \\\end{array}$$ 
	$$\begin{array}{lll}
	&= \begin{pmatrix}
	\bar{\eta}_0 t_0 u_0 & 0 
	\end{pmatrix}\begin{pmatrix} \xi_0^*\\\alpha_0^T \end{pmatrix} \vspace{2ex} = \bar{\eta}_0t_0 u_0 \xi_0^*\vspace{2ex} \\
	&= t_0 \displaystyle\frac{zy_0}{\bar{h}_0} \frac{\bar{z} \bar{h_0}}{h_0} \frac{x_0^*}{\bar{h}_0}\vspace{2ex} 
             = \displaystyle\frac{t_0 y_0 x_0^*}{|h_0|^2}.

\end{array}$$
Let $j$ be a non-negative integer such that $T_j=0$ and $T_i \neq 0$ for $ 1\le i < j$. By the commutativity of the diagrams \eqref{diagr1122} and \eqref{hfj+1}, $H_{F_j} =0$, and therefore $F_j \in H^\infty(\Dd,\Cc^{(m-j)\times (n-j)})$.
By Proposition \ref{epsilonjj}, the superoptimal analytic approximant $\mathcal{A}G$ satisfies equation \eqref{epsilonnj}, that is,
\begin{equation}\label{AG}
G-\mathcal{A}G= W_0^* W_1^* \cdots W_{j-1}^* 
\begin{pmatrix}
t_0 u_0 & 0 &\dots &0 \\
0&t_1u_1 &\dots &0 \\
\vdots & \hspace{10ex}\ddots &  &\vdots \\
0  &\cdots &t_{j-1}u_{j-1} & 0\\ 
0 & \dots & \dots & 0
\end{pmatrix}
V_{j-1}^* \cdots V_1^* V_0^*,
\end{equation}
where, for $i=0,1,\dots, j-1$,
$$ \tilde{V}_{i} = \begin{pmatrix}
\alpha_{i-1}^T \cdots \alpha_0^T \xi_{i} & \bar{ \alpha}_{i} 
\end{pmatrix}, \quad \tilde{W}_{i}^T = \begin{pmatrix}
\beta_{i-1}^T  \cdots \beta_0^T \eta_{i} & \bar{\beta}_{i}
\end{pmatrix} $$ are unitary-valued functions, as described in Proposition \ref{tildevjwj}, 
  $u_i = \frac{\bar{z} \bar{h}_i}{h_i}$ are quasi-continuous unimodular functions, and $$V_{i} = \begin{pmatrix}
I_{i} & 0 \\
0 & \tilde{V}_{i} 
\end{pmatrix}, \quad W_{i} = \begin{pmatrix}
I_{i} & 0 \\
0 & \tilde{W}_{i}
\end{pmatrix} .$$
Recall that, by equations \eqref{xij+1etaj+1}, for $i =0, \dots, j-1$,
\begin{equation}
\xi_{i} = \frac{x_{i}}{h_{i}}, \quad \eta_{i}=\frac{\bar{z}\bar{y}_{i}}{h_{i}}.\end{equation}
By Proposition \ref{xjwevjetajwewj}, for $i =0, \dots, j-1$,
$$ |h_i(z)| = \|x_i(z)\|_{\Cc^n} = \|y_i(z)\|_{\Cc^m}$$
almost everywhere on $\Tt$.

With the aid of the formulae \eqref{WjW0*} and \eqref{V0Vj}, equation \eqref{AG}
simplifies to 
\begin{align}\label{G-AG-main}
G-\mathcal{A}G =& \displaystyle\frac{t_0 y_0 x_0^*}{|h_0|^2} + 
t_1\frac{1}{|h_1|^2}B_0 B_0^*y_1 x_1^* \bar{A}_0 A_0^T+\dots \nonumber\\
       & \; + t_{j-1} \frac{1}{|h_{j-1}|^2}B_{j-1} B_{j-1}^*y_{j-1} x_{j-1}^* \bar{A}_{j-1} A_{j-1}^T.
\end{align}
By equations \eqref{x_j+1} and \eqref{y_j+1}, for $i =0, \dots, j-1$,
$$ x_i^* = x_i^* \bar{A}_{i-1} A_{i-1}^T  \; \; \text{and} \;\; 
y_i = B_{i-1} B_{i-1}^* y_i.$$
Thus 
$$ G -\mathcal{A}G = \sum\limits_{i=0}^{r-1}\frac{t_i y_i x_i^*} {|h_i|^2} $$ and the assertion has been proved. \end{proof}


%% file: application.tex
\section{Application of the algorithm}\label{applic}

Let us now apply the new algorithm to the example Peller and Young solved in \cite{Constr}. 
\begin{problem}\label{applic-ex}
	Let $G=B^{-1}A \in L^\infty(\mathbb{C}^{2\times 2})$ where 
	
	$$A(z)= \begin{pmatrix}
	
	\sqrt{3}+2z \hspace{2ex}    &0\\
	0         \hspace{2ex}       &1
	
	\end{pmatrix},    \;\;                          B(z) =  \displaystyle\frac{1}{\sqrt{2}}  \begin{pmatrix}
	z^2 \hspace{2ex} z\\
	z  \hspace{2ex} -1
	
	\end{pmatrix} \qquad \mbox{ for all  } z\in\T.
$$
	Find the superoptimal singular values of $G$ and its superoptimal approximant $\mathcal{A}G \in H^{\infty},$ that is, the unique $\mathcal{A}G \in H^\infty(\Dd, \mathbb{C}^{2\times 2})$ such that the sequence 
	$$s^\infty(G-\mathcal{A}G) = (s_0^\infty ( G-\mathcal{A}G) , s_1^\infty(G-\mathcal{A}G), \cdots ) $$ is lexicographically minimized.
\end{problem}
On $\mathbb{T}$ $G$ is 
$$
G(z) = \displaystyle\frac{1}{\sqrt{2}}\begin{pmatrix}
\sqrt{3}\bar{z}^2 +2\bar{z} \hspace{2ex} &\bar{z}\\
\sqrt{3}\bar{z}+2   \hspace{2ex}   &-1
\end{pmatrix}.
$$
\textbf{Step 0:}
 The operator $H_G^* H_G$ with the respect to the orthonormal basis 
$$ B=\left\{ \begin{pmatrix}
1\\
0
\end{pmatrix},\begin{pmatrix}
z\\
0
\end{pmatrix},  \begin{pmatrix}
0\\
1
\end{pmatrix},\begin{pmatrix}
0\\
z
\end{pmatrix} \right\} $$of $(z^2H^2)^\perp,$
has matrix representation

$$H^*_G H_G \sim \displaystyle\frac{1}{\sqrt{2}} \begin{pmatrix} 
&10 \hspace{2ex} &2\sqrt{3} \hspace{2ex} &2 \hspace{2ex} &0\\
&2\sqrt{3} \hspace{2ex} &3 \hspace{2ex} &\sqrt{3} \hspace{2ex} &0\\
&2 \hspace{2ex} &\sqrt{3} \hspace{2ex}    &1 \hspace{2ex}  &0\\
&0 \hspace{2ex} &0 \hspace{2ex} &0 \hspace{2ex} &0\\
\end{pmatrix}.$$	
		
	 Then	$\|H_G\| = \sqrt{6}$ and a non-zero vector $ x_0 \in H^2(\Dd,\Cc^2)$ such that $$\|H_Gx_0\|_{H^{2}(\Dd,\Cc^2)^\perp} = \|H_G\| \|x_0\|_{H^2(\Dd,\Cc^2)} $$  is 
			$$x_0(z)= \begin{pmatrix}
			4+\sqrt{3}z\\
			1
			\end{pmatrix}.$$  
			
			 For $(x_0,y_0)$ to be a Schmidt pair for $H_G$ corresponding to $\|H_G\|,$ the vector $y_0 \in H^2(\Dd,\Cc^2)^\perp$ can be calculated by
			
			$$y_0 (z)= \frac{H_G x_0(z)}{\|H_G\|} = \displaystyle 2\bar{z} \begin{pmatrix}
			\bar{z} + \sqrt{3}\\
			1
			\end{pmatrix} \in H^2(\Dd,\Cc^2)^\perp.
$$
Perform the inner-outer factorizations $$x_0=\xi_{0}h_0,\quad  \bar{z}\bar{y}_0=\eta_{0}h_0$$ 
for some inner $\xi_{0}, \eta_{0}  \in H^{\infty}(\Dd,\Cc^2)$ and some scalar outer $h_0 \in H^2(\Dd,\Cc).$   

In this example
$$\xi_0(z) = \frac{x_0}{h_0}= \frac{a}{4\sqrt{3}(1-\gamma z)}\begin{pmatrix}
4 +\sqrt{3}z \\1 
\end{pmatrix},$$ $$\bar{\eta}_0(z) =\frac{zy_0}{\bar{h}_0} = \frac{2a}{4\sqrt{3}(1-\gamma \bar{z})}  \begin{pmatrix}
\bar{z} + \sqrt{3}\\
1
\end{pmatrix} ,$$ where $$h_0(z) = \frac{4\sqrt{3}}{a}(1-\gamma z) ,$$ $a= \sqrt{10-2\sqrt{13}}$ and $\gamma= -\frac{a^2}{4\sqrt{3}}.$

  A function $Q_1 \in H^\infty(\Dd,\Cc^{2\times 2})$ that satisfies 
$$(G-Q_1)x_0 = t_0 y_0 ,\quad (G-Q_1)^*y_0 = t_0 x_0 $$ is 

$$ Q_1(z) =  \begin{pmatrix}
0 & \sqrt{6} \\ 2\sqrt{2} & -\sqrt{6}(z +\sqrt{3})
\end{pmatrix}    .$$  
\textbf{Step 1:}
 Let $X_1= \xi_0 \telwe H^2(\Dd,\Cc^2)$ and $Y_1 = \bar{\eta}_0 \telwe H^2(\Dd,\Cc^2)^\perp.$  

 Let $T_1:X_1 \to Y_1$ be given by $$T_1 (\xi_0 \telwe x) = P_{Y_1}(\bar{\eta}_0 \telwe (G-Q_1)x)$$ for all $x \in H^2(\Dd,\Cc^2).$

 Note that
$$\begin{array}{cll}X_1&=\left\{\xi_0 \telwe \begin{pmatrix}
f_1 \\f_2 
\end{pmatrix}\;:\; f_i \in H^2(\Dd,\Cc) \right\} \vspace{2ex} \\ & = \left\{\displaystyle\frac{a}{4\sqrt{3}} \frac{(4+\sqrt{3}z)f_2-f_1}{1-\gamma z}\; : \; f_i \in H^2(\Dd,\Cc)\right\} .\end{array}$$  
If we choose $$ f_1 = -\frac{4\sqrt{3}}{a}(1-\gamma z)g\quad\text{and}\quad f_2 =0$$ for some $g\in H^2(\Dd,\Cc),$ we obtain $X_1 = H^2(\Dd,\Cc).$

Also $$\begin{array}{cll} Y_1 &= \left\{\bar{\eta}_0 \telwe \begin{pmatrix}
\bar{z} \bar{\phi}_1\\
\bar{z} \bar{\phi}_2
\end{pmatrix}\; : \; \phi_i \in H^2(\Dd,\Cc)\right\}\vspace{2ex}\\& = \left\{\displaystyle\frac{a\bar{z}}{2\sqrt{3}} \frac{(\bar{z}+\sqrt{3})\bar{\phi}_2 - \bar{\phi}_1}{1-\gamma \bar{z}}\; : \; \phi_i \in H^2(\Dd,\Cc) \right\}.\end{array}
$$  
If we choose  $$\phi_1 =-\frac{2\sqrt{3}}{a}(1-\gamma z)\psi, \quad \text{and}\quad \phi_2=0 $$ for some $\psi \in H^2(\Dd,\Cc),$ we find that $Y_1 = H^2(\Dd,\Cc)^\perp.$
 
 We have  
$$T_1\left(\xi_0 \telwe \begin{pmatrix} f_1 \\ f_2\end{pmatrix} \right)=T_1\left(\xi_0 \telwe \begin{pmatrix} -\frac{4\sqrt{3}}{a} (1-\gamma z)g \\ 0\end{pmatrix} \right) = \frac{u(\gamma)}{z-\gamma}$$ where $$u(\gamma) = \sqrt{2}(1-\gamma^2)(2\sqrt{3}\gamma +1)g(\gamma).$$
Then $t_1= \|T_1\| = \sqrt{2}(4-\sqrt{13}) .$ 
Since $T_1 $ is a compact operator, there exist $v_1 \in H^2(\Dd,\Cc^2),$ $w_1 \in H^2(\Dd,\Cc^2)^\perp$ such that 
$$T_1(\xi_0\telwe v_1) = t_1 (\bar{\eta}_0 \telwe w_1) , \quad T_1^* (\bar{\eta}_0 \telwe w_1) = t_1 (\xi_0 \telwe v_1). $$   
 
 Here we can choose $$v_1(z) = \frac{4\sqrt{3}}{a}\begin{pmatrix}-1 \\ 0 \end{pmatrix}, \quad w_1(z)= \frac{2\sqrt{3}}{a}\bar{z}\begin{pmatrix}
1 \\ 0
\end{pmatrix}.$$
\noindent  Perform the inner-outer factorisation of $\xi_0 \telwe v_1\in H^2(\Dd,\we^2\Cc^2).$ The function $h_1(z)= \frac{1}{1-\gamma z}$ is the scalar outer factor of $\xi_0 \telwe v_1.$   

\noindent  Let $$x_1 = (I - \xi_0(z)\xi_0^*(z))v_1(z), \quad y_1(z)= (I-\bar{\eta}_0(z) \eta_0^T(z))w_1(z).$$

\noindent  Then $$x_1 = \frac{\gamma}{\alpha} \frac{1}{(1-\gamma z)(1-\gamma \bar{z})} \begin{pmatrix}
\frac{-4\sqrt{3}}{\gamma} (1-\gamma z)(1-\gamma \bar{z})-19-4\sqrt{3}(z+\bar{z})\\
-4-\sqrt{3}\bar{z}
\end{pmatrix}
$$ and
$$y_1 = \frac{2\gamma \bar{z}}{\alpha} \frac{1}{(1-\gamma z)(1-\gamma \bar{z})} \begin{pmatrix}
\frac{\sqrt{3}}{\gamma} (1-\gamma z)(1-\gamma \bar{z})+4+\sqrt{3}(z+\bar{z})\\
z+\sqrt{3} 
\end{pmatrix} .$$ Calculations yield 
$$x_1 = \frac{\gamma}{\alpha} \frac{1}{(1-\gamma z)(1-\gamma \bar{z})} \begin{pmatrix}
1\\
-4-\sqrt{3}\bar{z}
\end{pmatrix}, \quad  y_1 = \frac{2\gamma \bar{z}}{\alpha} \frac{1}{(1-\gamma z)(1-\gamma \bar{z})} \begin{pmatrix}
-1\\
z+\sqrt{3} 
\end{pmatrix}.$$ 
 
\noindent The algorithm stops after at most $ \min(m,n)$ steps, hence in this case after $2$ steps. Then, by Theorem \ref{Tkismultipleofhankel}, the unique analytic superoptimal approximant $\mathcal{A}G$ is given by the formula $$\mathcal{A}G = G - \frac{t_0 y_0 x_0^*}{|h_0|^2}-  \frac{t_1 y_1 x_1^*}{|h_1|^2} .$$ 
All terms of $\mathcal{A}G$ can be calculated now, to give
$$\mathcal{A}G= \frac{\sqrt{2}}{1-\gamma z} \begin{pmatrix}
-\gamma & \sqrt{3}+4\gamma \\
2+\gamma \sqrt{3} -\gamma z & -(\sqrt{3}+4\gamma)(\sqrt{3}+z)
\end{pmatrix} ,$$  which is the \textit{unique superoptimal analytic approximant} for the given $G$  in Problem \ref{applic-ex}.